\def\E{\end{document}}
\documentclass[10pt]{article}
\usepackage{amssymb,amsmath}

\begin{document}
\title{
Critical exponents line, scattering theories for a weighted gradient system of semilinear Schr\"{o}dinger equations
}
  \author{Xianfa Song{\thanks{E-mail: songxianfa2004@163.com(X.F. Song)
 }}\\
\small Department of Mathematics, School of Mathematics, Tianjin University,\\
\small Tianjin, 300072, P. R. China
}

\maketitle
\date{}

\newtheorem{theorem}{Theorem}[section]
\newtheorem{definition}{Definition}[section]
\newtheorem{lemma}{Lemma}[section]
\newtheorem{proposition}{Proposition}[section]
\newtheorem{corollary}{Corollary}[section]
\newtheorem{remark}{Remark}[section]
\renewcommand{\theequation}{\thesection.\arabic{equation}}
\catcode`@=11 \@addtoreset{equation}{section} \catcode`@=12

\begin{abstract}

In this paper, we consider the following Cauchy problem of a weighted(or essential) gradient system of semilinear Schr\"{o}dinger equations
\begin{equation*}
\left\{
\begin{array}{lll}
iu_t+\Delta u=f(|u|^2,|v|^2)u,\quad iv_t+\Delta v=g(|u|^2,|v|^2)v,\quad x\in \mathbb{R}^d,\ t\in \mathbb{R},\\
u(x,0)=u_0(x),\quad v(x,0)=v_0(x),\quad x\in \mathbb{R}^d.
\end{array}\right.
\end{equation*}
Here $d\geq 3$, $f(\cdot,\cdot)$ and $g(\cdot,\cdot)$ are real-valued functions, $(u_0,v_0)$ belongs to $H^1(\mathbb{R}^d)\times H^1(\mathbb{R}^d)$ or $\dot{H}^1(\mathbb{R}^d)\times \dot{H}^1(\mathbb{R}^d)$
or $H^s(\mathbb{R}^d)\times H^s(\mathbb{R}^d)$. Under certain assumptions, we establish the local wellposedness of the $H^1\times H^1$-solution, $\dot{H}^1\times \dot{H}^1$-solution and $H^s\times H^s$-solution of the system with different types of initial data.

If $f(|u|^2,|v|^2)u=\lambda |u|^{\alpha}|v|^{\beta+2}u$ and $g(|u|^2,|v|^2)v=\mu |u|^{\alpha+2}|v|^{\beta}v$ with $\lambda>0$, $\mu>0$, $\alpha\geq 0$ and $\beta\geq 0$, it is surprised that there exists a critical exponents line $\alpha+\beta=2$ when $d=3$ in the following sense: The system always has a unique bounded $H^1\times H^1$-solution for any initial data $(u_0,v_0)\in H^1(\mathbb{R}^3)\times H^1(\mathbb{R}^3)$ if $\alpha+\beta\leq 2$, yet we can find some initial data $(u_0,v_0)\in H^1(\mathbb{R}^3)\times H^1(\mathbb{R}^3)$ such that it doesn't possess the global $H^1\times H^1$-solution if $\alpha+\beta>2$ and $\alpha=\beta$. While when $d=4$, we call $(\alpha,\beta)=(0,0)$ is the critical exponents point in the following sense: The system always has a unique bounded $H^1\times H^1$-solution for any initial data $(u_0,v_0)\in H^1(\mathbb{R}^4)\times H^1(\mathbb{R}^4)$ if $(\alpha,\beta)=(0,0)$, yet there exist some initial data $(u_0,v_0)\in H^1(\mathbb{R}^4)\times H^1(\mathbb{R}^4)$ such that it doesn't have the global $H^1\times H^1$-solution if $\alpha+\beta>0$ and $\alpha=\beta$.

Moreover, we establish the $H^1\times H^1$ and $\Sigma\times \Sigma$ scattering theories for the solution if $\alpha+\beta<2$(i.e., $(\alpha,\beta)$ is below the critical exponents line $\alpha+\beta=2$) when $d=3$, $\Sigma\times \Sigma$ scattering theory for the solution if $(\alpha,\beta)$ is on the critical exponents line $\alpha+\beta=2$ excluding the endpoints when $d=3$, $\dot{H}^1\times \dot{H}^1$ scattering theory for the solution if $(\alpha,\beta)$ is the endpoint of the critical exponents line $\alpha+\beta=2$ when $d=3$ and $(\alpha,\beta)=(0,0)$ when $d=4$. We also establish $\dot{H}^{s_c}\times \dot{H}^{s_c}$ scattering theories for the corresponding solutions if $(\alpha,\beta)$ is above the critical exponents line $\alpha+\beta=2$ when $d=3$ and $(\alpha,\beta)>(0,0)$ when $d\geq 4$. Here $s_c=\frac{d}{2}-\frac{2}{\alpha+\beta+2}$.



\end{abstract}

\newpage

{\bf Contents.}

{\bf 1 Introduction.}

{\bf 2 Preliminaries.}

2.1 Weighted(or essential) gradient system of Schr\"{o}dinger equations.

2.2 Some notations and lemmas about basic harmonic analysis.

{\bf 3 Local wellpossedness results on (\ref{1}).}

3.1 Local wellpossedness of $H^1\times H^1$-solution, $\Sigma\times \Sigma$-solution to (\ref{1}).

3.2 Existence of critical exponents line for (\ref{826x1}) when $d=3$ and critical exponents point for (\ref{826x1})) when $d=4$

3.3 Local wellpossedness of $H^s\times H^s$-solution to (\ref{826x1}).

{\bf 4 $H^1\times H^1$ and $\Sigma\times \Sigma$ scattering theories for (\ref{826x1}) with defocusing nonlinearities when $d=3$.}

4.1 $H^1\times H^1$ and $\Sigma\times \Sigma$ scattering theories for (\ref{826x1}) with defocusing nonlinearities and $\alpha+\beta<2$ when $d=3$.

4.2 $\Sigma\times \Sigma$ scattering theory for (\ref{826x1}) with defocusing nonlinearities and $\alpha+\beta=2$ excluding the endpoints when $d=3$.

{\bf 5 $\dot{H}^1\times \dot{H}^1$ scattering theory for (\ref{826x1}) with defocusing nonlinearities and $(\alpha,\beta)$ is the endpoint of the critical exponents line when $d=3$.}

5.1 Stability result about (\ref{826x1}).

5.2 Concentration compactness and reduction to almost periodic solution.

5.3 Strichartz estimates.

5.4 Impossibility of rapid frequency cascades.

5.5 The frequency-localized weight-coupled interaction Morawetz inequality.

5.6 Impossibility of quasi-solution.

{\bf 6 $\dot{H}^1\times \dot{H}^1$ scattering theory for (\ref{826x1}) with defocusing nonlinearities when $d=4$.}

6.1 Long-time Strichartz estimates.

6.2 The rapid frequency-cascade scenario.

6.3 The quasi-solution scenario.

{\bf 7 $\dot{H}^{s_c}\times \dot{H}^{s_c}$ scattering theory for (\ref{826x1}) in energy-supercritical case when $d\geq 3$.}

 7.1 $\dot{H}^{s_c}\times \dot{H}^{s_c}$ scattering theory for (\ref{826x1}) with special radial initial data.

 7.2 $\dot{H}^{s_c}\times \dot{H}^{s_c}$ scattering theory for (\ref{826x1}) with initial data has many bubbles.

{\bf 8 Some results on (\ref{826x1}) with focusing nonlinearities.}

{\bf References.}

\section{Introduction}
\qquad In this paper, we consider the following Cauchy problem:
\begin{equation}
\label{1}
\left\{
\begin{array}{lll}
iu_t+\Delta u=f(|u|^2,|v|^2)u,\quad iv_t+\Delta v=g(|u|^2,|v|^2)v,\quad x\in \mathbb{R}^d,\ t\in \mathbb{R},\\
u(x,0)=u_0(x),\quad v(x,0)=v_0(x),\quad x\in \mathbb{R}^d.
\end{array}\right.
\end{equation}
Here $d\geq 3$, $f(\cdot,\cdot)$ and $g(\cdot,\cdot)$ are real-valued functions, $(u_0,v_0)$ belongs to $H^1(\mathbb{R}^d)\times H^1(\mathbb{R}^d)$ or $\dot{H}^1(\mathbb{R}^d)\times \dot{H}^1(\mathbb{R}^d)$ or $H^s(\mathbb{R}^d)\times H^s(\mathbb{R}^d)$. It is well known that Schr\"{o}dinger equation often appears in quantum mechanics, in nonlinear optics,  in plasma physics, in the theory of Heisenberg ferromagnet and magnons, and in condensed matter theory.
There are many interesting topics on Schr\"{o}dinger equation, a attractive one is scattering theory.

There are many classical results on the Cauchy problem of the scalar Sch\"{o}dinger equation $iu_t+\Delta u=f(|u|^2)u$, we can refer to the books \cite{ Bourgain19992, Cazenave2003, Tao2006} and the numerous references therein to see more details. Here we would like to briefly review some scattering results on the Cauchy problem of
\begin{equation}
\label{1104w1}
\left\{
\begin{array}{lll}
iu_t+\Delta u=\lambda|u|^pu,\quad x\in \mathbb{R}^d,\ t\in \mathbb{R},\\
u(x,0)=u_0(x),\quad x\in \mathbb{R}^d,
\end{array}\right.
\end{equation}
where $p>0$, $d\geq 3$ and $\lambda\geq 0$, i.e., the equation in (\ref{1104w1}) has the defocusing nonlinearity. By the results of \cite{Barab1984, Strauss1981, Tsutsumi1984}, there no nontrivial solution to (\ref{1104w1}) has scattering states when $p\leq \frac{2}{d}$ even for the $L^2(\mathbb{R}^d)$ topology, while the solution to (\ref{1104w1}) has scattering states for the $L^2(\mathbb{R}^d)$ topology when $\frac{2}{d}<p<\frac{4}{d-2}$. Using pseudoconformal conservation law and decay of solutions in the weighted $L^2(\mathbb{R}^d)$ space, the authors respectively in \cite{Cazenave1992,Ginibre1979,Tsutsumi1985} showed that the solution to (\ref{1104w1}) has scattering states for the weighted $L^2(\mathbb{R}^d)$ topology when $\frac{4}{d-2}>p\geq \alpha_0=\frac{2-d+\sqrt{d^2+12d+4}}{2d}$. Based on Morawetz's estimate and decay of solutions in the energy space, many authors proved that the solution to (\ref{1104w1}) has scattering states for the $H^1(\mathbb{R}^d)$ topology when $\frac{4}{d}<p<\frac{4}{d-2}$, see \cite{Dodson20190, Dodson2019, Dodson2017, Duyckaerts2008, Ginibre19851,Ginibre19852, Murphy2015, Tao2004}. When $p=\frac{4}{d}$, (\ref{1104w1}) is in the mass critical case. We can refer to \cite{Dodson2015, Dodson20193,  Killip2013, Killip20082, Tao20072, Tao2008, Tsutsumi1990} and the references therein to see the scattering results on (\ref{1104w1}). When $p=\frac{4}{d-2}$, (\ref{1104w1}) is in the energy critical case. To establish the scattering theory, the key point is to estimate the spacetime bound. We can refer to \cite{Bourgain19991, Bourgain19992, Colliander2004, Colliander2008, Dodson2012, Grillakis2000, Kenig2006, Killip2019, Killip20102, Killip2012, Killip2013, Ryckman2007, Tao20051, Tao20052, Visan2006, Visan2007, Visan2011} and see the results on $H^1$-scattering  theory or $\dot{H}^1$-scattering theory for (\ref{1104w1}). When $p>\frac{4}{d-2}$, (\ref{1104w1}) is in the energy supercritical case. To study it, the main difficulty is the lack of conservation laws beyond $H^1(\mathbb{R}^d)$ Sobolev space. Recently, the results on $\dot{H}^s$ scattering for (\ref{1104w1}) are established by many authors, see  \cite{Gao2019, Kenig2010, Killip20102, Murphy2014,Murphy20142, Murphy2015, Xie2013} and the references therein. Very recently, in \cite{Beceanu??}, Beceanu et. al. constructed some classes of initial data, which could be arbitrarily large in critical Sobolev space $\dot{H}^{s_c}(\mathbb{R}^d)$ such that the corresponding solutions are globally in time and scatter.

In addition, some authors established the scattering theory for the Schr\"{o}dinger equation with combined power-type nonlinearities. We can refer to \cite{Killip2017, Killip20172, Miao2013, Tao20071, Zhang2006} and the references therein.

Some authors studied the stable manifolds for an orbitally unstable nonlinear Schr\"{o}dinger equation in $\mathbb{R}^3$ and obtained many interesting results including the scattering results, see \cite{ Krieger2006, Nakanishi2012,Schlag2009, Schlag20101, Schlag20102}.

Although there are some results on the scattering theory for a system of Schr\"{o}dinger equations(see\cite{Cassano2015, Farah2017, Xu2016}), some new difficulties arising in the study of the system of Schr\"{o}dinger equations with more general nonlinearities. Besides the local wellposedness results, differing to the scalar case, some Hamiltonian may be difficult to get. For example, although we can establish the local wellposedness result, we even wonder how to define the energy of the coupled system $iu_t+\Delta u=|v|^pu$, $iv_t+\Delta v=|u|^qv$ for general $p,q>0$. However, in this paper, we will show that, if there exist some constants $a$, $b$ and function $G(w,z)$ such that $\frac{\partial G}{\partial w}=af(w,z)$ and $\frac{\partial G}{\partial z}=bg(w,z)$, we call the system $iu_t+\Delta u=f(|u|^2,|v|^2)u$, $iv_t+\Delta v=g(|u|^2,|v|^2)v$ as a weighted(or essential) gradient one. We can define the weighted energy, the weighted mass and the weighted momentum as follows:
\begin{align}
&E_w(u(t),v(t)):=\int_{\mathbb{R}^d}\left[\frac{a}{2}|\nabla u(t,x)|^2+\frac{b}{2}|\nabla v(t,x)|^2+\frac{1}{2}G(|u(t,x)|^2,|v(t,x)|^2)\right]dx,\label{1102x1}\displaybreak\\
&M_w(u(t),v(t)):=\int_{\mathbb{R}^d}[a|u(t,x)|^2+b|v(t,x)|^2]dx,\label{1106x1}\\
&P_w(u(t),v(t)):=\Im\int_{\mathbb{R}^d}[au(t,x)\nabla \bar{u}(t,x)+bv(t,x)\nabla \bar{v}(t,x)]dx.\label{1106x2}
\end{align}
Under certain assumptions, each is Hamiltonian.

Although there are many difficulties, we find some very interesting phenomena on the system of Schr\"{o}dinger equations. We would like to say our special contributions below.

Under certain assumptions on $f(\cdot,\cdot)$, $g(\cdot,\cdot)$ and the initial data $(u_0,v_0)$, we will establish the local wellposedness of the $H^1\times H^1$-solution, $\dot{H}^1\times \dot{H}^1$-solution and $H^s\times H^s$-solution to (\ref{1}) with different types of initial data.

We will mainly consider the asymptotic behavior for the solution to the following special case of (\ref{1})
\begin{equation}
\label{826x1}
\left\{
\begin{array}{lll}
iu_t+\Delta u=\lambda |u|^{\alpha}|v|^{\beta+2}u,\quad iv_t+\Delta v=\mu|u|^{\alpha+2}|v|^{\beta}v,\quad x\in \mathbb{R}^d,\ t>0,\\
u(x,0)=u_0(x),\quad v(x,0)=v_0(x),\quad x\in \mathbb{R}^d,
\end{array}\right.
\end{equation}
where $d\geq 3$, $\lambda\in \mathbb{R}$, $\mu\in \mathbb{R}$, $\alpha\geq 0$ and $\beta\geq 0$. If $\lambda>0$, $\mu>0$, $\alpha\geq 0$ and $\beta\geq 0$, i.e., the nonlinearities are defocusing, it is surprised that there exists a critical exponents line $\alpha+\beta=2$ when $d=3$ in the following sense: The system always has a unique bounded $H^1\times H^1$-solution for any initial data $(u_0,v_0)\in H^1(\mathbb{R}^3)\times H^1(\mathbb{R}^3)$ if $\alpha+\beta\leq 2$, yet we can find some initial data $(u_0,v_0)\in H^1(\mathbb{R}^3)\times H^1(\mathbb{R}^3)$ such that it doesn't possess the global $H^1\times H^1$-solution if $\alpha+\beta>2$ and $\alpha=\beta$. While when $d=4$, we call $(\alpha,\beta)=(0,0)$ is the critical exponents point in the following sense: The system always has a unique bounded $H^1\times H^1$-solution for any initial data $(u_0,v_0)\in H^1(\mathbb{R}^4)\times H^1(\mathbb{R}^4)$ if $(\alpha,\beta)=(0,0)$, yet we can find some initial data $(u_0,v_0)\in H^1(\mathbb{R}^4)\times H^1(\mathbb{R}^4)$ such that it doesn't have the global $H^1\times H^1$-solution if $\alpha+\beta>0$ and $\alpha=\beta$. Note that in \cite{Christ20032}, Christ, Colliander and Tao showed that the Cauchy problem of the scalar equation $iu_t+\Delta u=|u|^{2\sigma}u$ is not well posed in $H^1$ in the energy supercritical case and in \cite{Alazard20091}, Alazard and Carles prove that there exist a sequence of initial data such that $\|u(t)\|_{H^1_x}\rightarrow \infty$ in finite time. Recently, in \cite{Tao2018}, Tao showed that there exist a class of defocusing nonlinear Schr\"{o}dinger systems which the solutions can blow up in finite time by suitably constructing the initial datum. Meanwhile, by the results of \cite{Alazard20091, Alazard20092, Birnir1996, Burq2005, Carles20071, Carles20072, Carles2012, Christ2003,  Kenig2001}, the solution map is highly unstable in energy supercritical case if it exists at all, although one can at least construct global weak solutions, the uniqueness isn't known, see \cite{Ginibre19852, Tao2009}. We think, finding the critical exponents line and critical exponents point, it is a new discovery in the study of a system of Schr\"{o}dinger equations.

Moreover, we define the weight-coupled interaction Morawetz identity. Based on it and decay of the solutions, we establish the $H^1\times H^1$ and $\Sigma\times \Sigma$ scattering theories for the solution if $\alpha+\beta<2$(i.e., $(\alpha,\beta)$ is below the critical exponents line $\alpha+\beta=2$) when $d=3$, $\Sigma\times \Sigma$ scattering theory for the solution if $(\alpha,\beta)$ is on the critical exponents line $\alpha+\beta=2$ excluding the endpoints when $d=3$, $\dot{H}^1\times \dot{H}^1$ scattering theory for the solution if $(\alpha,\beta)$ is the endpoint of the critical exponents line $\alpha+\beta=2$ when $d=3$ and $(\alpha,\beta)=(0,0)$ when $d=4$. We also establish $\dot{H}^{s_c}\times \dot{H}^{s_c}$ scattering theories for the corresponding solutions if $(\alpha,\beta)$ is above the critical exponents line $\alpha+\beta=2$ when $d=3$ and $(\alpha,\beta)>(0,0)$ when $d\geq 4$. Here $s_c=\frac{d}{2}-\frac{2}{\alpha+\beta+2}$.

Now we would like to say something about the techniques applied in this paper.

To establish $H^1\times H^1$ and $\Sigma\times \Sigma$ scattering theories for the solution of (\ref{826x1}) with defocusing nonlinearities and $\alpha+\beta<2$, $(u_0,v_0)\in H^1\times H^1$ when $d=3$, the main step is to obtain the spacetime bound for the solution. We define the weight-coupled interaction Morawetz potential in the form of (\ref{829w1}) and get the bound for
$$\int_{\mathbb{R}}\int_{\mathbb{R}^3}[|u(t,x)|^4+|v(t,x)|^4]dxdt.$$
Using this bound, we obtain the finite global Strichartz norms, which will be applied to establish the $H^1\times H^1$ and $\Sigma\times \Sigma$ scattering theories. Thanks to his famous paper \cite{Morawetz1968}, Morawetz estimate becomes an
important tool to construct scattering operator on the energy space. Although our technical route is inspired by others(see Chapter 7 in \cite{Cazenave2003} and \cite{Tao20071}), we define the weight-coupled interaction Morawetz identity which is a very useful tool in the scattering theory for the system of Schr\"{o}dinger equations.

To establish $\Sigma\times \Sigma$ scattering theory for the solution of (\ref{826x1}) with defocusing nonlinearities and $\alpha+\beta=2$ excluding $(\alpha,\beta)=(0,2)$ and $(\alpha,\beta)=(2,0)$, $(u_0,v_0)\in \Sigma\times \Sigma$ when $d=3$, we establish the weight-coupled pseudoconformal conservation law and obtain decay of the solution, which can be applied to establish the $\Sigma\times \Sigma$ scattering theory. Although our technique route follows that dealing with the scalar Sch\"{o}dinger equation(see Chapter 7 in \cite{Cazenave2003} and  \cite{Ginibre1979,Ginibre19851, Ginibre19852, Ginibre1992}), the weight-coupled pseudoconformal conservation law is essential for the study of the asymptotic behavior for the weight-gradient system of Schr\"{o}dinger equations with initial data belonging to $\Sigma\times \Sigma$.

It is very difficult to obtain $\dot{H}^1\times \dot{H}^1$ scattering theory for the solution of (\ref{826x1}) with defocusing nonlinearities and $\alpha+\beta=2$, i.e., $(\alpha,\beta)$ is on the critical exponents line $\alpha+\beta=2$ when $d=3$. The main reason is that we cannot establish nonlinear estimates of the solution for all $(\alpha,\beta)$ satisfying $\alpha+\beta=2$, which consequently leads the lack of the corresponding $\dot{H}^1\times \dot{H}^1$ stability result. However, if $(\alpha,\beta)=(0,2)$ or $(\alpha,\beta)=(2,0)$, i.e., $(\alpha,\beta)$ is the endpoint of the critical exponents line $\alpha+\beta=2$, we can use the technique of reduction to almost periodic solutions, which is now fairly standard in the field of
nonlinear dispersive PDE, especially in the setting of NLS. In \cite{Keraani2006}, Keraani originally
established the existence of minimal blowup solutions to NLS, while in \cite{Killip20101}, Kenig and
Merle were the first to use them as a tool to prove global well-posedness. Such argument had been used by many authors in a variety of settings and has proven to be extremely effective, see \cite{Dodson2012, Holmer2008, Kenig2010, Killip2009, Killip20101, Killip20102, Killip2012, Killip2013, Killip20082, Miao2014, Murphy2014, Murphy20142, Ryckman2007, Tao20051, Tao2008, Visan2011}. It may be outlined as follows: Step 1, we assume contradictorily a spacetime bound of the solution doesn't hold, then there must be a minimal almost-counterexample which is a minimal-energy solution with enormous spacetime norm. Step 2, we show that it must have good tightness and equicontinuity properties because it has minimal-energy. Step 3, we define the weight-coupled interaction Morawetz inequality and prove that a solution must undergo a dramatic change of spatial scale in a short span of time. Step 4, We show that such a rapid change is inconsistent with simultaneous conservation of mass and energy.

The $\dot{H}^1\times \dot{H}^1$ scattering theory for the solution of (\ref{826x1}) with defocusing nonlinearities and $(\alpha,\beta)=(0,0)$ when $d=4$ is also established by using the concentration-compactness approach and argument by contradiction. The weight-coupled interaction Morawetz inequality plays an important role in the course of discussions.

The main difficulty in the study of energy supercritical Schr\"{o}dinger equation is the lack of the conservation
laws beyond $H^1(\mathbb{R}^d)$ Sobolev space. In \cite{Killip20102}, Killip and Visan developed the technique treating energy-critical NLS and applied in the  energy supercritical case when $d\geq 5$, and showed that critical $\dot{H}^s$-bounds imply scattering. The work \cite{Miao2014} treated energy supercritical case when $d=4$ by a similar approach. In \cite{Murphy2015}, Murphy applied the method of reduction to almost periodic solutions and Lin-Strauss Morawetz inequality, established the scattering theory for Schr\"{o}dinger equation in energy supercritical case. Some authors borrowed the ideas from the study of wave equations(see \cite{Beceanu??, Kenig2006, Kenig2008, Kenig20111, Kenig20112, Killip20111,Killip20112, Merle2015,  Roy2011, Wang2016}), they used concentration-compactness-rigidity method or constructing the explicit formula of the outgoing and incoming of the radial linear flow, and obtained some priori estimates for the solution, then established the $\dot{H}^{s_c}$-scattering theory. In this paper, we use the method of constructing the explicit formula of the outgoing and incoming of the radial linear flow(inspired by \cite{Beceanu??, Tao2004}), and establish the $\dot{H}^{s_c}\times \dot{H}^{s_c}$-scattering theory.

Our first result is about the local wellposedness of $H^1\times H^1$ solution to (\ref{1}).

{\bf Theorem 1(Local wellposedness of the $H^1\times H^1$ solution to (\ref{1})).} {\it Assume that $u_0(x)\in H^1(\mathbb{R}^d)$, $v_0(x)\in H^1(\mathbb{R}^d)$, $d\geq 3$. Then there exists a unique, strong $H^1\times H^1$-solution of (\ref{1}), defined on a maximal time interval $(-T_{\min}, T_{\max})$ in one of the following cases:

Case 1.  $d\geq 3$,
\begin{align*}
&f(|u|^2, |v|^2)u=f_1(|u|^2,|v|^2)u+...+f_m(|u|^2,|v|^2)u,\\
&g(|u|^2, |v|^2)v=g_1(|u|^2,|v|^2)v+...+g_m(|u|^2,|v|^2)v.
\end{align*}
Suppose that $f_k\in C(\mathbb{R}\times \mathbb{R}; \mathbb{R})$, $g_k\in C(\mathbb{R}\times \mathbb{R}; \mathbb{R})$, there exist $2\leq r_k<\frac{2d}{d-2}$ and $2\leq \rho_k<\frac{2d}{d-2}$ such that
\begin{align}
\|f_k(|u|^2,|v|^2)u\|_{W^{1,\rho'_k}}+\|g_k(|u|^2,|v|^2)v\|_{W^{1,\rho'_k}}\leq C(M)(1+\|u\|_{W^{1,r_k}}+\|v\|_{W^{1,r_k}})\label{8241}
\end{align}
 for all $u, v\in H^1(\mathbb{R}^d)\cap W^{1,r_k}(\mathbb{R}^d)$ such that $\|u\|_{H^1}\leq M$, $\|v\|_{H^1}\leq M$, and
\begin{align}
\|f_k(|u_1|^2,|v_1|^2)u_1-f_k(|u_2|^2,|v_2|^2)u_2\|_{L^{\rho'_k}}&\leq C(M)(\|u_1-u_2\|_{L^{r_k}}+\|v_1-v_2\|_{L^{r_k}})\label{8242}\\
\|g_k(|u_1|^2,|v_1|^2)v_1-g_k(|u_2|^2,|v_2|^2)v_2\|_{L^{\rho'_k}}&\leq C(M)(\|u_1-u_2\|_{L^{r_k}}+\|v_1-v_2\|_{L^{r_k}})\label{8243}
\end{align}
 for all $u_1, u_2, v_1, v_2\in H^1(\mathbb{R}^d)$ such that $\|u_1\|_{H^1}\leq M$, $\|u_2\|_{H^1}\leq M$, $\|v_1\|_{H^1}\leq M$ and $\|v_2\|_{H^1}\leq M$, $k=1,...,m$.

Case 2. $d=3$, $f(|u|^2,|v|^2)u=\lambda|u|^{\alpha}|v|^{\beta+2}u$, $g(|u|^2,|v|^2)v=\mu |u|^{\alpha+2}|v|^{\beta}v$, $\alpha\geq 0$, $\beta\geq 0$, $\alpha+\beta=2$, $\lambda, \mu\in \mathbb{R}$ and $\lambda\mu>0$.

Case 3. $d=4$, $f(|u|^2,|v|^2)u=\lambda |v|^2u$, $g(|u|^2,|v|^2)v=\mu |u|^2v$, $\lambda, \mu\in \mathbb{R}$ and $\lambda\mu>0$.

Moreover, for any admissible pair $(q,r)$,
$$
u \in L^q_{loc}((-T_{\min}, T_{\max}), W^{1,r}(\mathbb{R}^d),\quad  v \in L^q_{loc}((-T_{\min}, T_{\max}), W^{1,r}(\mathbb{R}^d).
$$
And the following properties hold:

(i) There is the blowup alternative in the following sense: $\|u(t,x)\|_{H^1_x}+\|v(t,x)\|_{H^1_x}\rightarrow \infty$
as $t\uparrow T_{\max}$ if $T_{\max}<\infty$ and as $t\downarrow T_{\min}$ if $T_{\min}<\infty$.

(ii) $(u,v)$ depends continuously on $(u_0,v_0)$, i.e., there exists $T>0$ which depends on $(u_0,v_0)$ and satisfies: if $(u_n,v_n)$ is the solution of (\ref{1}) with the corresponding initial data $(u_{n0},v_{n0})\rightarrow (u_0,v_0)$, then $(u_n,v_n)$ is defined on $[-T,T]$ for $n$ large enough and $(u_n,v_n)\rightarrow (u,v)$ in $C([-T,T], L^p(\mathbb{R}^d))\times C([-T,T], L^p(\mathbb{R}^d))$ for all $2\leq p<\frac{2d}{d-2}$.

(iii) The conservation of mass for each component of $(u,v)$, i.e., $\|u\|_{L^2}=\|u_0\|_{L^2}$ and $\|v\|_{L^2}=\|v_0\|_{L^2}$.

(iv) If there exist $a, b\in \mathbb{R}$ and real value functions $G_k(w,z)$ such that $\frac{\partial G_k}{\partial w}=a f_k(w,z)$, $\frac{\partial G_k}{\partial z}=b g_k(w,z)$ for $w\geq 0$ and $z\geq 0$, then there are conservation laws of weighted mometum and weighted energy for the weighted gradient system of Schr\"{o}dinger equations. For example, if we take $(c_1,c_2)=(a,b)$ in Case 1, $(c_1,c_2)=(\frac{\alpha+2}{2\lambda},\frac{\beta+2}{2\mu})$ in Case 2 and in Case 3$((\alpha,\beta)=(0,0))$, then
\begin{align*}
 &\Im\int_{\mathbb{R}^d}[c_1u(t,x)\nabla \bar{u}(t,x)c_2v(t,x)\nabla \bar{v}(t,x)]dx=\Im\int_{\mathbb{R}^d}[c_1u_0\nabla \bar{u}_0+c_2v_0\nabla \bar{v}_0]dx,\\
 &E_w(u,v)=E_w(u_0,v_0).
\end{align*}
Here
\begin{align*}
&E_w(u,v):=\int_{\mathbb{R}^d}[\frac{c_1}{2}|\nabla u|^2+\frac{c_2}{2}|\nabla v|^2+\frac{1}{2}\sum_{k=1}^mG_k(|u|^2,|v|^2)]dx\quad {\rm in \ Case}\ 1,\\
&E_w(u,v):=\int_{\mathbb{R}^3}[\frac{\alpha+2}{2\lambda}|\nabla u|^2+\frac{\beta+2}{2\mu}|\nabla v|^2+|u|^{\alpha+2}|v|^{\beta+2}]dx\quad {\rm in \ Case}\ 2\ {\rm and \ in \ Case}\ 3.
\end{align*}
Moreover, if $(u_0,v_0)\in H^1(\mathbb{R}^d)\times H^1(\mathbb{R}^d)$ and $(xu_0(x), xv_0(x))\in L^2(\mathbb{R}^d)\times L^2(\mathbb{R}^d)$, then the vector-valued function $t\rightarrow (|\cdot|u(t,\cdot), |\cdot|v(t,\cdot)|)$ belongs to $C((-T_{\min},T_{\max}), L^2(\mathbb{R}^d)\times L^2(\mathbb{R}^d))$ and
$(u,v)\in \Sigma\times \Sigma$, where
\begin{align}
\Sigma=\{w: w\in H^1(\mathbb{R}^d),\quad |xw|\in L^2(\mathbb{R}^d).\label{1208x1}
\end{align}
Furthermore,
\begin{align}
Y(t)=\int_{\mathbb{R}^d}|x|^2[c_1|u(t,x)|^2+c_2|v(t,x)|^2]dx\label{1113w1}
\end{align}
is in $C^2(-T_{\min},T_{\max})$ with
\begin{align}
Y'(t)&=4\Im\int_{\mathbb{R}^d}[c_1\bar{u}x\cdot \nabla u+c_2\bar{v}x\cdot \nabla v]dx,\label{1113w2}\\
Y''(t)&=4\int_{\mathbb{R}^d}[d(c_1f(|u|^2,|v|^2)|u|^2+c_2g(|u|^2,|v|^2)|v|^2)-(d+2)G(|u|^2,|v|^2)]dx\nonumber\\
&\quad +16E_w(u,v)\quad {\rm in\  Case\  1},\label{1113w3}\\
Y''(t)&=2[3(\alpha+\beta)+2]sgn(\lambda)\int_{\mathbb{R}^d}|u|^{\alpha+2}|v|^{\beta+2}dx+8E_w(u,v)\quad {\rm in\  Cases\  2\ and\ 3}\label{1113w31}
\end{align}
for all $t\in (-T_{\min},T_{\max})$.

}

We find some more interesting phenomena below.

{\bf Theorem 2.} {\it 1. There exists the {\bf critical exponents line} $\alpha+\beta=2$ for (\ref{826x1}) with $\lambda>0$, $\mu>0$, $\alpha\geq 0$, $\beta\geq 0$  when $d=3$ in the following sense:
If $\alpha+\beta\leq 2$, then for every $(u_0,v_0)\in H^1(\mathbb{R}^3)\times H^1(\mathbb{R}^3)$, (\ref{826x1}) always has a unique bounded $H^1\times H^1$-solution; If $\alpha+\beta>2$ and $\alpha=\beta$, there exists $(u_0,v_0)$ such that (\ref{826x1}) doesn't possess the global $H^1\times H^1$-solution.

2. $(\alpha,\beta)=(0,0)$ is the {\bf critical exponents point} for (\ref{826x1}) with $\lambda>0$, $\mu>0$, $(\alpha,\beta)=(0,0)$ when $d=4$ in the following sense:
If $(\alpha,\beta)=(0,0)$, then for every $(u_0,v_0)\in H^1(\mathbb{R}^4)\times H^1(\mathbb{R}^4)$, (\ref{826x1}) always has a unique bounded $H^1\times H^1$-solution; If $\alpha+\beta>0$ and $\alpha=\beta$, there exists $(u_0,v_0)$ such that (\ref{826x1}) doesn't possess the global $H^1\times H^1$-solution.
}

We will compare the critical exponents line $\alpha+\beta=2$ for (\ref{826x1}) with the critical Fujita exponents line $p_1+q=\frac{2}{d}$ when $p_1>1$ and the critical Fujita exponents curve
$pq=\frac{2}{d}+p+c(p_1,q_2)$  or $pq=\frac{2}{d}+q+c(p_1,q_2)$ when $0\leq p_1\leq 1$ for the following system of parabolic equations
\begin{equation}
\label{1108w1}
\left\{
\begin{array}{lll}
w_t-\Delta w= w^{p_1}z^q,\quad z_t-\Delta z=w^pz^{q_2},\quad x\in \mathbb{R}^d,\ t>0,\\
w(x,0)=w_0(x),\quad z(x,0)=z_0(x),\quad x\in \mathbb{R}^d,
\end{array}\right.
\end{equation}
Escobedo and Levine, in \cite{Escobedo1995}(we also can refer to \cite{Escobedo1991} and see the related results), proved that

{\it A. Suppose that $p_1>1$.

(Ai). If $(p_1+q-1)>\frac{2}{d}$, then there are both nontrivial global solutions and
nonglobal solutions.

(Aii). If $(p_1+q-1)\leq \frac{2}{d}$, then every nontrivial solution is nonglobal.

B. Suppose that $0\leq p_1\leq 1$ and $pq>0$.

(Bi). If $\nu=\max(\lambda,\mu)<0$, then all solutions are global.

(Bii). If $0\leq \nu<\frac{2}{d}$, then there are both global nontrivial solutions and nonglobal
solutions.

(Biii). If $\nu\geq \frac{2}{d}$, then all nontrivial solutions are nonglobal.

Here
$$
\lambda=\frac{q_2-q-1}{p_1q_2-pq-(p_1+q_2)+1},\quad \mu=\frac{p_1-p-1}{p_1q_2-pq-(p_1+q_2)+1}
$$
and $p_1q_2-pq-(p_1+q_2)+1\neq 0$.
}

In another word, comparing their results with ours, we find that, there exists critical Fujita exponents line or critical Fujita exponents curve as the watershed to judge whether the system of parabolic equations (\ref{1108w1}) has $L^{\infty}_{t,x}\times L^{\infty}_{t,x}$-solution or not, while  there exists critical exponents line(when $d=3$) or point(when $d=4$) as the watershed to judge whether the weighted gradient system of Schr\"{o}dinger equations (\ref{826x1}) always has a unique $L^{\infty}_tH^1_x\times L^{\infty}_tH^1_x$-solution or not. In this sense, finding the critical exponents line and critical exponents point, it is a subtle result in the direction of studying a system of Schr\"{o}ding equations.

Our third theorem is about the local wellposedness of $H^s\times H^s$-solution to (\ref{826x1}).

{\bf Theorem 3(Local wellposedness of the $H^s\times H^s$-solution to (\ref{826x1})).} {\it
Suppose that $d\geq 3$, $\alpha\geq 0$, $\beta\geq 0$ and $s>0$ satisfy
\begin{align}
0<s<\frac{d}{2},\quad \alpha+\beta+2\leq \frac{4}{d-2s},\label{8271}
\end{align}
and if $\alpha+\beta+2$ is not an even integer,
\begin{align}
[s]<\alpha+\beta+2.\label{8272}
\end{align}
Let $(\gamma, \rho)$ be the admissible pair, defined by
\begin{align}
\gamma=\frac{4(\alpha+\beta+4)}{(\alpha+\beta+2)(d-2s)},\quad \rho=\frac{(\alpha+\beta+4)d}{d+(\alpha+\beta+2) s}.
\end{align}
Then, for every $(u_0,v_0)\in H^s(\mathbb{R}^d)\times H^s(\mathbb{R}^d)$, there exist $T^*_{\min}=T^*_{\min}(u_0,v_0)$, $T^*_{\max}=T^*_{\max}(u_0,v_0)$ and a solution $(u,v)\in [C([-T^*_{\min},T^*_{\max});H^s(\mathbb{R}^d))\cap L^{\gamma}_{loc}((-T^*_{\min},T^*_{\max});B^s_{\rho,2}(\mathbb{R}^d))]^2$ of (\ref{826x1}) which satisfies the following properties:

(i) $(u,v)\in [L^q((-T,T); B^s_{r,2}(\mathbb{R}^d))]^2$ for every admissible pair $(q,r)$ and every $T<T^*_{\min}$ and $T<T^*_{\max}$;

(ii) $(u,v)$ is unique in $[L^{\gamma}((-T,T);B^s_{\rho,2}(\mathbb{R}^d))]^2$ for every $T<T^*_{\min}$ and $T<T^*_{\max}$;

(iii) If $s\geq 1$, $E(u,v)=E(u_0,v_0)$ for every $-T^*_{\min}<t<T^*_{\max}$;

(iv) If $T^*_{\max}<\infty$, then $\|u\|_{L^{\gamma}((0,T^*_{\max}); B^s_{\rho,2}(\mathbb{R}^d))}+\|v\|_{L^{\gamma}((0,T^*_{\max}); B^s_{\rho,2}(\mathbb{R}^d))}=\infty$; If $T^*_{\min}<\infty$, then $\|u\|_{L^{\gamma}((-T^*_{\min},0); B^s_{\rho,2}(\mathbb{R}^d))}+\|v\|_{L^{\gamma}((-T^*_{\min},0); B^s_{\rho,2}(\mathbb{R}^d))}=\infty$;

(v) If $\alpha+\beta+2<\frac{4}{d-2s}$ and $ T^*_{\max}<\infty$, then
\begin{align}
&\lim_{t\uparrow T^*_{\max}}[\|(-\Delta)^{\frac{s}{2}}u(t)\|_{L^2(\mathbb{R}^d)}+\|(-\Delta)^{\frac{s}{2}}v(t)\|_{L^2(\mathbb{R}^d)}]=\infty,\label{827x1}\\
&[\|(-\Delta)^{\frac{s}{2}}u(t)\|_{L^2(\mathbb{R}^d)}+\|(-\Delta)^{\frac{s}{2}}v(t)\|_{L^2(\mathbb{R}^d)}]\geq \frac{C}{(T^*_{\max}-t)^{\frac{1}{\alpha+\beta+2}-\frac{d-2s}{4}}};\label{827x2}
\end{align}
If $T^*_{\min}<\infty$, then
\begin{align}
&\lim_{t\downarrow -T^*_{\min}}[\|(-\Delta)^{\frac{s}{2}}u(t)\|_{L^2(\mathbb{R}^d)}+\|(-\Delta)^{\frac{s}{2}}v(t)\|_{L^2(\mathbb{R}^d)}]=\infty,\label{827x3}\\
&[\|(-\Delta)^{\frac{s}{2}}u(t)\|_{L^2(\mathbb{R}^d)}+\|(-\Delta)^{\frac{s}{2}}v(t)\|_{L^2(\mathbb{R}^d)}]\geq \frac{C}{(T^*_{\min}+t)^{\frac{1}{\alpha+\beta+2}-\frac{d-2s}{4}}};\label{827x4}
\end{align}

(vi) If $\alpha+\beta+2=\frac{4}{d-2s}$ and $\|(-\Delta)^{\frac{s}{2}}u_0\|_{L^2(\mathbb{R}^d)}+\|(-\Delta)^{\frac{s}{2}}v_0\|_{L^2(\mathbb{R}^d)}$ is sufficiently small, then $T^*_{\min}=\infty$, $T^*_{\max}=\infty$ and $(u,v)\in [L^q((-\infty,+\infty); B^s_{r,2}(\mathbb{R}^d))]^2$ for every admissible pair $(q,r)$.

(vii) There exists $T>0$ satisfying $T<T^*_{\min}$ and $T<T^*_{\max}$, such that if $(u_{0n},v_{0n})$ is a sequence in $H^s(\mathbb{R}^d)\times H^s(\mathbb{R}^d)$ and $(u_{0n},v_{0n})\rightarrow (u,v)\in H^s(\mathbb{R}^d)\times H^s(\mathbb{R}^d)$, then for $n$ large enough, $T<T^*_{\min}(u_{0n},v_{0n})$, $T<T^*_{\max}(u_{0n},v_{0n})$ and the corresponding solutions $(u_n,v_n)$ of (\ref{826x1}) form a bounded sequence in $[L^q((-T,T);B^s_{r,2}(\mathbb{R}^d))]^2$. And for every admissible pair $(q,r)$, $(u_n,v_n)\rightarrow (u,v)$ in $[L^q((-T,T);L^r(\mathbb{R}^d))]^2$ and
$[C((-T,T);H^{s-\epsilon}(\mathbb{R}^d))]^2$ for every $\epsilon>0$.
}

The fourth theorem will establish $H^1\times H^1$ and $\Sigma\times \Sigma$ scattering theories for (\ref{826x1}) with defocusing nonlinearities and $\alpha+\beta<2$ when $d=3$. That is, the exponents pair $(\alpha,\beta)$ is below the critical exponents line $\alpha+\beta=2$. The results can be stated as follows.

{\bf Theorem 4($H^1\times H^1$  and $\Sigma\times \Sigma$ scattering theories for (\ref{826x1}) when $(\alpha,\beta)$ is below the critical exponents line).} {\it Let $(u,v)$ be the global solution of (\ref{826x1}), $\lambda>0$, $\mu>0$, $\alpha\geq 0$, $\beta\geq 0$ and $\alpha+\beta<2$ when $d=3$. Assume that $(u_0(x), v_0(x))\in H^1(\mathbb{R}^3)\times H^1(\mathbb{R}^3)$. There exist unique scattering states $(u_+,v_+)$ and $(u_-,v_-)$ such that
\begin{align}
&\|e^{-it\Delta}u(t)-u_+\|_{H^1_x}+\|e^{-it\Delta}v(t)-v_+\|_{H^1_x}\longrightarrow 0\quad {\rm as} \quad t\rightarrow +\infty,\label{1023w1}\\
&\|e^{-it\Delta}u(t)-u_-\|_{H^1_x}+\|e^{-it\Delta}v(t)-v_-\|_{H^1_x}\longrightarrow 0\quad {\rm as} \quad t\rightarrow -\infty.\label{1023w2}
\end{align}
Moreover, if $(u_0(x), v_0(x))\in \Sigma\times \Sigma$, then
\begin{align}
&\|e^{-it\Delta}u(t)-u_+\|_{\Sigma}+\|e^{-it\Delta}v(t)-v_+\|_{\Sigma}\longrightarrow 0\quad {\rm as} \quad t\rightarrow +\infty,\label{11201}\\
&\|e^{-it\Delta}u(t)-u_-\|_{\Sigma}+\|e^{-it\Delta}v(t)-v_-\|_{\Sigma}\longrightarrow 0\quad {\rm as} \quad t\rightarrow -\infty.\label{11202}
\end{align}
  }

The fifth theorem is about $\Sigma\times \Sigma$ scattering theory for (\ref{826x1}) when $(\alpha,\beta)$ is on the critical exponents line excluding the endpoints $(\alpha,\beta)=(0,2)$ and $(\alpha,\beta)=(2,0)$.

{\bf Theorem 5($\Sigma\times \Sigma$ scattering theory for (\ref{826x1}) when $d=3$ and $(\alpha,\beta)$ is on the critical exponents line excluding the endpoints).} {\it Let $(u,v)$ be the global solution of (\ref{826x1}), $\lambda>0$, $\mu>0$, $\alpha>0$, $\beta>0$ and $\alpha+\beta=2$ when $d=3$ and $(u_0(x), v_0(x))\in \Sigma\times \Sigma$. Then there exist scattering states $(u_+,v_+)$ and $(u_-,v_-)$ such that
\begin{align}
&\|e^{-it\Delta}u(t)-u_+\|_{\Sigma}+\|e^{-it\Delta}v(t)-v_+\|_{\Sigma}\longrightarrow 0\quad {\rm as} \quad t\rightarrow +\infty,\label{1121x5}\\
&\|e^{-it\Delta}u(t)-u_-\|_{\Sigma}+\|e^{-it\Delta}v(t)-v_-\|_{\Sigma}\longrightarrow 0\quad {\rm as} \quad t\rightarrow -\infty.\label{1121x6}
\end{align}
  }

If $(\alpha,\beta)$ is the endpoint of the critical exponents line $\alpha+\beta=2$, i.e., $(\alpha,\beta)=(0,2)$ or $(\alpha,\beta)=(2,0)$, our results can be stated as follows.

{\bf Theorem 6($\dot{H}^1\times \dot{H}^1$ scattering theory for (\ref{826x1}) when $d=3$ and $(\alpha,\beta)$ is the endpoint of the critical exponents line).} {\it Let $\lambda>0$, $\mu>0$ and $(u_0,v_0)\in \dot{H}^1(\mathbb{R}^3)\times \dot{H}^1(\mathbb{R}^3)$. Then there exists unique global strong solution $(u,v)\in [C^0_t\dot{H}^1_x(\mathbb{R}\times \mathbb{R}^3)]^2$ of (\ref{826x1}) if $(\alpha,\beta)=(0,2)$ or $(\alpha,\beta)=(2,0)$. And
\begin{align}
\int_{\mathbb{R}}\int_{\mathbb{R}^3}[|u(t,x)|^{10}+|v(t,x)|^{10}]dxdt\leq C(\|u_0\|_{\dot{H}_x}+\|v_0\|_{\dot{H}_x}).\label{931}
\end{align}
Moreover, there exist $(u_{+},v_{+})\in \dot{H}^1\times \dot{H}^1$ and $(u_{-},v_{-})\in \dot{H}^1\times \dot{H}^1$ such that
\begin{align}
&\|u(t)-e^{it\Delta}u_{+}\|_{\dot{H}^1_x}+\|v(t)-e^{it\Delta}v_{+}\|_{\dot{H}^1_x}\rightarrow 0\quad {\rm as}\quad t\rightarrow +\infty,\label{932}\\
&\|u(t)-e^{it\Delta}u_{-}\|_{\dot{H}^1_x}+\|v(t)-e^{it\Delta}v_{-}\|_{\dot{H}^1_x}\rightarrow 0\quad {\rm as}\quad t\rightarrow -\infty.\label{933}
\end{align}
 On the other hand, for any $(u_+,v_+)\in \dot{H}^1\times \dot{H}^1$ and $(u_-,v_-)\in \dot{H}^1\times \dot{H}^1$, there exists a unique global solution $(u,v)$ of (\ref{826x1}) with $(\alpha,\beta)=(0,2)$ or $(\alpha,\beta)=(2,0)$ such that (\ref{932}) and (\ref{933}) are true.

 }

Similar, if $(\alpha,\beta)$ is the critical exponents point when $d=4$, i.e., $(\alpha,\beta)=(0,0)$, we will establish $\dot{H}^1\times \dot{H}^1$ scattering theory for (\ref{826x1}) below.

{\bf Theorem 7($\dot{H}^1\times \dot{H}^1$ scattering theory for (\ref{826x1}) when $d=4$).} {\it Let $\lambda>0$, $\mu>0$ and $(u_0,v_0)\in \dot{H}^1(\mathbb{R}^4)\times \dot{H}^1(\mathbb{R}^4)$. Then there exists unique global strong solution $(u,v)\in [C^0_t\dot{H}^1_x(\mathbb{R}\times \mathbb{R}^4)]^2$ of (\ref{826x1}) if $(\alpha,\beta)=(0,0)$. And
\begin{align}
\int_{\mathbb{R}}\int_{\mathbb{R}^4}[|u(t,x)|^6+|v(t,x)|^6]dxdt\leq C(\|u_0\|_{\dot{H}^1_x}+\|v_0\|_{\dot{H}^1_x}).\label{931'}
\end{align}
Moreover, there exist $(u_{+},v_{+})\in \dot{H}^1\times \dot{H}^1$ and $(u_{-},v_{-})\in \dot{H}^1\times \dot{H}^1$ such that
\begin{align}
&\|u(t)-e^{it\Delta}u_{+}\|_{\dot{H}^1_x}+\|v(t)-e^{it\Delta}v_{+}\|_{\dot{H}^1_x}\rightarrow 0\quad {\rm as}\quad t\rightarrow +\infty,\label{932'}\\
&\|u(t)-e^{it\Delta}u_{-}\|_{\dot{H}^1_x}+\|v(t)-e^{it\Delta}v_{-}\|_{\dot{H}^1_x}\rightarrow 0\quad {\rm as}\quad t\rightarrow -\infty.\label{933'}
\end{align}
 On the other hand, for any $(u_+,v_+)\in \dot{H}^1\times \dot{H}^1$ and $(u_-,v_-)\in \dot{H}^1\times \dot{H}^1$, there exists a unique global solution $(u,v)$ of (\ref{826x1}) with $(\alpha,\beta)=(0,2)$ or $(\alpha,\beta)=(2,0)$ such that (\ref{932'}) and (\ref{933'}) are true.
 }

The eighth theorem will establish the $\dot{H}^{s_c}\times \dot{H}^{s_c}$ scattering theory for the solution of (\ref{826x1}) with special radial initial data when $d=3,4,5$. In order to state the results, we need some notations below. For a small fixed constant $\epsilon_0>0$, denote
\begin{align*}
s_1=\max\{\frac{s_c}{d}-\epsilon_0, s_c-\frac{4d-1}{4d-2}+\frac{2}{d}+\epsilon_0\}
\end{align*}
and
\begin{align*}
s_2=\max\{-\epsilon_0, s_c-\frac{d-2}{2(d-1)}+\frac{d-1}{2d-1}+\epsilon_0\}.
\end{align*}
Let $\hat{W}^{s,1}(\mathbb{R}^d)$ be the space of functions
$$
\|h\|_{\hat{W}^{s,1}(\mathbb{R}^d)}:=\|\langle\xi\rangle^s\hat{h}\|_{L^1(\mathbb{R}^d)}<+\infty.
$$

Now the results can be stated as follows.

{\bf Theorem 8($\dot{H}^{s_c}\times \dot{H}^{s_c}$ scattering theory for the solution of (\ref{826x1}) with special radial initial data).}  {\it
Let $\alpha>0$, $\beta>0$ and $d=3,4,5$. Assume that there exists a small constant $\delta_0$ such that
the radial functions $w_1$, $w_2$, $z_1$ and $z_2$ satisfy $supp z_1\in\{x:|x|\leq 1\}$ and $supp z_2\in\{x:|x|\leq 1\}$,
\begin{align}
&\|\chi_{\leq 1}w_1\|_{\dot{H}^{s_c}_x(\mathbb{R}^d)}+\|\chi_{\leq 1}w_2\|_{\dot{H}^{s_c}_x(\mathbb{R}^d)}+\|\chi_{\geq 1}w_1\|_{\dot{H}^{s_1}_x(\mathbb{R}^d)}+\|\chi_{\geq 1}w_2\|_{\dot{H}^{s_1}_x(\mathbb{R}^d)}\leq \delta_0,\label{1011w1}\\
& \|\langle\xi\rangle^{s_2}\hat{z_1}\|_{L^1(\mathbb{R}^d)}+
\|\langle\xi\rangle^{s_2}\hat{z_2}\|_{L^1(\mathbb{R}^d)}\leq \delta_0.\label{1011w2}
\end{align}
Then if the initial data $(u_0,v_0)$ is of form
$$
u_0=w_{1+}+z_1,\quad v_0=w_{2+}+z_2,\quad ({\rm or}\quad u_0=w_{1-}+z_1,\quad v_0=w_{2-}+z_2 ),
$$
the corresponding solution $(u,v)$ exists globally forward ( or backward) in time and
$$
(u,v)\in [C(\mathbb{R}^+;H^{s_1}(\mathbb{R}^d)+\hat{W}^{s_2,1}(\mathbb{R}^d))]^2\quad ( {\rm or}\quad
(u,v)\in [C(\mathbb{R}^-;H^{s_1}(\mathbb{R}^d)+\hat{W}^{s_2,1}(\mathbb{R}^d))]^2 ).
$$
Here $w_{1+}$, $w_{1-}$, $w_{2+}$ and $w_{2-}$ are the modified outgoing and incoming components of $w_1$ and $w_2$.
Furthermore, there exists $(u_+,v_+)\in [H^{s_1}(\mathbb{R}^d)+\hat{W}^{s_2,1}(\mathbb{R}^d)]^2$(or $(u_-,v_-)\in [H^{s_1}(\mathbb{R}^d)+\hat{W}^{s_2,1}(\mathbb{R}^d)]^2$) such that
\begin{align}
&\lim_{t\rightarrow +\infty}[\|u(t)-e^{it\Delta}u_+\|_{\dot{H}^{s_c}_x(\mathbb{R}^d)}+\|v(t)-e^{it\Delta}v_+\|_{\dot{H}^{s_c}_x(\mathbb{R}^d)}]=0\label{1011w3}\\
({\rm or}\quad &\lim_{t\rightarrow -\infty}[\|u(t)-e^{it\Delta}u_-\|_{\dot{H}^{s_c}_x(\mathbb{R}^d)}+\|v(t)-e^{it\Delta}v_-\|_{\dot{H}^{s_c}_x(\mathbb{R}^d)}]=0.\quad )\label{1011w4}
\end{align}
}

Last, we state $\dot{H}^{s_c}\times \dot{H}^{s_c}$ scattering theory for the solution of (\ref{826x1}) with initial data has many bubbles. We introduce the following condition

{\bf (C9)} {\it Given a constant $\epsilon\in (0,1]$. Suppose that $w$ is of the form
\begin{align*}
w=\sum_{k=0}^{+\infty}w_k, \quad w_k\in \dot{H}^{s_c}(\mathbb{R}^d),\quad supp \hat{w}_k=\{\xi:2^k\leq |\xi|\leq (1+\epsilon)2^k\}
\end{align*}
and
\begin{align*}
\|w\|_{\dot{H}^{s_c}_x(\mathbb{R}^d)}\leq \epsilon^{-\alpha_0}
\end{align*}
for an absolute constant $\alpha_0>0$.
}

{\bf Theorem 9($\dot{H}^{s_c}\times \dot{H}^{s_c}$ scattering theory for the solution of (\ref{826x1}) with initial data has many bubbles).} {\it Assume that $\alpha\geq 0$, $\beta\geq 0$, $\alpha+\beta>2$ when $d=3$ and $\alpha+\beta>0$ when $d\geq 4$.
Then the solution $(u,v)$ of (\ref{826x1}) with initial data $(u_0,v_0)$ satisfying
\begin{align}
u_0=w_1+z_1,\quad v_0=w_2+z_2\label{1031x1}
\end{align}
is global existence in time and $(u,v)\in [C_t\dot{H}^{s_c}_x(\mathbb{R}\times \mathbb{R}^d)]^2$ if there exists some constant $\epsilon_0\in (0,1]$ such that $w_1$ and $w_2$ satisfy {\bf (C9)} with
\begin{align*}
w_1=\sum_{k=0}^{+\infty}w_{1k},\quad w_2=\sum_{k=0}^{+\infty}w_{2k}
\end{align*}
and
\begin{align}
\|z_1\|_{\dot{H}^{s_c}_x(\mathbb{R}^d)}+\|z_2\|_{\dot{H}^{s_c}_x(\mathbb{R}^d)}\leq \epsilon \label{10312}
\end{align}
for $\epsilon\in (0,\epsilon_0]$. Moreover, there exist scattering states $(u_{0+},v_{0+})$ and $(u_{0-},v_{0-})$ such that
\begin{align}
&\lim_{t\rightarrow + \infty} [\|u(t)-e^{it\Delta}u_{0+}\|_{\dot{H}^{s_c}_x(\mathbb{R}^d)}+\|v(t)-e^{it\Delta}v_{0+}\|_{\dot{H}^{s_c}_x(\mathbb{R}^d)}]=0,
\label{10313}\\
&\lim_{t\rightarrow - \infty} [\|u(t)-e^{it\Delta}u_{0-}\|_{\dot{H}^{s_c}_x(\mathbb{R}^d)}+\|v(t)-e^{it\Delta}v_{0-}\|_{\dot{H}^{s_c}_x(\mathbb{R}^d)}]=0.
\label{10313'}
\end{align}
}

We summarize our results in the following table

\begin{tabular}{llll}
  \hline
  dimension & conditions on $(\alpha,\beta)$ & initial data & scattering  \cr\\
 $d=3$ & $\alpha\geq 0$, $\beta\geq 0$, $\alpha+\beta<2$  & $H^1\times H^1$, $\Sigma\times \Sigma$  & $H^1\times H^1$, $\Sigma\times \Sigma$ \cr\\
 $d=3$ & $\alpha>0$ and $\beta>0$, $\alpha+\beta=2$  & $\Sigma\times \Sigma$ & $\Sigma\times \Sigma$\cr \\
 $d=3$ & $\alpha=0$ or $\beta=0$, $\alpha+\beta=2$  & $\dot{H}^1\times \dot{H}^1$ & $\dot{H}^1\times \dot{H}^1$\cr \\
 $d=3$ & $\alpha\geq 0$, $\beta\geq 0$, $\alpha+\beta>2$ & $\dot{H}^{s_c}\times \dot{H}^{s_c}$ & $\dot{H}^{s_c}\times \dot{H}^{s_c}$\cr \\
 $d=4$ & $(\alpha,\beta)=(0,0)$, & $\dot{H}^1\times \dot{H}^1$ & $\dot{H}^1\times \dot{H}^1$ \cr\\
 $d\geq 4$ & $\alpha\geq 0$, $\beta\geq 0$, $\alpha+\beta>0$ & $\dot{H}^{s_c}\times \dot{H}^{s_c}$ & $\dot{H}^{s_c}\times \dot{H}^{s_c}$\cr \\
  \hline
\end{tabular}\\

~~~~~~~~~~~~~~~~~~~~~~~~~~~~~~~~~~~~\\

The rest of this paper is organized as follows. In Section 2, we give some notations and useful lemmas. In Section 3, we establish the local wellposedness of
$H^1\times H^1$-solution, $\Sigma\times \Sigma$-solution and $H^s\times H^s$-solution of (\ref{1}) under different assumptions. In Section 4, we establish $H^1\times H^1$ and $\Sigma\times \Sigma$ scattering theories for the solution if $\alpha+\beta<2$(i.e., $(\alpha,\beta)$ is below the critical exponents line $\alpha+\beta=2$) when $d=3$, $\Sigma\times \Sigma$ scattering theory for the solution if $(\alpha,\beta)$ is on the critical exponents line $\alpha+\beta=2$ excluding the endpoint when $d=3$. In Section 5, we will be concerned with $\dot{H}^1\times \dot{H}^1$ scattering theory for the solution of (\ref{826x1}) containing defocusing nonlinearities with $(\alpha, \beta)$ is the endpoint of the critical exponents line $\alpha+\beta=2$ when $d=3$. In Section 6, we will get $\dot{H}^1\times \dot{H}^1$ scattering theory for the solution of (\ref{826x1}) containing defocusing nonlinearities with $(\alpha, \beta)=(0,0)$ when $d=4$. In Section 7, we discuss the $\dot{H}^{s_c}\times \dot{H}^{s_c}$ scattering theory for the solution of (\ref{826x1}) containing defocusing nonlinearities with $\alpha+\beta>2$ when $d=3$ and $\alpha+\beta>0$ when $d\geq 4$. In Section 8, we will deal with (\ref{826x1}) containing focusing nonlinearities and give some discussions.

\section{Preliminaries}
\qquad In this section, we will give some notations and useful lemmas.

\subsection{Weighted(or essential) gradient system of Schr\"{o}dinger equations}
\qquad First, we give the definition of ``weighted(or essential) gradient system of Schr\"{o}dinger equations" as follows.

{\bf Definition 2.1 (Weighted(or essential) gradient system of Schr\"{o}dinger equations).} {\it We say that
$$
iu_t+\Delta u=f(|u|^2,|v|^2)u,\quad iv_t+\Delta v=g(|u|^2,|v|^2)v
$$
is a weighted(or essential) gradient system of Schr\"{o}dinger equations if there exist $a\in\mathbb{R}$, $b\in\mathbb{R}$ and real value function $G(w,z)$ such that $\frac{\partial G}{\partial w}=af(w,z)$, $\frac{\partial G}{\partial z}=bg(w,z)$ for $w\geq 0$ and $z\geq 0$. $(a,b)$ is called as the weighted coefficients pair. Especially, if $(a,b)=(1,1)$, then the system is a gradient one.
}

By this definition, if $\lambda\neq 0$ and $\mu\neq 0$, then the system
$$
iu_t+\Delta u=\lambda |u|^{\alpha}|v|^{\beta+2}u,\quad iv_t+\Delta v=\mu|u|^{\alpha+2}|v|^{\beta}v
$$
is a weighted(or essential) gradient system of Schr\"{o}dinger equations. In fact, taking $a=\frac{\alpha+2}{2\lambda}$, $b=\frac{\beta+2}{2\mu}$ and $G(w,z)=w^{\frac{\alpha+2}{2}}z^{\frac{\beta+2}{2}}$ for $w\geq 0$ and $z\geq 0$, then $\frac{\partial G(w,z)}{\partial w}=a\lambda w^{\frac{\alpha}{2}}z^{\frac{\beta+2}{2}}$ and
$\frac{\partial G(w,z)}{\partial z}=b\mu w^{\frac{\alpha+2}{2}}z^{\frac{\beta}{2}}$ for $w\geq 0$ and $z\geq 0$.

But the following system isn't a weighted(or essential) gradient system of Schr\"{o}dinger equations:
$$
iu_t+\Delta u=\lambda |v|^pu,\quad iv_t+\Delta v=\mu|u|^qv
$$
if $(p,q)\neq (2,2)$.

In Section 3, we will show that, there exist the weighted conservation laws for the weighted(or essential) gradient system of Sch\"{o}dinger equations, which are similar to those for the scalar Sch\"{o}dinger equation.

\subsection{Some notations and lemmas about basic harmonic analysis}
\qquad We give some notations below. The notation $X\lesssim Y$ means that $X\leq CY$ for some constant $C$, while $X\thicksim Y$ implies that
 $Y\lesssim X\lesssim Y$. If $C$ depends up some additional parameters, we will indicate this with subscripts, for example, $X\lesssim_u Y$ means that
 $X\leq C_u Y$ for some constant $C_u$ depending on $u$. And we use $X\pm$ to denote any quantity of the form $X\pm \epsilon$ for any $\epsilon>0$.
In convenience, we will use $C$, $C'$, and so on, to denote some constants in the sequels, the values of it may vary line to line.

Let $L^r(\mathbb{R}^d)$ be the Banach space of functions $f:\mathbb{R}^d\rightarrow \mathbb{C}$ with the norm
$$
\|f\|_{L^r_x}:=\left(\int_{\mathbb{R}^d} |f(x)|^rdx\right)^{\frac{1}{r}}<+\infty,
$$
with the usual modifications when $r=\infty$. For any non-negative integer $k$, let $W^{k,r}(\mathbb{R}^d)$ be the Sobolev space where the norm of $f$ is defined as
$$
\|f\|_{W^{k,r}(\mathbb{R}^d)}:=\sum_{|\alpha|\leq k}\|\frac{\partial^{\alpha}}{\partial x^{\alpha}}f\|_{L^r_x}.
$$

Denote the spacetime norm $L^q_tL^r_x$ by
$$
\|u\|_{L^q_tL^r_x(\mathbb{R}\times \mathbb{R}^d)}:=\left(\int_{\mathbb{R}}\left(\int_{\mathbb{R}^d} |f(x)|^rdx\right)^{\frac{q}{r}}dt\right)^{\frac{1}{q}},
$$
with the usual modifications when $q$ or $r$ is infinity, or when $\mathbb{R}\times \mathbb{R}^d$ is replaced by some smaller region. Especially, we abbreviate $L^q_tL^r_x$ by $L^q_{t,x}$ if $q=r$.

If $(q,r)$ satisfies $\frac{2}{q}=d(\frac{1}{2}-\frac{1}{r})$ and $2\leq q, r\leq \infty$, we call it admissible pair. We define the $\dot{S}^0(I\times \mathbb{R}^d)$ and $\dot{S}^1(I\times \mathbb{R}^d)$ Strichartz norms by
\begin{align}
\|u\|_{\dot{S}^0(I\times \mathbb{R}^d)}:=sup \|u\|_{L^q_tL^r_x(I\times \mathbb{R}^d)},\quad \|u\|_{\dot{S}^1(I\times \mathbb{R}^d)}:=\|\nabla u\|_{\dot{S}^0(I\times \mathbb{R}^d)}
\end{align}
if $I\times \mathbb{R}^d$ is a spacetime slab, where the sup is taken over all admissible pairs $(q,r)$. While denote the dual space of $\dot{S}^0(I\times \mathbb{R}^d)$ by $\dot{N}^0(I\times \mathbb{R}^d)$ and
$$
\dot{N}^1(I\times \mathbb{R}^d):=\{u; \nabla u\in \dot{N}^0(I\times \mathbb{R}^d)\}.
$$

The Fourier transform on $\mathbb{R}^d$ and the fractional differential operators $|\nabla |^s$ are defined by
$$
\hat{f}(\xi):=\int_{\mathbb{R}^d}e^{-2\pi ix\cdot \xi}f(x)dx,\quad \widehat{|\nabla|^sf}(\xi)=|\xi|^s\hat{f}(\xi),
$$
which define the homogeneous Sobolev norms
$$
\|f\|_{\dot{H}^s_x}:=\||\nabla|^sf\|_{L^2_x}.
$$

Let $\varphi(\xi)$ be a cut-off function defined as
\begin{equation}
\varphi(\xi)=
\left\{
\begin{array}{lll}
1,\quad |\xi|\leq 1,\\
{\rm smooth\ bump\ supported},\quad 1\leq |\xi|\leq 2,\\
0,\quad |\xi|\geq 2.
\end{array}\right.
\end{equation}
We define the Littlewood-Paley operators
\begin{align*}
&\widehat{P_{\leq N}f}(\xi):=\varphi(\frac{\xi}{N})\hat{f}(\xi),\quad \widehat{P_{> N}f}(\xi):=(1-\varphi(\frac{\xi}{N}))\hat{f}(\xi),\\
&\widehat{P_N f}(\xi):=(\varphi(\frac{\xi}{N})-\varphi(\frac{2\xi}{N}))\hat{f}(\xi),
\end{align*}
similarly define $P_{<N}$, $P_{\geq N}$ and $P_{M<\cdot\leq N}:=P_{\leq N}-P_{\leq M}$, whenever $M$ and $N$ are dyadic numbers. Frequently, we will write $f_{\leq N}$ for $P_{\leq N}f$ and for other operators similarly.

Let $e^{it\Delta}$ be the free Schr\"{o}dinger propagator and the generated group of isometries $(\mathcal{J}(t))_{t\in \mathcal{R}}$. Then
\begin{align}
&\mathcal{J}(t)f(x):=e^{it\Delta}f(x)=\frac{1}{(4\pi it)^{\frac{d}{2}}}\int_{\mathbb{R}^d}e^{i|x-y|^2/4t}f(y)dy,\quad t\neq 0,\label{829x1} \\
&\widehat{e^{it\Delta}f}(\xi)=e^{-4\pi^2it|\xi|^2}\hat{f}(\xi).\label{829x3}
\end{align}
And
\begin{align}
& \|e^{it\Delta}f\|_{L^{\infty}_x(\mathbb{R}^d)}\lesssim |t|^{-\frac{d}{2}}\|f\|_{L^1_x(\mathbb{R}^d)},\quad   t\neq 0,\label{829x2}\\
&\|e^{it\Delta}f\|_{L^p_x(\mathbb{R}^d)}\lesssim |t|^{-d(\frac{1}{2}-\frac{1}{p}}\|f\|_{L^{p'}_x(\mathbb{R}^d)}, \quad   t\neq 0,\label{11101}
\end{align}
where $2\leq p\leq \infty$ and $\frac{1}{p}+\frac{1}{p'}=1$.

Note that Duhamel's formula
\begin{align}
u(t)=e^{i(t-t_0)\Delta}u(t_0)-i\int_{t_0}^te^{i(t-s)\Delta}(iu_t+\Delta u)(s)ds.\label{829x4}
\end{align}
We have the following Strichartz estimate.

{\bf Lemma 2.1(Strichartz estimate).} {\it Assume that $I$ is a compact time interval, and $u:I\times \mathbb{R}^d\rightarrow \mathbb{C}$ is a solution to the forced Schr\"{o}dinger equation
$$
iu_t+\Delta u=f.
$$
Then for given $s\geq 0$, we have
$$
\||\nabla|^su\|_{S^0(I)}\lesssim \||\nabla|^su(t_0)\|_{L^2_x}+\||\nabla|^sf\|_{N^0(I)}\quad {\rm for \ any}\quad t_0\in I.
$$
}

The following special Strichartz estimate for radial data was first proved in \cite{Shao20091,Shao20092} and developed in \cite{Cho2013, Guo2014}, it was presented as Lemma 2.9 in \cite{Beceanu??}.

{\bf Lemma 2.2(Radial Strichartz estimtates)} {\it
Assume that $f\in L^2(\mathbb{R}^d)$ is a radial function, and $q$, $r$, $\gamma$ satisfy
\begin{align}
\gamma\in \mathbb{R},\quad q\geq 2,\quad r>2,\quad \frac{2}{q}+\frac{2d-1}{r}<\frac{2d-1}{2},\quad \frac{2}{q}+\frac{d}{r}=\frac{d}{2}+\gamma.\label{1227x1}
\end{align}
Then
$$
\||\nabla|^{\gamma}e^{it\Delta}f\|_{L^q_tL^r_x(\mathbb{R}\times \mathbb{R}^d)}\lesssim \|f\|_{L^2_x(\mathbb{R}^d)}
$$
Moreover, if $G\in L^{\tilde{q}'}_tL^{\tilde{r}'}_x(\mathbb{R}\times \mathbb{R}^d)$ is a radial function in $x$, then
\begin{align}
\|\int_0^te^{i(t-s)\Delta}G(s)ds\|_{L^q_tL^r_x}+\||\nabla|^{-\gamma}\int_0^te^{i(t-s)\Delta}G(s)ds\|_{L^{\infty}_tL^2_x}\lesssim \|G\|_{L^{\tilde{q}'}_tL^{\tilde{r}'}_x}.\label{1227x2}
\end{align}
All the norms above are on $\mathbb{R}\times \mathbb{R}^d$. Here $(q,r,\gamma)$ and $(\tilde{q},\tilde{r},-\gamma)$ satisfy (\ref{1227x1}).

}

The following lemmas are about basic harmonic analysis.

{\bf Lemma 2.3(Bernstein estimates).} {\it
\begin{align*}
&\||\nabla|^{\pm s}P_Nf\|_{L^r_x(\mathbb{R}^d)}\thicksim N^{\pm s}\|P_Nf\|_{L^r_x(\mathbb{R}^d)},\quad
\|P_Nf\|_{L^q_x(\mathbb{R}^d)}\lesssim N^{\frac{d}{r}-\frac{d}{q}}\|P_Nf\|_{L^r_x(\mathbb{R}^d)},\\
&\|P_{\leq N}f\|_{L^q_x(\mathbb{R}^d)}\lesssim N^{\frac{d}{r}-\frac{d}{q}}\|P_{\leq N}f\|_{L^r_x(\mathbb{R}^d)}\quad {\rm
for}\quad  1\leq r\leq q\leq \infty.
\end{align*}
}

{\bf Lemma 2.4(Product rule,\cite{Christ1991}).} {\it Assume that $s\in (0,1]$, $1<r, r_1,r_2,q_1,q_2<+\infty$ and satisfy $\frac{1}{r}=\frac{1}{r_1}+\frac{1}{q_1}=\frac{1}{r_2}+\frac{1}{q_2}$. Then
\begin{align*}
\||\nabla|^s(fg)\|_{L^r_x}\lesssim \|f\|_{L^{r_1}_x}\||\nabla|^sg\|_{L^{q_1}_x}+\||\nabla|^sf\|_{L^{r_2}_x}\|g\|_{L^{q_2}_x}.
\end{align*}
}

{\bf Lemma 2.5(Fractional chain rule,\cite{Christ1991}).} {\it Assume that $s\in (0,1]$, $1<q,q_1,q_2<+\infty$ and satisfy $\frac{1}{q}=\frac{1}{q_1}+\frac{1}{q_2}$, $G\in C^1(\mathbb{C})$. Then
\begin{align*}
\||\nabla|^sG(u)\|_{L^q_x}\lesssim \|G'(u)\|_{L^{q_1}_x}\||\nabla|^su\|_{L^{q_2}_x}.
\end{align*}
}

{\bf Lemma 2.6(Fractional chain rule for a H\"{o}lder continuous function,\cite{Visan2007}).} {\it
Assume that $G$ is a H\"{o}lder continuous function of order $0<p<1$. Then
\begin{align*}
\||\nabla|^sG(u)\|_{L^q_x}\lesssim \||u|^{p-\frac{s}{\sigma}}\|_{L^{q_1}_x}\||\nabla|^{\sigma}u\|_{L^{\frac{sq_2}{\sigma}}_x}^{\frac{s}{\sigma}}
\end{align*}
for every $0<s<p$, $1<q<\infty$ and $\frac{s}{p}<\sigma<1$ satisfying $\frac{1}{q}=\frac{1}{q_1}+\frac{1}{q_2}$ and $(1-\frac{s}{p\sigma})q_1>1$.
}

{\bf Lemma 2.7(Another fractional chain rule for a H\"{o}lder continuous function,\cite{Killip20101}).} {\it
Assume that $G$ is a H\"{o}lder continuous function of order $0<p\leq 1$, $0<s<\sigma p<p$, $1<q,q_1,q_2,r_1,r_2,r_3<\infty$ and satisfy
\begin{align*}
(1-p)r_1>1,\quad (p-\frac{s}{\sigma})r_2>1,\\
\frac{1}{q}=\frac{1}{q_1}+\frac{1}{q_2}=\frac{1}{r_1}+\frac{1}{r_2}+\frac{1}{r_3}.
\end{align*}
 Then
\begin{align*}
&\quad\||\nabla|^s[w\cdot(G(u+v)-G(u))]\|_{L^q_x}\\
&\lesssim\||\nabla|^{\sigma}w\|_{L^{q_1}_x}\|v\|_{L^{pq_2}_x}^p
+\|w\|_{L^{r_1}}
\|v\|_{L^{(p-\frac{s}{\sigma})r_2}}^{p-\frac{s}{\sigma}}
\left(\||\nabla|^{\sigma}u\|_{L^{\frac{sr_3}{\sigma}}_x}+\||\nabla|^{\sigma}v\|_{L^{\frac{sr_3}{\sigma}}_x}\right)^{\frac{s}{\sigma}}.
\end{align*}
}

{\bf Lemma 2.8(Nonlinear Bernstein,\cite{Killip2013}).} {\it
Assume that $G:\mathbb{C}\rightarrow \mathbb{C}$ is a H\"{o}lder continuous function of order $0<p\leq 1$. Then
\begin{align*}
\|P_NG(u)\|_{L^{\frac{q}{p}}_x(\mathbb{R}^d)}\lesssim N^{-p}\|\nabla u\|^p_{L^q_x(\mathbb{R}^d)}\quad {\rm for \ any } \quad 1\leq q<\infty.
\end{align*}
}

\section{Local wellposedness results on (\ref{1})}

\subsection{Local wellposedness of $H^1\times H^1$-solution, $\Sigma\times \Sigma$-solution to (\ref{1})}
\qquad In this subsection, using Kato's method(see \cite{Kato1987}), we will prove Theorem 1 and establish the local well-posedness of $H^1\times H^1$-solution to (\ref{1}) if $(u_0,v_0)\in H^1(\mathbb{R}^d)\times H^1(\mathbb{R}^d)$, $\Sigma\times \Sigma$-solution to (\ref{1}) if $(u_0,v_0)\in \Sigma\times \Sigma$.

{\bf Proof of Theorem 1:}

{\bf Case 1.} Let $r=\max\{r_1,..., r_m, \rho'_1,..., \rho'_m\}$ and $(q,r)$ be the corresponding admissible pair.  Giving $M, T>0$ to be chosen later, we
define the following complete metric space
\begin{align}
S&=\left\{w=(w_1,w_2)\in \left[ L^{\infty}((-T,T), H^1(\mathbb{R}^d))\cap L^q((-T,T), W^{1,r}(\mathbb{R}^d))\right]^2, \right.\nonumber\\
&\qquad  \left. \sum_{j=1}^2\|w_j\|_{L^{\infty}((-T,T),H^1)}\leq M,\quad  \sum_{j=1}^2\|w_j\|_{L^q((-T,T),W^{1,r})}\leq M\right\}.\label{8251}
\end{align}
subject to the distance
\begin{align}
d(w,\tilde{w})=\sum_{j=1}^2[\|w_j-\tilde{w}_j\|_{L^q((-T,T),L^r)}+\|w_j-\tilde{w}_j\|_{L^{\infty}((-T,T),L^2)}].\label{8252}
\end{align}

{\bf Existence.}

For given $1\leq k\leq m$, let $(q_k,r_k)$ and $(\gamma_k,\rho_k)$ be the corresponding admissible pairs.
Since
\begin{align*}
\|w_j\|_{L^{q_k}((-T,T), W^{1,r_k})}\leq \|w_j\|_{L^{\infty}((-T,T),H^1)}^{\frac{2(r-r_k)}{r_k(r-2)}}\|w_j\|_{L^q((-T,T),W^{1,r})}^{\frac{r(r_k-2)}{r_k(r-2)}}
\end{align*}
for any $w_j\in W^{1,r_k}$, we know that if $(u,v)\in S$, then $u, v \in L^{q_k}((-T,T), W^{1,r_k}(\mathbb{R}^d))$ and
\begin{align}
\|u\|_{L^{q_k}((-T,T), W^{1,r_k})}\leq M^{\frac{2(r-r_k)}{r_k(r-2)}}M^{\frac{r(r_k-2)}{r_k(r-2)}}=M,\quad \|v\|_{L^{q_k}((-T,T), W^{1,r_k})}\leq M\label{8253}
\end{align}
for all $1\leq k\leq m$.

Note that $f_k\in C(\mathbb{R}\times \mathbb{R}; \mathbb{R})$ and $g_k\in C(\mathbb{R}\times \mathbb{R}; \mathbb{R})$. By the assumptions of (\ref{8241})--(\ref{8243}), if $(u,v)\in E$, then $f_k(|u|^2,|v|^2)u:(-T,T)\rightarrow L^{\rho'_k}(\mathbb{R}^d)$ and $g_k(|u|^2,|v|^2)v:(-T,T)\rightarrow L^{\rho'_k}(\mathbb{R}^d)$ are measurable and $f_k(|u|^2,|v|^2)u\in L^{\infty}((-T,T), L^{\rho'_k}(\mathbb{R}^d))$ and $g_k(|u|^2,|v|^2)v\in L^{\infty}((-T,T), L^{\rho'_k}(\mathbb{R}^d))$. Moreover,
\begin{align}
\|f_k(|u|^2,|v|^2)u\|_{L^{q_k}((-T,T), W^{1,\rho'_k})}&\leq C_M(T^{\frac{1}{q_k}}+\|u\|_{L^{q_k}((-T,T),W^{1,r_k})}+\|v\|_{L^{q_k}((-T,T),W^{1,r_k})})\nonumber\\
&\leq C_M(T^{\frac{1}{q_k}}+M),\label{8254}\\
\|f_k(|u|^2,|v|^2)u\|_{L^{\gamma'_k}((-T,T), W^{1,\rho'_k})}&\leq C_M(T^{\frac{1}{q_k}}+M)T^{\frac{q_k-\gamma'_k}{q_k\gamma'_k}},\label{8255}\\
\|g_k(|u|^2,|v|^2)v\|_{L^{q_k}((-T,T), W^{1,\rho'_k})}&\leq C_M(T^{\frac{1}{q_k}}+M),\label{8256}\\
\|g_k(|u|^2,|v|^2)v\|_{L^{\gamma'_k}((-T,T), W^{1,\rho'_k})}&\leq C_M(T^{\frac{1}{q_k}}+M)T^{\frac{q_k-\gamma'_k}{q_k\gamma'_k}}\label{8257}
\end{align}
for all $1\leq k\leq m$.

Given $(u_0,v_0)\in H^1(\mathbb{R}^d)\times H^1(\mathbb{R}^d)$ and $(u,v)\in S$, let $\mathcal{H}(u)$ and $\mathcal{H}(v)$ be defined by
\begin{align}
\mathcal{H}(u)(t)&=\mathcal{J}(t)u_0-i\int_0^t\mathcal{J}(t-s)f(|u|^2,|v|^2)u(s)ds,\label{8258}\\
\mathcal{H}(v)(t)&=\mathcal{J}(t)v_0-i\int_0^t\mathcal{J}(t-s)g(|u|^2,|v|^2)v(s)ds.\label{8259}
\end{align}
Using Strichartz inequalities and (\ref{8254})--(\ref{8257}), we know that, if $T\leq 1$, then
$\mathcal{H}(u), \mathcal{H}(v)\in L^q((-T,T), W^{1,r}(\mathbb{R}^d))\cap C([-T,T], H^1(\mathbb{R}^d))$, and
\begin{align}
\|\mathcal{H}(u)\|_{L^q((-T,T),W^{1,r})}+\|\mathcal{H}(u)\|_{L^{\infty}((-T,T), H^1)}&\leq C\|u_0\|_{H^1_x}+CC_M(1+M)T^{\delta},\label{825w1}\\
\|\mathcal{H}(v)\|_{L^q((-T,T),W^{1,r})}+\|\mathcal{H}(v)\|_{L^{\infty}((-T,T), H^1)}&\leq C\|v_0\|_{H^1_x}+CC_M(1+M)T^{\delta},\label{825w2}
\end{align}
where
$$\delta=\min_{1\leq k\leq m} \frac{q_k-\gamma'_k}{q_k\gamma'_k}>0.$$

Taking $M>2C(\|u_0\|_{H^1_x}+\|v_0\|_{H^1_x})\in S$ and $T$ such that  $2CC_M(1+M)T^{\delta}\leq M$, if $(u,v)\in S$, then
${\bf{H}(u,v)}:=(\mathcal{H}(u),\mathcal{H}(v))\in S$. Moreover, if $T$ small enough,
\begin{align}
d({\bf{H}(u_1,v_1)},{\bf{H}(u_2,v_2)})\leq \frac{1}{2} d((u_1,v_1),(u_2,v_2)) \quad {\rm for \ all}\ (u_1,v_1),(u_2,v_2)\in S.\label{825w3}
\end{align}
Therefore, $\bf{H}$ at least has a fixed point $(u,v)$ which satisfies (\ref{8258}) and (\ref{8259}).

{\bf Uniqueness.}

Assume that $(u_1,v_1)$ and $(u_2,v_2)$ are two solutions of (\ref{1}). Let $w_1=u_1-u_2$ and $w_2=v_1-v_2$. Then
\begin{align}
w_1(t)&=i\int_0^t\mathcal{J}(t-s)\sum_{k=1}^m[f_k(|u_1|^2,|v_1|^2)u_1(s)-f_k(|u_2|^2,|v_2|^2)u_2(s)]ds,\label{825w4}\\
w_2(t)&=i\int_0^t\mathcal{J}(t-s)\sum_{k=1}^m[g_k(|u_1|^2,|v_1|^2)v_1(s)-g_k(|u_2|^2,|v_2|^2)v_2(s)]ds.\label{825w5}
\end{align}
Using Strichartz estimates, we have
\begin{align*}
\sum_{k=1}^m\|w_1(t)\|_{L^{q_k}(J,L^{r_k})}&\leq \sum_{k=1}^m\|[f_k(|u_1|^2,|v_1|^2)u_1-f_k(|u_2|^2,|v_2|^2)u_2]\|_{L^{\gamma'_k}(J,L^{\rho'_k})},\\
\sum_{k=1}^m\|w_2(t)\|_{L^{q_k}(J,L^{r_k})}&\leq \sum_{k=1}^m\|[g_k(|u_1|^2,|v_1|^2)v_1-g_k(|u_2|^2,|v_2|^2)v_2]\|_{L^{\gamma'_k}(J,L^{\rho'_k})}
\end{align*}
for every interval $J$ such that $0\in J\subset [-T,T]$. Recalling (\ref{8242}) and (\ref{8243}), we get
\begin{align}
&\quad \sum_{k=1}^m[\|w_1(t)\|_{L^{q_k}(J,L^{r_k})}+\|w_2(t)\|_{L^{q_k}(J,L^{r_k})}]\nonumber\\
&\leq C\sum_{k=1}^m[\|w_1(t)\|_{L^{\gamma'_k}(J,L^{r_k})}+\|w_2(t)\|_{L^{\gamma'_k}(J,L^{r_k})}].\label{825w6}
\end{align}
Applying the results of Lemma 4.2.2 and Lemma 4.2.4 in \cite{Cazenave2003} to (\ref{825w6}), we obtain
 $$\|w_1(t)\|_{L^{r_k}}+\|w_2(t)\|_{L^{r_k}}=0.$$
 That is,
$(u_1,v_1)=(u_2,v_2)$.

{\bf Maximality and blowup alternative.}

Consider $u_0\in H^1(\mathbb{R}^d)$, $v_0\in H^1(\mathbb{R}^d)$ and let
\begin{align*}
T_{\max}(u_0,v_0)&=\sup \{T>0: {\rm there \ exists \ a \ solution \ of} \ (\ref{1})\ {\rm on} \ [0,T],\\
T_{\min}(u_0,v_0)&=\sup \{T>0: {\rm there \ exists \ a \ solution \ of} \ (\ref{1})\ {\rm on} \ [-T,0].
\end{align*}
Note that the unique solution $(u,v)$ of (\ref{1}) satisfies
$$
u, v\in C((-T_{\min},T_{\max}), H^1(\mathbb{R}^d))\cap C^1((-T_{\min},T_{\max}),H^{-1}(\mathbb{R}^d)).
$$
Assume that $T_{\max}<+\infty$ and there exist a sequence $t_j\uparrow T_{\max}$ such that $\|u(t_j)\|_{H^1_x}+\|v(t_j)\|_{H^1_x}\leq M$ for some $M<+\infty$. Let $l$ be
such that $t_l+T(M)>T_{\max}(u_0,v_0)$. Then we can look $(u(t_l),v(t_l))$ as the initial data, and extend $(u,v)$ to $t_l+T(M)$, which is a contradiction to maximality. Hence
$$\|u(t)\|_{H^1_x}+\|v(t)\|_{H^1_x}\rightarrow +\infty\quad {\rm as}\quad t\uparrow T_{\max}.$$
Similarly, if $T_{\min}(u_0,v_0)<+\infty$, then
$$\|u(t)\|_{H^1_x}+\|v(t)\|_{H^1_x}\rightarrow +\infty\quad {\rm as}\quad t\downarrow -T_{\min}.$$

{\bf Continuous dependence.}

Suppose that $(u_{0n},v_{0n})\rightarrow (u_0,v_0)$ in $H^1(\mathbb{R}^d)\times H^1(\mathbb{R}^d)$. Since $\|u_{0n}\|_{H^1_x}+\|v_{0n}\|_{H^1_x}\leq 2(\|u_0\|_{H^1_x}+\|v_0\|_{H^1_x})$ for $n$ large enough, there exist $T=T(\|u_0\|_{H^1_x},\|v_0\|_{H^1_x})$ and $n_0$ such that
$(u,v)$ and $(u_n,v_n)$ are defined on $[-T,T]$ for $n\geq n_0$ and
\begin{align*}
&\quad \|u\|_{L^{\infty}((-T,T), H^1)}+\|v\|_{L^{\infty}((-T,T), H^1)}+sup_{n\geq n_0}[\|u_n\|_{L^{\infty}((-T,T), H^1)}+\|v_n\|_{L^{\infty}((-T,T), H^1)}]\\
&\leq C[\|u_0\|_{H^1_x}+\|v_0\|_{H^1_x}].
\end{align*}
Note that
\begin{align*}
&u_n(t)-u(t)=\mathcal{J}(t)(u_{0n}-u_0)+\mathcal{H}(u_n)(t)-\mathcal{H}(u)(t),\\
&v_n(t)-v(t)=\mathcal{J}(t)(v_{0n}-v_0)+\mathcal{H}(v_n)(t)-\mathcal{H}(v)(t).
\end{align*}
We can use Strichartz estimates, Gagliardo-Nirenberg's inequality, and a covering argument as Step 3 in the proof of Theorem 4.4.1 in \cite{Cazenave2003} to establish the continuous dependence result. We omit the details here.

Let $[0,t]$ be the time interval when the solution exists.

{\bf Mass conservation law for each component of $(u,v)$}.

Multiplying the first equation and the second one of (\ref{1}) by $2\bar{u}$ and $2\bar{v}$ respectively, taking the imaginary parts of the results, then integrating them over $\mathbb{R}^d\times [0,t]$ respectively, we have
\begin{align*}
\int_{\mathbb{R}^d}|u|^2dx=\int_{\mathbb{R}^d}|u_0|^2dx,\quad \int_{\mathbb{R}^d}|v|^2dx=\int_{\mathbb{R}^d}|v_0|^2dx,
\end{align*}
which implies mass conservation law for each component of $(u,v)$.

As a byproduct, for any $c_1,c_2\in \mathbb{R}$, we have
\begin{align*}
\int_{\mathbb{R}^d}[c_1|u|^2+c_2|v|^2]dx=\int_{\mathbb{R}^d}[c_1|u_0|^2+c_2|v_0|^2]dx.
\end{align*}

{\bf Momentum and energy conservations laws.}

Multiplying the first equation and the second one of (\ref{1}) by $c_1\nabla \bar{u}$ and $c_2\nabla\bar{v}$ respectively, taking the real parts of the results, then integrating them over $\mathbb{R}^d\times [0,t]$, and sum them up, we obtain
\begin{align*}
\int_{\mathbb{R}^d}[c_1u(t,x)\nabla \bar{u}(t,x)+c_2v(t,x)\nabla \bar{v}(t,x)]dx=\int_{\mathbb{R}^d}[c_1\bar{u}_0\nabla u_0+c_2v_0\nabla \bar{v}_0]dx,
\end{align*}
which implies momentum conservation law.

Multiplying the first equation and the second one of (\ref{1}) by $2c_1\bar{u}_t$ and $2c_2\bar{v}_t$ respectively, taking the real part of the result, then integrating them over $\mathbb{R}^d\times [0,t]$, and sum them up, we obtain
\begin{align*}
\int_{\mathbb{R}^d}[c_1|\nabla u|^2+c_2|\nabla v|^2+G(|u|^2,|v|^2)]dx=\int_{\mathbb{R}^d}[c_1|\nabla u_0|^2+c_2|\nabla v_0|^2+G(|u_0|^2,|v_0|^2)]dx,
\end{align*}
which implies energy conservation law. Here $G=\sum_{k=1}^mG_k$.

{\bf Case 2.} It is a critical case in $H^1(\mathbb{R}^3)\times H^1(\mathbb{R}^3)$. Given $\Lambda\in \mathbb{N}$, let
\begin{equation}
\label{8261} f_{\Lambda}(|u|^2,|v|^2)u=\left\{
\begin{array}{lll}
\lambda|u|^{\alpha}|v|^{\beta+2}u\quad & {\rm if} \quad |u|\leq \Lambda, \ |v|\leq \Lambda,\\
\lambda \Lambda^{\alpha}|v|^{\beta+2}u \quad & {\rm if} \quad |u|\geq \Lambda, \ |v|\leq \Lambda,\\
\lambda \Lambda^{\beta+2}|u|^{\alpha}u \quad & {\rm if} \quad |u|\leq \Lambda, \ |v|\geq \Lambda,\\
\lambda \Lambda^4u \quad & {\rm if} \quad |u|\geq \Lambda, \ |v|\geq \Lambda,
\end{array}\right.
\end{equation}
\begin{equation}
\label{8262} g_{\Lambda}(|u|^2,|v|^2)v=\left\{
\begin{array}{lll}
\mu |u|^{\alpha+2}|v|^{\beta}v\quad & {\rm if} \quad |u|\leq \Lambda, \ |v|\leq \Lambda,\\
\mu \Lambda^{\alpha+2}|v|^{\beta}v \quad & {\rm if} \quad |u|\geq \Lambda, \ |v|\leq \Lambda,\\
\mu \Lambda^{\beta}|u|^{\alpha+2}v \quad & {\rm if} \quad |u|\leq \Lambda, \ |v|\geq \Lambda,\\
\mu \Lambda^4v \quad & {\rm if} \quad |u|\geq \Lambda, \ |v|\geq \Lambda.
\end{array}\right.
\end{equation}
Consider the following truncated problem:
\begin{align}
\mathcal{H}_{\Lambda}(u)(t)&=\mathcal{J}(t)u_0-i\int_0^t\mathcal{J}(t-s)f_{\Lambda}(|u|^2,|v|^2)u(s)ds,\label{8263}\\
\mathcal{H}_{\Lambda}(v)(t)&=\mathcal{J}(t)v_0-i\int_0^t\mathcal{J}(t-s)g_{\Lambda}(|u|^2,|v|^2)v(s)ds.\label{8264}
\end{align}
Similar to Step 2 in the proof of Theorem 4.5.1. in \cite{Cazenave2003}, we can prove that there exists a unique, global solution $(u_{\Lambda},v_{\Lambda})$ satisfying (\ref{8263}) and (\ref{8264}). Similar to Step 3 there, passing to to the limit by letting $\Lambda\rightarrow \infty$, we can prove that (\ref{1}) has a unique strong solution $(u,v)\in H^1(\mathbb{R}^3)\times H^1(\mathbb{R}^3)$. Next, similar to Step 4, Step 5 and Step 6 there, we can prove the conservation laws, the blowup alternative, continuous dependence respectively.

{\bf Case 3.} It is a critical case in $H^1(\mathbb{R}^4)\times H^1(\mathbb{R}^4)$. Given $\Lambda\in \mathbb{N}$, let
\begin{equation}
\label{8265} f_{\Lambda}(u,v)=\left\{
\begin{array}{lll}
\lambda |v|^2u\quad & {\rm if} \quad  |v|\leq \Lambda,\\
\lambda \Lambda^2u \quad & {\rm if} \quad  |v|\geq \Lambda,
\end{array}\right.\quad
 g_{\Lambda}(u,v)=\left\{
\begin{array}{lll}
\mu |u|^2v\quad & {\rm if} \quad  |u|\leq \Lambda,\\
\mu \Lambda^2v \quad & {\rm if} \quad  |u|\geq \Lambda.
\end{array}\right.
\end{equation}
Consider the following truncated problem:
\begin{align}
\mathcal{H}_{\Lambda}(u)(t)&=\mathcal{J}(t)u_0-i\int_0^t\mathcal{J}(t-s)f_{\Lambda}(u,v)(s)ds,\label{8266}\\
\mathcal{H}_{\Lambda}(v)(t)&=\mathcal{J}(t)v_0-i\int_0^t\mathcal{J}(t-s)g_{\Lambda}(u,v)v(s)ds.\label{8267}
\end{align}
Similar to the proof of Theorem 4.5.1. in \cite{Cazenave2003}, we also can prove the existence and uniqueness, the conservation laws, the blowup alternative, continuous dependence.

{\bf The main steps of proving $(xu, xv)\in L^2(\mathbb{R}^d)\times L^2(\mathbb{R}^d)$.}

 If $(u_0,v_0)\in H^1(\mathbb{R}^d)\times H^1(\mathbb{R}^d)$ and $(xu_0, xv_0)\in L^2(\mathbb{R}^d)\times L^2(\mathbb{R}^d)$, in order to prove $(xu, xv)\in L^2(\mathbb{R}^d)\times L^2(\mathbb{R}^d)$, we can follow the idea of the proof of Proposition 6.5.1 in \cite{Cazenave2003}. Step 1, we use auxiliary function $\|e^{-\epsilon |x|^2}|x|u(t)\|^2_{L^2_x}+\|e^{-\epsilon |x|^2}|x|v(t)\|^2_{L^2_x}$ to get
the vector-valued function $t\rightarrow (|\cdot|u(t,\cdot), |\cdot|v(t,\cdot))$ is continuous $(-T_{\min},T_{\max})\rightarrow L^2(\mathbb{R}^d)\times L^2(\mathbb{R}^d)$. Step 2, we can prove that if the sequence of initial data $(u_{0n},v_{0n})\rightarrow (u_0,v_0)$ in $H^1(\mathbb{R}^d)\times H^1(\mathbb{R}^2)$ and $(xu_{0n},xv_{0n})\rightarrow (xu_0,xv_0)$ in $L^2(\mathbb{R}^d)\times L^2(\mathbb{R}^2)$, then the corresponding sequence of solutions to (\ref{1}) $(u_n,v_n)\rightarrow (u,v)$ in $L^{\infty}_tH^1_x((-T_{\min},T_{\max})\times\mathbb{R}^d)\times L^{\infty}_tH^1_x((-T_{\min},T_{\max})\times\mathbb{R}^d)$ and $(xu_n,xv_n)\rightarrow (xu,xv)$ in $L^{\infty}_tL^2_x((-T_{\min},T_{\max})\times\mathbb{R}^d)\times L^{\infty}_tL^2_x((-T_{\min},T_{\max})\times\mathbb{R}^d)$ as $n\rightarrow +\infty$. Step 3, in order to prove (\ref{1113w1})--(\ref{1113w3}), choosing $(u_{0n},v_{0n})\in H^2(\mathbb{R}^d)\times H^2(\mathbb{R}^2)$, using the $H^2\times H^2$ regularity, we can prove that (\ref{1113w1})--(\ref{1113w3})
holds for the corresponding solutions $(u_n,v_n)$, then let $n\rightarrow +\infty$ and get the conclusions.

By the way, we denote the space $\Sigma\times \Sigma$ by
\begin{align}
\Sigma\times \Sigma=\{(u,v)|(u,v)\in H^1(\mathbb{R}^d)\times H^1(\mathbb{R}^d),\quad (xu,xv)\in L^2(\mathbb{R}^d)\times L^2(\mathbb{R}^d)\}.\label{1113w4}
\end{align}

Using the equations in (\ref{1}), integrating by parts, after some elementary computations and recalling (\ref{1113w1}), we can obtain (\ref{1113w2}) and (\ref{1113w3}).\hfill$\Box$

{\bf Remark 3.1.} By the proof of Theorem 1, we find that the local existence, uniqueness, continuous dependence, blowup alternative and the mass of each component conservation law hold for general system of Sch\"{o}dinger equations if the nonlinearities satisfy the assumptions of this theorem. However, we only prove that the conservation laws of weighted mass, weighted momentum and weighted energy are true for the weighted(or essential) gradient system of Sch\"{o}dinger equations.

\subsection{Existence of critical exponents line for (\ref{826x1}) when $d=3$ and critical exponents point for (\ref{826x1}) when $d=4$}
\qquad In this subsection, we will prove Theorem 2 and show the existence of critical exponents line for (\ref{826x1}) when $d=3$ and critical exponents point for (\ref{826x1}) when $d=4$.

{\bf Proof of Theorem 2:} 1. If $\alpha+\beta\leq 2$, taking $\rho=r=\alpha+\beta+2$, we can verify the assumptions of Theorem 1 when $d=3$, and prove that for every $(u_0,v_0)\in H^1(\mathbb{R}^3)\times H^1(\mathbb{R}^3)$, (\ref{826x1}) has a unique strong $H^1\times H^1$-solution. Moreover, if $\lambda>0$, $\mu>0$, $\alpha\geq 0$ and $\beta\geq 0$, since
\begin{align}
E_w(u,v)=\int_{\mathbb{R}^3}[\frac{\alpha+2}{2\lambda}|\nabla u|^2+\frac{\beta+2}{2\mu}|\nabla v|^2+|u|^{\alpha+2}|v|^{\beta+2}]dx,\label{1141}
\end{align}
by the conservation laws of mass and energy, obviously, we have
\begin{align*}
\frac{1}{2}\int_{\mathbb{R}^3}[\frac{\alpha+2}{\lambda}|\nabla u|^2+\frac{\beta+2}{\mu}|\nabla v|^2]dx\leq E_w(u_0,v_0)<+\infty,
\end{align*}
which implies that $(u,v)$ is global existence and uniformly bounded in the norm $L^{\infty}_tH^1_x(\mathbb{R}\times \mathbb{R}^3)\times L^{\infty}_tH^1_x(\mathbb{R}\times \mathbb{R}^3)$ if $\lambda>0$, $\mu>0$, $\alpha\geq 0$, $\beta\geq 0$ and $\alpha+\beta\leq 2$ when $d=3$.

If $\alpha+\beta>2$ and $\alpha=\beta$, (\ref{826x1}) meets with the case of $\sigma=2(\alpha+1)$ in \cite{Alazard20091, Christ20032}. By their results, the Cauchy problem of $iu_t+\Delta u=|u|^{2\sigma}u$ is not well posed and there exist a sequence of initial data such that $\|u(t)\|_{H^1_x}\rightarrow \infty$ in finite time.

Consequently, if $\lambda>0$, $\mu>0$, $\alpha\geq 0$ and $\beta\geq 0$ when $d=3$, the critical exponents line $\alpha+\beta=2$ is the watershed for whether (\ref{826x1}) always has a unique bounded $H^1\times H^1$-solution or not.

2. The proof is similar to the discussions above, we omit the details. \hfill$\Box$

{\bf Remark 3.2.} If $d=3$, $u_0(x)\equiv 0$ or $v_0(x)\equiv 0$, then (\ref{826x1}) always has a semi-trivial solution $(\tilde{u}(t,x),0)$ or $(0, \tilde{v}(t,x))$. Here
$\tilde{u}(t,x)$ and $\tilde{v}(t,x)$ are the solutions of
\begin{equation}
\label{826w4}
\left\{
\begin{array}{lll}
iu_t+\Delta u=0,\quad  x\in \mathbb{R}^3,\ t\in \mathbb{R},\\
u(x,0)=u_0(x),\quad x\in \mathbb{R}^3.
\end{array}\right.
\end{equation}
and
\begin{equation}
\label{826w5}
\left\{
\begin{array}{lll}
iv_t+\Delta v=0,\quad  x\in \mathbb{R}^3,\ t\in \mathbb{R},\\
v(x,0)=v_0(x),\quad x\in \mathbb{R}^3.
\end{array}\right.
\end{equation}
Therefore, a very interesting open question is: What assumptions on $(\alpha,\beta), (\lambda,\mu)\in \mathbb{R}^2$ and $(u_0,v_0)\in H^1(\mathbb{R}^3)\times H^1(\mathbb{R}^3)$ can guarantee (\ref{826x1}) having a unique bounded solution $(u(t),v(t))\in H^1(\mathbb{R}^3)\times H^1(\mathbb{R}^3)$ if $\alpha+\beta>2$ and $d=3$? Another opposite conjecture is: (\ref{826x1}) doesn't have the bounded $H^1\times H^1$-solution for any $\lambda,\mu\in\mathbb{R}$, $\lambda\neq 0$, $\mu\neq 0$ and any initial data $(u_0,v_0)\in H^1(\mathbb{R}^3)\times H^1(\mathbb{R}^3)$ with $u_0\not\equiv 0$, $v_0\not\equiv 0$ if $\alpha+\beta>2$ and $d=3$.

(\ref{826x1}) with $\alpha+\beta>0$ when $d=4$ also has similar open question and opposite conjecture.

Although we cannot solve the open problems above, in the next subsection, we will show that there exists $H^s\times H^s$-solution to (\ref{826x1}) with initial data $(u_0,v_0)\in H^s(\mathbb{R}^d)\times H^s(\mathbb{R}^d)$, $d\geq 3$.

\subsection{Local wellposedness of the $H^s\times H^s$-solution to (\ref{826x1})}
\qquad In this subsection,  we will prove Theorem 3 and establish the local well-posedness of $H^s\times H^s$-solution to (\ref{1}) if $(u_0,v_0)\in H^s(\mathbb{R}^d)\times H^s(\mathbb{R}^d)$

{\bf Proof of Theorem 3:} Similar to the proof of Theorem 1.1 in \cite{Cazenave1990}, we proceed in several steps.

{\bf Step 1.} First we prove a fact below: Let $(u_1,v_1), (u_2, v_2)\in [L^{\gamma}((-T,T); B^s_{\rho,2})]^2$ and $\delta=1-\frac{\alpha+\beta+4}{\gamma}$.
Then for any admissible pair $(q,r)$,
\begin{align}
&\quad \|\mathcal{J}(|u_1|^{\alpha}|v_1|^{\beta+2}u_1)-\mathcal{J}(|u_2|^{\alpha}|v_2|^{\beta+2}u_2)\|_{L^q((-T,T);L^r)}\nonumber\\
&\quad +\|\mathcal{J}(|u_1|^{\alpha+2}|v_1|^{\beta}v_1)-\mathcal{J}(|u_2|^{\alpha+2}|v_2|^{\beta}v_2)\|_{L^q((-T,T);L^r)}\nonumber\\
&\leq CT^{\delta}\{\|u_1\|^{\alpha+\beta+2}_{L^{\gamma}((-T,T);\dot{B}^s_{\rho,2})}+\|u_2\|^{\alpha+\beta+2}_{L^{\gamma}((-T,T);\dot{B}^s_{\rho,2})}
+\|v_1\|^{\alpha+\beta+2}_{L^{\gamma}((-T,T);\dot{B}^s_{\rho,2})}
+\|v_2\|^{\alpha+\beta+2}_{L^{\gamma}((-T,T);\dot{B}^s_{\rho,2})}\}\nonumber\\
&\qquad\quad \times \{\|u_1-u_2\|_{L^{\gamma}((-T,T);L^{\rho})}+\|v_1-v_2\|_{L^{\gamma}((-T,T);L^{\rho})}\},\label{827x5}
\end{align}
and
\begin{align}
&\quad \|\mathcal{J}(|u_1|^{\alpha}|v_1|^{\beta+2}u_1)\|_{L^q((-T,T);\dot{B}^s_{r,2})}+\|\mathcal{J}(|u_1|^{\alpha+2}|v_1|^{\beta}v_1)\|_{L^q((-T,T);\dot{B}^s_{r,2})}
\nonumber\\
&\leq CT^{\delta}\{\|u_1\|^{\alpha+\beta+3}_{L^{\gamma}((-T,T);\dot{B}^s_{\rho,2})}+\|v_1\|^{\alpha+\beta+3}_{L^{\gamma}((-T,T);\dot{B}^s_{\rho,2})}
\},\label{827x6}\\
&\quad \|\mathcal{J}(|u_2|^{\alpha}|v_2|^{\beta+2}u_2)\|_{L^q((-T,T);\dot{B}^s_{r,2})}+\|\mathcal{J}(|u_2|^{\alpha+2}|v_2|^{\beta}v_2)\|_{L^q((-T,T);\dot{B}^s_{r,2})}\nonumber\\
&\leq CT^{\delta}\{\|u_2\|^{\alpha+\beta+3}_{L^{\gamma}((-T,T);\dot{B}^s_{\rho,2})}
+\|v_2\|^{\alpha+\beta+3}_{L^{\gamma}((-T,T);\dot{B}^s_{\rho,2})}\}.\label{827x7}
\end{align}

Indeed,
\begin{align*}
&\quad \|\mathcal{J}(|u_1|^{\alpha}|v_1|^{\beta+2}u_1)-\mathcal{J}(|u_2|^{\alpha}|v_2|^{\beta+2}u_2)\|_{L^q((-T,T);L^r)}\nonumber\\
&\quad +\|\mathcal{J}(|u_1|^{\alpha+2}|v_1|^{\beta}v_1)-\mathcal{J}(|u_2|^{\alpha+2}|v_2|^{\beta}v_2)\|_{L^q((-T,T);L^r)}\nonumber\\
&\leq  \|(|u_1|^{\alpha}|v_1|^{\beta+2}u_1-|u_2|^{\alpha}|v_2|^{\beta+2}u_2)\|_{L^{\gamma'}((-T,T);L^{\rho'})}\nonumber\\
&\quad +\|(|u_1|^{\alpha+2}|v_1|^{\beta}v_1-|u_2|^{\alpha+2}|v_2|^{\beta}v_2)\|_{L^{\gamma'}((-T,T);L^{\rho'})}\nonumber\\
&\leq C\|[|u_1|^{\alpha+\beta+2}+|u_2|^{\alpha+\beta+2}+|v_1|^{\alpha+\beta+2}+|v_2|^{\alpha+\beta+2}][|u_1-u_2|+|v_1-v_2|]\|_{L^{\gamma'}((-T,T);L^{\rho'})}\nonumber\\
&\leq C\left(\int_{-T}^T \|[|u_1|+|u_2|+|v_1|+|v_2|]\|^{(\alpha+\beta+2)\gamma'}_{\dot{B}^s_{\rho,2}}\|[|u_1-u_2|+|v_1-v_2|]\|^{\gamma'}_{L^{\rho}}dt\right)\nonumber\\
&\leq CT^{\delta}\{\|u_1\|^{\alpha+\beta+2}_{L^{\gamma}((-T,T);\dot{B}^s_{\rho,2})}+\|u_2\|^{\alpha+\beta+2}_{L^{\gamma}((-T,T);\dot{B}^s_{\rho,2})}
+\|v_1\|^{\alpha+\beta+2}_{L^{\gamma}((-T,T);\dot{B}^s_{\rho,2})}
+\|v_2\|^{\alpha+\beta+2}_{L^{\gamma}((-T,T);\dot{B}^s_{\rho,2})}\}\nonumber\\
&\qquad\quad \times \{\|u_1-u_2\|_{L^{\gamma}((-T,T);L^{\rho})}+\|v_1-v_2\|_{L^{\gamma}((-T,T);L^{\rho})}\},
\end{align*}
and
\begin{align*}
&\quad \|\mathcal{J}(|u_1|^{\alpha}|v_1|^{\beta+2}u_1)\|_{L^q((-T,T);\dot{B}^s_{r,2})}+\|\mathcal{J}(|u_1|^{\alpha+2}|v_1|^{\beta}v_1)\|_{L^q((-T,T);\dot{B}^s_{r,2})}
\nonumber\\
&\leq \|\mathcal{J}(|u_1|^{\alpha}|v_1|^{\beta+2}u_1)\|_{L^{\gamma'}((-T,T);\dot{B}^s_{\rho',2})}
+\|\mathcal{J}(|u_1|^{\alpha+2}|v_1|^{\beta}v_1)\|_{L^{\gamma'}((-T,T);\dot{B}^s_{\rho',2})}
\nonumber\\
&\leq C\left(\int_{-T}^T\|u_1\|^{(\alpha+\beta+3)\gamma'}_{\dot{B}^s_{\rho',2}}\right)^{ \frac{1}{\gamma'}}
+C\left(\int_{-T}^T\|v_1\|^{(\alpha+\beta+3)\gamma'}_{\dot{B}^s_{\rho',2}}\right)^{ \frac{1}{\gamma'}}\nonumber\\
&\leq CT^{\delta}[\|u_1\|^{(\alpha+\beta+3)}_{L^{\gamma}((0,T);\dot{B}^s_{\rho,2})}+\|v_1\|^{(\alpha+\beta+3)}_{L^{\gamma}((0,T);\dot{B}^s_{\rho,2})}].
\end{align*}

{\bf Step 2.} We show another fact below: Let $(u_0,v_0)\in H^s\times H^s$ and $(u,v)\in [L^{\gamma}((-T,T);B^s_{\rho,2})]^2$ be a solution of (\ref{826x1}).
Then $(u,v)\in [L^q((0,T);B^s_{r,2})\cap C([-T,T]; H^s)]^2$ for every admissible pair $(q,r)$. If $(w,z)\in [L^{\gamma}((-T,T);B^s_{\rho,2})]^2$ is also a solution of (\ref{826x1}), then $(u,v)=(w,z)$.

Indeed, we can write $u(t)=\mathcal{J}(\cdot)u_0-i\mathcal{J}(|u|^{\alpha}|v|^{\beta+2}u)$ and $v(t)=\mathcal{J}(\cdot)v_0-i\mathcal{J}(|u|^{\alpha+2}|v|^{\beta}v)$, then using Strichartz estimates, we can prove the first statement.

About the uniqueness, we assume that $(u(t),v(t))\neq (w(t),z(t))$ for some $t\in [-T,T]$. If
$t_0=\inf\{t\in [-T,T], (u(t),v(t))\neq (w(t),z(t))\}\geq 0$,
then $(U(t), V(t)):=(u(t+t_0),v(t+t_0))$ and $(W(t),Z(t)):=(w(t+t_0),z(t+t_0))$ both satisfy
$u(t)=\mathcal{J}(\cdot)\phi-i\mathcal{J}(|u|^{\alpha}|v|^{\beta+2}u)$ and $v(t)=\mathcal{J}(\cdot)\psi-i\mathcal{J}(|u|^{\alpha+2}|v|^{\beta}v)$ on $[-T, T-t_0]$, where $(\phi,\psi):=(u(t_0),v(t_0))=(w(t_0),z(t_0))$. Taking the admissible pairs $(q,r)=(\gamma,\rho)$, for all $t\in [t_0,T]$, we have
\begin{align*}
&\quad \|u-w\|_{L^{\gamma}((t_0,t); L^{\rho})}+\|v-z\|_{L^{\gamma}((t_0,t); L^{\rho})}\nonumber\\
&=\|\mathcal{J}(|u|^{\alpha}|v|^{\beta+2}u)-\mathcal{J}(|w|^{\alpha}|z|^{\beta+2}w)\|_{L^{\gamma}((t_0,t); L^{\rho})}
+\|\mathcal{J}(|u|^{\alpha+2}|v|^{\beta}v)-\mathcal{J}(|w|^{\alpha+2}|z|^{\beta}z)\|_{L^{\gamma}((t_0,t); L^{\rho})}\nonumber\\
&\leq C(t-t_0)^{\delta}\left\{\|u\|^{\alpha+\beta+2}_{L^{\gamma}((t_0,t);B^2_{\rho,2})}+\|v\|^{\alpha+\beta+2}_{L^{\gamma}((t_0,t);B^2_{\rho,2})}
+\|w\|^{\alpha+\beta+2}_{L^{\gamma}((t_0,t);B^2_{\rho,2})}+\|z\|^{\alpha+\beta+2}_{L^{\gamma}((t_0,t);B^2_{\rho,2})}\right\}\nonumber\\
&\qquad \qquad \qquad\times\left\{\|u-w\|_{L^{\gamma}((t_0,t); L^{\rho})}+\|v-z\|_{L^{\gamma}((t_0,t); L^{\rho})}\right\}.
\end{align*}
For $t>t_0$ but close to $t_0$ enough such that
$$
C(t-t_0)^{\delta}\left\{\|u\|^{\alpha+\beta+2}_{L^{\gamma}((t_0,t);B^2_{\rho,2})}+\|v\|^{\alpha+\beta+2}_{L^{\gamma}((t_0,t);B^2_{\rho,2})}
+\|w\|^{\alpha+\beta+2}_{L^{\gamma}((t_0,t);B^2_{\rho,2})}+\|z\|^{\alpha+\beta+2}_{L^{\gamma}((t_0,t);B^2_{\rho,2})}\right\}<1.
$$
Then
$$
\|u-w\|_{L^{\gamma}((t_0,t); L^{\rho})}+\|v-z\|_{L^{\gamma}((t_0,t); L^{\rho})}=0,
$$
which is a contradiction to the definition of $t_0$, so $(u(t),v(t))\equiv(w(t),z(t))$ for all $t\in [-T,T]$.

Similarly, if $t_0=\sup\{t\in [-T,T], (u(t),v(t))\neq (w(t),z(t))\}\leq 0$, we can obtain the parallel conclusions for $t\in [-T, t_0]$ and prove that $(u(t),v(t))\equiv(w(t),z(t))$ for all $t\in [-T,T]$.

{\bf Step 3.} Existence of the solution to (\ref{826x1}). Let $M>0$ be finite and
$$
\mathcal{X}=\mathcal{X}(T,M)=\{(u, v)\in [L^{\gamma}((-T,T); B^s_{\rho,2})]^2:\|u\|_{L^{\gamma}((-T,T); \dot{B}^s_{\rho,2})}+\|v\|_{L^{\gamma}((-T,T); \dot{B}^s_{\rho,2})}\leq M\}.
$$
By the results of Strichartz estimates, this space is never empty because $(\mathcal{J}(t)u_0, \mathcal{J}(t)v_0)$ is in $\mathcal{X}(\infty,M)$ if $(u_0,v_0)\in H^s\times H^s$ and $\|u_0\|_{\dot{H}^s_x}+\|v_0\|_{\dot{H}^s_x}$ is small enough. Endowed with the metric
$$
d((u,v),(w,z))=\|u-w\|_{L^{\gamma}((-T,T);L^{\rho})}+\|v-z\|_{L^{\gamma}((-T,T);L^{\rho})},
$$
$\mathcal{X}$ is a complete metric space.

Consider
\begin{align}
\mathcal{F}_1u=\mathcal{J}(\cdot)u_0-i\mathcal{J}(|u|^{\alpha}|v|^{\beta+2}u),\ \mathcal{F}_2v=\mathcal{J}(\cdot)v_0-i\mathcal{J}(|u|^{\alpha+2}|v|^{\beta}v),
\end{align}
and denote ${\bf F}(u,v)=(\mathcal{F}_1u,\mathcal{F}_2v)$. We will seek for the conditions on $T$ and $M$ to let ${\bf F}$ be a strict contraction on $\mathcal{X}$.

By Strichartz estimates and the results in Step 1, we know that if $(u,v)\in \mathcal{X}$, then $(\mathcal{F}_1u,\mathcal{F}_2v)\in [L^{\gamma}((-T,T); \dot{B}^s_{\rho,2})]^2$. Moreover, if
\begin{align}
\|\mathcal{J}(\cdot)u_0\|_{L^{\gamma}((-T,T); \dot{B}^s_{\rho,2})}+\|\mathcal{J}(\cdot)v_0\|_{L^{\gamma}((-T,T); \dot{B}^s_{\rho,2})}+CT^{\delta}M^{\alpha+\beta+3}\leq M,\label{8281}
\end{align}
then ${\bf F}(u,v)=(\mathcal{F}_1u,\mathcal{F}_2v)\in \mathcal{X}$. Furthermore, by (\ref{827x5}), if
\begin{align}
CT^{\delta}M^{\alpha+\beta+2}<1,\label{8282}
\end{align}
then ${\bf F}$ is a strict contraction on $\mathcal{X}$. Therefore, if we take $M>\|u_0\|_{H^s_x}+\|v_0\|_{H^s_x}$ and $T$ is small enough such that (\ref{8281}) and (\ref{8282}) hold, then ${\bf F}$ is a strict contraction on $\mathcal{X}$ and at least has a fixed point, which is the solution of (\ref{826x1}).

{\bf Step 4.} By the results of Step 1 to Step 3, we have proved that there exist an unique solution $(u,v)$ of (\ref{826x1}) satisfying (i) and (ii). Now the properties (iii)--(vii) can be obtained by following the same standard steps as in the proof of Theorem 1.1 in \cite{Cazenave1990},  we just need to replace $u$ by $(u,v)$ and the norm $\|\phi\|_{\dot{H}^s_x}$ by $\|u_0\|_{\dot{H}^s_x}+\|v_0\|_{\dot{H}^s_x}$, and so on. We omit the details here.\hfill $\Box$

\section{$H^1\times H^1$  and $\Sigma\times \Sigma$ scattering theories  for (\ref{826x1}) with defocusing nonlinearities when $d=3$}

\qquad In this section, we consider the global solution of (\ref{826x1}) with defocusing nonlinearities($\lambda>0$ and $\mu>0$) when $d=3$.

\subsection{$H^1\times H^1$  and $\Sigma\times \Sigma$ scattering theories for (\ref{826x1}) with defocusing nonlinearities and $\alpha+\beta<2$ when $d=3$ }
\qquad In this subsection, we will establish $H^1\times H^1$ and $\Sigma\times \Sigma$ scattering theories for (\ref{826x1}) with defocusing nonlinearities, $\alpha\geq 0$, $\beta\geq 0$ and $\alpha+\beta<2$ when $d=3$. That is, the exponents pair $(\alpha,\beta)$ is below the critical exponents line $\alpha+\beta=2$.

The main steps of the proof Theorem 4 are as follows:

Step 1. We give the weight-coupled interaction Morawetz estimates and obtain
\begin{align}
\int_I\int_{\mathbb{R}^3}[|u(t,x)|^4+|v(t,x)|^4]dxdt\leq C.\label{830w5}
\end{align}
Here $C$ is a uniform constant independent of the time interval $I$.

Step 2. We will prove that
\begin{align}
\|u(t,x)\|_{S^1(\mathbb{R}\times \mathbb{R}^3)}+\|v(t,x)\|_{S^1(\mathbb{R}\times \mathbb{R}^3)}<\infty.\label{830w6}
\end{align}

Step 3. We use Duhamel formulae and Strichartz estimates to establish $H^1\times H^1$ and $\Sigma\times \Sigma$ scattering theories for (\ref{826x1}) under the assumptions of Theorem 4.

{\bf Step 1. Weight-coupled interaction Morawetz estimate.}

Let $(u,v)$ be the solution of (\ref{826x1}).  For a given smooth real function $a(x,y)$, we define the following weight-coupled interaction Morawetz potential:
\begin{align}
M^{\otimes_2}_a(t)&=2A\int_{\mathbb{R}^d}\int_{\mathbb{R}^d}\widetilde{\nabla} a(x,y)\Im[\bar{u}(t,x)\bar{u}(t,y)\widetilde{\nabla}(u(t,x)u(t,y))]dxdy\nonumber\\
&\quad+2B\int_{\mathbb{R}^d}\int_{\mathbb{R}^d}\widetilde{\nabla} a(x,y)\Im[\bar{v}(t,x)\bar{v}(t,y)\widetilde{\nabla}(v(t,x)v(t,y))]dxdy\nonumber\\
&\quad+2C\int_{\mathbb{R}^d}\int_{\mathbb{R}^d}\widetilde{\nabla} a(x,y)\Im[\bar{u}(t,x)\bar{v}(t,y)\widetilde{\nabla}(u(t,x)v(t,y))]dxdy\nonumber\\
&\quad+2D\int_{\mathbb{R}^d}\int_{\mathbb{R}^d}\widetilde{\nabla} a(x,y)\Im[\bar{v}(t,x)\bar{u}(t,y)\widetilde{\nabla}(v(t,x)u(t,y))]dxdy,\label{829w1}
\end{align}
where $\widetilde{\nabla}=(\nabla_x, \nabla_y)$, $x\in \mathbb{R}^d$ and $y\in \mathbb{R}^d$, and
\begin{align}
&A=\mu^2(\alpha+2)^2,\quad B=\lambda^2(\beta+2)^2,\quad C=D=\lambda\mu(\alpha+2)(\beta+2),\label{1112w1}\\
&L_1=L_2=2(\alpha+\beta)\lambda \mu^2(\alpha+2),\quad L_3=L_4=2(\alpha+\beta)\mu\lambda^2(\beta+2).\label{11241}
\end{align}
 Then
\begin{align}
&\qquad \frac{d}{dt}M^{\otimes_2}_a(t)\nonumber\\
&=\int_{\mathbb{R}^d}\int_{\mathbb{R}^d}\left(-[A|u|^2+D|v|^2]\Delta_x\Delta_xa+4a_{jk}Re(Au_j\bar{u}_k+Dv_j\bar{v}_k)
+L_1|u|^{\alpha+2}|v|^{\beta+2}\Delta_xa \right)dx|u(t,y)|^2dy\nonumber\\
&+\int_{\mathbb{R}^d}\int_{\mathbb{R}^d}\left(-[C|u|^2+B|v|^2]\Delta_x\Delta_xa+4a_{jk}Re(Cu_j\bar{u}_k+Bv_j\bar{v}_k)
+L_2|u|^{\alpha+2}|v|^{\beta+2}\Delta_xa \right)dx|v(t,y)|^2dy\nonumber\displaybreak\\
&+\int_{\mathbb{R}^d}\int_{\mathbb{R}^d}\left(-[A|u|^2+C|v|^2]\Delta_y\Delta_ya+4a_{jk}Re(Au_j\bar{u}_k+Cv_j\bar{v}_k)
+L_3|u|^{\alpha+2}|v|^{\beta+2}\Delta_ya \right)dy|u(t,x)|^2dx\nonumber\\
&+\int_{\mathbb{R}^d}\int_{\mathbb{R}^d}\left(-[D|u|^2+B|v|^2]\Delta_y\Delta_ya+4a_{jk}Re(Du_j\bar{u}_k+Bv_j\bar{v}_k)
+L_4|u|^{\alpha+2}|v|^{\beta+2}\Delta_ya \right)dy|v(t,x)|^2dx\nonumber\\
&-4\int_{\mathbb{R}^d}\int_{\mathbb{R}^d}\nabla_xa\cdot\Im[\bar{u}(t,x)\nabla_xu(t,x)]dx
[\nabla_y\cdot(\Im[A\bar{u}(t,y)\nabla_yu(t,y)+C\bar{v}(t,y)\nabla_yv(t,y)])]dy\nonumber\\
&-4\int_{\mathbb{R}^d}\int_{\mathbb{R}^d}\nabla_xa\cdot\Im[\bar{v}(t,x)\nabla_xv(t,x)]dx
[\nabla_y\cdot(\Im[B\bar{v}(t,y)\nabla_yv(t,y)+D\bar{u}(t,y)\nabla_yu(t,y)])]dy\nonumber\\
&-4\int_{\mathbb{R}^d}\int_{\mathbb{R}^d}\nabla_ya\cdot\Im[\bar{u}(t,y)\nabla_yu(t,y)]dy
[\nabla_x\cdot(\Im[A\bar{u}(t,x)\nabla_xu(t,x)+D\bar{v}(t,x)\nabla_xv(t,x)])])]dx\nonumber\\
&-4\int_{\mathbb{R}^d}\int_{\mathbb{R}^d}\nabla_ya\cdot\Im[\bar{v}(t,y)\nabla_yv(t,y)]dy
[\nabla_x\cdot(\Im[B\bar{v}(t,x)\nabla_xv(t,x)+C\bar{u}(t,x)\nabla_xu(t,x)])]dx.\label{829w2}
\end{align}

If $d=3$ and $a(x,y)=|x-y|$, denoting
$$
\overline{\nabla}_y:=\nabla_x-\frac{x-y}{|x-y|}(\frac{x-y}{|x-y|}\cdot \nabla_x),\quad
\overline{\nabla}_x:=\nabla_y-\frac{y-x}{|x-y|}(\frac{y-x}{|x-y|}\cdot \nabla_y),$$
 then we have
\begin{align}
&\quad\frac{d}{dt}M^{\otimes_2}_a(t)\nonumber\\
&=16\pi\int_{\mathbb{R}^3}[A|u(t,x)|^4+(C+D)|u|^2|v|^2+B|v(t,x)|^4]dx\nonumber\\
&\quad +2\int_{\mathbb{R}^3}\int_{\mathbb{R}^3}\frac{|u(t,x)|^{\alpha+2}|v(t,x)|^{\beta+2}[(L_1+L_3)|u(t,y)|^2+(L_2+L_4)|v(t,y)|^2]}{|x-y|}dxdy\nonumber\\
&\quad +4\int_{\mathbb{R}^3}\int_{\mathbb{R}^3}\frac{
[A|\overline{\nabla}_yu(t,x)|^2+D|\overline{\nabla}_yv(t,x)|^2]|u(t,y)|^2}{|x-y|}dxdy\nonumber\\
&\quad +4\int_{\mathbb{R}^3}\int_{\mathbb{R}^3}\frac{
[A|\overline{\nabla}_xu(t,y)|^2+C|\overline{\nabla}_xv(t,y)|^2]|u(t,x)|^2}{|x-y|}dxdy\nonumber\\
&\quad +4\int_{\mathbb{R}^3}\int_{\mathbb{R}^3}\frac{
[B|\overline{\nabla}_yv(t,x)|^2+C|\overline{\nabla}_yu(t,x)|^2]|v(t,y)|^2}{|x-y|}dxdy\nonumber\\
&\quad +4\int_{\mathbb{R}^3}\int_{\mathbb{R}^3}\frac{
[B|\overline{\nabla}_xv(t,y)|^2+D|\overline{\nabla}_xu(t,y)|^2]|v(t,x)|^2}{|x-y|}dxdy\nonumber\\
&\quad-4\int_{\mathbb{R}^3}\int_{\mathbb{R}^3}\frac{x-y}{|x-y|}\cdot\Im[\bar{u}(t,x)\nabla_xu(t,x)]dx
[\nabla_y\cdot(\Im[A\bar{u}(t,y)\nabla_yu(t,y)+C\bar{v}(t,y)\nabla_yv(t,y)])]dy\nonumber\\
&\quad-4\int_{\mathbb{R}^3}\int_{\mathbb{R}^3}\frac{x-y}{|x-y|}\cdot\Im[\bar{v}(t,x)\nabla_xv(t,x)]dx
[\nabla_y\cdot(\Im[B\bar{v}(t,y)\nabla_yv(t,y)+D\bar{u}(t,y)\nabla_yu(t,y)])]dy\nonumber\\
&\quad-4\int_{\mathbb{R}^3}\int_{\mathbb{R}^3}\frac{y-x}{|x-y|}\cdot\Im[\bar{u}(t,y)\nabla_yu(t,y)]dy
[\nabla_x\cdot(\Im[A\bar{u}(t,x)\nabla_xu(t,x)+D\bar{v}(t,x)\nabla_xv(t,x)])])]dx\nonumber\\
&\quad-4\int_{\mathbb{R}^3}\int_{\mathbb{R}^3}\frac{y-x}{|x-y|}\cdot\Im[\bar{v}(t,y)\nabla_yv(t,y)]dy
[\nabla_x\cdot(\Im[B\bar{v}(t,x)\nabla_xv(t,x)+C\bar{u}(t,x)\nabla_xu(t,x)])]dx\nonumber\\
&:=(I)+(II)+4(III)+4(IV)+4(V)+4(VI)-4(VII)-4(VIII)-4(X)-4(XI).\label{830x1}
\end{align}
We will prove that
\begin{align}
4(III)+4(IV)+4(V)+4(VI)-4(VII)-4(VIII)-4(X)-4(XI)\geq 0.\label{830w1}
\end{align}
We only estimate the term
$$
-4\int_{\mathbb{R}^3}\int_{\mathbb{R}^3}\frac{x-y}{|x-y|}\cdot\Im[\bar{v}(t,x)\nabla_xv(t,x)]dx
[\nabla_y\cdot(\Im[D\bar{u}(t,y)\nabla_yu(t,y)])]dy,
$$
the other terms can be estimated similarly.

Denote $p_1(x)=2\Im [\bar{u}(t,x)\nabla_x u(t,x)]$, $p_1(y)=2\Im [\bar{u}(t,y)\nabla_y u(t,y)]$, $p_2(x)=2\Im [\bar{v}(t,x)\nabla_x v(t,x)]$ and $p_2(y)=2\Im [\bar{v}(t,y)\nabla_y v(t,y)]$. Since
\begin{align*}
&\partial_{x_k}\left(\frac{x_j-y_j}{|x-y|}\right)=\frac{\delta_{jk}}{|x-y|}-\frac{(x_k-y_k)(x_j-y_j)}{|x-y|^3},\quad j,k=1,2,3,\\
&\partial_{y_k}\left(\frac{x_j-y_j}{|x-y|}\right)=-\frac{\delta_{jk}}{|x-y|}+\frac{(x_k-y_k)(x_j-y_j)}{|x-y|^3},\quad j,k=1,2,3,
\end{align*}
integrating by parts, we have
\begin{align}
&\quad-4\int_{\mathbb{R}^3}\int_{\mathbb{R}^3}\frac{x-y}{|x-y|}\cdot\Im[\bar{v}(t,x)\nabla_xv(t,x)]dx
[\nabla_y\cdot(\Im[D\bar{u}(t,y)\nabla_yu(t,y)])]dy\nonumber\\
&=D\int_{\mathbb{R}^3}\int_{\mathbb{R}^3}[p_1(y)p_2(x)-(p_1(y)\frac{x-y}{|x-y|})(p_2(x)\frac{x-y}{|x-y|})]\frac{dxdy}{|x-y|}.\label{830w2}
\end{align}
Noticing that
\begin{align*}
&|\pi_{(x-y)^{\perp}}p_1(y)|=|p_1(y)-\frac{x-y}{|x-y|}(\frac{x-y}{|x-y|}p_1(y))|\leq 2|\Im[\bar{u}(t,y)\overline{\nabla}_xu(t,y)]|\leq 2|u(t,y)||\overline{\nabla}_xu(t,y)|,\nonumber\\
&|\pi_{(y-x)^{\perp}}p_2(x)|=|p_2(x)-\frac{x-y}{|x-y|}(\frac{x-y}{|x-y|}p_2(x))|\leq 2|\Im[\bar{v}(t,x)\overline{\nabla}_yv(t,x)]|\leq 2|v(t,x)||\overline{\nabla}_yv(t,x)|,
\end{align*}
(\ref{830w2}) implies that
\begin{align*}
&\quad-4\int_{\mathbb{R}^3}\int_{\mathbb{R}^3}\frac{x-y}{|x-y|}\cdot\Im[\bar{v}(t,x)\nabla_xv(t,x)]dx
[\nabla_y\cdot(\Im[D\bar{u}(t,y)\nabla_yu(t,y)])]dy\nonumber\\
&\geq -4D\int_{\mathbb{R}^3}\int_{\mathbb{R}^3}|u(t,y)||\overline{\nabla}_xu(t,y)||v(t,x)||\overline{\nabla}_yv(t,x)|\frac{dxdy}{|x-y|}\nonumber\\
&\geq -2D\int_{\mathbb{R}^3}\int_{\mathbb{R}^3}\frac{|u(t,y)|^2|\overline{\nabla}_yv(t,x)|^2}{|x-y|}dxdy
-2D\int_{\mathbb{R}^3}\int_{\mathbb{R}^3}\frac{|v(t,x)|^2|\overline{\nabla}_xu(t,y)|^2}{|x-y|}dxdy.
\end{align*}
Similarly,  we also can obtain the results on -4(VI), -4(VII) and -4(VIII) as above. Summing them up, we get (\ref{830w1}).

(\ref{830x1}) and (\ref{830w1}) mean that
\begin{align}
&\int_I\int_{\mathbb{R}^3}[|u(t,x)|^4+|v(t,x)|^4]dxdt\lesssim |M^{\otimes_2}_a(t)|\lesssim \|u\|^4_{L^{\infty}_tH^1_x(I\times \mathbb{R}^3)}+\|v\|^4_{L^{\infty}_tH^1_x(I\times \mathbb{R}^3)},\label{830w3}\\
&\int_I\int_{\mathbb{R}^3}\int_{\mathbb{R}^3}\frac{|u(t,x)|^{\alpha+2}|v(t,x)|^{\beta+2}[|u(t,y)|^2+|v(t,y)|^2]}{|x-y|}dxdy\nonumber\\
&\lesssim \|u\|^4_{L^{\infty}_tH^1_x(I\times \mathbb{R}^3)}+\|v\|^4_{L^{\infty}_tH^1_x(I\times \mathbb{R}^3)}.\label{830w4}
\end{align}

{\bf Step 2. Estimate for $\|u(t,x)\|_{S^1(\mathbb{R}\times \mathbb{R}^3)}+\|v(t,x)\|_{S^1(\mathbb{R}\times \mathbb{R}^3)}$.}

Since $\int_{\mathbb{R}}\int_{\mathbb{R}^3}[|u(t,x)|^4+|v(t,x)|^4]dxdt<+\infty$, for any $0<\delta<1$, we can subdivide $(-\infty, +\infty)$ into $J=J(E,M,\delta)$ subintervals $J_k$ such that
\begin{align}
\int_{J_k}\int_{\mathbb{R}^3}[|u(t,x)|^4+|v(t,x)|^4]dxdt<\delta.\label{8311}
\end{align}
Note that
\begin{align}
u(t)=\mathcal{J}(t)u_0-i\lambda\int_0^t\mathcal{J}(t-s)[|u|^{\alpha}|v|^{\beta+2}u](s)ds,\label{8312}\\
v(t)=\mathcal{J}(t)v_0-i\mu\int_0^t\mathcal{J}(t-s)[|u|^{\alpha+2}|u|^{\beta}v](s)ds.\label{8313}
\end{align}
Using Strichartz estimates, for any admissible pair $(q,r)$, we have
\begin{align*}
\|u(t)\|_{L^q((-\infty, \infty);W^{1,r}(\mathbb{R}^3))}&\lesssim \|u_0\|_{H^1}
+\sum_{k=1}^J\|[|u|^{\alpha}|v|^{\beta+2}u]\|_{L^2(J_k;W^{1,\frac{6}{5}})},\nonumber\\
\|v(t)\|_{L^q((-\infty, \infty);W^{1,r}(\mathbb{R}^3))}&\lesssim \|v_0\|_{H^1}
+\sum_{k=1}^J\|[|u|^{\alpha+2}|v|^{\beta}v]\|_{L^2(J_k;W^{1,\frac{6}{5}})}.
\end{align*}
Taking over all admissible pair $(q,r)$, we get
\begin{align}
&\qquad \|u(t)\|_{S^1((-\infty, +\infty);\mathbb{R}^3)}+\|v(t)\|_{S^1((-\infty, +\infty);\mathbb{R}^3)}\lesssim \|u_0\|_{H^1}+\|v_0\|_{H^1}
\nonumber\\
&\quad +\sum_{k=1}^J\left(\|[|u|^{\alpha}|v|^{\beta+2}u]\|_{L^2(J_k;W^{1,\frac{6}{5}})}+\|[|u|^{\alpha+2}|v|^{\beta}v]\|_{L^2(J_k;W^{1,\frac{6}{5}})}\right)\nonumber\\
&\lesssim \sum_{k=1}^J\left\{\int_{J_k}\left[\int_{\mathbb{R}^3}(|u|^{\alpha+\beta+2}+|v|^{\alpha+\beta+2})^{\frac{6}{5}}
(|u|+|\nabla u|+|v|+|\nabla v|)^{\frac{6}{5}}dx\right]^{\frac{5}{3}}dt\right\}^{\frac{1}{2}}+\|u_0\|_{H^1}+\|v_0\|_{H^1}\nonumber\\
&\lesssim \sum_{k=1}^J\left\{\int_{J_k}\left[\left(\int_{\mathbb{R}^3}(|u|^2+|v|^2)dx\right)^{\frac{1}{\tau_5}}\left(\int_{\mathbb{R}^3}(|u|^4+|v|^4)dx\right)^{\frac{1}{\tau_6}}
\left(\int_{\mathbb{R}^3}(|u|^6+|v|^6)dx\right)^{\frac{1}{\tau_7}}\right]^{\frac{5\tau_4}{3\tau_2}}dt\right\}^{\frac{1}{2\tau_4}}\nonumber\\
&\quad \times\left\{\int_{J_k}\left(\int_{\mathbb{R}^3}
(|u|+|\nabla u|+|v|+|\nabla v|)^{\frac{6\tau_1}{5}}dx\right)^{\frac{5\tau_3}{3\tau_1}}dt\right\}^{\frac{1}{2\tau_3}}+\|u_0\|_{H^1}+\|v_0\|_{H^1}\nonumber\\
&\lesssim \delta^{\frac{5-\tau_1}{4\tau_1}}\left(\|u(t)\|_{S^1((-\infty, +\infty);\mathbb{R}^3)}+\|v(t)\|_{S^1((-\infty, +\infty);\mathbb{R}^3)}\right)+\|u_0\|_{H^1}+\|v_0\|_{H^1},\label{911}
\end{align}
which implies that
\begin{align}
\|u(t)\|_{S^1((-\infty, +\infty);\mathbb{R}^3)}+\|v(t)\|_{S^1((-\infty, +\infty);\mathbb{R}^3)}<+\infty. \label{912}
\end{align}
Here
\begin{align*}
&\max\left(\frac{25}{13},\frac{5}{5-2(\alpha+\beta)}, \frac{15}{7-2(\alpha+\beta)}\right)<\tau_1<5, \quad  \tau_2=\frac{\tau_1}{\tau_1-1},\quad  \tau_3=\frac{2\tau_1}{3\tau_1-5},\quad
\tau_4=\frac{2\tau_1}{5-\tau_1},\\
&\tau_5=\frac{20(\tau_1-1)}{[21-6(\alpha+\beta)]\tau_1-45},\quad \tau_6=\frac{10(\tau_1-1)}{3(5-\tau_1)},\quad \tau_7=\frac{20(\tau_1-1)}{[5+6(\alpha+\beta)]\tau_1-5}.
\end{align*}

{\bf Step 3. $H^1\times H^1$ and $\Sigma\times \Sigma$ scattering theories for (\ref{826x1}).}

For $0<t<+\infty$, we define
\begin{align}
u_{+}(t)&=u_0-i\lambda\int_0^t\mathcal{J}(-s)[|u|^{\alpha}|v|^{\beta+2}u](s)ds,\label{913}\\
v_{+}(t)&=v_0-i\mu\int_0^t\mathcal{J}(-s)[|u|^{\alpha+2}|u|^{\beta}v](s)ds.\label{914}
\end{align}
Then for any $0<\tau<t$, we have
\begin{align*}
\|u_{+}(t)-u_{+}(\tau)\|_{H^1_x}&\lesssim \|\int_{\tau}^t\mathcal{J}(t-s)[|u|^{\alpha}|v|^{\beta+2}u](s)ds\|_{L^{\infty}_t((\tau,t);H^1_x(\mathbb{R}^3))},\\
\|v_{+}(t)-v_{+}(\tau)\|_{H^1_x}&\lesssim \|\int_{\tau}^t\mathcal{J}(t-s)[|u|^{\alpha+2}|v|^{\beta}v](s)ds\|_{L^{\infty}_t((\tau,t);H^1_x(\mathbb{R}^3))}.
\end{align*}
Using Strichartz estimates, and similar to (\ref{911}), we obtain
\begin{align*}
&\quad \|u_{+}(t)-u_{+}(\tau)\|_{H^1_x}+\|v_{+}(t)-v_{+}(\tau)\|_{H^1_x}\nonumber\\
&\lesssim \|[|u|^{\alpha}|v|^{\beta+2}u]\|_{L^2((\tau,t);W^{1,\frac{6}{5}})}+\|[|u|^{\alpha+2}|v|^{\beta}v]\|_{L^2((\tau,t);W^{1,\frac{6}{5}})}\nonumber\\
&\leq \|u(t)\|_{S^1((\tau, t);\mathbb{R}^3)}+\|v(t)\|_{S^1((\tau, t);\mathbb{R}^3)}.
\end{align*}
Recalling (\ref{912}), we know that for any $\epsilon>0$, there exists $T=T(\epsilon)>0$ such that
$$
\|u_{+}(t)-u_{+}(\tau)\|_{H^1_x}+\|v_{+}(t)-v_{+}(\tau)\|_{H^1_x}\leq \epsilon \quad {\rm for \ any}\quad \tau,\ t>T(\epsilon),
$$
which means that both $u_+(t)$ and $v_+(t)$ respectively converge in $H^1_x$ as $t\rightarrow +\infty$ to $u_+$ and $v_+$ as follows:
$$
u_{+}:=u_0-i\lambda\int_0^{+\infty}\mathcal{J}(-s)[|u|^{\alpha}|v|^{\beta+2}u](s)ds,\quad
v_{+}:=v_0-i\mu\int_0^{+\infty}\mathcal{J}(-s)[|u|^{\alpha+2}|u|^{\beta}v](s)ds.
$$
Moreover, by Strichartz estimates, (\ref{911}) and (\ref{912}), we can get
\begin{align*}
&\quad \|e^{-it\Delta}u(t)-u_{+}\|_{H^1_x}+\|e^{-it\Delta}v(t)-v_{+}\|_{H^1_x}\nonumber\\
&=\|\lambda\int_t^{+\infty}\mathcal{J}(-s)[|u|^{\alpha}|v|^{\beta+2}u](s)ds\|_{H^1_x}+
\|\mu\int_t^{+\infty}\mathcal{J}(-s)[|u|^{\alpha+2}|u|^{\beta}v](s)ds\|_{H^1_x}\nonumber\\
&\lesssim \|u(t)\|_{S^1((t,+\infty);\mathbb{R}^3)}+\|v(t)\|_{S^1((t,+\infty);\mathbb{R}^3)}\rightarrow 0 \quad {\rm as}\ t\rightarrow +\infty.
\end{align*}

Similarly, we can prove that there exist $u_{-}$ and $v_{-}$ such that
\begin{align*}
\|e^{-it\Delta}u(t)-u_{-}\|_{H^1_x}+\|e^{-it\Delta}v(t)-v_{-}\|_{H^1_x}\rightarrow 0 \quad {\rm as}\ t\rightarrow -\infty.
\end{align*}

(\ref{1023w1}) and (\ref{1023w2}) are obtained.

Now we will prove (\ref{11201}) and (\ref{11202}). Letting $w(t)=(x+2it\nabla)u$, $z(t)=(x+2it\nabla)v$, we get the following Duhamel formulae
\begin{align}
w(t)=\mathcal{J}(t)(xu_0)-i\lambda\int_0^t\mathcal{J}(t-s)\tilde{f}_1(u,v,w,z,)(s)ds,\label{11203}\\
z(t)=\mathcal{J}(t)(xv_0)-i\mu\int_0^t\mathcal{J}(t-s)\tilde{f}_2(u,v,w,,z,)(s)ds.\label{11204}
\end{align}
Here
\begin{align}
&\tilde{f}_1(u,v,w,z)=(\frac{\alpha}{2}+1)|u|^{\alpha}|v|^{\beta+2}w-\frac{\alpha u}{2\bar{u}}|u|^{\alpha}|v|^{\beta+2}\bar{w}+(\frac{\beta}{2}+1)|u|^{\alpha}u|v|^{\beta}(\bar{v}z-v\bar{z}),\label{1120w1}\\
&\tilde{f}_2(u,v,w,z)=(\frac{\beta}{2}+1)|v|^{\beta}|u|^{\alpha+2}z-\frac{\beta v}{2\bar{v}}|v|^{\beta}|u|^{\alpha+2}\bar{z}+(\frac{\alpha}{2}+1)|v|^{\beta}v|u|^{\alpha}(\bar{u}w-u\bar{w}).\label{1120w2}
\end{align}

Noticing that $|\bar{u}|=|u|$, $|\bar{v}|=|v|$, $|\bar{w}|=|w|$, $|\bar{z}|=|z|$ and
$$
|u|^{\alpha+1}|v|^{\beta+1}\lesssim |u|^{\alpha}|v|^{\beta+2}+|u|^{\alpha+2}|v|^{\beta},
$$
using Strichartz estimates and similar to (\ref{911}), we can obtain
\begin{align}
\|w(t)\|_{S^0((-\infty, +\infty);\mathbb{R}^3)}+\|z(t)\|_{S^0((-\infty, +\infty);\mathbb{R}^3)}<+\infty,\label{1120w3}
\end{align}
and
\begin{align}
\|\tilde{f}_1(u,v,w,z,)\|_{L^2((-\infty, +\infty),L^{\frac{6}{5}}(\mathbb{R}^3))}+\|\tilde{f}_2(u,v,w,z)\|_{L^2((-\infty, +\infty),L^{\frac{6}{5}}(\mathbb{R}^3))}<+\infty.\label{1120w4}
\end{align}
Letting $\tilde{u}(t)=e^{-it\Delta}u(t)$ and $\tilde{v}(t)=e^{-it\Delta}v(t)$, we get
\begin{align}
x\tilde{u}(t)-x\tilde{u}(\tau)=i\lambda\int^t_{\tau}e^{-is\Delta}\tilde{f}_1(u,v,w,z)(s)ds,\label{1120w5}\\ x\tilde{v}(t)-x\tilde{v}(\tau)=i\lambda\int^t_{\tau}e^{-is\Delta}\tilde{f}_2(u,v,w,z)(s)ds.\label{1120w6}
\end{align}
Consequently, using (\ref{1120w3})--(\ref{1120w6}) and Strichartz estimates, we have
\begin{align}
&\quad \|x\tilde{u}(t)-x\tilde{u}(\tau)\|_{L^2}+\|x\tilde{v}(t)-x\tilde{v}(\tau)\|_{L^2}\nonumber\\
&\lesssim \|\tilde{f}_1(u,v,w,z)\|_{L^2((t,\tau),L^{\frac{6}{5}})}+\|\tilde{f}_1(u,v,w,z)\|_{L^2((t,\tau),L^{\frac{6}{5}})}\longrightarrow 0\label{1121xj1}
\end{align}
as $t,\tau \rightarrow +\infty$, which implies that
\begin{align}
\|x\tilde{u}(t)-xu_+\|_{L^2}+\|x\tilde{v}(t)-xv_+\|_{L^2} \rightarrow 0\quad {\rm as}\quad t\rightarrow +\infty.\label{1120w7}
\end{align}
Combining (\ref{1023w1}) and (\ref{1120w7}), we can get (\ref{11201}).

Similarly, (\ref{11202}) can be obtained. Theorem 4 is proved.\hfill $\Box$

\subsection{$\Sigma\times \Sigma$ scattering theory for (\ref{826x1}) with defocusing nonlinearities and $\alpha+\beta=2$ excluding the endpoints when $d=3$}
\qquad In this subsection, we will establish $\Sigma\times \Sigma$ scattering theory for (\ref{826x1}) with defocusing nonlinearities(i.e., $\lambda>0$, $\mu>0$) and $\alpha+\beta=2$ excluding the endpoints when $d=3$. That is, the exponents pair $(\alpha,\beta)$ is on the critical exponents line $\alpha+\beta=2$ excluding the endpoints $(\alpha,\beta)=(0,2)$ and $(\alpha,\beta)=(2,0)$. In this case, the method used in the subcritical case is invalid, the essential difficulty in technique is that we cannot obtain the expected estimates when we use H\"{o}lder's inequality. We will use the weight-coupled pseudo-conformal conservation law to give some estimates below.

{\bf Lemma 4.1(Weight-coupled pseudo-conformal conservation law).} {\it  Assume that $(u,v)$ is the global solution of (\ref{826x1}) with defocusing nonlinearities when $d=3$, $\alpha>0$, $\beta>0$, $\alpha+\beta=2$ and $(u_0,v_0)\in\Sigma\times \Sigma$. Then
\begin{align}
P(t)&=\int_{\mathbb{R}^3}[c_1|(x+2it\nabla)u|^2+c_2|(x+2it\nabla)v|^2]dx+4t^2\int_{\mathbb{R}^3} |u|^{\alpha+2}|v|^{\beta+2}dx\nonumber\\
&=\int_{\mathbb{R}^3}[c_1|xu_0|^2+c_2|xv_0|^2]dx-16\int_0^t\tau \int_{\mathbb{R}^3} |u|^{\alpha+2}|v|^{\beta+2}dx d\tau.\label{1113w5}
\end{align}
Here $(c_1,c_2)=(\frac{\alpha+2}{2\lambda},\frac{\beta+2}{2\mu})$.
}

{\bf Proof:} Using $E_w(u,v)=E_w(u_0,v_0)$, we get
\begin{align}
P(t)&=\int_{\mathbb{R}^3}[c_1|xu|^2+c_2|xv|^2]dx-4t\Im \int_{\mathbb{R}^3}[c_1\bar{u}(x\cdot \nabla u)+c_2\bar{v}(x\cdot \nabla v)]dx+4t^2E_w(u_0,v_0).\label{11142}
\end{align}
Recalling that
$$\frac{d}{dt} \int_{\mathbb{R}^3}|x|^2|u|^2dx=4\Im \int_{\mathbb{R}^3}\bar{u}(x\cdot \nabla u)dx,\quad \frac{d}{dt} \int_{\mathbb{R}^3}|x|^2|v|^2dx=4\Im \int_{\mathbb{R}^3}\bar{v}(x\cdot \nabla v)dx,$$
we obtain
\begin{align}
P'(t)&=\frac{d}{dt}\int_{\mathbb{R}^3}[c_1|xu|^2+c_2|xv|^2]dx-4\Im \int_{\mathbb{R}^3}[c_1\bar{u}(x\cdot \nabla u)+c_2\bar{v}(x\cdot \nabla v)]dx\nonumber\\
&\quad -4t\frac{d}{dt}\Im \int_{\mathbb{R}^3}[c_1\bar{u}(x\cdot \nabla u)+c_2\bar{v}(x\cdot \nabla v) ]dx+8tE_w(u_0,v_0)\nonumber\\
&=-4t\frac{d}{dt}\Im \int_{\mathbb{R}^3}[c_1\bar{u}(x\cdot \nabla u)+c_2\bar{v}(x\cdot \nabla v)]dx+8tE_w(u_0,v_0)\nonumber\\
&=-16t\int_{\mathbb{R}^3} |u|^{\alpha+2}|v|^{\beta+2}dx.\label{11143}
\end{align}
Integrating (\ref{11143}) from $0$ to $t$, we have (\ref{1113w5}).\hfill$\Box$

As the direct result of (\ref{1113w5}), we give the decay rate of the global solution $(u,v)$ of (\ref{826x1}).

{\bf Lemma 4.2(Decay rate of the global solution).} {\it Assume that $(u,v)$ is the global solution of (\ref{826x1}) with defocusing nonlinearities when $d=3$, $\alpha>0$, $\beta>0$, $\alpha+\beta=2$ and $(u_0,v_0)\in\Sigma\times \Sigma$.
Then
\begin{align}
\int_{\mathbb{R}^3}|u|^{\alpha+2}|v|^{\beta+2}dx\leq \frac{C}{t^2}.\label{11144}
\end{align}
}

{\bf Proof:} By (\ref{1113w5}), we have
\begin{align}
t^2\int_{\mathbb{R}^N}|u|^{\alpha+2}|v|^{\beta+2}dx\leq \int_{\mathbb{R}^3}[c_1|xu_0|^2+c_2|xv_0|^2]dx\leq C,\label{11212}
\end{align}
which implies (\ref{11144}).\hfill $\Box$

Now we will establish the $\Sigma\times \Sigma$ scattering result for the solution of (\ref{826x1}).

Note that $(2,6)$ is an admissible pair $(q,r)$ when $d=3$.

 We first prove that if $\alpha>0$, $\beta>0$ and $\alpha+\beta=2$, i.e., $(\alpha,\beta)$ is on the critical exponents line $\alpha+\beta=2$ excluding the endpoints $(0,2)$ and $(2,0)$, then
 \begin{align}
\|u\|_{L^2((0,t), W^{1,6})}+\|v\|_{L^2((0,t), W^{1,6})}\leq C \quad {\rm for }\quad t>0. \label{11216}
 \end{align}
Using Duhamel's principle and Strichartz estimates with $(q,r)=(2,6)$ and $(q',r')=(2,\frac{6}{5})$, and noticing that
$$|u|^{\alpha+1}|v|^{\beta+1}\leq C[|u|^{\alpha}|v|^{\beta+2}+|u|^{\alpha+2}|v|^{\beta}],$$
 we have
\begin{align}
&\quad \|u\|_{L^2((0,t), W^{1,6})}+\|v\|_{L^2((0,t), W^{1,6})}\nonumber\\
&\leq C\|u_0\|_{H^1}+C\|v_0\|_{H^1}+C\||u|^{\alpha}|v|^{\beta+2}u\|_{L^2((0,t),W^{1,\frac{6}{5}})}+C\||u|^{\alpha+2}|v|^{\beta}v\|_{L^2((0,t),W^{1,\frac{6}{5}})}\nonumber\\
&\leq C+C\left(\int_0^t\left(\int_{\mathbb{R}^3}[|u|^{\alpha}|v|^{\beta+2}+|u|^{\alpha+2}|v|^{\beta}]^{\frac{3}{2}}dx\right)^{\frac{4}{5}}
\left(\int_{\mathbb{R}^3}[|u|+|v|+|\nabla u|+|\nabla v|]^6dx\right)^{\frac{1}{3}}dt\right)^{\frac{1}{2}}\nonumber\displaybreak\\
&\leq C+C\left(\int_0^T\left(\int_{\mathbb{R}^3}[|u|^{\alpha}|v|^{\beta+2}+|u|^{\alpha+2}|v|^{\beta}]^{\frac{3}{2}}dx\right)^{\frac{4}{5}}
\left(\int_{\mathbb{R}^3}[|u|+|v|+|\nabla u|+|\nabla v|]^6dx\right)^{\frac{1}{3}}dt\right)^{\frac{1}{2}}\nonumber\\
&\qquad+C\left(\int_T^t\left(\int_{\mathbb{R}^3}[|u|^{\alpha}|v|^{\beta+2}+|u|^{\alpha+2}|v|^{\beta}]^{\frac{3}{2}}dx\right)^{\frac{4}{5}}
\left(\int_{\mathbb{R}^3}[|u|+|v|+|\nabla u|+|\nabla v|]^6dx\right)^{\frac{1}{3}}dt\right)^{\frac{1}{2}}\nonumber\\
&\leq C'+C\left(\int_T^t\left(\int_{\mathbb{R}^3}|u|^{\alpha+2}|v|^{\beta+2}dx\right)^{\frac{1}{\tau_1}}
\left(\int_{\mathbb{R}^3}|u|^6dx\right)^{\frac{1}{\tau_2}}
\left(\int_{\mathbb{R}^3}[|u|+|v|+|\nabla u|+|\nabla v|]^6dx\right)^{\frac{1}{3}}dt\right)^{\frac{1}{2}}\nonumber\\
&\qquad +C\left(\int_T^t\left(\int_{\mathbb{R}^3}|u|^{\alpha+2}|v|^{\beta+2}dx\right)^{\frac{1}{\tau_3}}
\left(\int_{\mathbb{R}^3}|v|^6dx\right)^{\frac{1}{\tau_4}}\left(\int_{\mathbb{R}^3}[|u|+|v|+|\nabla u|+|\nabla v|]^6dx\right)^{\frac{1}{3}}dt\right)^{\frac{1}{2}}\nonumber\\
&\leq C'+\max_{(T,t)} \left[\left(\int_{\mathbb{R}^3}|u|^{\alpha+2}|v|^{\beta+2}dx\right)^{\frac{1}{\tau_1}}
+\left(\int_{\mathbb{R}^3}|u|^{\alpha+2}|v|^{\beta+2}dx\right)^{\frac{1}{\tau_3}}\right]\nonumber\\
 &\qquad \qquad \times [\|u\|_{L^2((T,t), W^{1,6})}+\|v\|_{L^2((T,t), W^{1,6})}]\nonumber\\
&\leq C+\frac{1}{2}[\|u\|_{L^2((0,t), W^{1,6})}+\|v\|_{L^2((0,t), W^{1,6})}]\label{11217}
\end{align}
 because $\int_{\mathbb{R}^3}|u|^{\alpha+2}|v|^{\beta+2}dx$ can be small enough if $\alpha>0$, $\beta>0$, $\alpha+\beta=2$  when $t$ and $T$ are large while $\int_{\mathbb{R}^3}|u|^6dx$ and $\int_{\mathbb{R}^3}|v|^6dx$ are bounded for all $t\geq 0$. Here
 $$
 \frac{1}{\tau_1}=\frac{3\alpha}{2\alpha+4},\quad \frac{1}{\tau_2}=\frac{4-\alpha}{2\alpha+4},\quad \frac{1}{\tau_3}=\frac{3\beta}{2\beta+4},\quad \frac{1}{\tau_4}=\frac{4-\beta}{2\beta+4}.
 $$
Consequently,
 \begin{align}
\|u\|_{L^2((0,t), W^{1,6})}+\|v\|_{L^2((0,t), W^{1,6})}\leq C \quad {\rm for \ all} \ t\geq 0.\label{11218}
\end{align}

As a byproduct of (\ref{11217}), we get
\begin{align*}
\||u|^{\alpha}|v|^{\beta+2}u\|_{L^2((0,t),W^{1,\frac{6}{5}})}+\||u|^{\alpha+2}|v|^{\beta}v\|_{L^2((0,t),W^{1,\frac{6}{5}})}\leq C,\quad t>0,
\end{align*}
which implies that
\begin{align}
\||u|^{\alpha}|v|^{\beta+2}u\|_{L^2((t,\tau),W^{1,\frac{6}{5}})}+\||u|^{\alpha+2}|v|^{\beta}v\|_{L^2((t,\tau),W^{1,\frac{6}{5}})}\longrightarrow 0 \quad {\rm as}\quad t,\ \tau\rightarrow +\infty. \label{11219'}
\end{align}
Therefore, we obtain
\begin{align}
&\quad \|e^{-it\Delta}u(t)-e^{-i\tau\Delta}u(\tau)\|_{H^1_x}+\|e^{-it\Delta}v(t)-e^{-i\tau\Delta}v(\tau)\|_{H^1_x}\nonumber\\
&\leq \|\int_t^{\tau} e^{-is\Delta}|u|^{\alpha}|v|^{\beta+2}u(s)ds\|_{H^1_x}+\|\int_t^{\tau} e^{-is\Delta}|u|^{\alpha+2}|v|^{\beta}v(s)ds\|_{H^1_x}\nonumber\\
&\leq C\||u|^{\alpha}|v|^{\beta+2}u\|_{L^2((t,\tau),W^{1,\frac{6}{5}})}+C\||u|^{\alpha+2}|v|^{\beta}v\|_{L^2((t,\tau),W^{1,\frac{6}{5}})}
\nonumber\\
&\quad \longrightarrow 0 \quad {\rm as}\quad t,\ \tau\rightarrow +\infty,\label{11219}
\end{align}
which implies that there exists $(u_+,v_+)$ such that
\begin{align}
\|e^{-it\Delta}u(t)-u_+\|_{H^1_x}+\|e^{-it\Delta}v(t)-v_+\|_{H^1_x}\longrightarrow 0 \quad {\rm as}\quad t,\ \tau\rightarrow +\infty.\label{1121x1}
\end{align}

Letting $w$, $z$, $\tilde{f}_1(u,v,w,z)$ and $\tilde{f}_2(u,v,w,z)$ be defined as in (\ref{11203})--(\ref{1120w3}), the difference is that $\alpha+\beta=2$
now. Similar to (\ref{11217})--(\ref{11219'}), we can prove that
\begin{align}
\|w(t)\|_{L^2((0,t), W^{1,6})}+\|z(t)\|_{L^2((0,t), W^{1,6})}\leq C \quad {\rm for }\quad t>0.,\label{1121x2}
\end{align}
and
\begin{align}
\|\tilde{f}_1(u,v,w,z,)\|_{L^2((-\infty, +\infty),L^{\frac{6}{5}}(\mathbb{R}^3))}+\|\tilde{f}_2(u,v,w,z)\|_{L^2((-\infty, +\infty),L^{\frac{6}{5}}(\mathbb{R}^3))}<+\infty.\label{1121x3}
\end{align}
Similar to (\ref{1120w5})--(\ref{1120w7}) and (\ref{11219}), denoting $\tilde{u}(t)=e^{-it\Delta}u(t)$ and $\tilde{v}(t)=e^{-it\Delta}v(t)$, we can obtain
\begin{align}
\|x\tilde{u}(t)-xu_+\|_{L^2_x}+\|x\tilde{v}(t)-xv_+\|_{L^2_x} \rightarrow 0\quad {\rm as}\quad t\rightarrow +\infty.\label{1121x4}
\end{align}
Hence, the solution of (\ref{826x1}) has scattering state in $\Sigma\times \Sigma$.

Therefore, we have proved Theorem 5.\hfill $\Box$

{\bf Remark 4.3.} We cannot establish the $H^1\times H^1$ or $\Sigma\times \Sigma$ scattering theory for (\ref{826x1}) when $d=3$ and $(\alpha,\beta)$ is the endpoint of the critical exponents line $\alpha+\beta=2$, i.e., $(\alpha,\beta)=(0,2)$ or $(\alpha,\beta)=(2,0)$. In the next section, we will establish $\dot{H}^1\times \dot{H}^1$ scattering theory for (\ref{826x1}) with defocusing nonlinearities and $(\alpha,\beta)$ is the endpoint on the critical exponents line when $d=3$.

\section{$\dot{H}^1\times \dot{H}^1$ scattering theory for (\ref{826x1}) with defocusing nonlinearities and $(\alpha,\beta)$ is the endpoint on the critical exponents line when $d=3$}
\qquad In this section, we establish the $\dot{H}^1\times \dot{H}^1$ scattering theory for (\ref{826x1}) when $d=3$ and $(\alpha,\beta)$ is the endpoint on the critical exponents line $\alpha+\beta=2$, i.e., $(\alpha,\beta)=(0,2)$ or $(\alpha,\beta)=(2,0)$.

In the sequels, we only prove Theorem 6 in the case of $(\alpha,\beta)=(0,2)$, the conclusion in the case of $(\alpha,\beta)=(2,0)$ can be proved similarly.

We will use the method of reduction to almost periodic solutions and argue by contradiction. We will prove that Theorem 6 holds for solution with small energy by simple contraction mapping arguments. If the conclusions of this theorem are not true, we can find the solution with a transition energy above which the energy no longer controls the spacetime norm. Therefore, first we will show that there is a minimal counterexample which has good compactness properties.

{\bf Definition 5.1(Almost periodicity).} {\it A solution $(u,v)\in [L^{\infty}_t\dot{H}^1_x(I\times \mathbb{R}^3)]^2$ of (\ref{826x1}) is said to be almost periodic(modulo symmetries) if there exist functions $\widetilde{N}:I\rightarrow \mathbb{R}^+$, $\tilde{x}:I\rightarrow \mathbb{R}^3$, and $C:\mathbb{R}^+\rightarrow \mathbb{R}^+$
such that for all $t\in I$ and $\eta>0$
\begin{align}
\int_{|x-\tilde{x}(t)|\geq \frac{C(\eta)}{\widetilde{N}(t)}}[|\nabla u(t,x)|^2+|\nabla v(t,x)|^2]dx
+\int_{|\xi|\geq C(\eta)\widetilde{N}(t)}|\xi|^2[|\hat{u}(t,\xi)|^2+|\hat{v}(t,\xi)|^2]d\xi\leq \eta.\label{934}
\end{align}
The functions $\widetilde{N}(t)$, $\tilde{x}(t)$ and $C(\eta)$ are called as the frequency scale function for the solution $(u,v)$, the spatial center function and the modulus of compactness respectively. }

By compactness, there exists $c(\eta)>0$ such that
\begin{align}
\int_{|x-\tilde{x}(t)|\geq \frac{C(\eta)}{\widetilde{N}(t)}}[|\nabla u(t,x)|^2+|\nabla v(t,x)|^2]dx
+\int_{|\xi|\geq C(\eta)\widetilde{N}(t)}|\xi|^2[|\hat{u}(t,\xi)|^2+|\hat{v}(t,\xi)|^2]d\xi\leq \eta,\label{12251}
\end{align}
and
\begin{align}
\int_{\mathbb{R}^3}[|\nabla u(t,x)|^2+|\nabla v(t,x)|^2]dx\lesssim_{(u,v)} \int_{\mathbb{R}^3}[|u(t,x)|^6+|v(t,x)|^6]dx\label{12252}
\end{align}
for all $t\in I$.

The following proposition is very important to the proof of Theorem 6.

{\bf Proposition 5.2(Reduction to almost periodic solution).} {\it
Assume that Theorem 6 failed. Then there exists a maximal-lifespan solution $(u,v):I\times \mathbb{R}^3\rightarrow \mathbb{C}\times \mathbb{C}$
of (\ref{826x1}) with $(\alpha,\beta)=(0,2)$ which is almost periodic and blows up both forward and backward in time in the sense of
$$
\int_{t_0}^{\sup I}\int_{\mathbb{R}^3}[|u(t,x)|^{10}+|v(t,x)|^{10}]dxdt=\int^{t_0}_{\inf I}\int_{\mathbb{R}^3}[|u(t,x)|^{10}+|v(t,x)|^{10}]dxdt=+\infty
$$
for all $t_0\in I$. }

The proof of this proposition will be given in Subsection 5.2 by some conclusions.

 By continuity, the modulation parameters $\tilde{x}(t)$ and $\widetilde{N}(t)$ cannot change rapidly,  similar to Lemma 5.18 in \cite{Killip2013}, we obtain

{\bf Lemma 5.3(Local constancy property).} {\it
 Let $(u,v):I\times \mathbb{R}^3\rightarrow \mathbb{C}\times \mathbb{C}$ be a maximal-lifespan almost period solution of (\ref{826x1}) with $(\alpha,\beta)=(0,2)$. Then there exists a small number $\delta$ which depends only on $(u,v)$, such that
 \begin{align}
 [t_0-\delta\widetilde{N}(t_0)^{-2}, t_0+\delta\widetilde{N}(t_0)^{-2}]\subset I\quad {\rm if }\quad t_0\in I,\label{12171}
 \end{align}
 and
 \begin{align}
 \widetilde{N}(t)\thicksim_{(u,v)}\widetilde{N}(t_0)\quad {\rm whenever}\quad |t-t_0|\leq \delta\widetilde{N}(t_0)^{-2}.\label{12172}
\end{align}
 }

As the direct result of (\ref{12171}) and (\ref{12172}), we give a corollary of Lemma 5.3 below

{\bf Corollary 5.4($\widetilde{N}(t)$ blows up).} {\it
Let $(u,v):I\times \mathbb{R}^3\rightarrow \mathbb{C}\times \mathbb{C}$ be a maximal-lifespan almost period solution of (\ref{826x1}) with $(\alpha,\beta)=(0,2)$. If $T$ is any finite endpoint of $I$, then $\widetilde{N}(t)\geq_{(u,v)}|T-t|^{-\frac{1}{2}}$. Consequently, $\lim_{t\rightarrow T} \tilde{N}(t)=+\infty$.
}

The following lemma reveals the relation between $\widetilde{N}(t)$ of an almost periodic solution $(u,v)$ and its Strichartz norms:

{\bf Lemma 5.5(Spacetime bounds).} {\it Let $(u,v)$ be an almost periodic solution of (\ref{826x1}) with $(\alpha,\beta)=(0,2)$ on a time interval $I$. Then for any admissible pair $(q,r)$ with $2\leq q<+\infty$,
\begin{align}
\int_I\widetilde{N}(t)^2dt\lesssim_{(u,v)}\|\nabla u\|^q_{L^q_tL^r_x(I\times \mathbb{R}^3)}+\|\nabla v\|^q_{L^q_tL^r_x(I\times \mathbb{R}^3)}\lesssim_{(u,v)}1+\int_I\widetilde{N}(t)^2dt\label{10241}
\end{align}
}

{\bf Proof:} Similar to the proof of Lemma 5.21 in \cite{Killip2013}, we first prove
$$
\|\nabla u\|^q_{L^q_tL^r_x(I\times \mathbb{R}^3)}+\|\nabla v\|^q_{L^q_tL^r_x(I\times \mathbb{R}^3)}\lesssim_{(u,v)} 1+\int_I\widetilde{N}(t)^2dt.
$$
Let $0<\epsilon<1$ be a small parameter to be chosen later and subdivide $I$ into
subintervals $I_j$ such that
\begin{align}
\int_{I_j}\widetilde{N}^2(t)dt\leq \epsilon, \label{1218x1}
\end{align}
the number of such intervals is at most $\epsilon^{-1}[1+\int_I\widetilde{N}(t)^2dt]$.

For each $j$, we may choose $t_j\in I_j$ such that
\begin{align}
 |I_j|^{\frac{d-2}{2(d+2)}}[\widetilde{N}^2(t_j)]^{\frac{2d}{d+2}}\leq 2\epsilon.\label{1218x2}
 \end{align}
Using Strichartz estimate, for any admissible pair $(q,r)$, we have
\begin{align*}
&\quad \|\nabla u\|^q_{L^q_tL^r_x}+\|\nabla v\|^q_{L^q_tL^r_x}\nonumber\\
&\lesssim_{(u,v)} \|\nabla u_{\geq M\widetilde{N}^2(t_j)}(t_j)\|^q_{L^2_x}+|I_j|^{\frac{q(d-2)}{2(d+2)}}[M\widetilde{N}^2(t_j)]^{\frac{2qd}{d+2}}\|\nabla u(t_j)\|^q_{L^2_x}+\|\nabla v_{\geq M\widetilde{N}^2(t_j)}(t_j)\|^q_{L^2_x}\nonumber\\
&\quad +|I_j|^{\frac{q(d-2)}{2(d+2)}}[M\widetilde{N}^2(t_j)]^{\frac{2qd}{d+2}}\|\nabla v(t_j)\|^q_{L^2_x}+[\|\nabla u\|^q_{L^q_tL^r_x}+\|\nabla v\|^q_{L^q_tL^r_x}]^{\frac{(d+2)}{d}},
\end{align*}
where all spacetime norms are taken on the slab $I_j\times \mathbb{R}^d$. If $M$ is large, by Definition 5.1,  the first term can be small. By (\ref{1218x2}), the second term also can small. Hence, we can use the bootstrap argument to get
\begin{align}
 \|\nabla u\|^q_{L^q_tL^r_x}+\|\nabla v\|^q_{L^q_tL^r_x}\lesssim_{(u,v)} 1+\int_I\widetilde{N}(t)^2dt.\label{1218w1}
\end{align}
On the other hand, since $(u,v)\not\equiv (0,0)$, $[\widetilde{N}(t)]^{-\frac{2}{q}}[\|\nabla u\|_{L^r_x}+\|\nabla v \|_{L^r_x}]$ will never vanish, almost periodicity implies that it is bounded away from zero, i.e.,
$$[\|\nabla u\|_{L^r_x}+\|\nabla v \|_{L^r_x}]\gtrsim_{(u,v)} [\widetilde{N}(t)]^{\frac{2}{q}}.$$
Integrating it on $I$, we obtain
\begin{align}
\int_I\widetilde{N}(t)^2dt\lesssim_{(u,v)}\|\nabla u\|^q_{L^q_tL^r_x(I\times \mathbb{R}^3)}+\|\nabla v\|^q_{L^q_tL^r_x(I\times \mathbb{R}^3)}.\label{1218w2}
\end{align}
(\ref{10241}) is proved.\hfill $\Box$

Without loss of generality, we can assume that $\widetilde{N}(t)\geq 1$ on $[0,T_{\max})$ because it can be realized by a simple rescaling argument. Finally, we get

{\bf Proposition 5.6(Two special scenarios for blowup).} {\it Assume that Theorem 6 failed. Then there exists an almost periodic solution $(u,v):[0,T_{\max})\times \mathbb{R}^3\rightarrow \mathbb{C}\times \mathbb{C}$ such that
$$
\|u\|_{L_{t,x}^{10}([0,T_{\max})\times \mathbb{R}^3)}+\|v\|_{L_{t,x}^{10}([0,T_{\max})\times \mathbb{R}^3)}=+\infty
$$
and $[0,T_{\max})=\cup_{k}J_k$, where $J_k$ are characteristic intervals on which $\widetilde{N}(t)\equiv N_k\geq 1$.
Moreover,
$$
{\rm either}\quad \int_0^{T_{\max}}\widetilde{N}(t)^{-1}dt<+\infty\quad {\rm or}\quad  \int_0^{T_{\max}}\widetilde{N}(t)^{-1}dt=+\infty.
$$
}

Consequently, we only need to preclude the existence of the two types of almost periodic solution described in Proposition 5.6, then we can give the proof of Theorem 6. The proof of Proposition 5.6 will be given later.

To prove the no-existence of cascade solutions, we need the following proposition.

{\bf Proposition 5.7(No-waste Duhamel formulae).} {\it
Let $(u,v):[0,T_{\max})\times \mathbb{R}^3\rightarrow \mathbb{C}\times \mathbb{C}$ be defined as in Proposition 5.6. Then for all $t\in [0,T_{\max})$
\begin{align}
u(t)=i\lim_{T\rightarrow T_{\max}} \int_t^Te^{i(t-s)\Delta}|v|^4u(s)ds,\quad v(t)=i\lim_{T\rightarrow T_{\max}} \int_t^Te^{i(t-s)\Delta}|u|^2|v|^2v(s)ds\label{1125s1}
\end{align}
in the weak $\dot{H}^1_x$ topology.
}

{\bf Proof:} If $T_{\max}<+\infty$, then by Corollary 5.4,
$$
\lim_{t\rightarrow T_{\max}} \widetilde{N}(t)=+\infty.
$$
By Definition 5.1, this implies that $(u(t),v(t))$ converges weakly to $(0,0)$ as $t\rightarrow T_{\max}$. Since
$T_{\max}<+\infty$ and the map $t\rightarrow e^{-it\Delta}$ is continuous in the strong operator topology on
$L^2_x$, we know that $(e^{-it\Delta}u(t),e^{-it\Delta}v(t))$ converges weakly to $(0,0)$.

If $T_{\max}=+\infty$, then for any test functions pair $(\varphi,\psi)\in C^{\infty}_0(\mathbb{R}^d)\times C^{\infty}_0(\mathbb{R}^d)$ and $\eta>0$, we have
\begin{align}
&\quad |<u(t),e^{it\Delta}\varphi>_{L^2_x(\mathbb{R}^d)}|+|<v(t),e^{it\Delta}\psi>_{L^2_x(\mathbb{R}^d)}|\nonumber\\
&\lesssim \int_{|x-x(t)|\leq \frac{C(\eta)}{\widetilde{N}(t)}}|e^{it\Delta}\varphi|^2dx+\eta\|\varphi\|^2_{L^2_x}+\int_{|x-x(t)|\leq \frac{C(\eta)}{\widetilde{N}(t)}}|e^{it\Delta}\psi|^2dx+\eta\|\psi\|^2_{L^2_x}.\label{1218w3}
\end{align}
Using Fraunhofer formula, by the results of Corollary 5.4 and changing variables, we can see that $(e^{-it\Delta}u(t),e^{-it\Delta}v(t))$ converges weakly to $(0,0)$.

In any case, $(e^{-it\Delta}u(t),e^{-it\Delta}v(t))$ converges weakly to $(0,0)$ as $t\rightarrow T_{\max}$, then we can use Duhamel formulae to get (\ref{1125s1}). \hfill $\Box$

\subsection{Stability result about (\ref{826x1})}
\qquad Since we have established the local wellposedness of $H^1\times H^1$-solution to (\ref{826x1}) with initial data $(u_0,v_0)\in H^1\times H^1$, to remove the constraint conditions $(u_0,v_0)\in L^2\times L^2$, in this subsection, we show a stability result about (\ref{826x1}) when $d=3$ and $(\alpha,\beta)$ is the endpoint of critical exponents line, which will prove the local wellposedness of $\dot{H}^1\times \dot{H}^1$-solution to (\ref{826x1}) with initial data $(u_0,v_0)\in\dot{H}^1\times \dot{H}^1$. Since there is the similar stability result about (\ref{826x1}) when $d=4$ and $(\alpha,\beta)=(0,0)$, we would like to put them together and state below. Although the proofs of some propositions and lemmas in this subsection are very similar to those in \cite{Killip2013} for the scalar equation and seem standard, however, to our best knowledge, the present paper is the first one discussing $\dot{H}^1\times \dot{H}^1$ scattering phenomenon on a coupled system of Schr\"{o}dinger equations in energy critical case, for completeness, we would like to give the details.

{\bf Proposition 5.8(Stability result).} {\it  Let $I$ be a compact time interval
and let $(\tilde{u},\tilde{v})$ be an approximate solution of (\ref{826x1}) on $I\times \mathbb{R}^d$ in the sense of
\begin{align}
i\tilde{u}_t+\Delta \tilde{u}=\lambda|\tilde{u}|^{\alpha}|\tilde{v}|^{\beta+2}\tilde{u}+e_1,\quad
i\tilde{v}_t+\Delta \tilde{v}=\mu|\tilde{u}|^{\alpha+2}|\tilde{v}|^{\beta}\tilde{u}+e_2.\label{9151}
\end{align}
for some functions $e_1$ and $e_2$. Here $(\alpha,\beta)=(0,2)$ or $(\alpha,\beta)=(2,0)$ when $d=3$, while $(\alpha,\beta)=(0,0)$ when $d=4$. Suppose that
\begin{align}
\|\tilde{u}\|_{L^{\infty}_t\dot{H}^1_x}+\|\tilde{v}\|_{L^{\infty}_t\dot{H}^1_x}\leq E,\quad \|\tilde{u}\|_{L^{\frac{2(d+2)}{d-2}}_{t,x}}+\|\tilde{v}\|_{L^{\frac{2(d+2)}{d-2}}_{t,x}}\leq L\label{9152}
\end{align}
for some positive constants $E$ and $L$. Assume that $t_0\in I$ and
\begin{align}
&\|u(t_0)-\tilde{u}(t_0)\|_{\dot{H}^1_x}+\|v(t_0)-\tilde{v}(t_0)\|_{\dot{H}^1_x}\leq E',\label{9153}\\
&\|e^{i(t-t_0)\Delta} [u(t_0)-\tilde{u}(t_0)]\|_{L^{\frac{2(d+2)}{d-2}}_{t,x}}+\|e^{i(t-t_0)\Delta} [v(t_0)-\tilde{v}(t_0)]\|_{L^{\frac{2(d+2)}{d-2}}_{t,x}}\leq\epsilon,\label{9154}\\
&\|\nabla e_1\|_{N^0(I)}+\|\nabla e_2\|_{N^0(I)}\leq \epsilon\label{9155}
\end{align}
for some positive constants $E'$ and $0<\epsilon<\epsilon_1=\epsilon_1(E,E',L)$. Then there exists a unique strong solution
$(u,v):I\times \mathbb{R}^d\rightarrow \mathbb{C}\times \mathbb{C}$ with initial data $(u(t_0),v(t_0))$ at $t_0$ which satisfies
\begin{align}
&\|\nabla u(t)]\|_{\dot{S}^0(I)}+\|\nabla v]\|_{\dot{S}^0(I)}\leq C(E,E',L),\label{9156}\\
&\|\nabla [u(t)-\tilde{u}(t)]\|_{\dot{S}^0(I)}+\|\nabla [v(t)-\tilde{v}(t)]\|_{\dot{S}^0(I)}\leq C(E,E',L)E',\label{9157}\\
&\|[u(t)-\tilde{u}(t)]\|_{L^{\frac{2(d+2)}{d-2}}_{t,x}}+\|[v(t)-\tilde{v}(t)]\|_{L^{\frac{2(d+2)}{d-2}}_{t,x}}\leq C(E,E',L)\epsilon^c\label{9158}
\end{align}
for some $0<c=c(d)<1$.
}

As a consequences of Proposition 5.8, we have the following local well-posedness result.

{\bf Corollary 5.9(Local well-posednes). } {\it
Let $I$ be a compact time interval, $t_0\in I$. Assume that $(\alpha,\beta)=(0,2)$ or $(\alpha,\beta)=(2,0)$ when $d=3$, while $(\alpha,\beta)=(0,0)$ when $d=4$, $(u_0,v_0)\in \dot{H}^1\times \dot{H}^1$ and
$$
\|u_0\|_{\dot{H}^1_x}+\|v_0\|_{\dot{H}^1_x}\leq E.
$$
For any $\epsilon>0$, there exists $\delta=\delta(E,\epsilon)$ such that if
$$
\|e^{i(t-t_0)\Delta} u(t_0)\|_{L^{\frac{2(d+2)}{d-2}}_{t,x}}+\|e^{i(t-t_0)\Delta} v(t_0)\|_{L^{\frac{2(d+2)}{d-2}}_{t,x}}\leq\delta,
$$
then (\ref{826x1}) has  a unique solution $(u,v)$ with initial data $(u(t_0),v(t_0))=(u(0),v(0))$ and
\begin{align*}
\|u(t)\|_{L^{\frac{2(d+2)}{d-2}}_{t,x}(I\times \mathbb{R}^d)}+\|v(t)\|_{L^{\frac{2(d+2)}{d-2}}_{t,x}(I\times \mathbb{R}^d)}\leq\epsilon,\quad \|\nabla u\|_{S^0(I)}+\|\nabla v\|_{S^0(I)}\leq 2E.
\end{align*}
}

In order to prove stability Theorem, similar to \cite{Killip2013}, for any time interval $I$, we introduce some spaces $X^0(I)$, $X(I)$, $Y(I)$ and denote
\begin{align}
&\|w\|_{X^0(I)}=\|w\|_{L^{\frac{d(d+2)}{2(d-2)}}_tL^{\frac{2d^2(d+2)}{(d+4)(d-2)^2}}_x(I\times \mathbb{R}^d)},\label{9159}\\
&\|w\|_{X(I)}=\||\nabla|^{\frac{4}{d+2}} w\|_{L^{\frac{d(d+2)}{2(d-2)}}_tL^{\frac{2d^2(d+2)}{d^3-4d+16}}_x(I\times \mathbb{R}^d)},\label{91510}\\
&\|f\|_{Y(I)}=\||\nabla|^{\frac{4}{d+2}} f\|_{L^{\frac{d}{2}}_tL^{\frac{2d^2(d+2)}{d^3+4d^2+4d-16}}_x(I\times \mathbb{R}^d)}.\label{91510'}
\end{align}
Then by the results of Lemma 3.10 and Lemma 3.11 in \cite{Killip2013},
\begin{align}
&\|w\|_{X^0(I)}\lesssim \|w\|_{X(I)}\lesssim \|\nabla w\|_{S^0(I)}, \label{91511}\\
&\|w\|_{X(I)}\lesssim \|w\|^{\frac{1}{d+2}}_{L^{\frac{2(d+2)}{d-2}}_{t,x}(I\times \mathbb{R}^d)}\|\nabla w\|^{\frac{d+1}{d+2}}_{S^0(I)},\label{91512}\\
&\|w\|_{L^{\frac{2(d+2)}{d-2}}_{t,x}(I\times \mathbb{R}^d)}\lesssim \|w\|^c_{X(I)}\|\nabla w\|^{1-c}_{S^0(I)}\quad {\rm for \ some}\ 0<c=c(d)\leq 1,\label{91513}\\
&\|\int^t_{t_0}e^{i(t-s)\Delta}f(s)\|_{X(I)}\lesssim \|f\|_{Y(I)}.\label{91514}
\end{align}

Now we will prove that

{\bf Lemma 5.10(Nonlinear estimates).} {\it  Denote $f_1(u,v)=\lambda|v|^4u$ and $f_2(u,v)=\mu|u|^2|v|^2v$ if $d=3$, while $f_1(u,v)=\lambda|v|^2u$ and $f_2(u,v)=\mu|u|^2v$ if $d=4$. Let $I$ be a time interval. Then
\begin{align}
&\|f_1(u,v)\|_{Y(I)}+\|f_2(u,v)\|_{Y(I)}\lesssim \|u\|^{\frac{d+2}{d-2}}_{X(I)}+\|v\|^{\frac{d+2}{d-2}}_{X(I)},\label{91515}\displaybreak\\
&\|f_{1u}(u_1+u_2,v_1+v_2)w_1\|_{Y(I)}+\|f_{1\bar{u}}(u_1+u_2,v_1+v_2)\bar{w}_1\|_{Y(I)}\nonumber\\
&\ +\|f_{1v}(u_1+u_2,v_1+v_2)w_2\|_{Y(I)}+\|f_{1\bar{v}}(u_1+u_2,v_1+v_2)\bar{w}_2\|_{Y(I)}\nonumber\\
&\ +\|f_{2u}(u_1+u_2,v_1+v_2)w_1\|_{Y(I)}+\|f_{2\bar{u}}(u_1+u_2,v_1+v_2)\bar{w}_1\|_{Y(I)}\nonumber\\
&\ +\|f_{2v}(u_1+u_2,v_1+v_2)w_2\|_{Y(I)}+\|f_{2\bar{v}}(u_1+u_2,v_1+v_2)\bar{w}_2\|_{Y(I)}\nonumber\\
&\lesssim \left(\sum_{m=1}^2[\|u_m\|^{\frac{8}{d^2-4}}_{X(I)}\|\nabla u_m\|^{\frac{4d}{d^2-4}}_{S^0(I)}+\|v_m\|^{\frac{8}{d^2-4}}_{X(I)}\|\nabla v_m\|^{\frac{4d}{d^2-4}}_{S^0(I)}] \right)[\|w_1\|_{X(I)}+\|w_2\|_{X(I)}].\label{9191}
\end{align}
}

{\bf Proof:} If $d=3$, then
\begin{align*}
&\|w\|_{X^0(I)}=\|w\|_{L^{\frac{15}{2}}_tL^{\frac{90}{7}}_x(I\times \mathbb{R}^d)},\quad
\|w\|_{X(I)}=\||\nabla|^{\frac{4}{5}} w\|_{L^{\frac{15}{2}}_tL^{\frac{90}{31}}_x(I\times \mathbb{R}^d)},\\
&\|f\|_{Y(I)}=\||\nabla|^{\frac{4}{5}} f\|_{L^{\frac{3}{2}}_tL^{\frac{90}{59}}_x(I\times \mathbb{R}^d)}.
\end{align*}
Using H\"{o}lder's and Young's inequalities, and the interpolation (\ref{91511})--(\ref{91514}), we have
\begin{align}
&\quad\|f_1(u,v)\|_{Y(I)}=\||\nabla|^{\frac{4}{5}} (|v|^4u)\|_{L^{\frac{3}{2}}_tL^{\frac{90}{59}}_x(I\times \mathbb{R}^d)}\nonumber\\
&\lesssim \|[\|u\|_{L^{\frac{90}{7}}_x}\|\nabla|^{\frac{4}{5}} (|v|^4)\|_{L^{\frac{90}{52}}_x}+\||v|^4\|_{L^{\frac{90}{28}}_x}\|\nabla|^{\frac{4}{5}}u\|_{L^{\frac{90}{31}}_x}]\|_{L^{\frac{3}{2}}_t}\nonumber\\
&\lesssim \|[\|u\|_{L^{\frac{90}{7}}_x}\|(|v|^3)\|_{L^{\frac{30}{7}}_x}\|\nabla|^{\frac{4}{5}} v\|_{L^{\frac{90}{31}}_x}+\|v\|^4_{L^{\frac{90}{7}}_x}\|\nabla|^{\frac{4}{5}}u\|_{L^{\frac{90}{31}}_x}]\|_{L^{\frac{3}{2}}_t}\nonumber\\
&\lesssim \|[\|u\|_{L^{\frac{90}{7}}_x}\|v\|^3_{L^{\frac{90}{7}}_x}\|\nabla|^{\frac{4}{5}} v\|_{L^{\frac{90}{31}}_x}+\|v\|^4_{L^{\frac{90}{7}}_x}\|\nabla|^{\frac{4}{5}}u\|_{L^{\frac{90}{31}}_x}]\|_{L^{\frac{3}{2}}_t}\nonumber\\
&\lesssim \|u\|_{X^0(I)}\|v\|^3_{X^0(I)}\|v\|_{X(I)}+\|v\|^4_{X^0(I)}\|u\|_{X(I)}\nonumber\\
&\lesssim \|u\|^5_{X(I)}+\|v\|^5_{X(I)},\label{9201}\\
&\quad \|f_2(u,v)\|_{Y(I)}=\||\nabla|^{\frac{4}{5}} (|u|^2|v|^2v)|_{L^{\frac{3}{2}}_tL^{\frac{90}{59}}_x(I\times \mathbb{R}^d)}\nonumber\\
&\lesssim \|[\||u|^2\|_{L^{\frac{45}{7}}_x}\||\nabla|^{\frac{4}{5}}|v|^2v\|_{L^2_x}+\||\nabla|^{\frac{4}{5}}(|u|^2)\|_{L^{\frac{90}{38}}_x}
\||v|^2v\|_{L^{\frac{30}{7}}_x}]\|_{L^{\frac{3}{2}}_t}\nonumber\\
&\lesssim \|[\|u\|^2_{L^{\frac{90}{7}}_x}\|v\|^2_{L^{\frac{90}{7}}_x}\||\nabla|^{\frac{4}{5}}v\|_{L^{\frac{90}{31}}_x}
+\|u\|_{L^{\frac{90}{7}}_x}\||\nabla|^{\frac{4}{5}}u\|_{L^{\frac{90}{31}}_x}
\|v\|^3_{L^{\frac{90}{7}}_x}]\|_{L^{\frac{3}{2}}_t}\nonumber\\
&\lesssim \|u\|^2_{X^0(I)}\|v\|^2_{X^0(I)}\|v\|_{X(I)}+\|u\|_{X^0(I)}\|u\|_{X(I)}\|v\|^3_{X^0(I)}\nonumber\\
&\lesssim \|u\|^5_{X(I)}+\|v\|^5_{X(I)}.\label{9202}
\end{align}
(\ref{91515}) is proved when $d=3$.

Since
\begin{align*}
&f_{1u}=\lambda|v|^4,\qquad f_{1\bar{u}}=0,\qquad \quad f_{1v}=2\lambda|v|^2\bar{v}u,\quad f_{1\bar{v}}=2\lambda|v|^2vu,\\
&f_{2u}=\mu|v|^2\bar{u}v,\quad f_{2\bar{u}}=\mu|v|^2uv,\quad f_{2v}=\mu|u|^2|v|^2,\quad f_{2\bar{v}}=\mu|u|^2v^2,
\end{align*}
\begin{align}
LHS(\ref{9191})&\lesssim \||v_1+v_2|^4w_1\|_{Y(I)}+\||v_1+v_2|^2(\bar{v}_1+\bar{v}_2)(u_1+u_2)w_2\|_{Y(I)}\nonumber\\
&\quad +\||v_1+v_2|^2(u_1+u_2)(v_1+v_2)\bar{w}_2\|_{Y(I)}+\||v_1+v_2|^2(\bar{u}_1+\bar{u}_2)(v_1+v_2)w_1\|_{Y(I)}\nonumber\\
&\quad +\||v_1+v_2|^2(u_1+u_2)(v_1+v_2)w_1\|_{Y(I)}+\||u_1+u_2|^2|v_1+v_2|^2w_2\|_{Y(I)}\nonumber\\
&\quad +\||u_1+u_2|^2(v_1+v_2)^2w_2\|_{Y(I)}\nonumber\\
&:=(I)+(II)+(III)+(IV)+(V)+(VI)+(VII).\label{920x1}
\end{align}
We just give the details of estimating (IV) below,
\begin{align}
(IV)&\lesssim\|[\||v_1+v_2|^2(\bar{u}_1+\bar{u}_2)(v_1+v_2)\|_{L^{\frac{45}{14}}_x}\||\nabla|^{\frac{4}{5}}w_1\|_{L^{\frac{90}{31}}_x}\nonumber\\
&\quad +\||\nabla|^{\frac{4}{5}}(|v_1+v_2|^2(\bar{u}_1+\bar{u}_2)(v_1+v_2))\|_{L^{\frac{45}{26}}_x}\|w_1\|_{L^{\frac{90}{7}}_x}]\|_{L^{\frac{3}{2}}_t}\nonumber\\
&\lesssim \|\left(\|v_1+v_2\|^3_{L^{\frac{90}{7}}_x}\|u_1+u_2\|_{L^{\frac{90}{7}}_x}\||\nabla|^{\frac{4}{5}}w_1\|_{L^{\frac{90}{31}}_x}\right)\|_{L^{\frac{3}{2}}_t}\nonumber\\
&\quad +\|\left(\||v_1+v_2|^2(v_1+v_2)\|_{L^{\frac{30}{7}}_x}\||\nabla|^{\frac{4}{5}}(\bar{u}_1+\bar{u}_2)\|_{L^{\frac{90}{31}}_x}
\|w_1\|_{L^{\frac{90}{7}}_x}\right)\|_{L^{\frac{3}{2}}_t}\nonumber\\
&\quad+\|\left(\|v_1+v_2\|^2_{L^{\frac{90}{7}}_x}\||\nabla|^{\frac{4}{5}}(v_1+v_2)\|_{L^{\frac{90}{31}}_x}
\|\bar{u}_1+\bar{u}_2\|_{L^{\frac{90}{7}}_x}\|w_1\|_{L^{\frac{90}{7}}_x}\right)\|_{L^{\frac{3}{2}}_t}\nonumber\\
&\lesssim \|v_1+v_2\|^3_{X^0(I)}\|u_1+u_2\|_{X(I)}\|w_1\|_{X^0(I)}\nonumber\\
&\quad+\|v_1+v_2\|^2_{X^0(I)}\|v_1+v_2\|_{X(I)}\|u_1+u_2\|_{X^0(I)}\|w_1\|_{X^0(I)}\nonumber\\
&\lesssim \left[\|u_1+u_2\|^{\frac{8}{5}}_{X(I)}\|\nabla (u_1+u_2)\|^{\frac{12}{5}}_{S^0(I)}+\|v_1+v_2\|^{\frac{8}{5}}_{X(I)}\|\nabla (v_1+v_2)\|^{\frac{12}{5}}_{S^0(I)}\right]\|w_1\|_{X(I)}\nonumber\\
&\lesssim RHS(\ref{9191}).\label{920x2}
\end{align}
Other terms can be estimated similarly, just using H\"{o}lder's inequality and Young's inequality and (\ref{91511}). Hence (\ref{9191}) is proved when $d=3$.

If $d=4$,
\begin{align*}
&\|w\|_{X^0(I)}=\|w\|_{L^6_tL^6_x(I\times \mathbb{R}^4)},\
\|w\|_{X(I)}=\||\nabla|^{\frac{2}{3}} w\|_{L^6_tL^3_x(I\times \mathbb{R}^4)},\ \|f\|_{Y(I)}=\||\nabla|^{\frac{2}{3}} f\|_{L^2_tL^{\frac{3}{2}}_x(I\times \mathbb{R}^4)}.
\end{align*}
Similar to (\ref{9202}), we can use H\"{o}lder's and Young's inequalities, and the interpolation (\ref{91511})--(\ref{91514}) to obtain
\begin{align}
\|f_1(u,v)\|_{Y(I)}\lesssim \|u\|^3_{X(I)}+\|v\|^3_{X(I)}.\label{9202'}
\end{align}
Similar to (\ref{920x1}) and (\ref{920x2}), we can prove (\ref{9191}) when $d=4$, we omit the details here.

Putting the results above, Lemma 5.10 is proved.\hfill$\Box$

To prove the stability result, we still need the following short-time perturbation result.

{\bf Lemma 5.11(Short-time perturbation).} {\it Denote $f_1(u,v)=\lambda|v|^4u$ and $f_2(u,v)=\mu|u|^2|v|^2v$ if $d=3$, while $f_1(u,v)=\lambda|v|^2u$ and $f_2(u,v)=\mu|u|^2v$ if $d=4$. Let $I$ be a compact time interval and $(\tilde{u},\tilde{v})$ be an approximate solution of (\ref{826x1}) on $I\times \mathbb{R}^d$ in the sense that
\begin{align*}
i\tilde{u}_t+\Delta \tilde{u}=f_1(\tilde{u},\tilde{v})+e_1,\quad i\tilde{v}_t+\Delta \tilde{v}=f_2(\tilde{u},\tilde{v})+e_2
\end{align*}
for some $e_1$ and $e_2$. Suppose that
$$
\|\tilde{u}\|_{L^{\infty}_t\dot{H}^1_x(I\times \mathbb{R}^d)}+\|\tilde{v}\|_{L^{\infty}_t\dot{H}^1_x(I\times \mathbb{R}^d)}\leq E
$$
for some $E$. Let $t_0\in I$ and $(u(t_0),v(t_0))$ satisfy
$$
\|u(t_0)-\tilde{u}(t_0)\|_{\dot{H}^1_x}+\|v(t_0)-\tilde{v}(t_0)\|_{\dot{H}^1_x}\leq E'
$$
for some $E'$. Assume that there exist some small $0<\delta=\delta(E)$ and $0<\epsilon<\epsilon_0(E,E')$ such that
\begin{align}
&\|\tilde{u}\|_{X(I)}+\|\tilde{v}\|_{X(I)}\leq \delta,\label{920x3}\\
&\|e^{i(t-t_0)\Delta}[u(t_0)-\tilde{u}(t_0)]\|_{X(I)}+\|e^{i(t-t_0)\Delta}[v(t_0)-\tilde{v}(t_0)]\|_{X(I)}\leq \epsilon,\label{920x4}\\
&\|\nabla e_1\|_{N^0(I)}+\|\nabla e_2\|_{N^0(I)}\leq \epsilon.\label{920x5}
\end{align}
Then there exists a unique solution $(u,v):I\times \mathbb{R}^d\rightarrow \mathbb{C}\times \mathbb{C}$ of (\ref{826x1}) with initial data $(u(t_0),v(t_0))$
at time $t=t_0$ such that
\begin{align}
&\|u-\tilde{u}\|_{X(I)}+\|v-\tilde{v}\|_{X(I)}\lesssim\epsilon,\label{920x6}\\
&\|\nabla(u-\tilde{u})\|_{S^0(I)}+\|\nabla(v-\tilde{v})\|_{S^0(I)}\lesssim E',\label{920x7}\\
&\|\nabla u\|_{S^0(I)}+\|\nabla v\|_{S^0(I)}\lesssim E+E',\label{920x8}\\
&\|f_1(u,v)-f_1(\tilde{u},\tilde{v})\|_{Y(I)}+\|f_2(u,v)-f_2(\tilde{u},\tilde{v})\|_{Y(I)}\lesssim \epsilon,\label{920x9}\\
&\|\nabla [f_1(u,v)-f_1(\tilde{u},\tilde{v})]\|_{N^0(I)}+\|\nabla[f_2(u,v)-f_2(\tilde{u},\tilde{v})]\|_{N^0(I)}\lesssim E'.\label{920x9'}
\end{align}

}

{\bf Proof:} Inspired by the proof of Lemma 3.13 in \cite{Killip2013}, we prove it under the additional assumption $M(u)+M(v)<+\infty$,
which can guarantee that $(u,v)$ exists by the result of local well-posedness of $H^1\times H^1$-solution. The additional assumption
$M(u)+M(v)<+\infty$ can be removed a posteriori by the usual limiting argument: Approximating
$(u(t_0),v(t_0))$ in $\dot{H}^1_x\times \dot{H}^1_x$ by a sequence $\{(u_n(t_0),v_n(t_0))\}_n\subset H^1_x\times H^1_x$, and letting
$(\tilde{u}, \tilde{v})=(u_m,v_m)$, $(u,v)=(u_n,v_n)$ and $e_1=e_2=0$, we can prove that the sequence of solutions $\{(u_n,v_n)\}_n$
 with initial data $\{(u_n(t_0),v_n(t_0))\}_n$
is a Cauchy sequence in energy-critical norms and convergent to a solution $(u,v)$ with initial
data $(u(t_0),v(t_0))$ which satisfies $(\nabla u, \nabla v)\in S^0(I)\times S^0(I)$. Therefore, we assume that the solution $(u,v)$ exists and $(\nabla u, \nabla v)\in S^0(I)\times S^0(I)$, we only to prove (\ref{920x6})--(\ref{920x9}) as a priori estimates.

First, we estimate some bounds on $(\tilde{u},\tilde{v})$ and $(u,v)$. Using Strichartz, interpolation inequalities (\ref{91511})--(\ref{91514}), the assumptions (\ref{920x3}) and (\ref{920x5}),  we have
\begin{align*}
\|\nabla \tilde{u}\|_{S^0(I)}+\|\nabla \tilde{v}\|_{S^0(I)}&\lesssim \|\tilde{u}\|_{L^{\infty}_t\dot{H}^1_x}+\|\nabla f_1(\tilde{u},\tilde{v})\|_{N^0(I)}+\|\nabla e_1\|_{N^0(I)}\nonumber\\
&\quad+ \|\tilde{v}\|_{L^{\infty}_t\dot{H}^1_x}+\|\nabla f_2(\tilde{u},\tilde{v})\|_{N^0(I)}+\|\nabla e_2\|_{N^0(I)}\nonumber\\
&\lesssim E+[\|\tilde{u}\|^{\frac{4}{d-2}}_{L^{\frac{2(d+2)}{d-2}}_{t,x}}+\|\tilde{v}\|^{\frac{4}{d-2}}_{L^{\frac{2(d+2)}{d-2}}_{t,x}}][\|\nabla \tilde{u}\|_{S^0(I)}+\|\nabla \tilde{v}\|_{S^0(I)}]+\epsilon\nonumber\\
&\lesssim E+\delta^{\frac{4c}{d-2}}[\|\nabla \tilde{u}\|^{1+\frac{4(1-c)}{d-2}}_{S^0(I)}+\|\nabla \tilde{v}\|^{1+\frac{4(1-c)}{d-2}}_{S^0(I)}]+\epsilon.
\end{align*}
Here $c=c(d)$ is the same as that in interpolation inequalities (\ref{91511})--(\ref{91514}). Choosing $\delta$ small depending on $d$ and $E$, $\epsilon_0$ small enough depending on $E$, we get
\begin{align}
\|\nabla \tilde{u} \|_{S^0(I)}+\|\nabla \tilde{v}\|_{S^0(I)} \lesssim E.\label{1025w2}
\end{align}
Meanwhile, by Strichartz estimates, nonlinear estimates, (\ref{920x3}) and (\ref{920x5}), if $\delta$ and $\epsilon$ are chosen small enough, then
\begin{align*}
&\quad \|e^{i(t-t_0)\Delta}\tilde{u}(t_0)\|_{X(I)}+\|e^{i(t-t_0)\Delta}\tilde{v}(t_0)\|_{X(I)}\nonumber\\
&\lesssim \|\tilde{u}\|_{X(I)}+\|f_1(\tilde{u},\tilde{v})\|_{Y(I)}
+\|\nabla e_1\|_{N^0(I)}\nonumber\\
&\quad +\|\tilde{v}\|_{X(I)}+\|f_2(\tilde{u},\tilde{v})\|_{Y(I)}
+\|\nabla e_2\|_{N^0(I)}\nonumber\\
&\lesssim \delta+\delta^{\frac{d+2}{d-2}}+\epsilon.
\end{align*}
Recalling (\ref{920x4}), we have
\begin{align*}
\|e^{i(t-t_0)\Delta}u(t_0)\|_{X(I)}+\|e^{i(t-t_0)\Delta}v(t_0)\|_{X(I)}\lesssim \delta.
\end{align*}

Using interpolation inequalities (\ref{91511})--(\ref{91514}) and nonlinear estimates again, we get
\begin{align*}
\|u\|_{X(I)}+\|v\|_{X(I)}&\lesssim \|e^{i(t-t_0)\Delta}u(t_0)\|_{X(I)}+\|f_1(u,v)\|_{Y(I)}+\|e^{i(t-t_0)\Delta}v(t_0)\|_{X(I)}+\|f_2(u,v)\|_{Y(I)}\\
&\lesssim \delta +[\|u\|_{X(I)}+\|v\|_{X(I)}]^{\frac{d+2}{d-2}}.
\end{align*}
Choosing $\delta$ small enough, by the usual bootstrap argument, we have
\begin{align}
\|u\|_{X(I)}+\|v\|_{X(I)}&\lesssim \delta.\label{1025w3}
\end{align}

Next, let $(w,z):= (u-\tilde{u}, v-\tilde{v})$. Then
\begin{align*}
&iw_t+\Delta w=f_1(\tilde{u}+w, \tilde{v}+z)-f_1(\tilde{u},\tilde{v})-e_1,\\
&iz_t+\Delta z=f_2(\tilde{u}+w, \tilde{v}+z)-f_2(\tilde{u},\tilde{v})-e_2,\\
&w(t_0)=u(t_0)-\tilde{u}(t_0),\quad z(t_0)=v(t_0)-\tilde{v}(t_0).
\end{align*}
Using Strichartz estimates, interpolation inequalities, nonlinear estimates and (\ref{920x5}), we obtain
\begin{align}
&\quad \|w\|_{X(I)}+\|z\|_{X(I)}\nonumber\\
&\lesssim \|e^{it\Delta}[u(t_0)-\tilde{u}(t_0)]\|_{X(I)}+\|\nabla e_1\|_{N^0(I)}+\|f_1(\tilde{u}+w, \tilde{v}+z)-f_1(\tilde{u},\tilde{v})\|_{Y(I)}\nonumber\\
&\quad +\|e^{it\Delta}[v(t_0)-\tilde{v}(t_0)]\|_{X(I)}+\|\nabla e_2\|_{N^0(I)}+\|f_2(\tilde{u}+w, \tilde{v}+z)-f_2(\tilde{u},\tilde{v})\|_{Y(I)}\nonumber\\
&\lesssim \epsilon+\delta^{\frac{8}{d^2-4}}E^{\frac{4d}{d^2-4}}[\|w\|_{X(I)}+\|z\|_{X(I)}]+[\|\nabla w\|_{S^0(I)}^{\frac{4d}{d^2-4}}+\|\nabla z\|_{S^0(I)}^{\frac{4d}{d^2-4}}][\|w\|^{\frac{d^2+4}{d^2-4}}_{X(I)}+\|z\|^{\frac{d^2+4}{d^2-4}}_{X(I)}].\label{1025w1}
\end{align}
Using Strichartz estimate, nonlinear estimates, H\"{o}lder's inequality, we get
\begin{align}
&\quad \|\nabla w\|_{S^0(I)}+\|\nabla z\|_{S^0(I)}\nonumber\\
&\lesssim \|u(t_0)-\tilde{u}(t_0)\|_{\dot{H}^1_x}+\|\nabla e_1\|_{N^0(I)}+\|f_1(\tilde{u}+w, \tilde{v}+z)-f_1(\tilde{u},\tilde{v})\|_{N^0(I)}\nonumber\\
&\quad +\|v(t_0)-\tilde{v}(t_0)\|_{\dot{H}^1_x}+\|\nabla e_2\|_{N^0(I)}+\|f_2(\tilde{u}+w, \tilde{v}+z)-f_2(\tilde{u},\tilde{v})\|_{N^0(I)}\nonumber\\
&\lesssim E'+\epsilon+[\|\nabla \tilde{u}\|_{S^0(I)}+\|\nabla \tilde{v}\|_{S^0(I)}][\|u\|_{X^0(I)}+\|\tilde{u}\|_{X^0(I)}]^{\frac{6-d}{d-2}}[\|w\|_{X^0(I)}+\|z\|_{X^0(I)}]\nonumber\\
&\quad+ [\|u\|_{X^0(I)}+\|v\|_{X^0(I)}]^{\frac{4}{d-2}}[\|\nabla w\|_{S^0(I)}+\|\nabla z\|_{S^0(I)}]\nonumber\\
&\lesssim E'+\epsilon+[E\delta^{\frac{6-d}{d-2}}+\delta^{\frac{4}{d-2}}][\|\nabla w\|_{S^0(I)}+\|\nabla z\|_{S^0(I)}].\label{1025w4}
\end{align}
Combining (\ref{1025w1}) with (\ref{1025w4}), we obtain (\ref{920x6}) and (\ref{920x7}). (\ref{920x7}) and (\ref{1025w2}) imply (\ref{920x8}).
(\ref{920x6}), (\ref{920x7}), (\ref{1025w1}) and (\ref{1025w4}) deduce (\ref{920x9}) and (\ref{920x9'}).\hfill $\Box$

Now we give the proof of stability result below.

{\bf Proof of Proposition 5.8:} First, we show that
\begin{align}
\|\nabla \tilde{u}\|_{S^0(I)}+\|\nabla \tilde{v}\|_{S^0(I)}\leq C(E,L).\label{10261}
\end{align}
Under the assumption (\ref{9152}), we can divide $I$ into $J_0=J(L,\eta)$ subintervals $I_j=(t_j,t_{j+1})$ such that
$$
\|\tilde{u}\|_{L^{\frac{2(d+2)}{d-2}}_{t,x}}+\|\tilde{v}\|_{L^{\frac{2(d+2)}{d-2}}_{t,x}}\leq \eta
$$
for some small $\eta>0$ to be chosen later on each spacetime slab $I_j\times \mathbb{R}^d$. Using Strichartz estimate and (\ref{9155}),
we have
\begin{align*}
&\quad \|\nabla \tilde{u}\|_{S^0(I_j)}+\|\nabla \tilde{v}\|_{S^0(I_j)}\nonumber\\
&\lesssim \|\tilde{u}(t_j)\|_{\dot{H}^1_x}+\|\nabla e_1\|_{N^0(I)}+\|\nabla f_1(\tilde{u},\tilde{v})\|_{N^0(I)}\\
&\quad+\|\tilde{v}(t_j)\|_{\dot{H}^1_x}+\|\nabla e_2\|_{N^0(I)}+\|\nabla f_2(\tilde{u},\tilde{v})\|_{N^0(I)}\\
&\lesssim E+\epsilon+[\|\tilde{u}\|^{\frac{4}{d-2}}_{L^{\frac{2(d+2)}{d-2}}_{t,x}}+\|\tilde{v}\|^{\frac{4}{d-2}}_{L^{\frac{2(d+2)}{d-2}}_{t,x}}]
[\|\nabla \tilde{u}\|_{S^0(I_j)}+\|\nabla \tilde{v}\|_{S^0(I_j)}]\\
&\lesssim E+\epsilon+\eta^{\frac{4}{d-2}}_{L^{\frac{2(d+2)}{d-2}}_{t,x}}[\|\nabla \tilde{u}\|_{S^0(I_j)}+\|\nabla \tilde{v}\|_{S^0(I_j)}].
\end{align*}
If we choose $\eta$ and $\epsilon$ small enough depending on $E$, we can obtain
$$
 \|\nabla \tilde{u}\|_{S^0(I_j)}+\|\nabla \tilde{v}\|_{S^0(I_j)}\lesssim E.
$$
Sum them over all $I_j$, we get (\ref{10261}).

By interpolation inequalities (\ref{91511})--(\ref{91514}) and (\ref{10261}), under the assumptions (\ref{9153})--(\ref{9154}),
we have
\begin{align}
&\|\tilde{u}\|_{X(I)}+\|\tilde{v}\|_{X(I)}\leq C(E,L),\label{10262}\\
&\|e^{i(t-t_0)\Delta}[u(t_0)-\tilde{u}(t_0)]\|_{X(I)}+\|e^{i(t-t_0)\Delta}[v(t_0)-\tilde{v}(t_0)]\|_{X(I)}\lesssim \epsilon^{\frac{1}{d+2}}(E')^{\frac{d+2}{d-2}}.\label{10263}
\end{align}
(\ref{10261}) implies that we can divide $I$ into $J_1=J_1(E,L)$ subintervals $I_j=[t_j,t_{j+1}]$ such that
$$
\|\tilde{u}\|_{X(I_j)}+\|\tilde{v}\|_{X(I_j)}\lesssim \delta
$$
for the same small $\delta=\delta(E)$ in (\ref{920x3}), while (\ref{10263}) can guarantee (\ref{920x4}) if we take $\epsilon$ small enough.

By Strichartz estimate and the inductive hypothesis,
\begin{align*}
&\quad \|e^{i(t-t_j)\Delta}[u(t_j)-\tilde{u}(t_j)]\|_{X(I_j)}+\|e^{i(t-t_j)\Delta}[v(t_j)-\tilde{v}(t_j)]\|_{X(I_j)}\nonumber\\
&\lesssim \|e^{i(t-t_j)\Delta}[u(t_j)-\tilde{u}(t_j)]\|_{X(I_j)}+\|\nabla e_1\|_{N^0(I)}+\|f_1(u,v)-f_1(\tilde{u},\tilde{v})\|_{Y([t_0,t_j])}\nonumber\\
&\quad +\|e^{i(t-t_j)\Delta}[v(t_j)-\tilde{v}(t_j)]\|_{X(I_j)}+\|\nabla e_2\|_{N^0(I)}+\|f_2(u,v)-f_2(\tilde{u},\tilde{v})\|_{Y([t_0,t_j])}\nonumber\\
&\lesssim \epsilon^c+\epsilon+\sum_{k=0}^{j-1}C(k)\epsilon^c,\\
&\quad \|(u(t_j)-\tilde{u}(t_j))\|_{\dot{H}^1_x}+\|(v(t_j)-\tilde{v}(t_j))\|_{\dot{H}^1_x}\nonumber\\
&\lesssim \|(u(t_j)-\tilde{u}(t_j))\|_{\dot{H}^1_x}+\|\nabla e_1\|_{N^0([t_0,t_j])}+\|f_1(u,v)-f_1(\tilde{u},\tilde{v})\|_{N^0([t_0,t_j])}\nonumber\\
&\quad +\|(v(t_j)-\tilde{v}(t_j))\|_{\dot{H}^1_x}+\|\nabla e_2\|_{N^0([t_0,t_j])}+\|f_2(u,v)-f_2(\tilde{u},\tilde{v})\|_{N^0([t_0,t_j])}\nonumber\\
&\lesssim E'+\epsilon+\sum_{k=0}^{j-1}C(k)E'.
\end{align*}
Summing the bounds over all subintervals $I_j$, and using the interpolation inequalities (\ref{91511})--(\ref{91514}), we can obtain (\ref{920x6})--(\ref{920x9'}) and prove Proposition 5.8.\hfill $\Box$

\subsection{Concentration compactness and reduction to almost periodic solution}
\qquad In this subsection, we first give the linear profile decomposition which will lead to the reduction to almost periodic solution.

First we give the symmetry group for a pair of functions.

{\bf Definition 5.12(Symmetry group)}{\it

(i) For any $\theta\in \mathbb{R}/2\pi\mathbb{Z}$, position $x_0\in \mathbb{R}^d$ and scaling parameter $\tilde{\lambda}>0$, we define a unitary transformation
$g_{\theta,x_0,\tilde{\lambda}}(\varphi, \psi):\dot{H}^1(\mathbb{R}^d)\times \dot{H}^1(\mathbb{R}^d)\rightarrow \dot{H}^1(\mathbb{R}^d)\times \dot{H}^1(\mathbb{R}^d)$
by
$$
[g_{\theta,x_0,\tilde{\lambda}}(\varphi, \psi)](x):=\tilde{\lambda}^{-\frac{2}{\alpha+\beta+2}}e^{i\theta}(\varphi(\frac{x-x_0}{\tilde{\lambda}}), \psi(\frac{x-x_0}{\tilde{\lambda}}));
$$

(ii) Denote the collection of such transformations in (i) by $G$;

(iii) For any $(u,v):I\times \mathbb{R}^d\rightarrow \mathbb{C}\times \mathbb{C}$, we define $T_{\theta,x_0,\tilde{\lambda}}(u, v):\tilde{\lambda}^2I\times \mathbb{R}^d\rightarrow \mathbb{C}\times \mathbb{C}$ by
$$
[T_{\theta,x_0,\tilde{\lambda}}(u, v)](x):=\tilde{\lambda}^{-\frac{2}{\alpha+\beta+2}}e^{i\theta}(u(\frac{t}{\tilde{\lambda}^2},\frac{x-x_0}{\tilde{\lambda}}), v(\frac{t}{\tilde{\lambda}^2},\frac{x-x_0}{\tilde{\lambda}})).
$$
}

We state the the linear profile decomposition below, which was proved in \cite{Keraani2001,Shao20091}.

{\bf Lemma 5.13(Linear profile decomposition).} {\it Let $(u_n,v_n)_{n\geq 1}$ be a sequence of function pairs bounded in $\dot{H}^1(\mathbb{R}^d)\times \dot{H}^1(\mathbb{R}^d)$, $d\geq 3$. Then after passing to a subsequence if necessary, there exist a sequence of group elements $\{g^j_n\}_{j\geq 1,n\geq 1}\in G$, times $\{t^j_n\}_{j\geq 1,n\geq 1}\subset \mathbb{R}$ and function pairs $\{(\varphi^j,\psi^j)\}_{j\geq 1}\subset \dot{H}^1(\mathbb{R}^d)\times \dot{H}^1(\mathbb{R}^d) $ such that for all $J\geq 1$, the following decomposition is true
\begin{align*}
u_n=\sum_{j=1}^Jg_n^je^{it_n^j\Delta}\phi_1^j+\varphi_n^J,\quad v_n=\sum_{j=1}^Jg_n^je^{it_n^j\Delta}\phi_2^j+\psi_n^J.
\end{align*}
Moreover, the following properties hold:

(i) $(\varphi_n^J,\psi_n^J)\in \dot{H}^1(\mathbb{R}^d)\times \dot{H}^1(\mathbb{R}^d)$ and
\begin{align*}
\lim_{J\rightarrow +\infty}\lim_{n\rightarrow +\infty}\sup \left(\|e^{it\Delta}\varphi_n^J\|_{L^{\frac{(\alpha+\beta+2)(d+2)}{2}}_{t,x}(\mathbb{R}\times \mathbb{R}^d)}+\|e^{it\Delta}\psi_n^J\|_{L^{\frac{(\alpha+\beta+2)(d+2)}{2}}_{t,x}(\mathbb{R}\times \mathbb{R}^d)}\right)=0.
\end{align*}

(ii) For any $j\neq j'$,
\begin{align*}
\frac{\tilde{\lambda}_n^j}{\tilde{\lambda}_n^{j'}}+\frac{\tilde{\lambda}_n^{j'}}{\tilde{\lambda}_n^j}
+\frac{|x_n^j-x_n^{j'}|^2}{\tilde{\lambda}_n^j\tilde{\lambda}_n^{j'}}
+\frac{|t_n^j(\tilde{\lambda}_n^j)^2-t_n^{j'}(\tilde{\lambda}_n^{j'})^2|}{\tilde{\lambda}_n^j\tilde{\lambda}_n^{j'}}
\rightarrow +\infty\quad {\rm as}\quad n\rightarrow +\infty.
\end{align*}

(iii) For any $J\geq 1$,
\begin{align*}
\lim_{n\rightarrow +\infty}\left(\|\nabla u_n\|_2^2-\sum_{j=1}^J\|\nabla \phi_1^j\|_2^2-\|\nabla \varphi_n^J\|_2^2\right)=0,\\
\lim_{n\rightarrow +\infty}\left(\|\nabla v_n\|_2^2-\sum_{j=1}^J\|\nabla \phi_2^j\|_2^2-\|\nabla \psi_n^J\|_2^2\right)=0,
\end{align*}
and for any $1\leq j \leq J$,
\begin{align*}
e^{-t_n^j\Delta}\left((g_n^j)^{-1}(\varphi_n^J,\psi_n^J)\right)\rightarrow (0,0)\quad {\rm weakly \ in} \quad \dot{H}^1(\mathbb{R}^d)\times \dot{H}^1(\mathbb{R}^d)
\quad {\rm as}\quad n\rightarrow +\infty.
\end{align*}
}

We will prove the existence of a minimal kinetic energy blowup solution to (\ref{826x1}) below. Our idea borrows from \cite{Kenig2006, Killip20101, Killip2013}.

Define
\begin{align}
S_I(u,v)&=\int_I\int_{\mathbb{R}^d}[|u(t,x)|^{\frac{2(d+2)}{d-2}}+|v(t,x)|^{\frac{2(d+2)}{d-2}}]dxdt,\label{94w1}\\
L^+(E_0)&=\sup\{S_I(u,v)|(u,v):I\times \mathbb{R}^d\rightarrow \mathbb{C}\times \mathbb{C},\
s.t.\ \sup_{t\in I} [\|\nabla u(t)\|^2_2+\|\nabla v(t)\|^2_2]=E_0\},\label{9142}
\end{align}
where the supremum is taken over all solutions $(u,v)$ of (\ref{826x1}), $L^+(E_0):[0,+\infty)\rightarrow [0,+\infty)$ is non-decreasing and
$L^+(E_0)\lesssim E_0^{\frac{d+2}{d-2}}$ for $E_0\leq \eta_0$.

By the stability result, $L^+$ is continuous, and there exists a unique critical kinetic energy $0<E_c\leq +\infty$ such that $L^+(E_0)<+\infty$ for $E_0<E_c$ and
$L^+(E_0)=+\infty$ for $E_0\geq E_c$.

If $(u,v): I\times \mathbb{R}^d\rightarrow \mathbb{C}\times \mathbb{C}$ is a maximal-lifespan solution of (\ref{826x1}) such that
$$\sup_{t\in I}[\|\nabla u\|_2^2+\|\nabla v\|_2^2]<E^c,$$
then $(u,v)$ is global and
$$
S_{\mathbb{R}}(u,v)\leq L(\sup_{t\in I}[\|\nabla u\|_2^2+\|\nabla v\|_2^2]).
$$

We can get the following key compactness result.

{\bf Lemma 5.14(Palais-Smale condition modulo symmetries).} {\it Assume that $d=3$ or $d=4$. Let $(u_n,v_n):I_n\times\mathbb{R}^d\rightarrow \mathbb{C}\times \mathbb{C}$ be a sequence of solutions of (\ref{826x1}) with $(\alpha,\beta)=(0,2)$ or $(\alpha,\beta)=(2,0)$ when $d=3$ and $(\alpha,\beta)=(0,0)$ when $d=4$, satisfying
\begin{align}
\lim_{n\rightarrow +\infty} \sup \sup_{t\in I_n}[\|\nabla u_n\|_2^2+\|\nabla v_n\|_2^2]=E_c\label{10271}
\end{align}
and there exist some sequence of times $t_n\in I_n$ such that
\begin{align*}
\lim_{n\rightarrow +\infty}S_{\geq t_n}(u_n,v_n)=\lim_{n\rightarrow +\infty}S_{\leq t_n}(u_n,v_n)=+\infty.
\end{align*}
Then the sequence $(u_n(t_n),v_n(t_n))$ has a subsequence which converges in $\dot{H}^1(\mathbb{R}^d)\times \dot{H}^1(\mathbb{R}^d)$ modulo symmetries.
 }

{\bf Proof:} Since the equations in (\ref{826x1}) are the time--translation symmetric, we can
set $t_n=0$ for all $n\geq 1$. Hence,
\begin{align}
\lim_{n\rightarrow +\infty}S_{\geq 0}(u_n,v_n)=\lim_{n\rightarrow +\infty}S_{\leq 0}(u_n,v_n)=+\infty.\label{10281}
\end{align}
Using the linear profile decomposition in Lemma 5.13 to the sequence $(u_n(0),v_n(0))$, up to a subsequence,
we can get
\begin{align*}
u_n(0)=\sum_{j=1}^Jg_n^je^{it_n^j\Delta}\phi_1^j+\varphi_n^J,\quad v_n(0)=\sum_{j=1}^Jg_n^je^{it_n^j\Delta}\phi_2^j+\psi_n^J.
\end{align*}
Similar to the argument in the proof of Proposition 5.3 in \cite{Killip2013}, we may assume that for each $j\geq 1$, either $t_n^j\equiv 0$ or $t_n^j\rightarrow \pm \infty$ as $n\rightarrow +\infty$. Similarly, we can define the nonlinear
profiles $(u^j,v^j):I^j\times\mathbb{R}^d\rightarrow \mathbb{C}\times\mathbb{C}$ and $(u_n^j,v_n^j):I_n^j\times\mathbb{R}^d\rightarrow \mathbb{C}\times\mathbb{C}$.

As the decoupling of the kinetic energy is asymptotic, we can find $J_0\geq 1$ such
that
$$
\|\nabla \phi^j_1\|_2^2+\|\nabla \phi^j_2\|_2^2\leq \eta_0
$$
for all $j\geq J_0$ and $\eta_0=\eta_0(d)$ such that the local wellposedness of the solution to (\ref{826x1}) holds.
By Corollary 5.9, the solutions $(u_n^j,v_n^j)$ are global and
\begin{align}
\sup_{t\in \mathbb{R}}[\|\nabla u_n^j(t)\|_2^2+\|\nabla v_n^j(t)\|_2^2]+S_{\mathbb{R}}(u_n^j,v_n^j)\leq \|\nabla \phi_1\|_2^2+\|\nabla \phi_2\|_2^2.\label{10282}
\end{align}

Since the kinetic
energy is not a conserved quantity,
even if $(u_n^j(0),v_n^j(0))=g_n^j(u_n^j(t_n^j),v_n^j(t_n^j))$ posses kinetic energy less
than the critical value $E_c$, it cannot guarantee the same will hold throughout
the lifespan of $(u_n^j,v_n^j)$, especially, it cannot
make sure neither global existence nor
global spacetime bounds. Consequently, we have to seek for a profile
responsible for the asymptotic blowup (\ref{10281}) as follows.

{\bf Claim 5.15(At least one bad profile).} {\it
There exists $1\leq j_0<J_0$ such that
$$
\lim_{n\rightarrow +\infty}\sup S_{[0,\sup I_n^{j_0})}(u_n^{j_0},v_n^{j_0})=+\infty.
$$
}

{\bf Proof:} Assume contradictorily that
\begin{align}
\lim_{n\rightarrow +\infty}\sup S_{[0,\sup I_n^{j})}(u_n^{j},v_n^{j})<+\infty,\label{10283}
\end{align}
which implies that $\sup I_n^{j}=+\infty$ for all $1\leq j<j_0$ and all $n$ large enough. Using (\ref{10282}), (\ref{10283}) and (\ref{10271}),
we have for all $n$ large enough
\begin{align}
\sum_{j\geq 1}S_{[0,+\infty)}(u_n^j,v_n^j)\lesssim 1+\sum_{j\geq J_0}[\|\nabla \phi_1^j\|_2^2+\|\nabla \phi_2^j\|_2^2]\lesssim 1+E_c.\label{10284}
\end{align}

To obtain the contradiction to (\ref{10281}), we will use the stability result and (\ref{10284}) to get
a bound on the scattering size of $(u_n,v_n)$. Therefore, we define the approximate solution $(U_n^J(t),V_n^J(t))$ as
$$
U^J_n(t):=\sum^J_{j=1}u_n^j(t)+e^{it\Delta}\varphi_n^J,\quad V^J_n(t):=\sum^J_{j=1}v_n^j(t)+e^{it\Delta}\psi_n^J.
$$
Recalling (\ref{10284}) and the asymptotic vanishing of the scattering size of $(e^{it\Delta}\varphi_n^J, e^{it\Delta}\psi_n^J)$, we have
\begin{align}
&\quad \lim_{J\rightarrow J^*}\lim_{n\rightarrow} \sup S_{[0,+\infty)}(U_n^J,V_n^J)\nonumber\\
&\lesssim \lim_{J\rightarrow J^*}\lim_{n\rightarrow} \sup \left(S_{[0,+\infty)}(\sum_{j=1}^Ju_n^j,\sum_{j=1}^Jv_n^j)+S_{[0,+\infty)}(e^{it\Delta}\varphi_n^J, e^{it\Delta}\psi_n^J)\right)\nonumber\\
&\lesssim \lim_{J\rightarrow J^*}\lim_{n\rightarrow} \sup \sum_{j=1}^J S_{[0,+\infty)}(u_n^j,v_n^j)\lesssim 1+E_c.\label{1028x1}
\end{align}

Moreover, we will see that $(U_n^J,V_n^J)$ is a good approximation to $(u_n,v_n)$ below. In fact, by the way $(u_n^j,v_n^j)$ is constructed,
$$
u_n(0)-\sum_{j=1}^Ju_n^j(0)-\varphi_n^J\rightarrow 0,\quad v_n(0)-\sum_{j=1}^Jv_n^j(0)-\psi_n^J\rightarrow 0
$$
in $\dot{H}^1_x\times \dot{H}^1_x$ as $n\rightarrow +\infty$. And by the linear decomposition in Lemma 5.13 and the definition of $(U_n^J,V_n^J)$, we get

{\bf Claim 5.16(Asymptotic agreement with initial data).} {\it For any $J\geq 1$, there holds
$$
\lim_{n\rightarrow +\infty}[\|U^J_n(0)-u_n(0)\|_{\dot{H}^1_x(\mathbb{R}^d)}+\|V^J_n(0)-v_n(0)\|_{\dot{H}^1_x(\mathbb{R}^d)}]=0.
$$
}

Similar to the proof of Lemma 3.2 in \cite{Killip2013}, we also show that

{\bf Claim 5.17(Asymptotic solution to the equation).} {\it For the approximate solution $(U_n^J(t),V_n^J(t))$,
\begin{align*}
&\lim_{J\rightarrow J^*}\lim_{n\rightarrow +\infty} \sup \|\nabla [(i\partial_t+\Delta)U_n^J-f_1(U^J_n,V^J_n)]\|_{L^{\frac{2(d+2)}{d+4}}_{t,x}(\mathbb{R}\times \mathbb{R}^d)}=0\\
&\lim_{J\rightarrow J^*}\lim_{n\rightarrow +\infty} \sup\|\nabla [(i\partial_t+\Delta)V_n^J-f_2(U^J_n,V^J_n)]\|_{L^{\frac{2(d+2)}{d+4}}_{t,x}(\mathbb{R}\times \mathbb{R}^d)}]=0.
\end{align*}
}

{\bf Proof:} By the definition of $(U_n^J(t),V_n^J(t))$, we have
\begin{align*}
(i\partial_t+\Delta)U_n^J=\sum_{j=1}^Jf_1(u^j_n,v^j_n),\quad (i\partial_t+\Delta)V_n^J=\sum_{j=1}^Jf_2(u^j_n,v^j_n),
\end{align*}
by triangle inequality, it sufficient to prove that
\begin{align*}
&\lim_{J\rightarrow J^*}\lim_{n\rightarrow +\infty} \sup \|f_1(U_n^J-e^{it\Delta}\varphi_n^J, V_n^J-e^{it\Delta}\psi_n^J) -f_1(U^J_n,V^J_n)\|_{L^{\frac{2(d+2)}{d+4}}_{t,x}(\mathbb{R}\times \mathbb{R}^d)}=0,\\
&\lim_{J\rightarrow J^*}\lim_{n\rightarrow +\infty} \sup
\|f_2(U_n^J-e^{it\Delta}\varphi_n^J, V_n^J-e^{it\Delta}\psi_n^J) -f_2(U^J_n,V^J_n)\|_{L^{\frac{2(d+2)}{d+4}}_{t,x}(\mathbb{R}\times \mathbb{R}^d)}=0,\\
&\lim_{n\rightarrow +\infty}\|f_1(\sum_{j=1}^Ju_n^j, \sum_{j=1}^Jv_n^j)-\sum_{j=1}^Jf_1(u_n^j, v_n^j)\|_{L^{\frac{2(d+2)}{d+4}}_{t,x}(\mathbb{R}\times \mathbb{R}^d)}=0\quad {\rm for\ all}\quad J\geq 1,\\
&\lim_{n\rightarrow +\infty}\|f_2(\sum_{j=1}^Ju_n^j, \sum_{j=1}^Jv_n^j)-\sum_{j=1}^Jf_2(u_n^j, v_n^j)\|_{L^{\frac{2(d+2)}{d+4}}_{t,x}(\mathbb{R}\times \mathbb{R}^d)}=0\quad {\rm for\ all}\quad J\geq 1.
\end{align*}
Since
\begin{align}
&\quad |f_1(u_1,v_1)-f_1(u_2,v_2)|+|f_2(u_1,v_1)-f_2(u_2,v_2)|\\
&\lesssim [|u_1|^{\alpha+\beta+2}+|v_1|^{\alpha+\beta+2}+|u_2|^{\alpha+\beta+2}+|v_2|^{\alpha+\beta+2}]
[|u_1-u_2|+|v_1-v_2|],\label{1028x2}
\end{align}
by the asymptotically vanishing scattering size of $(e^{it\Delta}\varphi_n^J, e^{it\Delta}\psi_n^J)$ and (\ref{1028x1}), we can obtain the first limit
and the second one are zero.

Since
\begin{align*}
&\quad |f_1(\sum_{j=1}^Ju_n^j, \sum_{j=1}^Jv_n^j)-\sum_{j=1}^Jf_1(u_n^j, v_n^j)|+|f_2(\sum_{j=1}^Ju_n^j, \sum_{j=1}^Jv_n^j)-\sum_{j=1}^Jf_2(u_n^j, v_n^j)|\\
&\lesssim \sum^J_{j=1,j'=1,j\neq j'}[|u_n^j|+|v_n^j|][|u_n^{j'}|^{\alpha+\beta+2}+|v_n^{j'}|^{\alpha+\beta+2}],
\end{align*}
the third limit and the fourth one are zero by the asymptotic orthogonality of the $(u_n^j,v_n^j)$. The details are similar to the proof of (3.10) in Lemma 3.2 of \cite{Killip2013}.  \hfill $\Box$

By the two claims above and (\ref{1028x1}), we have for $n$ large enough,
$$
S_{[0,+\infty)}(u_n,v_n)\lesssim 1+E_c,
$$
which contradicts to (\ref{10281}). Claim 5.15 is proved. \hfill $\Box$

Let's come back to the proof of Lemma 5.14. Rearranging the indices, we can assume that there exists
$1\leq J_1\leq J_0$ such that
\begin{align*}
&\lim_{n\rightarrow +\infty}\sup S_{[0,\sup I_n^j)}(u^j_n,v^j_n)=+\infty \quad {\rm for}\quad 1\leq j\leq J_1,\\
&\lim_{n\rightarrow +\infty}\sup S_{[0,\sup I_n^j)}(u^j_n,v^j_n)=+\infty \quad {\rm for}\quad  j> J_1.
\end{align*}
Up to a subsequence in $n$ if necessary,  there holds $S_{[0,\sup I_n^1)}(u^1_n,v^1_n)=+\infty$.

Since there are maybe two or more profiles
that take turns at driving the scattering norm of $(u_n,v_n)$ to infinity, we have to prove that only one profile is
responsible for the asymptotic blowup (\ref{10281}). To do this, we must
prove that kinetic energy decoupling for the nonlinear profiles for large periods of time is sufficiently large such that that the kinetic energy of $(u_n^1,v^1_n)$ has achieved the critical one.

For each $m,n>1$, we can define a integer $j(m,n)\in \{1,...,J_1\}$ and a time interval $K_n^m$ of the form $[0,\tau]$ by
\begin{align}
\sup_{1\leq j\leq J} S_{K^m_n}(u_n^j,v_n^j)=S_{K^m_n}(u_n^{j(m,n)},v_n^{j(m,n)})=m.\label{1029x1}
\end{align}
Note that there is a $1\leq j_1\leq J_1$ such that for infinitely many $m$
one has $j(m,n)=j_1$ for infinitely many $n$ by the pigeonhole principle. Without loss of generality, we reorder these indices and let $j_1=1$. Then
\begin{align}
\lim_{m\rightarrow \infty}\sup \lim_{n\rightarrow \infty}\sup \sup_{t\in K^m_n} [\|\nabla u^1_n(t)\|^2_2+\|\nabla u^1_n(t)\|^2_2]\geq E_c,\label{1029x2}
\end{align}
where $E_c$ is the critical kinetic energy.

Meanwhile, (\ref{1029x1}) implies that all $(u^j_n,v^j_n)$ have finite scattering size on $K^m_n$
for each $m\geq 1$. Similar to the discussions in Claim 5.15, we can find that for
$n$ and $J$ large enough, $(U_n^J,V_n^J)$ is a good approximation to $(u_n,v_n)$ on each $K^m_n$ for each $m\geq 1$ in the following sense
\begin{align}
\lim_{J\rightarrow J^*}\lim_{n\rightarrow \infty}\sup [\|U^J_n-u_n\|_{L_t^{\infty}\dot{H}^1_x(K^m_n\times \mathbb{R}^d)}+\|V^J_n-v_n\|_{L_t^{\infty}\dot{H}^1_x(K^m_n\times \mathbb{R}^d)}=0.\label{1029x3}
\end{align}

Similar to the proof of Lemma 5.10 in \cite{Killip2013}, we can prove the following claim but we omit the details here.

{\bf Claim 5.18(Kinetic energy decoupling for $(U_n^J, V_n^J)$).} {\it For all $J\geq 1$ and $m\geq 1$,
\begin{align*}
&\quad\lim_{n\rightarrow \infty} \sup \sup_{t\in K^m_n}|\|\nabla U^J_n(t)\|^2_2-\sum^J_{j=1}\|\nabla u^j_n(t)\|^2_2-\|\nabla \varphi^J_n\|^2_2|\\
&+\lim_{n\rightarrow \infty} \sup \sup_{t\in K^m_n}|\|\nabla V^J_n(t)\|^2_2-\sum^J_{j=1}\|\nabla v^j_n(t)\|^2_2-\|\nabla \psi^J_n\|^2_2|=0.
\end{align*}
}

By this claim and using (\ref{10271}), (\ref{1029x3}), we find for each $m\geq 1$
\begin{align*}
E_c&\geq \lim_{n\rightarrow +\infty}\sup \sup_{t\in K^m_n}[\|\nabla u_n(t)\|^2_2+\|\nabla v_n(t)\|^2_2]\\
&=\lim_{J\rightarrow \infty}\lim_{n\rightarrow \infty}\sup \{\|\nabla \varphi_n^J\|^2_2+\|\nabla \psi_n^J\|^2_2+
\sup_{t\in K^m_n}\sum_{j=1}^J[\|\nabla u_n^j(t)\|^2_2+\|\nabla v_n^j(t)\|^2_2].\}
\end{align*}
Invoking (\ref{1029x2}), we can find some $g_n\in G$, $\tau_n\in\mathbb{R}$, some functions $\phi_1$, $\phi_2$, $\varphi_n$ and $\psi_n$ in
$\dot{H}^1(\mathbb{R}^d)$ with $(\varphi_n,\psi_n)\rightarrow (0,0)$ strongly in $\dot{H}^1(\mathbb{R}^d)\times \dot{H}^1(\mathbb{R}^d)$
such that
\begin{align}
u_n(0)=g_ne^{i\tau_n\Delta}\phi_1+\varphi_n,\quad v_n(0)=g_ne^{i\tau_n\Delta}\phi_2+\psi_n.\label{1029x4}
\end{align}
Moreover, similar to the proof of Proposition 5.3 in \cite{Killip2013}, we can prove $\tau_n\equiv0$. Lemma 5.14 is proved.\hfill $\Box$

{\bf Lemma 5.19(Reduction to almost periodic solution)} {\it Suppose that Theorem 6 failed when $d=3$ or Theorem 7 failed when $d=4$. Then there exists a maximal-lifespan solution $(u,v):I\times \mathbb{R}^d\rightarrow \mathbb{C}\times \mathbb{C}$ such that
$$
\sup_{t\in I} [\|\nabla u(t)\|^2_2+\|\nabla v(t)\|^2_2]=E_c,
$$
$(u,v)$ is almost periodic modulo symmetric and blows up both forward and backward in time.
}

{\bf Proof:} By the definition of the critical energy $E_c$(and the
continuity of $L$), there exist a sequence $(u_n,v_n):I_n\times \mathbb{R}^d\rightarrow \mathbb{C}\times \mathbb{C}$ of maximal-lifespan
solutions with $E(u_n,v_n)\leq E_c$ and $\lim_{n\rightarrow +\infty} S_{I_n}(u_n,v_n)=+\infty$. We can chose $t_n\in I_n$ to be
the median time of the $L^{\frac{2(d+2)}{d-2}}_{t,x}$ norm of $(u_n,v_n)$
such that (\ref{10281}) holds. Without loss of generality, we take
$t_n = 0$ because of the time-translation invariance.

Noticing Lemma 5.14 and up to a subsequence if necessary, there exist
group elements $g_n\in G$ such that $g_n(u_n(0),v_n(0))$ converges strongly in $\dot{H}^1(\mathbb{R}^d)\times \dot{H}^1(\mathbb{R}^d)$ to some
$(u_0,v_0)\in  \dot{H}^1(\mathbb{R}^d)\times \dot{H}^1(\mathbb{R}^d)$.
 By applying the group action $T_{g_n}$ to the solutions $(u_n,v_n)$, we
can take $g_n$ to all be the identity. Consequently, $(u_n(0),v_n(0))$ converge strongly in $\dot{H}^1(\mathbb{R}^d)\times \dot{H}^1(\mathbb{R}^d)$ to some
$(u_0,v_0)$, which particularly implies that $E(u_0,v_0)<E_c$.

Let $(u,v):I\times \mathbb{R}^d\rightarrow \mathbb{C}\times \mathbb{C}$ be the maximal-lifespan solution to (\ref{826x1})) with initial data
$u(0) =u_0$ and $v(0)=v_0$. We will show that $(u,v)$ blows up both forward
and backward in time. In fact, if $(u,v)$ does not blow up forward in time, then
$[0,+\infty)\subseteq I$ and $S_{\geq 0}(u,v)<\infty$. By energy-critical stability result, for $n$
large enough, we have $[0,+\infty)\subseteq I$ and
$$
\lim_{n\rightarrow +\infty} \sup S_{\geq 0}(u_n,v_n)<+\infty,
$$
which is a contradiction to (\ref{10281}). By the definition of $E_c$, it can derive that  $E(u_0,v_0)\geq E_c$ and $E(u_0,v_0)$ have to equal $E_c$.

Last, we will show that the solution u is almost periodic modulo $G$. Consider
an arbitrary sequence $(u(t'_n),v(t'_n))$ in the orbit $\{(u(t),v(t)):t\in I\}$. Although $(u,v)$ blows up both
forward and backward in time, it is locally in $L^{\frac{2(d+2)}{d-2}}_{t,x}$. Therefore,
\begin{align*}
S_{\geq t'_n}(u,v)=S_{\leq t'_n}(u,v)=\infty.
\end{align*}
By the Palais-Smale condition modulo symmetries, there exist a  a convergent subsequence of $G(u(t'_n),v(t'_n))$ in
$G\setminus [\dot{H}^1(\mathbb{R}^d)\times \dot{H}^1(\mathbb{R}^d)$. Hence, the orbit $\{G(u(t),v(t)):t\in I\}$ is precompact in $G\setminus [\dot{H}^1(\mathbb{R}^d)\times \dot{H}^1(\mathbb{R}^d)$. \hfill $\Box$

 {\bf Proof of Proposition 5.2:} Proposition 5.2 is a direct result of the conclusions above.\hfill$\Box$

Now we come back to consider (\ref{826x1}) when $d=3$ below.

\subsection{Strichartz estimates}
\qquad In this subsection, we give some Strichartz estimates, some of their proofs can be found in \cite{Colliander2008, Killip20101}.

{\bf Lemma 5.20(Strichartz inequality).} {\it Let $I$ be a compact time interval and $w:I\times \mathbb{R}^3\rightarrow \mathbb{C}$ be a solution to the forced Schr\"{o}dinger equation
$$
iw_t+\Delta w=f.
$$
Then
\begin{align}
\{\sum_{N\in 2^{\mathcal{Z}}}\|\nabla w_N\|_{L^q_tL^r_x(I\times \mathbb{R}^3)}\}^{\frac{1}{2}}\lesssim \|w(t_0)\|_{\dot{H}^1_x(\mathbb{R}^3)}+\|\nabla f\|_{L^{\gamma'}_tL^{\rho'}_x(I\times \mathbb{R}^3)}\label{93x1}
\end{align}
for any $t_0\in I$ and any admissible pairs $(q,r)$ and $(\gamma,\rho)$, where $\gamma'$ and $\rho'$ are the dual exponents to $\gamma$ and $\rho$ respectively.
}

{\bf Lemma 5.21(An endpoint estimate).} {\it
For any $w:I\times \mathbb{R}^3\rightarrow \mathbb{R}$, the following inequality hold
$$
\|w\|_{L^4_tL^{\infty}_x(I\times \mathbb{R}^3)}\lesssim \|\nabla w\|_{L^{\infty}_tL^2_x}^{\frac{1}{2}}\left\{\sum_{N\in 2^{\mathbb{Z}}}\|\nabla w_N\|_{L^2_tL^6_x(I\times \mathbb{R}^3)}\right\}^{\frac{1}{4}}.
$$
And
$$
\|w_{\leq N}\|_{L^4_tL^{\infty}_x(I\times \mathbb{R}^3)}\lesssim \|\nabla w_{\leq N}\|_{L^{\infty}_tL^2_x}^{\frac{1}{2}}\left\{\sum_{M\leq N}\|\nabla w_M\|_{L^2_tL^6_x(I\times \mathbb{R}^3)}\right\}^{\frac{1}{4}}.
$$
for any frequency $N>0$.
 }

{\bf Proposition 5.22.} {\it If $w$ satisfies $iv_t+\Delta v=f+g$ on a compact interval $[0,T]$, then for each $6<q\leq +\infty$,
$$
\|M(t)^{\frac{3}{q}-1}\|P_{M(t)}w(t)\|_{L^q_x}\|_{L^2_x}\lesssim \||\nabla|^{-\frac{1}{2}}w\|_{L^{\infty}_tL^2_x}+\||\nabla|^{-\frac{1}{2}}f\|_{L^2_tL^{\frac{6}{5}}}+\|g\|_{L^2_tL^1_x}
$$
uniformly for all functions $M:[0,T]\rightarrow 2^{\mathbb{Z}}$. Here the spacetime estimates are over $[0,T]\times \mathbb{R}^3$.
}

{\bf Proposition 5.23.} {\it If $w$ satisfies $iw_t+\Delta w=f+g$ on a compact interval $[0,T]$ and
\begin{align}
[\mathcal{S}_Rw](t,x):=\left((\pi R^2)^{-\frac{3}{2}}\int_{\mathbb{R}^3}|v(t,x+y)|^2e^{-\frac{|y|^2}{R^2}}dy\right)^{\frac{1}{2}},\label{93w1}
\end{align}
then for each $0<R<+\infty$ and $6<q\leq +\infty$,
\begin{align}
R^{\frac{1}{2}-\frac{3}{q}}\|\mathcal{S}_Rw\|_{L^2_tL^q_x}\leq \|w\|_{L^{\infty}_tL^2_x}+R^{-\frac{1}{2}}\|f\|_{L^2_tL^1_x}+\|g\|_{L^2_tL^{\frac{6}{5}}_x},\label{93w2}
\end{align}
where the spacetime norms are over $[0,T]\times \mathbb{R}^3$.
}

{\bf Lemma 5.24.} {\it
\begin{align}
\sup_{M>0}M^{\frac{3}{q}-1}\|f_M\|_{L^q_x}\lesssim \sup_{M>0}M^{\frac{3}{q}-1}\|\mathcal{S}_{M^{-1}}(f_M)\|_{L^q_x}\quad {\rm for \ fixed}\quad  6<q\leq +\infty.\label{93w3}
\end{align}
}

{\bf Lemma 5.25.} {\it The intergral kernel
$$
K_R(\tau,z;s,y;x):=(\pi R^2)^{-\frac{3}{2}}<\delta_z, e^{i\tau\Delta} e^{\frac{-|\cdot-x|^2}{R^2}}e^{is\Delta}\delta_y>
$$
satisfies
\begin{align}
\sup_{R>0}\int_0^{\infty}\int_0^{\infty}R^{2-\frac{6}{q}}\|K_R(\tau,z;s,y;x)\|_{L^{\infty}_{z,y}L^{\frac{q}{2}}_x}f(t+\tau)f(t-s)dsd\tau
\lesssim|[\mathcal{M}f](t)|^2,\label{11151}
\end{align}
for fixed $6<q\leq +\infty$, where $\mathcal{M}$ denotes the Hardy-Littlewood maximal operator and $f:\mathbb{R}\rightarrow [0,+\infty)$.
}

{\bf Proposition 5.26(Long-time Strichartz estimate).} {\it Assume that $(u,v):(-T_{\min},T_{\max})\times \mathbb{R}^3\rightarrow \mathbb{C}\times \mathbb{C}$
is a maximal-lifespan almost periodic solution to (\ref{826x1}) with $(\alpha,\beta)=(0,2)$ and $I\subset (-T_{\min},T_{\max})$ a time interval which is titled by finitely many characteristic interval $J_k$. For any fixed $6<q<+\infty$ and any frequency $N>0$, define
\begin{align}
&K:=\int_I\widetilde{N}(t)^{-1}dt,\label{94x2}\\
&A(N):=\{\sum_{M\leq N} \|\nabla u_M\|_{L^2_tL^6_x(I\times \mathbb{R}^3)}+\|\nabla v_M\|_{L^2_tL^6_x(I\times \mathbb{R}^3)}\}^{\frac{1}{2}},\label{941}\\
&\tilde{A}_q(N):=N^{\frac{3}{2}}\|\sup_{M\geq N} M^{\frac{3}{q}-1}[\|u_M(t)\|_{L^q_x(\mathbb{R}^3)}+\|v_M(t)\|_{L^q_x(\mathbb{R}^3)}]\|_{L^2_t(I)}.\label{94x1}
\end{align}
Then
\begin{align}
A(N)+\tilde{A}(N)\lesssim_{(u,v)} 1+N^{\frac{3}{2}}K^{\frac{1}{2}},\label{94x3}
\end{align}
where the implicit constant is independent of the interval $I$.}

The proof of Proposition 5.26 is based on the following facts and Proposition 5.27 below.

By Sobolev embedding, for a small parameter $\eta>0$ to be chosen later, there exists such that $c=c(\eta)$
\begin{align}
\|u_{\leq c\widetilde{N}(t)}\|_{L^{\infty}_tL^6_x}+\|v_{\leq c\widetilde{N}(t)}\|_{L^{\infty}_tL^6_x}
+\|\nabla u_{\leq c\widetilde{N}(t)}\|_{L^{\infty}_tL^2_x}+\|\nabla v_{\leq c\widetilde{N}(t)}\|_{L^{\infty}_tL^2_x}\leq \eta.\label{94x4}
\end{align}

By the results of Lemma 5.24 and 5.25, we have
\begin{align}
&\|\nabla u_{\leq N}\|_{L^2_tL^6_x}+\|\nabla v_{\leq N}\|_{L^2_tL^6_x}\lesssim A(N),\label{98xj1}\\
&\| u_{\leq N}\|_{L^4_tL^{\infty}_x}+\|v_{\leq N}\|_{L^4_tL^{\infty}_x}\lesssim A(N)^{\frac{1}{2}}[\|\nabla u_{\leq N}\|_{L^{\infty}_tL^2_x}+\|\nabla v_{\leq N}\|_{L^{\infty}_tL^2_x}]^{\frac{1}{2}}\lesssim_{(u,v)}A(N)^{\frac{1}{2}}.\label{98xj2}
\end{align}

By the results of Lemma 5.5 and Lemma 5.20, we know that
\begin{align}
[A(N)]^2\lesssim_{(u,v)} 1+\int_I[\widetilde{N}(t)]^2dt\lesssim_{(u,v)} \int_I[\widetilde{N}(t)]^2dt.\label{1219x1}
\end{align}
While by the results of Proposition 5.22 and using Berstein's estimates, we get
\begin{align}
\tilde{A}_q(N)&\lesssim N^{\frac{3}{2}}\left(\||\nabla|^{-\frac{1}{2}}u_{\geq N}\|_{L^{\infty}_tL^2_x}+\||\nabla|^{-\frac{1}{2}}P_{\geq N}(|v|^4u)\|_{L^2_tL^{\frac{6}{5}}_x}\right)\nonumber\\
&\quad+\left(\||\nabla|^{-\frac{1}{2}}v_{\geq N}\|_{L^{\infty}_tL^2_x}+\||\nabla|^{-\frac{1}{2}}P_{\geq N}(|u|^2|v|^2v)\|_{L^2_tL^{\frac{6}{5}}_x}\right)\nonumber\\
&\lesssim 1+[\|\nabla u\|_{L^2_tL^6_x}+\|\nabla v\|_{L^2_tL^6_x}][\|u\|^4_{L^{\infty}_tL^6_x}+\|v\|^4_{L^{\infty}_tL^6_x}]
\lesssim_{(u,v)} \left(\int_I[\widetilde{N}(t)]^2dt\right)^{\frac{1}{2}}.\label{12192}
\end{align}
Consequently,
\begin{align}
A(N)+\tilde{A}_q(N)\lesssim_{(u,v)} N^{\frac{3}{2}}K^{\frac{1}{2}}\quad {\rm whenever}\quad N\geq \left(\frac{\int_I[\widetilde{N}(t)]^2dt}{\int_I[\widetilde{N}]^{-1}dt}\right)^{\frac{1}{3}}.\label{1219x3}
\end{align}

{\bf Proposition 5.27(Recurrence relations for $A(N)$ and $\tilde{A}_q(N)$).} {\it For $\eta$ sufficient small and $c=c(\eta)$ is as in (\ref{94x4}),
\begin{align}
&A(N)\lesssim_{(u,v)} 1+c^{-\frac{3}{2}}N^{\frac{3}{2}}K^{\frac{1}{2}}+\eta^2\tilde{A}_q(2N),\label{94x5}\\
&\tilde{A}_q(N)\lesssim_{(u,v)} 1+c^{-\frac{3}{2}}N^{\frac{3}{2}}K^{\frac{1}{2}}+\eta A(N)+\eta^2\tilde{A}_q(2N),\label{94x5'}
\end{align}
uniformly in $N\in 2^{\mathbb{Z}}$.
}

{\bf Proof:} The recurrence relations for $A(N)$ and $\tilde{A}_q(N)$ rely on the Strichartz inequality in Lemma 5.20 and maximal Strichartz estimate in Proposition 5.22.
To deal with the contributions of the nonlinearities $|v|^4u$ and $|u|^2|v|^2v$, we use the notation $\varnothing(X)$ to denote a quantity that resembles $X$, and write
\begin{align}
&|v|^4u=\varnothing(v^2_{>c\widetilde{N}(t)}v^2u)+\varnothing(v^2_{\leq c\widetilde{N}(t)}v^2_{>N}u)
+\varnothing(v^2_{\leq c\widetilde{N}(t)}v^2_{\leq N}u),\label{94x6}\\
&|u|^2|v|^2v=\varnothing(u^2_{>c\widetilde{N}(t)}v^3)+\varnothing(u^2_{\leq c\widetilde{N}(t)}v^2_{>N}v)
+\varnothing(u^2_{\leq c\widetilde{N}(t)}v^2_{\leq N}u).\label{94x7}
\end{align}
Using the Strichartz inequality in Lemma 5.20, by(\ref{94x6}), (\ref{94x7}) and Bernstein's inequality, we have
\begin{align}
A(N)&\lesssim \|\nabla u_{\leq N}\|_{L^{\infty}_tL^2_x}+\|\nabla v_{\leq N}\|_{L^{\infty}_tL^2_x}+\|\nabla P_{\leq N}\varnothing(v^2_{>c\widetilde{N}(t)}v^2u) \|_{L^2_tL^{\frac{6}{5}}_x}\nonumber\\
&\quad+\|\nabla P_{\leq N}\varnothing(v^2_{\leq c\widetilde{N}(t)}v^2_{>N}u \|_{L^2_tL^{\frac{6}{5}}_x}
+\|\nabla P_{\leq N}\varnothing(v^2_{\leq c\widetilde{N}(t)}v^2_{\leq N}u)\|_{L^2_tL^{\frac{6}{5}}_x}\nonumber\\
&\quad +\|\nabla P_{\leq N}\varnothing(u^2_{>c\widetilde{N}(t)}v^3)\|_{L^2_tL^{\frac{6}{5}}_x}+\|\nabla P_{\leq N}\varnothing(u^2_{\leq c\widetilde{N}(t)}v^2_{>N}v)\|_{L^2_tL^{\frac{6}{5}}_x}\nonumber\\
&\quad +\|\nabla P_{\leq N}\varnothing(u^2_{\leq c\widetilde{N}(t)}v^2_{\leq N}v)\|_{L^2_tL^{\frac{6}{5}}_x}\nonumber\\
&\lesssim_{(u,v)} 1+N^{\frac{3}{2}}\|v^2_{>c\widetilde{N}(t)}v^2u\|_{L^2_tL^1_x}
+N^{\frac{3}{2}}\|v^2_{\leq c\widetilde{N}(t)}v^2_{>N}u\|_{L^2_tL^1_x} \nonumber\\
&\quad +\|\nabla \varnothing(v^2_{\leq c\widetilde{N}(t)}v^2_{\leq N}u)\|_{L^2_tL^{\frac{6}{5}}_x}+N^{\frac{3}{2}}\|u^2_{>c\widetilde{N}(t)}v^3\|_{L^2_tL^1_x}\nonumber\\
&\quad +N^{\frac{3}{2}}\|u^2_{\leq c\widetilde{N}(t)}v^2_{>N}v) \|_{L^2_tL^1_x}+\|\nabla \varnothing(u^2_{\leq c\widetilde{N}(t)}v^2_{\leq N}v)\|_{L^2_tL^{\frac{6}{5}}_x}\nonumber\\
&:=1+(I)+(II)+(III)+(IV)+(V)+(VI).\label{94x8}
\end{align}
Using maximal Strichartz estimate in Proposition 5.22, by (\ref{94x6}), (\ref{94x7}) and Bernstein's inequality, we get
\begin{align}
\tilde{A}_q(N)&\lesssim \||\nabla|^{-\frac{1}{2}} u_{\geq N}\|_{L^{\infty}_tL^2_x}+\||\nabla|^{-\frac{1}{2}} v_{\geq N}\|_{L^{\infty}_tL^2_x}+\|\varnothing(v^2_{>c\widetilde{N}(t)}v^2u) \|_{L^2_tL^1_x}\nonumber\\
&\quad+\|\varnothing(v^2_{\leq c\widetilde{N}(t)}v^2_{>N}u) \|_{L^2_tL^1_x}
+\||\nabla|^{-\frac{1}{2}} P_{\geq N}\varnothing(v^2_{\leq c\widetilde{N}(t)}v^2_{\leq N}u)\|_{L^2_tL^{\frac{6}{5}}_x}\nonumber\\
&\quad +\|\varnothing(u^2_{>c\widetilde{N}(t)}v^3)\|_{L^2_tL^1_x}+\|\varnothing(u^2_{\leq c\widetilde{N}(t)}v^2_{>N}v)\|_{L^2_tL^1_x}\nonumber\\
&\quad +\||\nabla|^{-\frac{1}{2}} P_{\geq N}\varnothing(u^2_{\leq c\widetilde{N}(t)}v^2_{\leq N}v)\|_{L^2_tL^{\frac{6}{5}}_x}\nonumber\\
&\lesssim_{(u,v)} 1+N^{\frac{3}{2}}\|v^2_{>c\widetilde{N}(t)}v^2u\|_{L^2_tL^1_x}
+N^{\frac{3}{2}}\|v^2_{\leq c\widetilde{N}(t)}v^2_{>N}u\|_{L^2_tL^1_x} \nonumber\\
&\quad +\|\nabla \varnothing(v^2_{\leq c\widetilde{N}(t)}v^2_{\leq N}u)\|_{L^2_tL^{\frac{6}{5}}_x}+N^{\frac{3}{2}}\|u^2_{>c\widetilde{N}(t)}v^3\|_{L^2_tL^1_x}\nonumber\\
&\quad +N^{\frac{3}{2}}\|u^2_{\leq c\widetilde{N}(t)}v^2_{>N}v) \|_{L^2_tL^1_x}+\|\nabla \varnothing(u^2_{\leq c\widetilde{N}(t)}v^2_{\leq N}v)\|_{L^2_tL^{\frac{6}{5}}_x}\nonumber\\
&:=1+(I)+(II)+(III)+(IV)+(V)+(VI).\label{94w1}
\end{align}

Before we consider these terms individually, we would like to point out that the terms (I)--(VI) in (\ref{94x8}) are the same as those in (\ref{94w1}), and we will show some facts below.

Fact 1: By Sobolev embedding,
\begin{align}
\|u(t,x)\|_{L^6_x}\lesssim \|u(t,x)\|_{\dot{H}^1_x},\quad \|v(t,x)\|_{L^6_x}\lesssim \|v(t,x)\|_{\dot{H}^1_x}\quad {\rm for\ all} \quad t\in I.\label{94w2}
\end{align}

Fact 2: Using Young's inequality, we have
\begin{align}
|u||v|^2\lesssim |u|^3+|v|^3.\label{94w3}
\end{align}

In order to estimate (I), we decompose the time interval $I$ into characteristic subintervals $J_k$ where $\widetilde{N}(t)\equiv \widetilde{N}_k$. On each these
subintervals, using H\"{o}lder's inequality, Bernstein's inequality, and Strichartz estimate for the admissible pair $(4,3)$, we obtain
\begin{align*}
&\quad \|v^2_{>c\widetilde{N}(t)}v^2u\|_{L^2_tL^1_x(J_k\times \mathbb{R}^3}\lesssim \|v^2_{>c\widetilde{N}(t)}(|u|^3+|v|^3)\|_{L^2_tL^1_x(J_k\times \mathbb{R}^3}\nonumber\\
&\lesssim \|v^2_{>c\widetilde{N}(t)}\|^2_{L^4_{t,x}(J_k\times \mathbb{R}^3)}[\|u\|^3_{L^{\infty}_tL^6_x}+\|v\|^3_{L^{\infty}_tL^6_x}]\nonumber\\
&\lesssim_{(u,v)} c^{-\frac{3}{2}}N^{-\frac{3}{2}}_k\|\nabla v^2_{>c\widetilde{N}(t)}\|^2_{L^4_tL^3_x(J_k\times \mathbb{R}^3)}\lesssim_{(u,v)} c^{-\frac{3}{2}}N^{-\frac{3}{2}}_k.
\end{align*}
Recalling that $\widetilde{N}(t)\equiv N_k$ on $J_k$ and $|J_k|\sim N_k^{-2}$, we can square and sum the estimates above over the subintervals $J_k$, and obtain
\begin{align}
(I)=N^{\frac{3}{2}}\|v^2_{>c\widetilde{N}(t)}v^2u\|_{L^2_tL^1_x}\lesssim_{(u,v)} c^{-\frac{3}{2}}N^{-\frac{3}{2}}K^{\frac{1}{2}}.\label{952}
\end{align}
Similarly,
\begin{align}
(IV)=N^{\frac{3}{2}}\|u^2_{>c\widetilde{N}(t)}v^3\|_{L^2_tL^1_x}\lesssim_{(u,v)} c^{-\frac{3}{2}}N^{-\frac{3}{2}}K^{\frac{1}{2}}.\label{953}
\end{align}

To estimate (II), we recall (\ref{94w3}) and point out
\begin{align}
|v^2_{\leq c\widetilde{N}(t)}v^2_{>N}u|\leq |v^2_{\leq c\widetilde{N}(t)}v_{>N}v_{>N}u|+|v^2_{\leq c\widetilde{N}(t)}v_{>N}v_{\leq N}u|:=(II1)+(II2).\label{95w1}
\end{align}
Using Berstein's inequality and Schur's test, for $6<q<\infty$, we can obtain
\begin{align}
\|v_{>N}v_{>N}u\|_{L^2_tL^{\frac{3}{2}}_x}&\lesssim_{(u,v)} \|\sum_{M_1\geq M_2\geq M_3, M_2>N}\|v_{M_1}\|_{L^2_x}\|v_{M_2}\|_{L^q_x}\|u_{M_3}\|_{L^{\frac{6q}{q-6}}_x}\|_{L^2_t}\nonumber\\
&\lesssim_{(u,v)} \|\sup_{M>N}\|M^{\frac{3}{q}-1}v_M(t)\|_{L^q_x}\sum_{M_1\geq M_3}(\frac{M_3}{M_1})^{\frac{3}{q}}\|\nabla v_{M_1}(t)\|_{L^2_x}\|\nabla u_{M_3}(t)\|_{L^2_x}\|_{L^2_t}\nonumber\\
&\lesssim_{(u,v)}N^{-\frac{3}{2}}\tilde{A}_q(2N).\label{95w2}
\end{align}
Using H\"{o}lder's inequality, by (\ref{94x4}) and (\ref{95w2}), we have
\begin{align}
(II1)\lesssim_{(u,v)} \|v_{\leq c\widetilde{N}(t)}\|^2_{L^{\infty}_2L^6_x}\|v_{>N}v_{>N}u\|_{L^2_tL^{\frac{3}{2}}_x}\lesssim_{(u,v)}N^{-\frac{3}{2}}\eta^2\tilde{A}_q(2N).\label{95w3}
\end{align}
About (II2), using H\"{o}lder's inequality and Berstein's inequality, we get
\begin{align}
\|v_{>N}v_{\leq N}u\|_{L^2_tL^{\frac{3}{2}}_x}&\leq\|\left(\int_{\mathbb{R}^3}|v_{>N}|^{\frac{3}{2}}[|v_{\leq N}|^3+|u|^3]dx\right)^{\frac{2}{3}}\|_{L^2_t}\nonumber\\
&\lesssim_{(u,v)} \|\left(\int_{\mathbb{R}^3}|v_{>N}|^3dx\right)^{\frac{1}{3}}\left(\int_{\mathbb{R}^3}[|v_{\leq N}|^6+|u|^6]dx\right)^{\frac{1}{3}}\|_{L^2_t}\nonumber\\
&\lesssim_{(u,v)}  \|\left(\int_{\mathbb{R}^3}|v_{>N}|^3dx\right)^{\frac{1}{3}}\|_{L^2_t}\|\left(\int_{\mathbb{R}^3}[|v_{\leq N}|^6+|u|^6]dx\right)^{\frac{1}{3}}\|_{L^{\infty}_t}\nonumber\\
&\lesssim_{(u,v)}\|\sup_{M> N}M^{\frac{3}{q}-1}\|u_M(t)\|_{L^q_x(\mathbb{R}^3)}\|_{L^2_t}\lesssim_{(u,v)} N^{-\frac{3}{2}}\tilde{A}_q(2N).\label{95w5}
\end{align}
Consequently, we get
\begin{align}
(II)\lesssim_{(u,v)}\eta^2\tilde{A}_q(2N).\label{95w8}
\end{align}
Similarly, we have
\begin{align}
(V)\lesssim_{(u,v)}\eta^2\tilde{A}_q(2N).\label{95w9}
\end{align}

Last, by (\ref{94x4}), (\ref{98xj1}) and (\ref{98xj2}), we have
\begin{align}
(III)&=\|\nabla\varnothing(v^2_{\leq c\widetilde{N}(t)}v^2_{\leq N}u)\|_{L^2_tL^{\frac{6}{5}}_x}\nonumber\\
&\lesssim \|\nabla v_{\leq N}\|_{L^2_tL^6_x}\|v_{\leq c\widetilde{N}(t)}\|_{L^{\infty}_tL^6_x}[\|u\|^3_{L^{\infty}_tL^6_x}+\|v\|^3_{L^{\infty}_tL^6_x}]\nonumber\\
&\quad+\|\nabla v\|_{L^{\infty}_tL^2_x}\|v_{\leq c\widetilde{N}(t)}\|_{L^{\infty}_tL^6_x}\|v_{\leq N}\|^2_{L^4_tL^{\infty}_x}[\|u\|_{L^{\infty}_tL^6_x}+\|v\|_{L^{\infty}_tL^6_x}]\nonumber\\
&\quad+ \|\nabla v_{\leq c\widetilde{N}(t)}\|_{L^{\infty}_tL^2_x}\|v_{\leq c\widetilde{N}(t)}\|_{L^{\infty}_tL^6_x}\|v_{\leq N}\|^2_{L^4_tL^{\infty}_x}\|u\|_{L^{\infty}_tL^6_x}\nonumber\\
&\lesssim_{(u,v)} \eta A(N).\label{11161}
\end{align}
Similarly,
\begin{align}
(IV)\lesssim_{(u,v)} \eta A(N).\label{11162}
\end{align}
Putting (\ref{952})--(\ref{11162}) together, we can get (\ref{94x5}) and (\ref{94x5'}).\hfill $\Box$

{\bf Proof of Proposition 5.26:} Taking $\eta$ small enough, using (\ref{1219x3}) and the results of Proposition 5.27, after a straightforward induction argument, we can obtain (\ref{94x3}). \hfill $\Box$

\subsection{Impossibility of rapid frequency cascades}
\qquad In this subsection, we prove that the first type of almost periodic solution described in Proposition 5.6 which satisfies $\int_0^{T_{\max}}\widetilde{N}(t)^{-1}dt<+\infty$ cannot exist.

{\bf Lemma 5.28(Finite mass).} {\it Assume that $(u,v):[0,T_{\max})\times \mathbb{R}^3\rightarrow \mathbb{C}\times \mathbb{C}$ is an almost periodic solution of (\ref{826x1}) with $(\alpha,\beta)=(0,2)$ satisfying
\begin{align}
\|u\|_{L^{10}_{t,x}([0,T_{\max})\times \mathbb{R}^3)}+\|v\|_{L^{10}_{t,x}([0,T_{\max})\times \mathbb{R}^3)}=+\infty\quad {\rm and} \quad K:=\int_0^{T_{\max}}N(t)^{-1}dt<+\infty.\label{961}
\end{align}
Then $(u,v)\in [L^{\infty}_tL^2_x]\times [L^{\infty}_tL^2_x]$ and for all $0<N<1$,
\begin{align}
&\quad \|u_{N\leq \cdot \leq 1}\|_{L^{\infty}_tL^2_x([0,T_{\max})\times \mathbb{R}^3)}+\|v_{N\leq \cdot \leq 1}\|_{L^{\infty}_tL^2_x([0,T_{\max})\times \mathbb{R}^3)}\nonumber\\
&+\frac{1}{N}\left\{\sum_{M<N}\|\nabla u_M\|_{L^2_tL^6_x([0,T_{\max})\times \mathbb{R}^3)}+\|\nabla v_M\|_{L^2_tL^6_x([0,T_{\max})\times \mathbb{R}^3)}\right\}^{\frac{1}{2}}\nonumber\\
&\lesssim_{(u,v)} 1.\label{962}
\end{align}
}

{\bf Proof:} First, we would like to point out that if (\ref{962}) holds, then letting $N\rightarrow 0$ in (\ref{962}), using $\nabla u, \nabla v\in L^{\infty}_tL^2_x$ and Berstein's inequality for high frequencies, we have
\begin{align}
\|u\|_{L^{\infty}_tL^2_x}+\|v\|_{L^{\infty}_tL^2_x}\leq \|u_{\leq 1}\|_{L^{\infty}_tL^2_x}+\|u_{>1}\|_{L^{\infty}_tL^2_x}+\|v_{\leq 1}\|_{L^{\infty}_tL^2_x}+\|v_{>1}\|_{L^{\infty}_tL^2_x}\lesssim_{(u,v)} 1,\label{963}
\end{align}
which implies the finiteness of the mass.

Since $K$ is finite, we extend the conclusions of Proposition 5.26 and Proposition 5.27 to the time interval $[0,T_{\max})$ and find that
the second summand in (\ref{962}) is $N^{-1}A(\frac{N}{2})$. Therefore, by Proposition 5.26 and Proposition 5.27 and Berstein's inequality,
\begin{align}
LHS(\ref{962})&\lesssim_{(u,v)} N^{-1}[\|\nabla u\|_{L^{\infty}_tL^2_x}+\|\nabla v\|_{L^{\infty}_tL^2_x}+A(\frac{N}{2})]\nonumber\\
&\lesssim_{(u,v)}
N^{-1}(1+N^3K)^{\frac{1}{2}}<+\infty.\label{964}
\end{align}

Fixing $0<N<1$, using Duhamel formulae and Strichartz's inequality, we have
\begin{align}
LHS(\ref{962})&\lesssim_{(u,v)} \frac{1}{N}\|\nabla P_{<N}(|v|^4u)\|_{L^2_tL^{\frac{6}{5}}_x}+\|P_{N\leq \cdot \leq 1}(|v|^4u)\|_{L^2_tL^{\frac{6}{5}}_x}\nonumber\\
&\qquad+\frac{1}{N}\|\nabla P_{<N}(|u|^2|v|^2v)\|_{L^2_tL^{\frac{6}{5}}_x}+\|P_{N\leq \cdot \leq 1}(|u|^2|v|^2v)\|_{L^2_tL^{\frac{6}{5}}_x}\nonumber\\
&:=(I)+(II)+(III)+(IV).\label{965}
\end{align}
Before we estimate the contributions of (I)--(IV) to (\ref{965}), we would like to point that
\begin{align}
||v|^4u|+||u|^2|v|^2u|&\lesssim |\varnothing(v^2_{>c\widetilde{N}(t)}v^2u)| +|\varnothing(v_{\leq c\widetilde{N}(t)}v^2_{<N}u_{<N}v_{>c\widetilde{N}(t)})|\nonumber\\
&\quad
+|\varnothing(v_{\leq c\widetilde{N}(t)}v^2_{<N}u_{>1}v_{>c\widetilde{N}(t)})|+|\varnothing(v_{\leq c\widetilde{N}(t)}v^2_{<N}u_{N\leq \cdot\leq 1}v_{>c\widetilde{N}(t)})|\nonumber\\
&\quad+|\varnothing(v_{\leq c\widetilde{N}(t)}v^2_{\leq 1}u_{N\leq \cdot\leq 1}v_{>c\widetilde{N}(t)})| +|\varnothing(v_{\leq c\widetilde{N}(t)}v^2_{>1}u_{\leq 1}v_{>c\widetilde{N}(t)})|\nonumber\\
&:=(i)+(ii1)+(ii2)+(iii1)+(iii2)+(iv).\label{971}
\end{align}

Using Bernstein's inequality and H\"{o}lder's inequality, similar to (\ref{952}), we have
\begin{align}
&\quad \frac{1}{N}\|\nabla P_{<N}\varnothing(v^2_{>c\widetilde{N}(t)}v^2u)\|_{L^2_tL^{\frac{6}{5}}_x}+\|P_{N\leq \cdot \leq 1}\varnothing(v^2_{>c\widetilde{N}(t)}v^2u)\|_{L^2_tL^{\frac{6}{5}}_x}\nonumber\\
&\lesssim_{(u,v)} (N^{\frac{1}{2}}+1)\|\varnothing(v^2_{>c\widetilde{N}(t)}v^2u)\|_{L^2_tL^1_x}\lesssim_{(u,v)} c^{-\frac{3}{2}}K^{\frac{1}{2}}.\label{972}
\end{align}

Similarly, we consider the contribution of (ii1) and get
\begin{align}
&\quad \frac{1}{N}\|\nabla P_{<N}\varnothing(v_{\leq c\widetilde{N}(t)}v^2_{<N}u_{<N}v_{>c\widetilde{N}(t)})\|_{L^2_tL^{\frac{6}{5}}_x}+\|P_{N\leq \cdot \leq 1}\varnothing(v_{\leq c\widetilde{N}(t)}v^2_{<N}u_{<N}v_{>c\widetilde{N}(t)})\|_{L^2_tL^{\frac{6}{5}}_x}\nonumber\\
&\lesssim_{(u,v)} \frac{1}{N}\|v_{\leq c\widetilde{N}(t)}\|_{L^{\infty}_tL^2_x}\|v_{<N}\|^2_{L^4_tL^{\infty}_x}[\|u_{<N}\|^2_{L^{\infty}_tL^6_x}+\|v_{>c\widetilde{N}(t)}\|^2_{L^{\infty}_tL^6_x}]\nonumber\\
&\quad+\frac{1}{N}\|v_{\leq c\widetilde{N}(t)}\|_{L^{\infty}_tL^6_x}\|\nabla v_{<N}\|_{L^2_tL^6_x}[\|u_{<N}\|^3_{L^{\infty}_tL^6_x}+\|v_{>c\widetilde{N}(t)}\|^3_{L^{\infty}_tL^6_x}]\nonumber\\
&\quad+ \frac{1}{N}\|v_{\leq c\widetilde{N}(t)}\|_{L^{\infty}_tL^6_x}\|v_{<N}\|^2_{L^4_tL^{\infty}_x}[\|u_{<N}\|_{L^{\infty}_tL^6_x}+\|v_{>c\widetilde{N}(t)}\|_{L^{\infty}_tL^6_x}]\nonumber\\
&\qquad \qquad\qquad\qquad \qquad\qquad\quad \times [\|\nabla u_{<N}\|_{L^{\infty}_tL^2_x}+\|\nabla v_{>c\widetilde{N}(t)}\|_{L^{\infty}_tL^2_x}]\nonumber\\
&\lesssim_{(u,v)} \eta LHS(\ref{962}).\label{972}
\end{align}
The contribution of (ii2) is entirely similar to that of (ii1).

Using Berstein's inequality and H\"{o}lder's inequality again, we consider the contribution of (iii1) and get
\begin{align}
&\frac{1}{N}\|\nabla P_{<N}\varnothing(v_{\leq c\widetilde{N}(t)}v^2_{<N}u_{N\leq \cdot\leq 1}v_{>c\widetilde{N}(t)})\|_{L^2_tL^{\frac{6}{5}}_x}+\|P_{N\leq \cdot \leq 1}\varnothing(v_{\leq c\widetilde{N}(t)}v^2_{<N}u_{N\leq \cdot\leq 1}v_{>c\widetilde{N}(t)})\|_{L^2_tL^{\frac{6}{5}}_x}\nonumber\\
&\lesssim_{(u,v)}\|v_{\leq c\widetilde{N}(t)}\|_{L^{\infty}_tL^6_x}\|v_{\leq 1}\|^2_{L^4_tL^{\infty}_x}\|v_{N\leq \cdot\leq 1}\|_{L^{\infty}_tL^2_x}\|v_{>c\widetilde{N}(t)}\|_{L^{\infty}_tL^6_x}\nonumber\\
&\lesssim_{(u,v)}\eta(1+K^{\frac{1}{2}})LHS(\ref{962}).\label{973}
\end{align}
The contribution of (iii2) is entirely similar to that of (iii1).

Similarly, we consider the contribution of (iv) and get
\begin{align}
&\quad \frac{1}{N}\|\nabla P_{<N}\varnothing(v_{\leq c\widetilde{N}(t)}v^2_{>1}u_{\leq 1}v_{>c\widetilde{N}(t)})\|_{L^2_tL^{\frac{6}{5}}_x}
+\|P_{N\leq \cdot \leq 1}\varnothing(v_{\leq c\widetilde{N}(t)}v^2_{>1}u_{\leq 1}v_{>c\widetilde{N}(t)})\|_{L^2_tL^{\frac{6}{5}}_x}\nonumber\\
&\lesssim_{(u,v)}(N^{\frac{1}{2}}+1)\|\varnothing(v_{\leq c\widetilde{N}(t)}v^2_{>1}u_{\leq 1}v_{>c\widetilde{N}(t)})\|_{L^2_tL^1_x}\nonumber\\
&\lesssim_{(u,v)} \|v_{\leq c\widetilde{N}(t)}\|_{L^{\infty}_tL^6_x}\|\varnothing(v^2_{>1}u_{\leq 1})\|_{L^2_tL^{\frac{3}{2}}_x}\|v_{>c\widetilde{N}(t)}\|_{L^{\infty}_tL^6_x}\lesssim_{(u,v)}\eta(1+K^{\frac{1}{2}}).\label{974}
\end{align}

Coming back to (\ref{965}), thanks to the results of (\ref{971})--(\ref{974}), we have
\begin{align}
(I)+(II)\lesssim_{(u,v)} 1+c^{-\frac{3}{2}}K^{\frac{1}{2}}+\eta(1+K^{\frac{1}{2}})LHS(\ref{962}).\label{975}
\end{align}
Similarly,
\begin{align}
(III)+(IV)\lesssim_{(u,v)} 1+c^{-\frac{3}{2}}K^{\frac{1}{2}}+\eta(1+K^{\frac{1}{2}})LHS(\ref{962}).\label{976}
\end{align}
(\ref{975}) and (\ref{976}) mean that
$$
LHS(\ref{962})\lesssim_{(u,v)} 1+c^{-\frac{3}{2}}K^{\frac{1}{2}}+\eta(1+K^{\frac{1}{2}})LHS(\ref{962}).
$$
If we chose $\eta$ small enough, we can get (\ref{962}).\hfill $\Box$

Now we are ready to prove the main result in this subsection.

{\bf Propostion 5.29(No rapid frequency-cascades).} {\it There are no almost periodic solution $(u,v):[0,T_{\max})\times \mathbb{R}^3\rightarrow \mathbb{C}\times \mathbb{C}$ of (\ref{826x1}) with $(\alpha,\beta)=(0,2)$ satisfying
\begin{align}
\|u\|_{L^{10}_{t,x}([0,T_{\max})\times \mathbb{R}^3)}+ \|v\|_{L^{10}_{t,x}([0,T_{\max})\times \mathbb{R}^3)}=+\infty,\quad
\int_0^{T_{\max}}\widetilde{N}(t)^{-1}dt<+\infty.\label{977}
\end{align}
}

{\bf Proof:} Assume contradictorily that $(u,v)$ is such a solution. Then
\begin{align}
\lim_{t\rightarrow T_{\max}}\widetilde{N}(t)=+\infty.\label{978}
\end{align}
Using Strichartz inequality and Bernstein's inequality
\begin{align}
 \|u_{\leq N}\|_{L^{\infty}_tL^2_x}+\|v_{\leq N}\|_{L^{\infty}_tL^2_x}&\lesssim_{(u,v)}\|P_{\leq N}(|v|^4u\|_{L^2_tL^{\frac{6}{5}}_x}
+\|P_{\leq N}(|u|^2|v|^2v\|_{L^2_tL^{\frac{6}{5}}_x}\nonumber\\
&\lesssim_{(u,v)}N^{\frac{1}{2}}\|(|v|^4u\|_{L^2_tL^1_x}
+N^{\frac{1}{2}}\|(|u|^2|v|^2v\|_{L^2_tL^1_x}.\label{979}
\end{align}

Young's inequality implies that
\begin{align}
||v|^4u|+|(|u|^2|v|^2v|\lesssim |u^5|+|v^5|.\label{97x1}
\end{align}

If we write
$$
u^5=\varnothing(u^3_{\leq 1}u^2)+\varnothing(u^3_{>1}u^2), \quad v^5=\varnothing(v^3_{\leq 1}v^2)+\varnothing(v^3_{>1}v^2),
$$
we have
\begin{align*}
&\quad \|u^5\|_{L^2_tL^1_x}+\|v^5\|_{L^2_tL^1_x}\nonumber\\
&\lesssim_{(u,v)} \|\varnothing(u^3_{\leq 1}u^2)\|_{L^2_tL^1_x}+\|\varnothing(u^3_{>1}u^2)\|_{L^2_tL^1_x}+
\|\varnothing(v^3_{\leq 1}v^2)\|_{L^2_tL^1_x}+\|\varnothing(v^3_{>1}v^2)\|_{L^2_tL^1_x}\nonumber\\
&\lesssim_{(u,v)}\|u_{\leq 1}\|_{L^4_tL^{\infty}_x}\|u_{\leq 1}\|_{L^{\infty}_{t,x}}\|u\|^2_{L^{\infty}_tL^2_x}+\|u\|^2_{L^{\infty}_tL^6_x}\|u^3_{>1}\|_{L^2_tL^{\frac{3}{2}}_x}\nonumber\\
&\quad +\|v_{\leq 1}\|_{L^4_tL^{\infty}_x}\|v_{\leq 1}\|_{L^{\infty}_{t,x}}\|v\|^2_{L^{\infty}_tL^2_x}+\|v\|^2_{L^{\infty}_tL^6_x}\|v^3_{>1}\|_{L^2_tL^{\frac{3}{2}}_x}\nonumber\\
&\lesssim_{(u,v)} 1.
\end{align*}
This estimate above, (\ref{979}) and (\ref{97x1}) imply that
\begin{align}
\|u_{\leq N}\|_{L^{\infty}_tL^2_x}+\|v_{\leq N}\|_{L^{\infty}_tL^2_x}\lesssim_{(u,v)} N^{\frac{1}{2}}.\label{97x2}
\end{align}

Recalling that for any $\eta>0$ there exists $c=c(u,v,\eta)>0$ such that
$$
\|\nabla u_{\leq c\widetilde{N}(t)}(t)\|_{L^2_x}+\|\nabla v_{\leq c\widetilde{N}(t)}(t)\|_{L^2_x}\lesssim_{(u,v)} \eta,
$$
we get
\begin{align*}
\|u(t)\|_{L^2_x}+\|v(t)\|_{L^2_x} &\lesssim_{(u,v)} \|u_{\leq c\widetilde{N}(t)}\|_{L^2_x}+\|P_{>N}u_{\leq c\widetilde{N}(t)}\|_{L^2_x}+\|u_{> c\widetilde{N}(t)}\|_{L^2_x}\nonumber\\
&\quad+\|v_{\leq c\widetilde{N}(t)}\|_{L^2_x}+\|P_{>N}v_{\leq c\widetilde{N}(t)}\|_{L^2_x}+\|v_{> c\widetilde{N}(t)}\|_{L^2_x}\nonumber\\
&\lesssim_{(u,v)} N^{\frac{1}{2}}+N^{-1}\|\nabla u_{\leq c\widetilde{N}(t)}\|_{L^2_x}+c^{-1}\widetilde{N}(t)^{-1}\|\nabla u\|_{L^{\infty}_tL^2_x}\nonumber\\
&\quad +N^{\frac{1}{2}}+N^{-1}\|\nabla v_{\leq c\widetilde{N}(t)}\|_{L^2_x}+c^{-1}\widetilde{N}(t)^{-1}\|\nabla u\|_{L^{\infty}_tL^2_x}\nonumber\\
&\lesssim_{(u,v)} N^{\frac{1}{2}}+N^{-1}\eta+c^{-1}\widetilde{N}(t)^{-1}.
\end{align*}
Using (\ref{978}), we can make the right-hand side here as small as we wish, which implies that $\|u\|_{L^{\infty}_tL^2_x}+\|v\|_{L^{\infty}_tL^2_x}=0$.
Consequently, $(u,v)\equiv(0,0)$, which is a contradiction to the hypothesis  $\|u\|_{L^{10}_{t,x}}+\|v\|_{L^{10}_{t,x}}=+\infty$.\hfill $\Box$

\subsection{The frequency-localized weight-coupled interaction Morawetz inequality}
\qquad In this subsection, we give a spacetime bound on the high-frequency portion of the solution.

{\bf Proposition 5.30(A frequency-localized interaction Morawetz estimate).} {\it
Assume that $(u,v):[0,T_{\max})\times \mathbb{R}^3\rightarrow \mathbb{C}\times \mathbb{C}$ is an almost periodic solution of (\ref{826x1}) with $(\alpha,\beta)=(0,2)$ such that $\widetilde{N}(t)\geq 1$. Let $I\subset [0,T_{\max})$ be a union of contiguous characteristic intervals $J_k$. Fix $0<\eta_0\leq 1$. Then for sufficiently small $N>0$ which depends on $\eta_0$ but not on $I$,
\begin{align}
\int_I\int_{\mathbb{R}^3}[|u_{>N}(t,x)|^4+|v_{>N}(t,x)|^4]dxdt\lesssim_{(u,v)} \eta_0(N^{-3}+K),\label{981}
\end{align}
where the implicit constant in the inequality above does not depend on $\eta_0$ or the interval $I$, $K:=\int_I\widetilde{N}(t)^{-1}dt$.
}

We also introduce the weight-coupled interaction Morawetz identity below.

{\bf Proposition 5.31.} {\it
Assume that
\begin{align}
i\phi_t+\Delta \phi=\lambda|\psi|^4\phi+\mathcal{F}_1,\quad i\psi_t+\Delta \psi=\mu|\phi|^2|\psi|^2\psi+\mathcal{F}_2,\label{982}
\end{align}
for some weight $a:\mathbb{R}^3\rightarrow \mathbb{R}$, we define the following weight-coupled Morawetz interaction:
\begin{align}
M^{\otimes_2}_a(t)&=2A\int_{\mathbb{R}^3}\int_{\mathbb{R}^3}\widetilde{\nabla} a(x-y)\Im[\bar{\phi}(t,x)\bar{\phi}(t,y)\widetilde{\nabla}(\phi(t,x)\phi(t,y))]dxdy\nonumber\\
&\quad+2B\int_{\mathbb{R}^3}\int_{\mathbb{R}^3}\widetilde{\nabla} a(x-y)\Im[\bar{\psi}(t,x)\bar{\psi}(t,y)\widetilde{\nabla}(\psi(t,x)\psi(t,y))]dxdy\nonumber\\
&\quad+2C\int_{\mathbb{R}^3}\int_{\mathbb{R}^3}\widetilde{\nabla} a(x-y)\Im[\bar{\phi}(t,x)\bar{\psi}(t,y)\widetilde{\nabla}(\phi(t,x)\psi(t,y))]dxdy\nonumber\\
&\quad+2D\int_{\mathbb{R}^3}\int_{\mathbb{R}^3}\widetilde{\nabla} a(x-y)\Im[\bar{\psi}(t,x)\bar{\phi}(t,y)\widetilde{\nabla}(\psi(t,x)\phi(t,y))]dxdy,\label{983}
\end{align}
where $\widetilde{\nabla}=(\nabla_x, \nabla_y)$, $x\in \mathbb{R}^3$ and $y\in \mathbb{R}^3$, the constants $A$, $B$, $C$ and $D$ are defined as in (\ref{1112w1}). Then

\begin{align}
&\qquad \frac{d}{dt}M^{\otimes_2}_a(t)\nonumber\\
&\thicksim \int_{\mathbb{R}^3}\int_{\mathbb{R}^3}\left(-[A|\phi|^2+D|\psi|^2]\Delta_x\Delta_xa+4a_{jk}Re(A\phi_j\bar{\phi}_k+D\psi_j\bar{\psi}_k)
+L_1|\phi|^2|\psi|^4\Delta_xa \right)dx|\phi(t,y)|^2dy\nonumber\\
&+\int_{\mathbb{R}^3}\int_{\mathbb{R}^3}\left(-[A|\phi|^2+C|\psi|^2]\Delta_y\Delta_ya+4a_{jk}Re(A\phi_j\bar{\phi}_k+C\psi_j\bar{\psi}_k)
+L_2|\phi|^2|\psi|^4\Delta_ya \right)dy|\phi(t,x)|^2dx\nonumber\\
&+\int_{\mathbb{R}^3}\int_{\mathbb{R}^3}\left(-[B|\psi|^2+C|\phi|^2]\Delta_x\Delta_xa+4a_{jk}Re(B\psi_j\bar{\psi}_k+C\phi_j\bar{\phi}_k)
+L_3|\phi|^2|\psi|^4\Delta_xa \right)dx|\psi(t,y)|^2dy\nonumber\\
&+\int_{\mathbb{R}^3}\int_{\mathbb{R}^3}\left(-[B|\psi|^2+D|\phi|^2]\Delta_y\Delta_ya+4a_{jk}Re(B\psi_j\bar{\psi}_k+D\phi_j\bar{\phi}_k)
+L_4|\phi|^2|\psi|^4\Delta_ya \right)dy|\psi(t,x)|^2dx\nonumber\\
&-4\sum_{m=1}^2\int_{\mathbb{R}^3}\int_{\mathbb{R}^3}\nabla_xa\cdot\Im[\bar{w}_m(t,x)\nabla_x w_m(t,x)]dx
[\nabla_y\cdot(\Im[\bar{\phi}(t,y)\nabla_y\phi(t,y)+\bar{\psi}(t,y)\nabla_y\psi(t,y)])]dy\nonumber\\
&-4\sum_{m=1}^2\int_{\mathbb{R}^3}\int_{\mathbb{R}^3}\nabla_ya\cdot\Im[\bar{w}_m(t,y)\nabla_y w_m(t,y)]dy
[\nabla_x\cdot(\Im[\bar{\phi}(t,x)\nabla_x\phi(t,x)+\bar{\psi}(t,x)\nabla_x\psi(t,x)])])]dx\nonumber\\
&+4\sum_{m,n=1}^2\int_{\mathbb{R}^3}\int_{\mathbb{R}^3}\Re(\overline{\mathcal{F}}_m\nabla_xw_m-\bar{w}_m\nabla_x \mathcal{F}_m)\cdot\nabla_xadx|w_n(t,y)|^2dy\nonumber\\
&+4\sum_{m,n=1}^2\int_{\mathbb{R}^3}\int_{\mathbb{R}^3}\Re(\overline{\mathcal{F}}_m\nabla_yw_m-\bar{w}_m\nabla_y \mathcal{F}_m)\cdot\nabla_yady|w_n(t,x)|^2dx\nonumber\\
&-4\sum_{m,n=1}^2\int_{\mathbb{R}^3}\int_{\mathbb{R}^3}\nabla_xa\cdot\Im[\bar{w}_m(t,x)\nabla_x w_m(t,x)]dx
\Im(\overline{\mathcal{F}}_nw_n)(t,y)dy\nonumber\\
&-4\sum_{m,n=1}^2\int_{\mathbb{R}^3}\int_{\mathbb{R}^3}\nabla_ya\cdot\Im[\bar{w}_m(t,y)\nabla_y w_m(t,y)]dy
\Im(\overline{\mathcal{F}}_nw_n)(t,x)dx.\label{984}
\end{align}
 Here $w_1$, $w_2$ can be taken as follows: (1) one is $\phi$, another is $\psi$; (2) both are $\phi$; (3) both are $\psi$. And the constants $L_1$, $L_2$, $L_3$ and $L_4$ are defined as in (\ref{11241}).
}

{\bf Proof:} (\ref{984}) can be obtained by (\ref{983}) after some elementary computations, we omit the details here. However,
we take ``$\frac{d}{dt}M^{\otimes_2}_a(t)\thicksim$" but not ``$\frac{d}{dt}M^{\otimes_2}_a(t)=$" because the coefficients in some terms on the righthand are different by the multiple of some constants.\hfill $\Box$

Similar to the discussions in \cite{Killip20101}, we choose $a$ be a smooth spherically symmetric function satisfying
\begin{align}
&a(0)=0,\quad a_r\geq 0,\quad a_{rr}\leq 0, \quad |\partial_r^ka_r|\lesssim_k J^{-1}r^{-k}\quad {\rm for\ each}\quad k\geq 1,\label{985}
\end{align}
and
\begin{equation}
\label{986} a_r=\left\{
\begin{array}{llll}
&1,\quad r\leq R,\\
&1-\frac{\log(\frac{r}{R})}{J},\quad eR\leq r \leq e^{J-J_0}R,\\
&0, \qquad e^JR\leq r,\\
\end{array}\right.
\end{equation}
where $J_0\geq 1$, $J\geq 2J_0$ and $R$ are parameters which will be determined in due course. Moreover, $a(x)=|x|$ if $|x|\leq R$,  $a(x)$ is a constant if $|x|\geq e^JR$, while
\begin{align}
\frac{2}{r}a_r\geq \frac{2J_0}{Jr},\quad |a_{rr}|\leq \frac{1}{Jr}\quad {\rm if}\quad eR\leq r \leq e^{J-J_0}R.\label{987}
\end{align}

For a parameter $N$ such that
\begin{align}
e^JRN=1,\label{988}
\end{align}
we will apply Proposition 5.31 with
\begin{align}
&\phi=u_{h_i}:=u_{>N},\quad \mathcal{F}_1=P_{hi}(|v|^4u)-|v_{hi}|^4u_{hi},\label{989}\\
&\psi=v_{h_i}:=v_{>N},\quad \mathcal{F}_2=P_{hi}(|u|^2|v|^2v)-|u_{hi}|^2|v_{hi}|^2v_{hi}.\label{989}
\end{align}
(By the way, we also write $u_{lo}:=u_{\leq N}$ and $v_{lo}:=v_{\leq N}$.) Since $(u,v)$ is almost periodic modulo symmetries and $\widetilde{N}(t)\geq 1$, if we take $R$ large enough, then $N$ can be small such that for given $\eta=\eta(\eta_0,u,v)$,
\begin{align}
&\quad \int_{\mathbb{R}^3}[|\nabla u_{lo}(t,x)|^2+|\nabla v_{lo}(t,x)|^2+|N u_{hi}(t,x)|^2+|N v_{hi}(t,x)|^2]dx\nonumber\\
&+\int_{|x-x(t)|>\frac{R}{2}}[|\nabla u_{hi}(t,x)|^2+|\nabla v_{hi}(t,x)|^2]dx<\eta^2\label{9810}
\end{align}
uniformly for $0\leq t<T_{\max}$.

{\bf Corollary 5.32(A priori bounds).} {\it
Under the assumption (\ref{9810}),
\begin{align}
\|u_{lo}\|_{L^4_tL^{\infty}_x(I\times \mathbb{R}^3)}+\|v_{lo}\|_{L^4_tL^{\infty}_x(I\times \mathbb{R}^3)}\lesssim_{(u,v)} \eta^{\frac{1}{2}}(1+N^3K)^{\frac{1}{4}}.\label{98x1}
\end{align}
For all admissible pair $(q,r)$, i.e., $\frac{2}{q}+\frac{3}{r}=\frac{3}{2}$, and any $s<1-\frac{3}{q}$,
\begin{align}
\|\nabla u_{lo}\|_{L^q_tL^r_x}+\|\nabla v_{lo}\|_{L^q_tL^r_x}+\|N^{1-s}|\nabla|^su_{hi}\|_{L^q_tL^r_x}
+\|N^{1-s}|\nabla|^sv_{hi}\|_{L^q_tL^r_x}\lesssim_{(u,v)}(1+N^3K)^{\frac{1}{q}}.\label{98x2}
\end{align}
For any $\rho\leq Re^J=N^{-1}$,
\begin{align}
\int_I\sup_{x\in \mathbb{R}^3}\int_{|x-y|\leq \rho}[|u_{hi}(t,y)|^2+|v_{hi}(t,y)|^2]dydt\lesssim_{(u,v)} \rho(K+N^{-3}).\label{98x3}
\end{align}
}

{\bf Proof:} (\ref{98x1}) is the direct results of Proposition 5.26 and (\ref{98xj2}). When $q=2$, the bound of $\|u_{lo}\|_{L^q_tL^r_x}+\|v_{lo}\|_{L^q_tL^r_x}$ in (\ref{98x2}) can be controlled similarly. If $q\neq 2$, using the conversation of energy, it can be deduced by interpolation.

To estimate $u_{hi}$ and $v_{hi}$, recalling that
$$
A(M):=\left\{\sum_{M'\leq M}[\|\nabla u_{M'}\|^2_{L^2_tL^6_x(I\times \mathbb{R}^3)}+\|\nabla v_{M'}\|^2_{L^2_tL^6_x(I\times \mathbb{R}^3)}]\right\}^{\frac{1}{2}}\lesssim_{(u,v)}(1+M^3K)^{\frac{1}{2}}
$$
uniformly in $M$, and by Bernstein's inequality, we have
\begin{align*}
&\quad M^{1-s}\||\nabla|^su_M\|_{L^q_tL^r_x}+M^{1-s}\||\nabla|^sv_M\|_{L^q_tL^r_x}\lesssim_{(u,v)}
\|\nabla u_M\|_{L^q_tL^r_x}+ \|\nabla v_M\|_{L^q_tL^r_x}\nonumber\\
&\lesssim_{(u,v)} A(M)^{\frac{2}{q}}[\|\nabla u\|^{\frac{(q-2)}{q}}_{L^{\infty}_tL^2_x}+ \|\nabla v\|^{\frac{(q-2)}{q}}_{L^{\infty}_tL^2_x}]\lesssim_{(u,v)}(1+M^3K)^{\frac{1}{q}},
\end{align*}
multiply through by $M^{s-1}$, and sum over $M\geq N$, then noticing $\frac{3}{q}+s<1$, we can guarantee the convergence of this sum and obtain the estimates for the parts containing $u_{hi}$ and $v_{hi}$ in (\ref{98x2}).

To establish (\ref{98x3}), we write
\begin{align*}
&|v|^4u=\varnothing(v^2_{>N}v^2u)+\varnothing(u_{> N}v_{>N}v^3)+\varnothing(u_{\leq N}v^3_{\leq N}v)+\varnothing(v^4_{\leq N}u)\\
&|u|^2|v|^2v=\varnothing(u^2_{>N}v^3)+\varnothing(u_{>N}v_{>N}v^2u)
+\varnothing(v^2_{>N}u^2v)+\varnothing(u_{\leq N}v^3_{\leq N}u)+\varnothing(u^2_{\leq N}v^2_{\leq N}v).
\end{align*}
Recalling Proposition 5.26, Proposition 5.27 and (\ref{94x8}), we have
\begin{align}
&\quad \|\varnothing(v^2_{>N}v^2u)\|_{L^2_tL^1_x}+\|\varnothing(u_{> N}v_{>N}v^3)\|_{L^2_tL^1_x}+\|\varnothing(u_{\leq N}v^3_{\leq N}v)\|_{L^2_tL^1_x}+\|\varnothing(v^4_{\leq N}u)\|_{L^2_tL^1_x}\nonumber\\
&\lesssim_{(u,v)}
[\|u\|^2_{L^{\infty}_tL^6_x}+\|v\|^2_{L^{\infty}_tL^6_x}][\|u^2_{>N}u\|_{L^2_tL^{\frac{3}{2}}_x}+\|u^2_{>N}v\|_{L^2_tL^{\frac{3}{2}}_x}
+\|v^2_{>N}v\|_{L^2_tL^{\frac{3}{2}}_x}+\|v^2_{>N}u\|_{L^2_tL^{\frac{3}{2}}_x}]\nonumber\\
&\lesssim_{(u,v)} N^{-\frac{3}{2}}+K^{\frac{1}{2}},\label{98w1}\\
&\quad \||u|^2|v|^2v\|_{L^2_tL^1_x}\lesssim_{(u,v)} N^{-\frac{3}{2}}+K^{\frac{1}{2}},\label{98w1'}
\end{align}
and by Bernstein's inequality, (\ref{98xj2}), we get
\begin{align}
&\quad\|\varnothing(u^2_{>N}v^3)\|_{L^2_tL^{\frac{6}{5}}_x}+\|\varnothing(u_{>N}v_{>N}v^2u)\|_{L^2_tL^{\frac{6}{5}}_x}
+\|\varnothing(v^2_{>N}u^2v)\|_{L^2_tL^{\frac{6}{5}}_x}\nonumber\\
&\quad +\|\varnothing(u_{\leq N}v^3_{\leq N}u)\|_{L^2_tL^{\frac{6}{5}}_x} +\|\varnothing(u^2_{\leq N}v^2_{\leq N}v)\|_{L^2_tL^{\frac{6}{5}}_x}\nonumber\\
&\lesssim_{(u,v)}N^{-1}[\|\nabla u\|_{L^{\infty}_tL^2_x}+\|\nabla v\|_{L^{\infty}_tL^2_x}]\|u\|^2_{L^{\infty}_tL^6_x}+\|v\|^2_{L^{\infty}_tL^6_x}][\|u_{\leq N}\|^2_{L^4_tL^{\infty}_x}+\|v_{\leq N}\|^2_{L^4_tL^{\infty}_x}]\nonumber\\
&\lesssim_{(u,v)} N^{-1}+N^{\frac{1}{2}}K^{\frac{1}{2}},\label{98w2}\\
&\quad\||v|^4u\|_{L^2_tL^{\frac{6}{5}}_x}\lesssim_{(u,v)} N^{-1}+N^{\frac{1}{2}}K^{\frac{1}{2}}.\label{98w2'}
\end{align}
Combing (\ref{98w1})--(\ref{98w2'}) with the results of Proposition 5.23, we obtain
\begin{align}
\rho^{\frac{1}{2}}[\|\mathcal{S}u_{>N}\|_{L^2_tL^{\infty}_x}+\|\mathcal{S}v_{>N}\|_{L^2_tL^{\infty}_x}]\lesssim_{(u,v)} N^{-1}+(N^{\frac{1}{2}}+\rho^{-\frac{1}{2}})(K+N^{-3})^{\frac{1}{2}}.\label{98w3}
\end{align}
Since $\rho\leq Re^J=N^{-1}$, modulo a factor of $\rho^{-\frac{3}{2}}$, $(\mathcal{S}u_{>N}(t,x), \mathcal{S}u_{>N}(t,x))$ controls the $L^2_x$-norm on the ball around $x$,
(\ref{98x3}) can be deduced by (\ref{98w3}).\hfill $\Box$

{\bf Lemma 5.33(Mass-mass interactions).} {\it
\begin{align}
&\quad 16\pi[\|u_{hi}\|^4_{L_{t,x}(I\times \mathbb{R}^3)}+\|v_{hi}\|^4_{L_{t,x}(I\times \mathbb{R}^3)}]\nonumber\\
&\quad +\int_I\int_{\mathbb{R}^3}\int_{\mathbb{R}^3}a_{jjkk}(x-y)[|u_{hi}(t,y)|^2+|v_{hi}(t,y)|^2][|u_{hi}(t,x)|^2+|v_{hi}(t,x)|^2]dxdydt\nonumber\\
&\lesssim_{(u,v)} \frac{\eta^2e^{2J}}{J}(K+N^{-3}).\label{991}
\end{align}
}

{\bf Proof:} Noticing that $\Delta |x|=2|x|^{-1}$ and $-\Delta (4\pi|x|)^{-1}=\delta(x)$, we only to estimate the error terms originating from the truncation of $a$ at radii $|x-y|\geq R$. By (\ref{985}), (\ref{988}), (\ref{9810}) and (\ref{98x3}), we have
\begin{align*}
&\quad\int_I\int_{\mathbb{R}^3}\int_{\mathbb{R}^3}a_{jjkk}(x-y)[|u_{hi}(t,y)|^2+|v_{hi}(t,y)|^2][|u_{hi}(t,x)|^2+|v_{hi}(t,x)|^2]dxdydt\nonumber\\
&\lesssim_{(u,v)} J^{-1}[\|u_{hi}\|^2_{L^{\infty}_tL^2_x}+\|v_{hi}\|^2_{L^{\infty}_tL^2_x}]\sum_{j=0}^J(Re^j)^{-2}(K+N^{-3})\nonumber\\
&\lesssim_{(u,v)}\frac{\eta^2e^{2J}}{J}(K+N^{-3}).\qquad\qquad\qquad\qquad\qquad\qquad\qquad\qquad\qquad\qquad\qquad\qquad\qquad\qquad\qquad \Box
\end{align*}

 Denote
\begin{align}
B_I&:=\int_I\int\int_{|x-y|\leq e^{J-J_0}R}|u_{hi}(t,x)|^2|v_{hi}(t,x)|^4\Delta_x a dx[L_1|u_{hi}|^2+L_3|v_{hi}|^2](t,y)dy\nonumber\\
&\qquad+\int_I\int\int_{|x-y|\leq e^{J-J_0}R}|u_{hi}(t,y)|^2|v_{hi}(t,y)|^4\Delta_y a dy[L_2|u_{hi}|^2+L_4|v_{hi}|^2](t,x)dx,\label{993}\\
B_{II}&:=\int_I\int\int_{|x-y|\geq e^{J-J_0}R}|u_{hi}(t,x)|^2|v_{hi}(t,x)|^4|\Delta_x a| dx[L_1|u_{hi}|^2+L_3|v_{hi}|^2](t,y)dy\nonumber\\
&\qquad+\int_I\int\int_{|x-y|\geq e^{J-J_0}R}|u_{hi}(t,y)|^2|v_{hi}(t,y)|^4|\Delta_y a| dy[L_2|u_{hi}|^2+L_4|v_{hi}|^2](t,x)dx,\label{994}
\end{align}

{\bf Lemma 5.34.} {\it Let $B_I$ and $B_{II}$ be defined as (\ref{993}) and (\ref{994}). Then
\begin{align}
B_I\geq 0,\quad B_{II}\lesssim_{(u,v)}\frac{J_0^2}{J}(K+N^{-3}).\label{995}
\end{align}
}

{\bf Proof:} $B_I\geq 0$ is a direct result of $\Delta_x a\geq 0$ and $\Delta_y a\geq 0$ by the construction of $a$. Noticing that $|\Delta_x a|\lesssim J_0(Jr)^{-1}$ and $|\Delta_y a|\lesssim J_0(Jr)^{-1}$ when $r\geq e^{J-J_0}R$ and $a_{kk}=0$ when $r>e^{JR}$, we get
\begin{align*}
B_{II}&\lesssim \int_I\int\int_{|x-y|\geq e^{J-J_0}R}\frac{J_0}{J|x-y|}|u_{hi}(t,x)|^2|v_{hi}(t,x)|^4dx[L_1|u_{hi}|^2+L_3|v_{hi}|^2](t,y)dy\nonumber\\
&\quad+\int_I\int\int_{|x-y|\geq e^{J-J_0}R}\frac{J_0}{J|x-y|}|u_{hi}(t,y)|^2|v_{hi}(t,y)|^4dy[L_2|u_{hi}|^2+L_4|v_{hi}|^2](t,x)dx\nonumber\\
&\lesssim \int_I\int\int_{e^{J-J_0}R\leq |x-y|\leq e^JR}\frac{J_0}{J|x-y|}[|u_{hi}|^6+|v_{hi}|^6](t,x)dx[L_1|u_{hi}|^2+L_3|v_{hi}|^2](t,y)dy\nonumber\\
&\quad+\int_I\int\int_{e^{J-J_0}R\leq |x-y|\leq e^JR}\frac{J_0}{J|x-y|}[|u_{hi}|^6+|v_{hi}|^6](t,y)dy[L_2|u_{hi}|^2+L_4|v_{hi}|^2](t,x)dx\nonumber\\
&\lesssim_{(u,v)} J_0[\|u_{hi}\|^6_{L^{\infty}_tL^6_x}+\|v_{hi}\|^6_{L^{\infty}_tL^6_x}]\sum_{j=J-J_0}^J(Je^jR)^{-1}(e^jR)(K+N^{-3})\nonumber\\
& \lesssim_{(u,v)}\frac{J^2_0}{J}(K+N^{-3}),
\end{align*}
which completes the proof of (\ref{995}).\hfill$\Box$

{\bf Lemma 5.35.} {\it
\begin{align}
&\quad \int_I\int_{\mathbb{R}^3}\int_{\mathbb{R}^3}a_{jk}Re[(u_{hi})_j(\bar{u}_{hi})_k+(v_{hi})_j(\bar{v}_{hi})_k](t,x)dx[|u_{hi}|^2+|v_{hi}|^2](t,y)dy\nonumber\\
&+\int_I\int_{\mathbb{R}^3}\int_{\mathbb{R}^3}a_{jk}Re[(u_{hi})_j(\bar{u}_{hi})_k+(v_{hi})_j(\bar{v}_{hi})_k](t,y)dy[|u_{hi}|^2+|v_{hi}|^2](t,x)dx\nonumber\\
&\lesssim_{(u,v)} (\eta^2+\frac{J_0}{J})(K+N^{-3})+\frac{1}{J_0}B_I.\label{996}
\end{align}
}

{\bf Proof:} Since the $3\times 3$-matrix $(a_{jk})$ is invariant under $x\leftrightarrow y$ and real symmetric for any $x$ and $y$, and $a_{jk}Re[(u_{hi})_j(\bar{u}_{hi})_k+(v_{hi})_j(\bar{v}_{hi})_k](t,x)[|u_{hi}|^2+|v_{hi}|^2](t,y)$ and $a_{jk}Re[(u_{hi})_j(\bar{u}_{hi})_k+(v_{hi})_j(\bar{v}_{hi})_k](t,y)[|u_{hi}|^2+|v_{hi}|^2](t,x)$ define two positive semi-definite quadratic forms on $\mathbb{R}^3$ for each $x$, $y$. Therefore, wherever $(a_{jk})$ is positive semi-definite, (\ref{996}) is true. By the choice of $a$, (\ref{986}) and (\ref{987}), we only to estimate
\begin{align}
(I):=\int_I\int\int_{R<|x-y|<e^JR}\frac{[|\nabla u_{hi}(t,x)|^2+|\nabla v_{hi}(t,x)|^2][u_{hi}(t,y)|^2+|v_{hi}(t,y)|^2]}{J|x-y|}dxdydt\label{997}
\end{align}
and
\begin{align}
(II):=\int_I\int\int_{R<|x-y|<e^JR}\frac{[|\nabla u_{hi}(t,y)|^2+|\nabla v_{hi}(t,y)|^2][u_{hi}(t,x)|^2+|v_{hi}(t,x)|^2]}{J|x-y|}dxdydt.\label{998}
\end{align}
We only estimate (\ref{997}). (\ref{998}) can be estimated similarly.  We break the integral into two regions: $|x-x(t)|>\frac{R}{2}$ and $|x-x(t)|\leq \frac{R}{2}$, the corresponding integrands are (I1) and (I2).

In the region $|x-x(t)|>\frac{R}{2}$, by (\ref{9810}) and (\ref{98x3}), we have
\begin{align}
(I1)&\lesssim_{(u,v)}[\|\nabla u_{hi}\|^2_{L^{\infty}_tL^2_x(|x-x(t)|>\frac{R}{2})}+\|\nabla v_{hi}\|^2_{L^{\infty}_tL^2_x(|x-x(t)|>\frac{R}{2})}]
\sum_{j=0}^J(Je^jR)^{-1}(e^jR)(K+N^{-3})\nonumber\\
&\lesssim_{(u,v)}\eta^2(K+N^{-3}).\label{999}
\end{align}

In the region $|x-x(t)|\leq \frac{R}{2}$, we further break the integral into two regions: $|x-y|>Re^{J-J_0}$ and $|x-y|\leq Re^{J-J_0}$, the corresponding integrands are (I21) and (I22). In the region $|x-y|>Re^{J-J_0}$, similar to (\ref{999}), we get
\begin{align}
(I21)&\lesssim_{(u,v)}[\|\nabla u_{hi}\|^2_{L^{\infty}_tL^2_x}+\|\nabla v_{hi}\|^2_{L^{\infty}_tL^2_x}]
\sum_{j=J-J_0}^J(Je^jR)^{-1}(e^jR)(K+N^{-3})\nonumber\\
&\lesssim_{(u,v)}\frac{J_0}{J}(K+N^{-3}).\label{9910}
\end{align}

Since $(u,v)$ is almost periodic,
\begin{align*}
&\quad\int_{\mathbb{R}^3}[|\nabla u_{hi}(t,x)|^2+|\nabla v_{hi}(t,x)|^2]dx\lesssim_{(u,v)}\int_{\mathbb{R}^3}[|u_{hi}(t,x)|^6+|v_{hi}(t,x)|^6]dx\\
&\lesssim_{(u,v)}\int_{|x-x(t)|\leq \frac{R}{2}}[|u_{hi}(t,x)|^6+|v_{hi}(t,x)|^6]dx\quad {\rm uniformly\ for}\quad t\in [0,T_{\max}).
\end{align*}
Since $\Delta a=a_{rr}+\frac{2}{r}a_r$ in space dimension 3, recalling (\ref{9810}), (\ref{98x3}) and $J_0\geq 1$, we know that the remaining integral is  $\lesssim_{(u,v)}\frac{1}{J}B_I$. \hfill$\Box$

Let
\begin{align}
\tilde{{\mathcal {M}}}(u_{hi},v_{hi})&:=4\sum_{m,n=1}^2\int_{\mathbb{R}^3}\int_{\mathbb{R}^3}\Re(\overline{\mathcal{F}}_m\nabla_xw_m-\bar{w}_m\nabla_x \mathcal{F}_m)\cdot\nabla_xadx|w_n(t,y)|^2dy\nonumber\\
&\qquad+4\sum_{m,n=1}^2\int_{\mathbb{R}^3}\int_{\mathbb{R}^3}\Re(\overline{\mathcal{F}}_m\nabla_yw_m-\bar{w}_m\nabla_y \mathcal{F}_m)\cdot\nabla_yady|w_n(t,x)|^2dx,\label{99w1}
\end{align}
where $w_m$ can be taken $u_{hi}$ or $v_{hi}$, the same does to $w_n$. We have

{\bf Lemma 5.36.} {\it For any $\epsilon\in (0,1]$,
\begin{align}
|\int_I\tilde{{\mathcal {M}}}(u_{hi},v_{hi})dt|\lesssim_{(u,v)}\epsilon B_I+\eta[\|u_{hi}\|^4_{L^4_{t,x}}+\|v_{hi}\|^4_{L^4_{t,x}}]+(\epsilon^{-1}\eta+\epsilon\frac{J^2_0}{J})(N^{-3}+K).\label{99w2}
\end{align}
}

{\bf Proof:} Denote
\begin{align*}
f_1(u,v)=|v|^4u,\quad f_2(u,v)=|u|^2|v|^2v,\quad \mathcal{F}_1=f_1(u,v)-P_{lo}f_1(u,v),\quad \mathcal{F}_2=f_2(u,v)-P_{lo}f_2(u,v).
\end{align*}
Then
\begin{align}
&\quad\sum_{m=1}^2\Re(\overline{\mathcal{F}}_m\nabla_xw_m-\bar{w}_m\nabla_x \mathcal{F}_m)\nonumber\\
&=\nabla_x\varnothing(u_{hi}f_1+v_{hi}f_2)+\varnothing(u_{hi}|v|^4\nabla_x u_{lo}+u_{hi}u|v|^2v\nabla_x v_{lo})\nonumber\\
&\quad +\varnothing(v_{hi}|u|^2|v|^2\nabla_x v_{lo}+v_{hi}u|v|^2v\nabla_x u_{lo})
+\nabla_x \varnothing(u_{hi}P_{lo}f_1+v_{hi}P_{lo}f_2)\nonumber\\
&\quad +\varnothing(u_{hi}\nabla_x P_{lo}f_1+v_{hi}\nabla_x P_{lo}f_2)\nonumber\\
&:=(I)+(II)+(III)+(IV)+(V).\label{99w3}
\end{align}
And
\begin{align}
&\quad \int_I\tilde{{\mathcal {M}}}(u_{hi},v_{hi})dt\nonumber\\
&=\int_I\int_{\mathbb{R}^3}\int_{\mathbb{R}^3}\nabla_xa\cdot \left[(I)+(II)+(III)+(IV)+(V)\right](t,x)dx|u_{hi}(t,y)|^2
dydt\nonumber\\
&+\int_I\int_{\mathbb{R}^3}\int_{\mathbb{R}^3}\nabla_xa\cdot \left[(I)+(II)+(III)+(IV)+(V)\right](t,x)dx|v_{hi}(t,y)|^2]
dydt\nonumber\\
&+\int_I\int_{\mathbb{R}^3}\int_{\mathbb{R}^3}\nabla_ya\cdot\left[(I)+(II)+(III)+(IV)+(V)\right](t,y)dy|u_{hi}(t,x)|^2
dxdt\nonumber\\
&+\int_I\int_{\mathbb{R}^3}\int_{\mathbb{R}^3}\nabla_ya\cdot\left[(I)+(II)+(III)+(IV)+(V)\right](t,y)dy|v_{hi}(t,x)|^2
dxdt\nonumber\\
&:=(1)+(2)+(3)+(4).\label{99w4}
\end{align}

We only estimate (1) in details, those of (2), (3) and (4) are similar to (1). First, we deal with
\begin{align}
&\quad |\int_I\int_{\mathbb{R}^3}\int_{\mathbb{R}^3}[\nabla_xa\cdot (I)](t,x)dx|u_{hi}(t,y)|^2dydt|\nonumber\\
&=|\int_I\int_{\mathbb{R}^3}\int_{\mathbb{R}^3}\Delta_xa\varnothing(u_{hi}f_1+v_{hi}f_2)(t,x)dx|u_{hi}(t,y)|^2dydt|\nonumber\\
&\lesssim\epsilon \int_I\int_{\mathbb{R}^3}\int_{\mathbb{R}^3}|\Delta_x a|[|u_{hi}(t,x)|^6+|v_{hi}(t,x)|^6]|u_{hi}(t,y)|^2dxdydt\nonumber\\
&\quad +\int_I\int_{\mathbb{R}^3}\int_{\mathbb{R}^3}\frac{|u_{hi}(t,y)|^2[|u_{hi}||v_{lo}|^2(|u_{hi}|+|v_{hi}|+|u_{lo}+|v_{lo}|)^3(t,x)]}{\epsilon |x-y|}dxdydt\nonumber\\
&\quad +\int_I\int_{\mathbb{R}^3}\int_{\mathbb{R}^3}\frac{|u_{hi}(t,y)|^2[|v_{hi}||u_{lo}|^2(|u_{hi}|+|v_{hi}|+|u_{lo}+|v_{lo}|)^3(t,x)]}{\epsilon |x-y|}dxdydt\nonumber\\
&\lesssim_{(u,v)} \frac{\epsilon J^2_0}{J}(K+N^{-3})+\frac{1}{\epsilon}\||x|^{-1}*|u_{hi}|^2\|_{L^4_tL^6_x}|u_{hi}|\|_{L^4_tL^3_x}\|v_{lo}\|^2_{L^4_tL^{\infty}_x}[\|u\|^3_{L^{\infty}_tL^6_x}
+\|v\|^3_{L^{\infty}_tL^6_x}]\nonumber\\
&\quad+\frac{1}{\epsilon}\||x|^{-1}*|u_{hi}|^2\|_{L^4_tL^6_x}|v_{hi}|\|_{L^4_tL^3_x}\|u_{lo}\|^2_{L^4_tL^{\infty}_x}[\|u\|^3_{L^{\infty}_tL^6_x}
+\|v\|^3_{L^{\infty}_tL^6_x}]\nonumber\displaybreak\\
&\lesssim_{(u,v)} \frac{\epsilon J^2_0}{J}(N^{-3}+K)+\frac{1}{\epsilon}[\|u_{hi}\|_{L^{\infty}_tL^2_x}+\|v_{hi}\|_{L^{\infty}_tL^2_x}]
[\|u_{hi}\|^2_{L^4_tL^3_x}+\|v_{hi}\|^2_{L^4_tL^3_x}]\nonumber\\
&\qquad \qquad \qquad \qquad \qquad \quad\times [\|u_{lo}\|^2_{L^4_tL^{\infty}_x}+\|v_{lo}\|^2_{L^4_tL^{\infty}_x}]
[\|u_{hi}\|^3_{L^{\infty}_tL^6_x}+\|v_{hi}\|^3_{L^{\infty}_tL^6_x}]\nonumber\\
&\lesssim_{(u,v)} [\frac{\epsilon J^2_0}{J}+\frac{\eta}{\epsilon}](N^{-3}+K).\label{99w5}
\end{align}

Since $|\nabla_x a|$ is uniformly bounded and
\begin{align}
|(II)|+|(III)&|\lesssim [|u_{hi}|^2+|v_{hi}|^2][|u|^3+|v|^3][|\nabla_x u_{lo}|+|\nabla_x v_{lo}|] \nonumber\\
&\quad+[|u_{hi}|+|v_{hi}|][|u|^2+|v|^2][|u_{lo}|^2+|v_{lo}|^2][|\nabla_x u_{lo}|+|\nabla_x v_{lo}|],\label{9101}
\end{align}
 we have
\begin{align}
&\quad |\int_I\int_{\mathbb{R}^3}\int_{\mathbb{R}^3}\nabla_xa\cdot [(II)+(III)](t,x)dx|u_{hi}(t,y)|^2dydt|\nonumber\\
&\lesssim \int_I\int_{\mathbb{R}^3}\int_{\mathbb{R}^3}|[(II)+(III)]|(t,x)dx|u_{hi}(t,y)|^2dydt\lesssim_{(u,v)}
\eta^2N^{-2}\int_I\int_{\mathbb{R}^3}[|(II)|+|(III)|](t,x)dxdt\nonumber\\
&\lesssim_{(u,v)}\eta^2N^{-2}[\|\nabla u_{lo}\|_{L^2_tL^{\infty}_x}+\|\nabla v_{lo}\|_{L^2_tL^{\infty}_x}][\|u_{hi}\|^2_{L^4_{t,x}}+\|v_{hi}\|^2_{L^4_{t,x}}]
[\|u\|^3_{L^{\infty}_tL^6_x}+\|v\|^3_{L^{\infty}_tL^6_x}]\nonumber\\
&\qquad \quad+\eta^2N^{-2}[\|u_{hi}\|_{L^{\infty}_tL^2_x}+\|v_{hi}\|_{L^{\infty}_tL^2_x}][\|\nabla u_{lo}\|_{L^2_tL^6_x}+\|\nabla v_{lo}\|_{L^2_tL^6_x}]
[\| u_{lo}\|^2_{L^4_tL^{\infty}_x}+\|v_{lo}\|^2_{L^4_tL^{\infty}_x}]\nonumber\\
&\qquad\qquad \qquad \qquad \times[\|u\|^2_{L^{\infty}_tL^6_x}+\|v\|^2_{L^{\infty}_tL^6_x}]\nonumber\\
&\lesssim_{(u,v)}\eta^2[\|u_{hi}\|^4_{L^4_{t,x}}+\|v_{hi}\|^4_{L^4_{t,x}}]+(\eta+\eta^2)(N^{-3}+K).\label{9102}
\end{align}
Here we use $\int_{\mathbb{R}^3}|u_{hi}|^2(t,y)dy\lesssim_{(u,v)}\eta^2N^{-2}$ by (\ref{98x3}).

Now we deal the term containing (IV) by integrating by parts
\begin{align}
&\quad \int_I\int_{\mathbb{R}^3}\int_{\mathbb{R}^3}[\nabla_xa\cdot (IV)](t,x)dx|u_{hi}(t,y)|^2dydt\nonumber\\
&=\int_I\int_{\mathbb{R}^3}\int_{\mathbb{R}^3}|u_{hi}|^2(t,y)\Delta_xa \varnothing(u_{hi}P_{lo}f_1+v_{hi}P_{lo}f_2)(t,x)dxdydt.\label{9103}
\end{align}
Writing $u_{hi}(t,x)=div (\nabla \Delta^{-1}u_{hi}(t,x)$ and $v_{hi}(t,x)=div (\nabla \Delta^{-1}v_{hi}(t,x)$, and integrating by parts once more, then applying the Mikhlin multiple theorem, we obtain
\begin{align}
&\quad |\int_I\int_{\mathbb{R}^3}\int_{\mathbb{R}^3}|u_{hi}|^2(t,y)\Delta_xa \varnothing(u_{hi}P_{lo}f_1+v_{hi}P_{lo}f_2)(t,x)dxdydt|\nonumber\\
&\lesssim \||x|^{-1}*|u_{hi}|^2\|_{L^4_tL^{12}_x}\||\nabla|^{-1}u_{hi}\|_{L^2_tL^6_x}\|\nabla P_{lo}f_1\|_{L^4_tL^{\frac{4}{3}}_x}\nonumber\\
&\quad+ \||x|^{-2}*|u_{hi}|^2\|_{L^4_tL^{\frac{12}{5}}_x}\||\nabla|^{-1}u_{hi}\|_{L^2_tL^6_x}\|\nabla P_{lo}f_1\|_{L^4_tL^{\frac{12}{5}}_x}\nonumber\\
&\quad+ \||x|^{-1}*|u_{hi}|^2\|_{L^4_tL^{12}_x}\||\nabla|^{-1}v_{hi}\|_{L^2_tL^6_x}\|\nabla P_{lo}f_2\|_{L^4_tL^{\frac{4}{3}}_x}\nonumber\\
&\quad+ \||x|^{-2}*|u_{hi}|^2\|_{L^4_tL^{\frac{12}{5}}_x}\||\nabla|^{-1}v_{hi}\|_{L^2_tL^6_x}\|\nabla P_{lo}f_2\|_{L^4_tL^{\frac{12}{5}}_x}.\label{9104}
\end{align}
Applying Hardy-Littlewood-Sobolev inequality to the first factor in each term above, and using Sobolev embedding on very last factors in the second term and the fourth one, we get
\begin{align}
RHS(\ref{9104})&\lesssim_{(u,v)} \||u_{hi}|^2\|_{L^4_tL^{\frac{4}{3}}_x}\||\nabla|^{-1}u_{hi}\|_{L^2_tL^6_x}\|\nabla P_{lo}f_1\|_{L^4_tL^{\frac{4}{3}}_x}\nonumber\\
&\qquad +\||u_{hi}|^2\|_{L^4_tL^{\frac{4}{3}}_x}\||\nabla|^{-1}v_{hi}\|_{L^2_tL^6_x}\|\nabla P_{lo}f_2\|_{L^4_tL^{\frac{4}{3}}_x}.\label{9105}
\end{align}

We need to estimate $\nabla P_{lo}f_1$ and $\nabla P_{lo}f_2$. Since $u=u_{hi}+u_{lo}$ and $v=v_{hi}+v_{lo}$, we have
\begin{align}
&|v_{hi}+v_{lo}|^4(u_{hi}+u_{lo})=[u_{lo}+u_{hi}][|v_{lo}|^4+4v_{hi}|\xi_1|^2\xi_1]\label{910x1}\\
&|u_{hi}+u_{lo}|^2|v_{hi}+v_{lo}|^2(v_{hi}+v_{lo})=[|u_{lo}|^2+2u_{hi}\xi_2][|v_{lo}|^2v_{lo}+3v_{hi}|\xi_3|^2],\label{910x2}
\end{align}
Here $\xi_2$ is some function ``between" $u_{lo}$ and $u_{hi}+u_{lo}$, while $\xi_1$ and $\xi_3$ are some functions ``between" $v_{lo}$ and $v_{hi}+v_{lo}$. Therefore, we can write
\begin{align}
f_1&=|v_{lo}|^4u_{lo}+\varnothing(u_{hi}|v|^4)+\varnothing(v_{hi}u|\xi_2|^2\xi_1),\label{910x3}\\
f_2&=|u_{lo}|^2|v_{lo}|^2v_{lo}+\varnothing(u_{hi}\xi_2|v|^2v)+\varnothing(v_{hi}|u|^2|\xi_3|^2).\label{910x4}
\end{align}
Using H\"{o}lder's, Bernetein's  and Young's inequalities, we can get
\begin{align}
&\quad \|\nabla P_{lo}(|v_{lo}|^4u_{lo})\|_{L^4_tL^{\frac{4}{3}}_x}\lesssim N^{\frac{3}{4}}\|\nabla (|v_{lo}|^4u_{lo})\|_{L^4_tL^1_x}\nonumber\\
&\lesssim_{(u,v)} N^{\frac{3}{4}}[\|\nabla u_{lo}\|_{L^4_tL^3_x}+\|\nabla v_{lo}\|_{L^4_tL^3_x}][\|u_{lo}\|^4_{L^{\infty}_tL^6_x}+\|u_{lo}\|^4_{L^{\infty}_tL^6_x}]\nonumber\\
&\lesssim_{(u,v)}N^{\frac{3}{4}}(1+N^3K)^{\frac{1}{4}},\label{910x5}\\
&\quad \|\nabla P_{lo}\varnothing(u_{hi}|v|^4)\|_{L^4_tL^{\frac{4}{3}}_x}\lesssim N^{\frac{3}{2}}
\|u_{hi}|v|^4\|_{L^4_tL^{\frac{12}{11}}_x}\nonumber\\
&\lesssim_{(u,v)} N^{\frac{3}{2}}\|u_{hi}\|_{L^4_{t,x}}[\|u\|^4_{L^{\infty}_tL^6_x}
\lesssim_{(u,v)} N^{\frac{3}{2}}\|u_{hi}\|_{L^4_{t,x}},\label{910x6}
\end{align}
and similarly,
\begin{align}
& \|\nabla P_{lo}\varnothing(v_{hi}u|\xi_1|^2\xi_1)\|_{L^4_tL^{\frac{4}{3}}_x}\lesssim_{(u,v)}N^{\frac{3}{2}}\|v_{hi}\|_{L^4_{t,x}},\label{910x7}\\
& \|\nabla P_{lo}|u_{lo}|^2|v_{lo}|^2v_{lo}\|_{L^4_tL^{\frac{4}{3}}_x}\lesssim_{(u,v)}N^{\frac{3}{4}}(1+N^3K)^{\frac{1}{4}},\label{910x8}\\
& \|\nabla P_{lo}\varnothing(u_{hi}\xi_2|v|^2v)\|_{L^4_tL^{\frac{4}{3}}_x}\lesssim_{(u,v)}N^{\frac{3}{2}}\|u_{hi}\|_{L^4_{t,x}},\label{910x9}\\
& \|\nabla P_{lo}\varnothing(v_{hi}|u|^2|\xi_3|^2)\|_{L^4_tL^{\frac{4}{3}}_x}\lesssim_{(u,v)}N^{\frac{3}{2}}\|v_{hi}\|_{L^4_{t,x}}.\label{910x10}
\end{align}

Putting (\ref{9103})--(\ref{910x10}) together, by the results of Corollary 5.32 and (\ref{98x3}), using H\"{o}lder's and Young's inequalities, we obtain
\begin{align}
&\quad |\int_I\int_{\mathbb{R}^3}\int_{\mathbb{R}^3}[\nabla_xa\cdot (IV)](t,x)dx|u_{hi}(t,y)|^2dydt|\nonumber\\
&\lesssim_{(u,v)} \|u_{hi}\|_{L^{\infty}_tL^2_x}\|u_{hi}\|_{L^4_{t,x}}\||\nabla|^{-1}u_{hi}\|_{L^2_tL^6_x}[N^{\frac{3}{2}}\|u_{hi}\|_{L^4_{t,x}}
+N^{\frac{3}{4}}(1+N^3K)^{\frac{1}{4}}]\nonumber\\
&\qquad +\|u_{hi}\|_{L^{\infty}_tL^2_x}\|u_{hi}\|_{L^4_{t,x}}\||\nabla|^{-1}v_{hi}\|_{L^2_tL^6_x}[N^{\frac{3}{2}}\|v_{hi}\|_{L^4_{t,x}}
+N^{\frac{3}{4}}(1+N^3K)^{\frac{1}{4}}]\nonumber\\
&\lesssim_{(u,v)} \eta N^{-1}\|u_{hi}\|_{L^4_{t,x}}N^{-2}(1+N^3K)^{\frac{1}{2}}[N^{\frac{3}{2}}\|u_{hi}\|_{L^4_{t,x}}
+N^{\frac{3}{4}}(1+N^3K)^{\frac{1}{4}}]\nonumber\\
&\qquad +\eta N^{-1}\|u_{hi}\|_{L^4_{t,x}}N^{-2}(1+N^3K)^{\frac{1}{2}}[N^{\frac{3}{2}}\|v_{hi}\|_{L^4_{t,x}}
+N^{\frac{3}{4}}(1+N^3K)^{\frac{1}{4}}]\nonumber\\
&\lesssim_{(u,v)} \eta [\|u_{hi}\|^4_{L^4_{t,x}}+\|v_{hi}\|^4_{L^4_{t,x}}+(N^{-3}+K)].\label{910w1}
\end{align}

Now we deal with the term containing (V). Similarly, we also write $u_{hi}(t,x)=div (\nabla \Delta^{-1}u_{hi}(t,x)$ and $v_{hi}(t,x)=div (\nabla \Delta^{-1}v_{hi}(t,x)$, and integrate by parts once, then
\begin{align}
&\quad |\int_I\int_{\mathbb{R}^3}\int_{\mathbb{R}^3}[\nabla_xa\cdot (V)](t,x)dx|u_{hi}(t,y)|^2dydt|\nonumber\\
&\lesssim_{(u,v)}  \||x|^{-1}*|u_{hi}|^2\|_{L^4_tL^{12}_x}\||\nabla|^{-1}u_{hi}\|_{L^2_tL^6_x}\|\nabla P_{lo}f_1\|_{L^4_tL^{\frac{4}{3}}_x}\nonumber\displaybreak\\
&\quad+ \|(\nabla a)*|u_{hi}|^2\|_{L^{\infty}_tL^{\infty}_x}\||\nabla|^{-1}u_{hi}\|_{L^2_tL^6_x}\|\Delta P_{lo}f_1\|_{L^2_tL^{\frac{6}{5}}_x}\nonumber\\
&\quad+ \||x|^{-1}*|u_{hi}|^2\|_{L^4_tL^{12}_x}\||\nabla|^{-1}v_{hi}\|_{L^2_tL^6_x}\|\nabla P_{lo}f_2\|_{L^4_tL^{\frac{4}{3}}_x}\nonumber\\
&\quad+ \|(\nabla a)*|u_{hi}|^2\|_{L^{\infty}_tL^{\infty}_x}\||\nabla|^{-1}v_{hi}\|_{L^2_tL^6_x}\|\Delta P_{lo}f_1\|_{L^2_tL^{\frac{6}{5}}_x}.\label{9111}
\end{align}
The first and third terms in (\ref{9111}) have been estimated above. We estimate the second and fourth ones. By (\ref{910x1}), we write
$$
f_1=|v_{lo}|^4u_{lo}+\varnothing(u_{hi}v^2_{lo}|v_{lo}|^2)+\varnothing(u_{lo}v_{hi}|\xi_1|^3)
+\varnothing(u_{hi}v_{hi}|\xi_1|^3).
$$
Using Bernstein, H\"{o}lder and Young's inequality, we can obtain
\begin{align*}
&\quad \|\Delta P_{lo}(|v_{lo}|^4u_{lo})\|_{L^2_tL^{\frac{6}{5}}_x}\lesssim_{(u,v)}N[\|\nabla u_{lo}\|_{L^2_tL^6_x}+\|\nabla v_{lo}\|_{L^2_tL^6_x}][\|u_{lo}\|^4_{L^{\infty}_tL^6_x}+\|v_{lo}\|^4_{L^{\infty}_tL^6_x}]\nonumber\\
&\lesssim_{(u,v)}N(1+N^3K)^{\frac{1}{2}},\\
&\quad \|\Delta P_{lo} \varnothing(u_{hi}v^2_{lo}|v_{lo}|^2)\|_{L^2_tL^{\frac{6}{5}}_x}\lesssim_{(u,v)}N^2\|u_{hi}\|_{L^{\infty}_tL^2_x}
\|v_{lo}\|^2_{L^4_tL^{\infty}_x}[\|u_{lo}\|^2_{L^{\infty}_tL^6_x}+\|v_{lo}\|^2_{L^{\infty}_tL^6_x}]\\
&\lesssim_{(u,v)}N(1+N^3K)^{\frac{1}{2}},\\
&\quad \|\Delta P_{lo} \varnothing(u_{lo}v_{hi}|\xi_1|^3)\|_{L^2_tL^{\frac{6}{5}}_x}\lesssim_{(u,v)} N^{\frac{5}{2}}\|v_{hi}\|^2_{L^4_{t,x}}[\|u_{lo}\|^3_{L^{\infty}_tL^6_x}+\|v_{lo}\|^3_{L^{\infty}_tL^6_x}]\nonumber\\
&\quad +N^2\|v_{hi}\|_{L^{\infty}_tL^2_x}
[\|u_{lo}\|^2_{L^4_tL^{\infty}_x}+\|v_{lo}\|^2_{L^4_tL^{\infty}_x}][\|u_{lo}\|^2_{L^{\infty}_tL^6_x}+\|v_{lo}\|^2_{L^{\infty}_tL^6_x}]\\
&\lesssim_{(u,v)}N^{\frac{5}{2}}\|v_{hi}\|^2_{L^4_{t,x}}+N(1+N^3K)^{\frac{1}{2}},\\
&\quad \|\Delta P_{lo} \varnothing(u_{hi}v_{hi}|\xi_1|^3)\|_{L^2_tL^{\frac{6}{5}}_x}\lesssim_{(u,v)}
N^{\frac{5}{2}}[\|u_{hi}\|^2_{L^4_{t,x}}+\|v_{hi}\|^2_{L^4_{t,x}}][\|u_{lo}\|^3_{L^{\infty}_tL^6_x}+\|\xi_1\|^3_{L^{\infty}_tL^6_x}]\nonumber\\
&\lesssim_{(u,v)}N^{\frac{5}{2}}[\|u_{hi}\|^2_{L^4_{t,x}}+\|v_{hi}\|^2_{L^4_{t,x}}].
\end{align*}
Here we use the fact that $|\xi_1|\lesssim |v|+|v_{hi}|+|v_{lo}|$. Consequently,
$$
\|\Delta P_{lo}f_1\|_{L^2_tL^{\frac{6}{5}}_x}\lesssim_{(u,v)} N(1+N^3K)^{\frac{1}{2}}+N^{\frac{5}{2}}[\|u_{hi}\|^2_{L^4_{t,x}}+\|v_{hi}\|^2_{L^4_{t,x}}]
$$
and
\begin{align*}
&\|(\nabla a)*|u_{hi}|^2\|_{L^{\infty}_tL^{\infty}_x}\||\nabla|^{-1}u_{hi}\|_{L^2_tL^6_x}\|\Delta P_{lo}f_1\|_{L^2_tL^{\frac{6}{5}}_x}\nonumber\\
&\lesssim_{(u,v)}\|u_{hi}\|^2_{L^{\infty}_tL^2_x}N^{-2}(1+N^3K)^{\frac{1}{2}}
\left(N(1+N^3K)^{\frac{1}{2}}+N^{\frac{5}{2}}[\|u_{hi}\|^2_{L^4_{t,x}}+\|v_{hi}\|^2_{L^4_{t,x}}]\right).
\end{align*}
Similarly,
$$
\|\Delta P_{lo}f_2\|_{L^2_tL^{\frac{6}{5}}_x}\lesssim_{(u,v)} N(1+N^3K)^{\frac{1}{2}}+N^{\frac{5}{2}}[\|u_{hi}\|^2_{L^4_{t,x}}+\|v_{hi}\|^2_{L^4_{t,x}}].
$$
and
\begin{align*}
&\|(\nabla a)*|u_{hi}|^2\|_{L^{\infty}_tL^{\infty}_x}\||\nabla|^{-1}u_{hi}\|_{L^2_tL^6_x}\|\Delta P_{lo}f_2\|_{L^2_tL^{\frac{6}{5}}_x}\nonumber\\
&\lesssim_{(u,v)}\|u_{hi}\|^2_{L^{\infty}_tL^2_x}N^{-2}(1+N^3K)^{\frac{1}{2}}
\left(N(1+N^3K)^{\frac{1}{2}}+N^{\frac{5}{2}}[\|u_{hi}\|^2_{L^4_{t,x}}+\|v_{hi}\|^2_{L^4_{t,x}}]\right).
\end{align*}
Combining all results above, we have estimated (1). Similar results on (2), (3) and (4) can also be obtained.  Putting all the results together, the proof of this lemma is finished.\hfill $\Box$

Now we deal with the mass bracket terms.

{\bf Lemma 5.37(Mass bracket terms).}  {\it
\begin{align}
&|\sum_{m,n=1}^2\int_I\int_{\mathbb{R}^d}\int_{\mathbb{R}^d}\nabla_xa\cdot\Im[\bar{w}_m(t,x)\nabla_xw_m(t,x)]dx
\Im(\overline{\mathcal{F}}_nw_n)(t,y)dydt|\nonumber\\
&+|\sum_{m,n=1}^2\int_I\int_{\mathbb{R}^d}\int_{\mathbb{R}^d}\nabla_ya\cdot\Im[\bar{w}_m(t,y)\nabla_yw_m(t,y)]dy
\Im(\overline{\mathcal{F}}_nw_n)(t,x)dxdt|\nonumber\\
&\lesssim_{(u,v)} \eta^{\frac{1}{4}}[\|u_{hi}\|^4_{L^4_{t,x}}+\|v_{hi}\|^4_{L^4_{t,x}}+(N^{-3}+K)].
\end{align}
Here $w_m$, $w_n$ can be taken as follows: (1) one is $u_{hi}$, another is $v_{hi}$; (2) both are $u_{hi}$; (3) both are $v_{hi}$.
}

{\bf Proof:} First, we write
\begin{align}
\Im(\overline{\mathcal{F}}_1u_{hi})&=\Im[u_{hi}P_{hi}(|v_{lo}|^4u_{lo})]-\Im[u_{hi}P_{lo}(|v_{hi}|^4u_{hi})]+\Im[\varnothing(v_{lo}u^2_{hi}|\zeta_1|^3)]
+\Im[\varnothing(u_{lo}u_{hi}v_{hi}|\zeta_1|^3)]\nonumber\\
&:=(I)+(II)+(III)+(IV)\label{9121}
\end{align}
and consequently
\begin{align}
&\quad \int_I\int_{\mathbb{R}^d}\int_{\mathbb{R}^d}\nabla_xa\cdot\Im[\bar{w}_m(t,x)\nabla_xw_m(t,x)]dx
\Im(\overline{\mathcal{F}}_1u_{hi})(t,y)dydt\nonumber\\
&=\int_I\int_{\mathbb{R}^d}\int_{\mathbb{R}^d}\nabla_xa\cdot\Im[\bar{w}_m(t,x)\nabla_xw_m(t,x)][(I)+(II)+(III)+(IV)]dxdydt.\label{11251}
\end{align}
We will treat their contributions in order.

Obviously,
\begin{align}
&\quad|\int_I\int_{\mathbb{R}^3}\int_{\mathbb{R}^3}\nabla_xa\cdot\Im[\bar{w}_m(t,x)\nabla_x w_m(t,x)][(I)(t,y)]dxdydt|\nonumber\\
&\lesssim \|\bar{w}_m\nabla w_m\|_{L^{\infty}_tL^1_x}\|u_{hi}P_{hi}(|v_{lo}|^4u_{lo})\|_{L^1_{t,x}}\lesssim_{(u,v)}\|\nabla w_m\|_{L^{\infty}_tL^2_x}\|w_m\|_{L^{\infty}_tL^2_x}N^{-1}\|\nabla (|v_{lo}|^4u_{lo}))\|_{L^1_tL^2_x}\nonumber\\
&\lesssim_{(u,v)}\eta^2N^{-3}[\|\nabla u_{lo}\|_{L^2_tL^6_x}+\|\nabla v_{lo}\|_{L^2_tL^6_x}]
[\|u_{lo}\|^2_{L^4_tL^{\infty}_x}+\|v_{lo}\|^2_{L^4_tL^{\infty}_x}][\|u_{lo}\|^2_{L^{\infty}_tL^6_x}+\|v_{lo}\|^2_{L^{\infty}_tL^6_x}]\nonumber\\
&\lesssim_{(u,v)}\eta^2(N^{-3}+K).\label{912w1}
\end{align}

For the term containing (II), we write $u_{hi}=div(\nabla \Delta^{-1}u_{hi})$ and integrate by parts. Then we get
\begin{align}
&\quad|\int_I\int_{\mathbb{R}^3}\int_{\mathbb{R}^3}\nabla_xa\cdot\Im[\bar{w}_m(t,x)\nabla_xw_m(t,x)][(II)(t,y)]dxdydt|\nonumber\\
&\lesssim \|\bar{w}_m\nabla w_m\|_{L^{\infty}_tL^1_x}\||\nabla|^{-1}u_{hi}\|_{L^2_tL^6_x}\|\nabla P_{lo}(|v_{hi}|^4u_{hi})\|_{L^2_tL^{\frac{6}{5}}_x}\nonumber\\
&\quad +\||x|^{-1}*|\bar{w}_m\nabla w_m|\|_{L^4_tL^{12}_x}\||\nabla|^{-1}u_{hi}\|_{L^2_tL^6_x}\| P_{lo}(|v_{hi}|^4u_{hi})\|_{L^4_tL^{\frac{4}{3}}_x}\nonumber\\
&\lesssim_{(u,v)} \|w_m\|_{L^{\infty}_tL^2_x}N^{-2}(1+N^3K)^{\frac{1}{2}}N^{\frac{3}{2}}\|(|v_{hi}|^4u_{hi})\|_{L^2_tL^1_x}\nonumber\\
&\quad+\|\nabla w_m\|_{L^{\infty}_tL^2_x}\|w_m\|_{L^4_{t,x}}N^{-2}(1+N^3K)^{\frac{1}{2}}N^{\frac{3}{4}}\|(|v_{hi}|^4u_{hi})\|_{L^4_tL^1_x}\nonumber\\
&\lesssim_{(u,v)}\eta(N^{-3}+K)^{\frac{1}{2}}[\|u_{hi}\|^2_{L^4_{t,x}}+\|v_{hi}\|^2_{L^4_{t,x}}]
[\|u_{hi}\|^3_{L^{\infty}_tL^6_x}+\|v_{hi}\|^3_{L^{\infty}_tL^6_x}]\nonumber\\
&\quad +\|w_m\|_{L^4_{t,x}}N^{\frac{1}{4}}(N^{-3}+K)^{\frac{1}{2}}[\|u_{hi}\|_{L^4_{t,x}}+\|v_{hi}\|_{L^4_{t,x}}]
[\|u_{hi}\|^{\frac{15}{4}}_{L^{\infty}_tL^6_x}+\|v_{hi}\|^{\frac{15}{4}}_{L^{\infty}_tL^6_x}]\nonumber\\
&\qquad\qquad \qquad\qquad\qquad\qquad \times
[\|u_{hi}\|^{\frac{1}{4}}_{L^{\infty}_tL^2_x}+\|v_{hi}\|^{\frac{1}{4}}_{L^{\infty}_tL^2_x}]\nonumber\\
&\lesssim_{(u,v)}(\eta+\eta^{\frac{1}{4}}) [\|u_{hi}\|^4_{L^4_{t,x}}+\|v_{hi}\|^4_{L^4_{t,x}}+N^{-3}+K].\label{912w2}
\end{align}

For the term containing (III), by the choice of $a$ and $Re^J=N^{-1}$, we have
\begin{align}
&\quad|\int_I\int_{\mathbb{R}^3}\int_{\mathbb{R}^3}\nabla_xa\cdot\Im[\bar{w}_m(t,x)\nabla_xw_m(t,x)][(III)(t,y)]dxdydt|\nonumber\\
&\lesssim_{(u,v)} \|\nabla w_m\|_{L^{\infty}_tL^2_x}\|w_m\|_{L^4_{t,x}}\|\nabla a\|_{L^{\infty}_tL^4_x}[\|u_{hi}\|^2_{L^4_{t,x}}+\|v_{hi}\|^2_{L^4_{t,x}}]
[\|u_{lo}\|_{L^4_tL^{\infty}_x}+\|v_{lo}\|_{L^4_tL^{\infty}_x}]\nonumber\\
&\qquad \qquad \qquad \qquad \qquad \qquad \times [\|u\|^3_{L^{\infty}_tL^6_x}+\|u_{hi}\|^3_{L^{\infty}_tL^6_x}+\|v\|^3_{L^{\infty}_tL^6_x}
+\|v_{hi}\|^3_{L^{\infty}_tL^6_x}]\nonumber\\
&\lesssim_{(u,v)} [\|u_{hi}\|^3_{L^4_{t,x}}+\|v_{hi}\|^3_{L^4_{t,x}}](e^JR)^{\frac{3}{4}}\eta^{\frac{1}{2}}(1+N^3K)^{\frac{1}{4}}\nonumber\\
&\lesssim_{(u,v)} \eta^{\frac{1}{2}}[\|u_{hi}\|^4_{L^4_{t,x}}+\|v_{hi}\|^4_{L^4_{t,x}}+N^{-3}+K].\label{912w3}
\end{align}

The term containing (IV) can be treated similarly, and
\begin{align}
&\quad|\int_I\int_{\mathbb{R}^3}\int_{\mathbb{R}^3}\nabla_xa\cdot\Im[\bar{w}_m(t,x)\nabla_xw_m(t,x)][(IV)(t,y)]dxdydt|\nonumber\\
&\lesssim_{(u,v)} \eta^{\frac{1}{2}}[\|u_{hi}\|^4_{L^4_{t,x}}+\|v_{hi}\|^4_{L^4_{t,x}}+N^{-3}+K].\label{912w4}
\end{align}
This completes the control of the mass bracket terms.\hfill $\Box$

{\bf Proof of Proposition 5.29:} If we take $\phi=u_{hi}$, $\psi=v_{hi}$ and the function $a$ as above, and choosing
$N$ is small enough such that (\ref{9810}) holds, using H\"{o}lder's inequality, we have
$$
|M^{\otimes_2}_a(t)|\lesssim_{(u,v)}[\|u_{hi}\|^3_{L^{\infty}_tL^2_x}+\|v_{hi}\|^3_{L^{\infty}_tL^2_x}]
[\|\nabla u_{hi}\|_{L^{\infty}_tL^2_x}+\|\nabla v_{hi}\|_{L^{\infty}_tL^2_x}]\lesssim_{(u,v)}\eta^3N^{-3}.
$$
All the integral are on $I\times \mathbb{R}^3$.

Putting all the conclusions of the above lemmas, we get
\begin{align}
16\pi[\|u_{hi}\|^4_{L^4_{t,x}}+\|v_{hi}\|^4_{L^4_{t,x}}]+B_I&\lesssim_{(u,v)}(\epsilon+\frac{1}{J_0})B_I
+\eta^{\frac{1}{4}}[\|u_{hi}\|^4_{L^4_{t,x}}+\|v_{hi}\|^4_{L^4_{t,x}}]\nonumber\\
&\lesssim_{(u,v)}(\eta^{\frac{1}{4}}+\frac{\eta}{\epsilon}+\frac{J^2_0}{J}+\eta^2\frac{e^{2J}}{J})(N^{-3}+K).\label{913x3}
\end{align}

Now we choose the parameters as follows: First, we choose $\epsilon+J^{-1}_0<<1$; Second, for these fixed $\epsilon$ and $J^{-1}_0$, we choose $\eta<<1$ and $J^{-}$ small enough such that $(\eta^{\frac{1}{4}}+\frac{\eta}{\epsilon}+\frac{J^2_0}{J}+\eta^2\frac{e^{2J}}{J})<\eta_0$; Last, for the fixed $J$, we choose $R$ and $N^{-1}$ large enough such that $NRe^J=1$, yet we can increase $N^{-1}$ or $R$ such that (\ref{9810}) holds.

Meanwhile,
\begin{align*}
&\quad [\|u_{hi}\|^4_{L^4_{t,x}}+\|v_{hi}\|^4_{L^4_{t,x}}]\lesssim_{(u,v)}[\||\nabla |^{\frac{1}{4}}u_{hi}\|^4_{L^4_tL^3_x}+\||\nabla |^{\frac{1}{4}}v_{hi}\|^4_{L^4_tL^3_x}]\nonumber\\
&\lesssim_{(u,v)} N^{-3}[\|u_{hi}\|^4_{L^4_tL^3_x}+\|v_{hi}\|^4_{L^4_tL^3_x}]\lesssim_{(u,v)}
N^{-3}+N^{-3}\int_I\widetilde{N}(t)^2dt,
\end{align*}
and
\begin{align*}
B_I&\lesssim_{(u,v)}[\||x|^{-1}*|u_{hi}|^2\|_{L^4_tL^{12}_x}+\||x|^{-1}*|v_{hi}|^2\|_{L^4_tL^{12}_x}][\|u_{hi}\|^{\frac{5}{4}}_{L^{\infty}_tL^2_x}
+\|v_{hi}\|^{\frac{5}{4}}_{L^{\infty}_tL^2_x}]\nonumber\\
&\qquad\qquad\qquad\qquad\qquad\qquad\times [\|u_{hi}\|^{\frac{19}{4}}_{L^{\frac{19}{3}}_tL^{\frac{114}{7}}_x}+\|v_{hi}\|^{\frac{19}{4}}_{L^{\frac{19}{3}}_tL^{\frac{114}{7}}_x}]\nonumber\\
&\lesssim_{(u,v)}[\|u_{hi}\|_{L^4_{t,x}}+\|v_{hi}\|_{L^4_{t,x}}][\|u_{hi}\|^{\frac{9}{4}}_{L^{\infty}_tL^2_x}+\|v_{hi}\|^{\frac{9}{4}}_{L^{\infty}_tL^2_x}]
[\|\nabla u_{hi}\|^{\frac{19}{4}}_{L^{\frac{19}{3}}_tL^{\frac{38}{15}}_x}+\|\nabla v_{hi}\|^{\frac{19}{4}}_{L^{\frac{19}{3}}_tL^{\frac{38}{15}}_x}]\nonumber\\
&\lesssim_{(u,v)}N^{-3}+N^{-3}\int_I\widetilde{N}(t)^2dt.
\end{align*}

By the choosing of the parameters, the terms containing $B_I$ and $[\|u_{hi}\|^4_{L^4_{t,x}}+\|v_{hi}\|^4_{L^4_{t,x}}]$ on the righthand side can be absorbed by the left ones, we get
$$
\qquad \qquad\qquad \qquad [\|u_{hi}\|^4_{L^4_{t,x}}+\|v_{hi}\|^4_{L^4_{t,x}}]\lesssim_{(u,v)}\eta_0(N^{-3}+K).\qquad \qquad\qquad \qquad\qquad  \qquad \hfill \Box
$$

\subsection{Impossibility of quasi-solution}

{\bf Proposition 5.38(No quasi-solutions).} {\it There are no almost periodic solutions $(u,v):[0,T_{\max})\times \mathbb{R}^3\rightarrow \mathbb{C}\times \mathbb{C}$ of (\ref{826x1})  with $\widetilde{N}(t)\equiv N_k\geq 1$ on each characteristic interval $J_k\subset [0,T_{\max})$
 which satisfy
\begin{align}
\|u\|_{L^{10}_{t,x}([0,T_{\max})\times \mathbb{R}^3)}+\|v\|_{L^{10}_{t,x}([0,T_{\max})\times \mathbb{R}^3)}=+\infty,\quad \int_0^{T_{\max}}\widetilde{N}(t)^{-1}(t)dt=+\infty.\label{913w1}
\end{align}
}

{\bf Proof:} Contradictorily, assume that such a solution $(u,v)$ exists.

Since $(u(t,x),v(t,x))$ is not identically zero for each times $t$ and $(u,v)$ is almost periodic, there exists $C(u,v)>0$ such that
\begin{align}
\widetilde{N}(t)\int_{|x-x(t)|\leq \frac{C(u,v)}{\widetilde{N}(t)}}[|u(t,x)|^4+|v(t,x)|^4]dx\geq \frac{1}{C(u,v)}\label{913w2}
\end{align}
uniformly for $t\in [0,T_{\max})$. Using H\"{o}lder's inequality, we have
$$
\widetilde{N}(t)\int_{|x-x(t)|\leq \frac{C(u,v)}{\widetilde{N}(t)}}[|u_{\leq N}(t,x)|^4+|v_{\leq N}(t,x)|^4]dx\lesssim_{(u,v)} [\|u_{\leq N}(t,\cdot)\|^4_{L^6_x}+\|v_{\leq N}(t,\cdot)\|^4_{L^6_x}]
$$
for any $N>0$ uniformly for $t\in [0,T_{\max})$. Recalling the results of Proposition 5.30  and (\ref{913w2}), we know that for each $\eta_0>0$, there exists some $N=N(\eta_0)$ sufficiently small such that
$$
\int_I[\widetilde{N}(t)]^{-1}dt\lesssim_{(u,v)}\eta_0N^{-3}+\eta_0\int_I[\widetilde{N}(t)]^{-1}dt
$$
uniformly for time intervals $I\subset [0,T_{\max})$ that are a union of characteristic subintervals $J_k$. Particularly, we can chose $\eta_0$ small enough such that the term on the righthand side containing $\int_I[\widetilde{N}(t)]^{-1}dt$ can be absorbed by the one on the left side. Consequently,
$$
\int_0^{T_{\max}}[\widetilde{N}(t)]^{-1}dt=\lim_{T\uparrow T_{\max}}\int_0^T[\widetilde{N}(t)]^{-1}dt\lesssim_{(u,v)} 1,
$$
which is a contradiction to (\ref{913w1}).\hfill $\Box$

Now we have precluded the existence of the two types of almost periodic solution described in Proposition 5.6 and proved Theorem 6.

\section{$\dot{H}^1\times \dot{H}^1$ scattering theories for the global solution of (\ref{826x1}) with defocusing nonlinearities when $d=4$ }
\qquad In this section, we will establish $\dot{H}^1\times \dot{H}^1$ scattering theories for the global solution of (\ref{826x1}) in defocusing case with $(\alpha,\beta)=(0,0)$ when $d=4$.

The organization of this section is similar to that of Section 5. We begin with the following definition, some propositions and lemmas below. The proofs of Proposition 6.2, Lemma 6.3 and Corollary 6.4 had been proved in Section 5 as well as those of Proposition 5.2, Lemma 5.3 and Corollary 5.4. The proofs of Lemma 6.5 and Proposition 6.7 are respectively similar to those of Lemma 5.5 and Proposition 5.7, we omit the details here.

{\bf Definition 6.1(Almost periodicity).} {\it A solution $(u,v)\in [L^{\infty}_t\dot{H}^1_x(I\times \mathbb{R}^4)]^2$ of (\ref{826x1}) is said to be almost periodic(modulo symmetries) if there exist functions $\widetilde{N}:I\rightarrow \mathbb{R}^+$, $\tilde{x}:I\rightarrow \mathbb{R}^4$, and $C:\mathbb{R}^+\rightarrow \mathbb{R}^+$
such that for all $t\in I$ and $\eta>0$
\begin{align}
\int_{|x-\tilde{x}(t)|\geq \frac{C(\eta)}{\widetilde{N}(t)}}[|\nabla u(t,x)|^2+|\nabla v(t,x)|^2]dx
+\int_{|\xi|\geq C(\eta)\widetilde{N}(t)}|\xi|^2[|\hat{u}(t,\xi)|^2+|\hat{v}(t,\xi)|^2]d\xi\leq \eta.\label{934}
\end{align}
The functions $\widetilde{N}(t)$, $\tilde{x}(t)$ and $C(\eta)$ are called as the frequency scale function for the solution $(u,v)$, the spatial center function and the modulus of compactness respectively. }

By compactness, there exists $c(\eta)>0$ such that
\begin{align}
\int_{|x-\tilde{x}(t)|\leq \frac{c(\eta)}{\widetilde{N}(t)}}[|\nabla u(t,x)|^2+|\nabla v(t,x)|^2]dx
+\int_{|\xi|\leq c(\eta)\widetilde{N}(t)}|\xi|^2[|\hat{u}(t,\xi)|^2+|\hat{v}(t,\xi)|^2]d\xi\leq \eta.\label{1225z1}
\end{align}

{\bf Proposition 6.2(Reduction to almost periodic solution).} {\it
Assume that Theorem 7 failed. Then there exists a maximal-lifespan solution $(u,v):I\times \mathbb{R}^4\rightarrow \mathbb{C}\times \mathbb{C}$
of (\ref{826x1}) with $(\alpha,\beta)=(0,0)$ which is almost periodic and blows up both forward and backward in time in the sense that for all $t_0\in I$,
$$
\int_{t_0}^{\sup I}\int_{\mathbb{R}^4}[|u(t,x)|^6+|v(t,x)|^6]dxdt=\int^{t_0}_{\inf I}\int_{\mathbb{R}^4}[|u(t,x)|^6+|v(t,x)|^6]dxdt=+\infty.
$$
 }

{\bf Lemma 6.3(Local constancy property).} {\it
 Let $(u,v):I\times \mathbb{R}^4\rightarrow \mathbb{C}\times \mathbb{C}$ be a maximal-lifespan almost period solution of (\ref{826x1}) with $(\alpha,\beta)=(0,0)$. Then there exists a small number $\delta$ which depends only on $(u,v)$, such that
 $$
 [t_0-\delta\widetilde{N}(t_0)^{-2}, t_0+\delta\widetilde{N}(t_0)^{-2}]\subset I\quad {\rm if }\quad t_0\in I,
 $$
 and
 $$
 \widetilde{N}(t)\thicksim_{(u,v)}\widetilde{N}(t_0)\quad {\rm whenever}\quad |t-t_0|\leq \delta\widetilde{N}(t_0)^{-2}.
 $$
 }

{\bf Corollary 6.4($\widetilde{N}(t)$ blows up).} {\it
Let $(u,v):I\times \mathbb{R}^4\rightarrow \mathbb{C}\times \mathbb{C}$ be a maximal-lifespan almost period solution of (\ref{826x1}) with $(\alpha,\beta)=(0,0)$. If $T$ is any finite endpoint of $I$, then $\widetilde{N}(t)\geq_{(u,v)}|T-t|^{-\frac{1}{2}}$. Consequently, $\lim_{t\rightarrow T} \tilde{N}(t)=+\infty$.
}

{\bf Lemma 6.5(Spacetime bounds).} {\it Let $(u,v)$ be an almost periodic solution of (\ref{826x1}) with $(\alpha,\beta)=(0,0)$ on a time interval $I$. Then for any admissible pair $(q,r)$ with $2\leq q<+\infty$,
\begin{align}
\int_I\widetilde{N}(t)^2dt\lesssim_{(u,v)}\|\nabla u\|^2_{L^2_tL^4_x(I\times \mathbb{R}^4)}+\|\nabla v\|^2_{L^2_tL^4_x(I\times \mathbb{R}^4)}\lesssim_{(u,v)}1+\int_I\widetilde{N}(t)^2dt\label{10241}
\end{align}
}

Since it can be realized by a simple rescaling argument, we can assume that $\widetilde{N}(t)\geq 1$ at least on half of $I$, say, on $[0,T_{\max})$. Finally, we get

{\bf Proposition 6.6(Two special scenarios for blowup).} {\it Assume that Theorem 7 failed. Then there exists an almost periodic solution $(u,v):[0,T_{\max})\times \mathbb{R}^4\rightarrow \mathbb{C}\times \mathbb{C}$ such that
$$
\|u\|_{L_{t,x}^6([0,T_{\max})\times \mathbb{R}^4)}+\|v\|_{L_{t,x}^6([0,T_{\max})\times \mathbb{R}^4)}=+\infty
$$
and $[0,T_{\max})=\cup_{k}J_k$, where $J_k$ are characteristic intervals on which $\widetilde{N}(t)\equiv N_k\geq 1$.
Moreover,
$$
{\rm either}\quad \int_0^{T_{\max}}\widetilde{N}(t)^{-1}dt<+\infty\quad {\rm or}\quad  \int_0^{T_{\max}}\widetilde{N}(t)^{-1}dt=+\infty.
$$
}

Similarly, we only need to preclude the existence of the two types of almost periodic solution described in Proposition 6.6, then we can give the proof of Theorem 7. To prove the no-existence of cascade solutions, we also need the following proposition.

{\bf Proposition 6.7(No-waste Duhamel formulae).} {\it
Let $(u,v):[0,T_{\max})\times \mathbb{R}^4\rightarrow \mathbb{C}\times \mathbb{C}$ be defined as in Proposition 6.6. Then for all $t\in [0,T_{\max})$
\begin{align}
u(t)=i\lim_{T\rightarrow T_{\max}} \int_t^Te^{i(t-s)\Delta}|v|^2u(s)ds,\quad v(t)=i\lim_{T\rightarrow T_{\max}} \int_t^Te^{i(t-s)\Delta}|u|^2v(s)ds\label{1125s3}
\end{align}
in the weak $\dot{H}^1_x$ topology.
}

\subsection{Long-time Strichartz estimates}
\qquad In this subsection, we give the Long-time Strichartz estimates.

The following lemmas are a consequence of dispersive estimates for the free Schr\"{o}dinger propagator $e^{it\Delta}$, which was proved in \cite{Ginibre1992, Keel1998, Strichartz1977}.

{\bf Lemma 6.8(Strichartz inequality).} {\it Assume that $I$ is a compact time interval and $w:I\times \mathbb{R}^4\rightarrow \mathbb{C}$ is a solution of the forced Sch\"{o}dinger equation
$$iw_t+\Delta w=G\quad {\rm for\ some\ function}\quad G.$$
Then for any time $t_0\in I$ and admissible pairs $(q,r)$ and $(\tilde{q},\tilde{r})$, i.e.,
\begin{align*}
\frac{1}{q}+\frac{2}{r}=\frac{1}{\tilde{q}}+\frac{2}{\tilde{r}}=1\quad {\rm and}\quad 2\leq q,\tilde{q}\leq \infty,
\end{align*}
the following inequality holds
\begin{align*}
\|\nabla w\|_{L^q_tL^r_x(I\times \mathbb{R}^4)}\lesssim \|w(t_0)\|_{\dot{H}^1_x(\mathbb{R}^4)}+\|\nabla G\|_{L^{\tilde{q}'}_tL^{\tilde{r}'}_x(I\times \mathbb{R}^4)},
\end{align*}
where $\tilde{q}'$ is the dual exponent to $\tilde{q}$ and satisfies $\frac{1}{\tilde{q}}+\frac{1}{\tilde{q}'}=1$.}

The following Lemmas can be found in \cite{Bourgain19992, Colliander2008, Ryckman2007, Visan2007, Visan2011}.

{\bf Lemma 6.9(Bilinear Strichartz).} {\it For any spacetime slab $I\times \mathbb{R}^4$ and any frequencies $M>0$ and $N>0$,
\begin{align*}
\|u_{\leq M}v_{\geq N}\|_{L^2_{t,x}(I\times \mathbb{R}^4)}\lesssim \frac{M^{\frac{1}{2}}}{N^{\frac{1}{2}}}\|\nabla u_{\leq M}\|_{S^*_0(I)}\|v_{\geq N}\|_{S^*_0(I)}.
\end{align*}
Here the norm $\|\cdot\|_{S^*_0(I)}$ is defined as
$$
\|u\|_{S^*_0(I)}:=\|u\|_{L^{\infty}_tL^2_x(I\times \mathbb{R}^4)}+\|(i\partial_t+\Delta)u\|_{L^{\frac{3}{2}}_{t,x}(I\times \mathbb{R}^4)}.
$$
 }

{\bf Lemma 6.10(Paraproduct estimate).} {\it
$$
\||\nabla|^{-\frac{2}{3}}(fg)\|_{L^{\frac{4}{3}}_x}\lesssim \||\nabla|^{-\frac{2}{3}}f\|_{L^p_x}\||\nabla|^{\frac{2}{3}}g\|_{L^q_x}
$$
for any $1<p,q<+\infty$ satisfying $\frac{1}{p}+\frac{1}{q}=\frac{11}{12}$.
}

Inspired by \cite{Dodson2012,Visan2011}, we give the long-time Strichartz estimates below.

{\bf Proposition 6.11(Long-time Strichartz estimates).} {\it Let $(u,v):[0,T_{\max})\times \mathbb{R}^4\rightarrow \mathbb{C}\times \mathbb{C}$
be an almost periodic solution of (\ref{826x1}) with $(\alpha,\beta)=(0,0)$ and $\widetilde{N}(t)\equiv N_k\geq 1$ on each characteristic interval $J_k\subset [0,T_{\max})$. Then
for any frequency $N>0$,
\begin{align}
\|\nabla u_{\leq N}\|_{L^2_tL^4_x(I\times \mathbb{R}^4)}+\|\nabla v_{\leq N}\|_{L^2_tL^4_x(I\times \mathbb{R}^4)}\lesssim_{(u,v)}
1+N^{\frac{3}{2}}K^{\frac{1}{2}}\quad {\rm with}\quad K:=\int_I[\widetilde{N}(t)]^{-1}dt\label{92101}
\end{align}
is true on any compact time interval $I\subset [0,T_{\max})$, which is a union of contiguous intervals $J_k$. And for any $\eta>0$ there exists $N_0=N_0(\eta)>0$ such that for all $N\leq N_0$,
\begin{align}
\|\nabla u_{\leq N}\|_{L^2_tL^4_x(I\times \mathbb{R}^4)}+\|\nabla v_{\leq N}\|_{L^2_tL^4_x(I\times \mathbb{R}^4)}\lesssim_{(u,v)}
\eta[1+N^{\frac{3}{2}}K^{\frac{1}{2}}]. \label{92102}
\end{align}
Here the constant $N_0$ and the implicit constants in (\ref{92101}) and (\ref{92102}) are independent of the interval $I$.
}

{\bf Proof:} Let $I$ be a fixed compact time interval $I\subset [0,T_{\max})$, which is a union of contiguous intervals $J_k$. Unless specifying, we let all spacetime norms be on $I\times \mathbb{R}^4$ throughout the proof. Let $\eta_0$ be a small parameter to be chosen later. By the definition of almost periodicity symmtries, for this $\eta_0$, there exists $c_0=c_0(\eta_0)$ such that
\begin{align}
\|\nabla u_{\leq c_0\widetilde{N}(t)}\|_{L^{\infty}_tL^2_x(I\times \mathbb{R}^4)}+\|\nabla v_{\leq c_0\widetilde{N}(t)}\|_{L^{\infty}_tL^2_x(I\times \mathbb{R}^4)}\leq \eta_0.\label{92103}
\end{align}

For $N>0$, we define
$$
A(N):=\|\nabla u_{\leq N}\|_{L^2_tL^4_x(I\times \mathbb{R}^4)}+\|\nabla v_{\leq N}\|_{L^2_tL^4_x(I\times \mathbb{R}^4)}.
$$
By the result of Lemma 6.5,
\begin{align}
A(N)\lesssim_{(u,v)} 1+N^{\frac{3}{2}}K^{\frac{1}{2}} \quad {\rm for\ any}\quad N\geq \left(\frac{\int_I \widetilde{N}(t)^2dt}{\int_I [\widetilde{N}(t)]^{-1}dt}  \right)^{\frac{1}{3}},\label{92104}
\end{align}
especially for $N\geq N_{\max}:=\sup_{t\in I}\widetilde{N}(t)$.

For arbitrary frequencies $N>0$, we will establish the result by induction. First, we will establish a recurrence relation for $A(N)$. Using Strichartz's inequality, we get
\begin{align}
A(N)&\lesssim_{(u,v)} \inf_{t\in I}[\|\nabla u_{\leq N}(t)\|_{L^2_x}+\|\nabla v_{\leq N}(t)\|_{L^2_x}]\nonumber\\
&\quad +\|\nabla P_{\leq N}f_1(u,v)\|_{L^2_tL^{\frac{4}{3}}_x}
+\|\nabla P_{\leq N}f_2(u,v)\|_{L^2_tL^{\frac{4}{3}}_x}.\label{921x1}
\end{align}

We only give the details of estimating $\|\nabla P_{\leq N}f_1(u,v)\|_{L^2_tL^{\frac{4}{3}}_x}$ below because the estimating $\|\nabla P_{\leq N}f_2(u,v)\|_{L^2_tL^{\frac{4}{3}}_x}$ is entirely similar. To do this, we decompose $u=u_{\leq \frac{N}{\eta_0}}+u_{>\frac{N}{\eta_0}}$ and $v=v_{\leq \frac{N}{\eta_0}}+v_{>\frac{N}{\eta_0}}$, then further decompose
$u=u_{\leq c_0\widetilde{N}(t)}+u_{>c_0\widetilde{N}(t)}$ and $v=v_{\leq c_0\widetilde{N}(t)}+v_{>c_0\widetilde{N}(t)}$ and write
\begin{align}
\nabla f_1(u,v)&=\varnothing[\nabla (u_{>\frac{N}{\eta_0}}v^2_{>\frac{N}{\eta_0}})]+\varnothing[\nabla u_{>\frac{N}{\eta_0}}v_{>\frac{N}{\eta_0}}v_{\leq \frac{N}{\eta_0}}]
+\varnothing[\nabla u_{>\frac{N}{\eta_0}}v_{\leq \frac{N}{\eta_0}}v_{\leq c_0\widetilde{N}(t)}]\nonumber\\
&\quad+\varnothing[\nabla u_{>\frac{N}{\eta_0}}v_{\leq \frac{N}{\eta_0}}v_{>c_0\widetilde{N}(t)}]+\varnothing[\nabla u_{\leq \frac{N}{\eta_0}}v^2_{\leq c_0\widetilde{N}(t)}]+\varnothing[\nabla u_{\leq \frac{N}{\eta_0}}v^2_{>c_0\widetilde{N}(t)}]\nonumber\\
&\quad +\varnothing[\nabla u_{\leq \frac{N}{\eta_0}}v_{>c_0\widetilde{N}(t)}v_{\leq c_0\widetilde{N}(t)}]+\varnothing[\nabla v_{>\frac{N}{\eta_0}}v_{>\frac{N}{\eta_0}}u_{>\frac{N}{\eta_0}}]+\varnothing[\nabla v_{>\frac{N}{\eta_0}}v_{>\frac{N}{\eta_0}}u_{\leq \frac{N}{\eta_0}}]\nonumber\\
&\quad +\varnothing[\nabla v_{>\frac{N}{\eta_0}}v_{\leq \frac{N}{\eta_0}}u_{\leq c_0\widetilde{N}(t)}]+\varnothing[\nabla v_{>\frac{N}{\eta_0}}v_{\leq \frac{N}{\eta_0}}u_{>c_0\widetilde{N}(t)}]+\varnothing[\nabla v_{\leq \frac{N}{\eta_0}}v_{\leq c_0\widetilde{N}(t)}u_{\leq c_0\widetilde{N}(t)}]\nonumber\\
&\quad +\varnothing[\nabla v_{\leq \frac{N}{\eta_0}}v_{>c_0\widetilde{N}(t)}u_{>c_0\widetilde{N}(t)}]+\varnothing[\nabla v_{\leq \frac{N}{\eta_0}}v_{>c_0\widetilde{N}(t)}u_{\leq c_0\widetilde{N}(t)}]\nonumber\\
&\quad +\varnothing[\nabla v_{\leq \frac{N}{\eta_0}}v_{\leq c_0\widetilde{N}(t)}u_{>c_0\widetilde{N}(t)}]\nonumber\\
&:=(I1)+(I2)+(I3)+(I4)+(I5)+(I6)+(I7)\nonumber\\
&\quad +(II1)+(II2)+(II3)+(II4)+(II5)+(II6)+(II7)+(II8).\label{921x2}
 \end{align}
We will estimate the contributions of each of these terms to (\ref{921x2}).

Using Berstein inequality, H\"{o}lder's inequality and Sobolev embedding, the contribution of (I1) can be estimated as follows:
\begin{align}
&\quad \|P_{\leq N} \varnothing[\nabla (u_{>\frac{N}{\eta_0}}v^2_{>\frac{N}{\eta_0}})] \|_{L^2_tL^{\frac{4}{3}}_x}\lesssim N^{\frac{5}{3}}\||\nabla|^{-\frac{2}{3}} (u_{>\frac{N}{\eta_0}}v^2_{>\frac{N}{\eta_0}})\|_{L^2_tL^{\frac{4}{3}}_x}\nonumber\\
&\lesssim N^{\frac{5}{3}}\||\nabla|^{-\frac{2}{3}}u_{>\frac{N}{\eta_0}}\|_{L^2_tL^4_x}\||\nabla|^{\frac{2}{3}}(v^2_{>\frac{N}{\eta_0}})\|_{L^{\infty}_tL^{\frac{3}{2}}_x}\nonumber\\
&\lesssim N^{\frac{5}{3}}\||\nabla|^{-\frac{2}{3}}u_{>\frac{N}{\eta_0}}\|_{L^2_tL^4_x}
\||\nabla|^{\frac{2}{3}}(v_{>\frac{N}{\eta_0}})\|_{L^{\infty}_tL^{\frac{12}{5}}_x}
\|v_{>\frac{N}{\eta_0}}\|_{L^{\infty}_tL^4_x}\nonumber\\
&\lesssim N^{\frac{5}{3}}\||\nabla|^{-\frac{2}{3}}u_{>\frac{N}{\eta_0}}\|_{L^2_tL^4_x}\|v_{>\frac{N}{\eta_0}}\|^2_{L^{\infty}_t\dot{H}^1_x}\nonumber\\
&\lesssim_{(u,v)} \sum_{M>\frac{N}{\eta_0}}\left(\frac{N}{M}\right)^{\frac{5}{3}}A(M).\label{922s1}
\end{align}
The contributions of (I2), (II1) and (II2) can be obtained similarly as follows:
\begin{align}
\|P_{\leq N} \varnothing[\nabla u_{>\frac{N}{\eta_0}}v_{>\frac{N}{\eta_0}}v_{\leq \frac{N}{\eta_0}}]\|_{L^2_tL^{\frac{4}{3}}_x}&\lesssim_{(u,v)} \sum_{M>\frac{N}{\eta_0}}\left(\frac{N}{M}\right)^{\frac{5}{3}}A(M),\label{922s2}\\
\|P_{\leq N} \varnothing[\nabla v_{>\frac{N}{\eta_0}}v_{>\frac{N}{\eta_0}}u_{>\frac{N}{\eta_0}}]\|_{L^2_tL^{\frac{4}{3}}_x}&\lesssim_{(u,v)} \sum_{M>\frac{N}{\eta_0}}\left(\frac{N}{M}\right)^{\frac{5}{3}}A(M),\label{922s3}\\
\|P_{\leq N} \varnothing[\nabla v_{>\frac{N}{\eta_0}}v_{>\frac{N}{\eta_0}}u_{\leq \frac{N}{\eta_0}}]\|_{L^2_tL^{\frac{4}{3}}_x}&\lesssim_{(u,v)} \sum_{M>\frac{N}{\eta_0}}\left(\frac{N}{M}\right)^{\frac{5}{3}}A(M),\label{922s4}
\end{align}

The contribution of (I3) can be estimated as follows:
\begin{align}
&\quad \|P_{\leq N} \varnothing[\nabla (u_{>\frac{N}{\eta_0}}v_{\leq \frac{N}{\eta_0}}v_{\leq c_0\widetilde{N}(t)})]\|_{L^2_tL^{\frac{4}{3}}_x}\lesssim N^{\frac{2}{3}}\||\nabla|^{-\frac{2}{3}} (u_{>\frac{N}{\eta_0}}v_{\leq \frac{N}{\eta_0}}v_{\leq c_0\widetilde{N}(t)})\|_{L^2_tL^{\frac{4}{3}}_x}\nonumber\\
&\lesssim N^{\frac{2}{3}}\||\nabla|^{\frac{2}{3}}v_{\leq \frac{N}{\eta_0}}\|_{L^2_tL^4_x}\||\nabla|^{-\frac{2}{3}}(\nabla u_{>\frac{N}{\eta_0}}v_{\leq c_0\widetilde{N}(t)})\|_{L^{\infty}_tL^{\frac{4}{3}}_x}\nonumber\\
&\lesssim N^{\frac{2}{3}}\||v_{\leq\frac{N}{\eta_0}}\|_{L^2_tL^4_x}
\||\nabla|^{\frac{2}{3}}v_{\leq c_0\widetilde{N}(t)}\|_{L^{\infty}_tL^{\frac{12}{5}}_x}
\||\nabla|^{\frac{1}{3}}(u_{>\frac{N}{\eta_0}})\|_{L^{\infty}_tL^2_x}\nonumber\\
&\lesssim_{(u,v)} N^{\frac{2}{3}}A(\frac{N}{\eta_0})\eta_0(\frac{N}{\eta_0})^{-\frac{2}{3}}\lesssim_{(u,v)} \eta_0^{\frac{5}{3}}A(\frac{N}{\eta_0}).\label{922s5}
\end{align}

The contribution of (I4) is
\begin{align}
\|P_{\leq N}\varnothing[\nabla u_{>\frac{N}{\eta_0}}v_{\leq \frac{N}{\eta_0}}v_{>c_0\widetilde{N}(t)}]\|_{L^2_xL^{\frac{4}{3}}_x}\lesssim_{(u,v)}
\frac{N^{\frac{3}{2}}K^{\frac{1}{2}}}{\eta_0^{\frac{1}{2}}c_0^{\frac{3}{2}}}\sup_{J_k\subset I}\|\nabla v_{\leq \frac{N}{\eta_0}}\|_{S^*_0(J_k)}.\label{922s6}
\end{align}

The contribution of (I5) is
\begin{align}
\|P_{\leq N}\varnothing[\nabla u_{\leq \frac{N}{\eta_0}}v^2_{\leq c_0\widetilde{N}(t)}]\|_{L^2_xL^{\frac{4}{3}}_x}\lesssim
\|\nabla u_{\leq \frac{N}{\eta_0}}\|_{L^2_tL^4_x}\|v_{\leq c_0\widetilde{N}(t)}\|^2_{L^{\infty}_tL^4_x}\lesssim_{(u,v)}\eta^2_0A(\frac{N}{\eta_0}).\label{922s7}
\end{align}

In convenience, we divide the time interval $I$ into subintervals $J_k$ where $\widetilde{N}(t)$ is constant and use the bilinear Strichartz estimate in Lemma 6.9 on each of these subintervals.  The contribution of (I6) is
\begin{align}
&\quad \|P_{\leq N}\varnothing[\nabla u_{\leq \frac{N}{\eta_0}}v^2_{>c_0\widetilde{N}(t)}]\|_{L^2_xL^{\frac{4}{3}}_x}\lesssim
\|v_{>\frac{N}{\eta_0}}\|_{L^{\infty}_tL^4_x}\|
\|\nabla u_{\leq \frac{N}{\eta_0}}v_{>c_0\widetilde{N}(t)}\|^2_{L_{t,x}}\nonumber\\
&\lesssim_{(u,v)}\frac{N^{\frac{1}{2}}K^{\frac{1}{2}}}{\eta_0^{\frac{1}{2}}c_0^{\frac{3}{2}}}\sup_{J_k\subset I}\|\Delta u_{\leq \frac{N}{\eta_0}}\|_{S^*_0(J_k)}\lesssim_{(u,v)}\frac{N^{\frac{3}{2}}K^{\frac{1}{2}}}{\eta_0^{\frac{3}{2}}c_0^{\frac{3}{2}}}\sup_{J_k\subset I}\|\Delta u_{\leq \frac{N}{\eta_0}}\|_{S^*_0(J_k)}.\label{922s8}
\end{align}

The contribution of (I7) can be controlled by the sum of those of (I5) and  (I6) because
$$|v_{>c_0\widetilde{N}(t)}v_{\leq c_0\widetilde{N}(t)}|\leq |v^2_{>c_0\widetilde{N}(t)}|+ |v^2_{\leq c_0\widetilde{N}(t)}|.$$
The contributions of (II3), (II4), (II5) and (II6) are similar to those of (I3), (I4), (I5) and (I6) respectively, while the contributions of (II7) and
(II8)  can be controlled by the sum of those of (II5) and  (II6). Therefore, we can put all of the results above together and get the following recurrence relation
\begin{align}
A(N)&\lesssim_{(u,v)}\inf_{t\in I}[\|\nabla u_{\leq N}(t)\|_{L^2_x}+\|\nabla v_{\leq N}(t)\|_{L^2_x}]+\sum_{M>\frac{N}{\eta_0}}\left(\frac{N}{M}\right)^{\frac{5}{3}}A(M)\nonumber\\
&\quad +\frac{N^{\frac{3}{2}}K^{\frac{1}{2}}}{\eta_0^{\frac{3}{2}}c_0^{\frac{3}{2}}}\sup_{J_k\subset I}[\|\nabla u_{\leq \frac{N}{\eta_0}}\|_{S^*_0(J_k)}+\|\nabla v_{\leq \frac{N}{\eta_0}}\|_{S^*_0(J_k)}].\label{922s9}
\end{align}

First, we address (\ref{92101}). Note that (\ref{92104}) holds for $N\geq N_{\max}$. That is, there exists some constant $C(u,v)>0$ such that
\begin{align}
A(N)\lesssim_{(u,v)}C(u,v)[1+N^{\frac{3}{2}}K^{\frac{1}{2}}]\quad {\rm for \ all}\quad N\geq N_{\max}.\label{92301}
\end{align}
Rewriting (\ref{922s9}) as
\begin{align}
A(N)\lesssim_{(u,v)} \tilde{C}(u,v)\left\{1+\frac{N^{\frac{3}{2}}K^{\frac{1}{2}}}{\eta_0^{\frac{3}{2}}c_0^{\frac{3}{2}}}
+\sum_{M>\frac{N}{\eta_0}}\left(\frac{N}{M}\right)^{\frac{5}{3}}A(M) \right\},\label{92302}
\end{align}
through halving the frequency $N$ at each step, we can prove (\ref{92101}) inductively. Without loss of generality, we assume that (\ref{92301})
holds for frequencies larger or equal to $N$, and apply (\ref{92302}) with $\eta_0\leq \frac{1}{2}$ such that
\begin{align*}
A(\frac{N}{2})&\leq \tilde{C}(u,v)\left\{1+\frac{(\frac{N}{2})^{\frac{3}{2}}K^{\frac{1}{2}}}{\eta_0^{\frac{3}{2}}c_0^{\frac{3}{2}}}
+\sum_{M>\frac{N}{2\eta_0}}\left(\frac{N}{2M}\right)^{\frac{5}{3}}A(M) \right\}\nonumber\\
&\leq \tilde{C}(u,v)\left\{1+\frac{(\frac{N}{2})^{\frac{3}{2}}K^{\frac{1}{2}}}{\eta_0^{\frac{3}{2}}c_0^{\frac{3}{2}}}
+\sum_{M>\frac{N}{2\eta_0}}\left(\frac{N}{2M}\right)^{\frac{5}{3}}A(M) \right\}.
\end{align*}
If we chose $\eta_0=\eta_0(u,v)$ small enough to satisfies $\eta_0^{\frac{1}{6}}\tilde{C}(u,v)\leq \frac{1}{2}$, then
\begin{align*}
A(\frac{N}{2})&\leq \frac{1}{2}C(u,v)\left\{1+(\frac{N}{2})^{\frac{3}{2}}K^{\frac{1}{2}}\right\}+ \tilde{C}(u,v)\left\{1+\frac{(\frac{N}{2})^{\frac{3}{2}}K^{\frac{1}{2}}}{\eta_0^{\frac{3}{2}}c_0^{\frac{3}{2}}}\right\}.
\end{align*}
Taking $C(u,v)\geq 2\tilde{C}(u,v)\eta_0^{-\frac{3}{2}}c_0^{-\frac{3}{2}}$, we can get (\ref{92101}).

Now we turn to (\ref{92102}). To show the small constant $\eta$, we will prove the following lemma.

{\bf Lemma 6.12(Vanishing of the small frequencies).} {\it Under the assumptions of Proposition 6.11, there holds
\begin{align*}
f(N)&:=\|\nabla u_{\leq N}\|_{L^{\infty}_tL^2_x([0,T_{\max})\times \mathbb{R}^4)}+\|\nabla v_{\leq N}\|_{L^{\infty}_tL^2_x([0,T_{\max})\times \mathbb{R}^4)}\\
&\qquad+\sup_{J_k\subset [0,T_{\max})}[\|\nabla u_{\leq N}\|_{S^*_0(J_k)}+\|\nabla v_{\leq N}\|_{S^*_0(J_k)}]\rightarrow 0\quad {\rm as}\quad N\rightarrow 0.
\end{align*}
}

{\bf Proof:} Under the assumption of Proposition 6.11, by compactness, we obtain
\begin{align}
\lim_{N\rightarrow 0}[\|\nabla u_{\leq N}\|_{L^{\infty}_tL^2_x([0,T_{\max})\times \mathbb{R}^4)}+\|\nabla v_{\leq N}\|_{L^{\infty}_tL^2_x([0,T_{\max})\times \mathbb{R}^4)}]=0.\label{923s1}
\end{align}

Fixing a characteristic interval $J_k\subset [0, T_{\max})$, by the spacetime bounds in Lemma 6.5, we have
\begin{align*}
\|\nabla u\|_{L^2_tL^4_x}+\|u\|_{L^3_tL^{12}_x}+\|u\|_{L^6_{t,x}}+\|\nabla v\|_{L^2_tL^4_x}+\|v\|_{L^3_tL^{12}_x}+\|v\|_{L^6_{t,x}}\lesssim_{(u,v)} 1.
\end{align*}
All spacetime norms in the estimates above are on $J_k\times \mathbb{R}^4$. For any frequency $N$, decomposing $u=u_{\leq N^{\frac{1}{2}}}+u_{>N^{\frac{1}{2}}}$ and $v=v_{\leq N^{\frac{1}{2}}}+v_{>N^{\frac{1}{2}}}$,  we get
\begin{align*}
&\quad\ \|\nabla u_{\leq N}\|_{S^*_0(J_k)}+\|\nabla v_{\leq N}\|_{S^*_0(J_k)}\nonumber\\
&=\|\nabla u_{\leq N}\|_{L^{\infty}_tL^2_x}+\|\nabla v_{\leq N}\|_{L^{\infty}_tL^2_x}+\|\nabla P_{\leq N}f_1(u,v)\|_{L^{\frac{3}{2}}_{t,x}}+\|\nabla P_{\leq N}f_2(u,v)\|_{L^{\frac{3}{2}}_{t,x}}\\
&\lesssim_{(u,v)}\|\nabla u_{\leq N}\|_{L^{\infty}_tL^2_x}+\|\nabla v_{\leq N}\|_{L^{\infty}_tL^2_x}
+\|\nabla P_{\leq N}f_1(u_{>N^{\frac{1}{2}}},v_{>N^{\frac{1}{2}}})\|_{L^{\frac{3}{2}}_{t,x}}\nonumber\\
&\quad+\|\varnothing(\nabla u_{>N^{\frac{1}{2}}}v_{\leq N^{\frac{1}{2}}}v\|_{L^{\frac{3}{2}}_{t,x}}+\|\varnothing(\nabla v_{>N^{\frac{1}{2}}}u_{\leq N^{\frac{1}{2}}}v\|_{L^{\frac{3}{2}}_{t,x}}+\|\varnothing(\nabla u_{>N^{\frac{1}{2}}}v_{\leq N^{\frac{1}{2}}}v_{>N^{\frac{1}{2}}}\|_{L^{\frac{3}{2}}_{t,x}}\nonumber\\
&\quad+\|\varnothing(\nabla v_{>N^{\frac{1}{2}}}u_{\leq N^{\frac{1}{2}}}v_{>N^{\frac{1}{2}}}\|_{L^{\frac{3}{2}}_{t,x}}+\|\varnothing(\nabla u_{\leq N^{\frac{1}{2}}}v^2\|_{L^{\frac{3}{2}}_{t,x}}+\|\varnothing(\nabla v_{\leq N^{\frac{1}{2}}}vu\|_{L^{\frac{3}{2}}_{t,x}}\nonumber\\
&\quad+\|\nabla P_{\leq N}f_2(u_{>N^{\frac{1}{2}}},v_{>N^{\frac{1}{2}}})\|_{L^{\frac{3}{2}}_{t,x}}+\|\varnothing(\nabla v_{>N^{\frac{1}{2}}}u_{\leq N^{\frac{1}{2}}}u\|_{L^{\frac{3}{2}}_{t,x}}+\|\varnothing(\nabla u_{>N^{\frac{1}{2}}}v_{\leq N^{\frac{1}{2}}}u\|_{L^{\frac{3}{2}}_{t,x}}\nonumber\\
&\quad+\|\varnothing(\nabla v_{>N^{\frac{1}{2}}}u_{\leq N^{\frac{1}{2}}}u_{>N^{\frac{1}{2}}}\|_{L^{\frac{3}{2}}_{t,x}}
+\|\varnothing(\nabla u_{>N^{\frac{1}{2}}}v_{\leq N^{\frac{1}{2}}}u_{>N^{\frac{1}{2}}}\|_{L^{\frac{3}{2}}_{t,x}}+\|\varnothing(\nabla v_{\leq N^{\frac{1}{2}}}u^2\|_{L^{\frac{3}{2}}_{t,x}}\nonumber\\
&\quad+\|\varnothing(\nabla u_{\leq N^{\frac{1}{2}}}vu\|_{L^{\frac{3}{2}}_{t,x}}\nonumber\displaybreak\\
&:=\|\nabla u_{\leq N}\|_{L^{\infty}_tL^2_x}+\|\nabla v_{\leq N}\|_{L^{\infty}_tL^2_x}+(I1)+(I2)+(I3)+(I4)+(I5)+(I6)+(I7)\nonumber\\
&\qquad +(II1)+(II2)+(II3)+(II4)+(II5)+(II6)+(II7).
\end{align*}
Using H\"{o}lder's inequality and Bernstein's inequality,
\begin{align*}
(I1)&\lesssim N\|u_{>N^{\frac{1}{2}}}\|_{L^2_tL^4_x}\|v_{>N^{\frac{1}{2}}}\|_{L^6_{t,x}}\|v_{>N^{\frac{1}{2}}}\|_{L^{\infty}_tL^4_x}\lesssim_{(u,v)} N^{\frac{1}{2}},\\
(I2)&\lesssim \|\nabla u_{>N^{\frac{1}{2}}}\|_{L^2_tL^4_x}\|v_{\leq N^{\frac{1}{2}}}\|_{L^{\infty}_tL^4_x}\|v\|_{L^6_{t,x}}\lesssim_{(u,v)} N^{\frac{1}{2}},\\
(I6)&\lesssim \|\nabla u_{\leq N^{\frac{1}{2}}}\|_{L^{\infty}_tL^2_x}\|v\|^2_{L^3_tL^{12}_x}\lesssim_{(u,v)}\|\nabla u_{\leq N^{\frac{1}{2}}}\|_{L^{\infty}_tL^2_x}.
\end{align*}
The estimate for (II1) is similar to that of (I1). The estimates for (I3), (I4), (I5), (II2), (II3), (II4) and (II5) are similar to that of (I2), they are all bounded by $N^{\frac{1}{2}}$. The estimate for (II6) is similar to that of (I6). While (I7) and (II7) can be controlled by the sum of (I6) and (II6) because
$|uv|\leq |u|^2+|v|^2$. Therefore,
\begin{align*}
&\quad\ \|\nabla u_{\leq N}\|_{S^*_0(J_k)}+\|\nabla v_{\leq N}\|_{S^*_0(J_k)}\nonumber\\
&\lesssim_{(u,v)} \|\nabla u_{\leq N}\|_{L^{\infty}_tL^2_x}+\|\nabla v_{\leq N}\|_{L^{\infty}_tL^2_x}+N^{\frac{1}{2}}+\|\nabla u_{\leq N^{\frac{1}{2}}}\|_{L^{\infty}_tL^2_x}+\|\nabla v_{\leq N^{\frac{1}{2}}}\|_{L^{\infty}_tL^2_x}
\end{align*}
All spacetime norms in the estimates above are also on $J_k\times \mathbb{R}^4$. Since $J_k\subset [0,T_{\max})$ is arbitrary, we get
\begin{align*}
&\sup_{J_k\subset [0,T_{\max})}[\|\nabla u_{\leq N}\|_{S^*_0(J_k)}+\|\nabla v_{\leq N}\|_{S^*_0(J_k)}\nonumber\\
&\lesssim_{(u,v)} N^{\frac{1}{2}}+\|\nabla u_{\leq N}\|_{L^{\infty}_tL^2_x([0,T_{\max})\times\mathbb{R}^4)}+\|\nabla v_{\leq N}\|_{L^{\infty}_tL^2_x([0,T_{\max})\times\mathbb{R}^4)}\nonumber\\
&\quad +\|\nabla u_{\leq N^{\frac{1}{2}}}\|_{L^{\infty}_tL^2_x([0,T_{\max})\times\mathbb{R}^4)}+\|\nabla v_{\leq N^{\frac{1}{2}}}\|_{L^{\infty}_tL^2_x([0,T_{\max})\times\mathbb{R}^4)}.
\end{align*}
This combines (\ref{923s1}) prove the claim of this lemma.\hfill $\Box$

Now we come back to the proof of (\ref{92102}). Using (\ref{92101}) and Lemma 6.12 , (\ref{922s9}) means that
\begin{align*}
A(N)&\lesssim_{(u,v)}f(N)+\frac{N^{\frac{3}{2}}K^{\frac{1}{2}}}{\eta_0^{\frac{3}{2}}c_0^{\frac{3}{2}}}f(N)
+\sum_{M>\frac{N}{\eta_0}}\left(\frac{N}{M}\right)^{\frac{5}{3}}A(M)\nonumber\\
&\lesssim_{(u,v)}f(N)+\eta^{\frac{5}{3}}_0+\left(\frac{f(N)}{\eta_0^{\frac{3}{2}}c_0^{\frac{3}{2}}}+\eta^{\frac{1}{6}}_0\right)N^{\frac{3}{2}}K^{\frac{1}{2}}.
\end{align*}
Hence, for any $\eta>0$, first we can chose $\eta_0=\eta_0(\eta)$ satisfying $\eta^{\frac{1}{6}}_0\leq \eta$, then $N_0=N_0(\eta)$ such that
$\frac{f(N)}{\eta_0^{\frac{3}{2}}c_0^{\frac{3}{2}}}\leq \eta$, and obtain
$$
A(N)\lesssim_{(u,v)} \eta(1+N^{\frac{3}{2}}K^{\frac{1}{2}})\quad {\rm for \ all}\quad N\leq N_0.
$$
Proposition 6.12 is proved.\hfill $\Box$

{\bf Corollary 6.13(Low and high frequencies control).}  {\it Let $(u,v):[0,T_{\max})\times \mathbb{R}^4\rightarrow \mathbb{C}\times \mathbb{C}$
be an almost periodic solution of (\ref{826x1}) with $\widetilde{N}(t)\equiv N_k\geq 1$ on each characteristic interval $J_k\subset [0,T_{\max})$. Then
for any frequency $N>0$ and the admissible pair $(q,r)$
\begin{align}
\| u_{\geq N}\|_{L^q_tL^r_x(I\times \mathbb{R}^4)}+\| v_{\geq N}\|_{L^q_tL^r_x(I\times \mathbb{R}^4)}\lesssim_{(u,v)}
N^{-1}\left(1+N^3K\right)^{\frac{1}{q}},\quad 3<q\leq +\infty\label{923z1}
\end{align}
is true on any compact time interval $I\subset [0,T_{\max})$, which is a union of contiguous intervals $J_k$. And for any $\eta>0$ there exists $N_0=N_0(\eta)>0$ such that for all $N\leq N_0$,
\begin{align}
\|\nabla u_{\leq N}\|_{L^q_tL^r_x(I\times \mathbb{R}^4)}+\|\nabla v_{\leq N}\|_{L^q_tL^r_x(I\times \mathbb{R}^4)}\lesssim_{(u,v)}
\eta\left(1+N^3K\right)^{\frac{1}{q}},\quad 2\leq q\leq +\infty. \label{923z2}
\end{align}
Here the constant $N_0$ and the implicit constants in (\ref{923z1}) and (\ref{923z2}) are independent of the interval $I$.
}

{\bf Proof:} First, we address (\ref{923z1}). Using (\ref{92101}) and Berstein's inequality, we obtain
\begin{align*}
&\quad \||\nabla|^{-\frac{1}{2}-\epsilon}u_{\geq N}\|_{L^2_tL^4_x(I\times \mathbb{R}^4)}+\||\nabla|^{-\frac{1}{2}-\epsilon}v_{\geq N}\|_{L^2_tL^4_x(I\times \mathbb{R}^4)}\nonumber\\
&\lesssim \sum_{M\geq N}M^{-\frac{3}{2}-\epsilon}[\|\nabla u_M\|_{L^2_tL^4_x(I\times \mathbb{R}^4)}+\|\nabla v_M\|_{L^2_tL^4_x(I\times \mathbb{R}^4)}]\nonumber\\
&\lesssim_{(u,v)} \sum_{M\geq N}M^{-\frac{3}{2}-\epsilon}(1+M^{\frac{3}{2}}K^{\frac{1}{2}})\lesssim_{(u,v)} N^{-\frac{3}{2}-\epsilon}(1+N^3K)^{\frac{1}{2}}
\end{align*}
for any $\epsilon>0$ and any frequencies $N>0$. Interpolating with the energy bound, we have
\begin{align*}
&\quad \|u_{\geq N}\|_{L^q_tL^r_x(I\times \mathbb{R}^4)}+\| v_{\geq N}\|_{L^q_tL^r_x(I\times \mathbb{R}^4)}\nonumber\\
&\lesssim_{(u,v)}\||\nabla|^{-\frac{1}{2}-\frac{q-3}{2}}u_{\geq N}\|^{\frac{2}{q}}_{L^2_tL^4_x(I\times \mathbb{R}^4)}\|\nabla u_{\geq N}\|^{1-\frac{2}{q}}_{L^{\infty}_tL^2_x(I\times \mathbb{R}^4)}\nonumber\\
&\quad +\||\nabla|^{-\frac{1}{2}-\frac{q-3}{2}}v_{\geq N}\|^{\frac{2}{q}}_{L^2_tL^4_x(I\times \mathbb{R}^4)}\|\nabla v_{\geq N}\|^{1-\frac{2}{q}}_{L^{\infty}_tL^2_x(I\times \mathbb{R}^4)}\nonumber\\
&\lesssim_{(u,v)} N^{-1}(1+N^3K)^{\frac{1}{q}},
\end{align*}
where $(q,r)$ is admissible pair with $3<q\leq +\infty$.

Now it is the turn of (\ref{923z2}). Since $\inf_{t\in I}\widetilde{N}(t)\geq 1$, by Remark, for any $\eta>0$, there exists $N_0(\eta)$ such that for all $N\leq N_0$,
$$
\|\nabla u_{\leq N}\|_{L^{\infty}_tL^2_x(I\times \mathbb{R}^4)}+\|\nabla v_{\leq N}\|_{L^{\infty}_tL^2_x(I\times \mathbb{R}^4)}\leq \eta.
$$
Interpolating with (\ref{92102}), we obtain the claim.\hfill $\Box$

\subsection{The rapid frequency-cascade scenario}
\qquad In this subsection, we preclude the existence of almost periodic solutions as in Proposition 6.6 for which $\int^{T_{\max}}_0\widetilde{N}(t)^{-1}dt<+\infty$.

{\bf Proposition 6.14(No rapid frequency-cascades).} {\it There are no almost periodic solutions $(u,v):[0,T_{\max})\times \mathbb{R}^4\rightarrow \mathbb{C}\times \mathbb{C}$ with $\widetilde{N}(t)\equiv N_k\geq 1$ on each characteristic interval $J_k\subset [0, T_{\max})$ satisfying $\|u\|_{L^6_{t,x}([0,T_{\max})\times \mathbb{R}^4)}+\|v\|_{L^6_{t,x}([0,T_{\max})\times \mathbb{R}^4)}=+\infty$ and
\begin{align}
\int^{T_{\max}}_0 \widetilde{N}(t)^{-1}dt<+\infty.\label{923x1}
\end{align}
}

{\bf Proof:} Contradictorily, assume that $(u,v)$ is such a solution. By Corollary 6.4, whether $T_{\max}$ is finite or infinite,
\begin{align}
\lim_{t\rightarrow T_{\max}}\widetilde{N}(t)=+\infty.\label{923x2}
\end{align}
By compactness, we have
\begin{align}
\lim_{t\rightarrow T_{\max}}[\|\nabla u_{\leq N}(t)\|_{L^2_x}+\|\nabla v_{\leq N}(t)\|_{L^2_x}]=0\quad {\rm for \ any} \quad N>0.\label{923x3}
\end{align}

Let $I_n$ be a nested sequence of compact subintervals of $[0,T_{\max})$ which are unions of contiguous characteristic subintervals $J_k$. Applying Proposition  6.11 on each $I_n$, noticing (\ref{922s9}) and (\ref{923x1}), we obtain
\begin{align*}
A_n(N)&:=\|\nabla u_{\leq N}\|_{L^2_tL^4_x(I_n\times \mathbb{R}^4)}+\|\nabla v_{\leq N}\|_{L^2_tL^4_x(I_n\times \mathbb{R}^4)}\nonumber\\
&\lesssim_{(u,v)}\inf_{t\in I}[\|\nabla u_{\leq N}(t)\|_{L^2_x}+\|\nabla v_{\leq N}(t)\|_{L^2_x}]+\sum_{M>\frac{N}{\eta_0}}\left(\frac{N}{M}\right)^{\frac{5}{3}}A_n(M)\nonumber\\
&\quad +\frac{N^{\frac{3}{2}}}{\eta_0^{\frac{3}{2}}c_0^{\frac{3}{2}}}\left(\int_0^{T_{\max}}[\widetilde{N}(t)]^{-1}dt\right)^{\frac{1}{2}}\nonumber\\
&\lesssim_{(u,v)}\inf_{t\in I}[\|\nabla u_{\leq N}(t)\|_{L^2_x}+\|\nabla v_{\leq N}(t)\|_{L^2_x}]+\sum_{M>\frac{N}{\eta_0}}\left(\frac{N}{M}\right)^{\frac{5}{3}}A_n(M)+\frac{N^{\frac{3}{2}}}{\eta_0^{\frac{3}{2}}c_0^{\frac{3}{2}}}
\end{align*}
for all frequencies $N>0$. Similar to the proof of (\ref{92101}), we get, for all $N>0$,
\begin{align*}
\|\nabla u_{\leq N}\|_{L^2_tL^4_x(I_n\times \mathbb{R}^4)}+\|\nabla v_{\leq N}\|_{L^2_tL^4_x(I_n\times \mathbb{R}^4)}
\lesssim_{(u,v)} \inf_{t\in I_n}[\|\nabla u_{\leq N}(t)\|_{L^2_x}+\|\nabla v_{\leq N}(t)\|_{L^2_x}]+N^{\frac{3}{2}}.
\end{align*}
Letting $n\rightarrow +\infty$ and recalling (\ref{923x3}), we have
\begin{align}
\|\nabla u_{\leq N}\|_{L^2_tL^4_x([0,T_{\max})\times \mathbb{R}^4)}+\|\nabla v_{\leq N}\|_{L^2_tL^4_x([0,T_{\max})\times \mathbb{R}^4)}
\lesssim_{(u,v)} N^{\frac{3}{2}}.\label{923x4}
\end{align}

Next, we hope to show that (\ref{923x4}) implies, for all $N>0$,
\begin{align}
\|\nabla u_{\leq N}\|_{L^{\infty}_tL^2_x([0,T_{\max})\times \mathbb{R}^4)}+\|\nabla v_{\leq N}\|_{L^{\infty}_tL^2_x([0,T_{\max})\times \mathbb{R}^4)}
\lesssim_{(u,v)} N^{\frac{3}{2}}.\label{923x5}
\end{align}

Decomposing $u=u_{\leq N}+u_{>N}$ and $v=v_{\leq N}+v_{>N}$, by the no-waste Duhamel formulae in Proposition 6.9, using Strichartz's inequality, Bernstein's inequality, H\"{o}lder's inequality and Sobolev embedding, we get
\begin{align*}
&\quad \|\nabla u_{\leq N}\|_{L^{\infty}_tL^2_x}+\|\nabla v_{\leq N}\|_{L^{\infty}_tL^2_x}\lesssim \|\nabla P_{\leq N}f_1(u,v)\|_{L^2_tL^{\frac{4}{3}}_x}+\|\nabla P_{\leq N}f_2(u,v)\|_{L^2_tL^{\frac{4}{3}}_x}\\
&\lesssim \|\nabla P_{\leq N}f_1(u_{\leq N},v_{\leq N})\|_{L^2_tL^{\frac{4}{3}}_x}
+\|\nabla P_{\leq N}\varnothing(u_{\leq N}v_{\leq N}v_{>N})\|_{L^2_tL^{\frac{4}{3}}_x}
+\|\nabla P_{\leq N}\varnothing(u_{\leq N}v_{>N}v_{>N})\|_{L^2_tL^{\frac{4}{3}}_x}\nonumber\\
&\quad+ \|\nabla P_{\leq N}u_{>N}v^2\|_{L^2_tL^{\frac{4}{3}}_x}+\|\nabla P_{\leq N}f_2(u_{\leq N},v_{\leq N})\|_{L^2_tL^{\frac{4}{3}}_x}
+\|\nabla P_{\leq N}\varnothing(v_{\leq N}u_{\leq N}u_{>N})\|_{L^2_tL^{\frac{4}{3}}_x}\nonumber\\
&\quad+\|\nabla P_{\leq N}\varnothing(v_{\leq N}u_{>N}u_{>N})\|_{L^2_tL^{\frac{4}{3}}_x}\nonumber\\
&\lesssim_{(u,v)} [\|\nabla u\|_{L^2_tL^4_x}+\|\nabla v\|_{L^2_tL^4_x}][\|u_{\leq N}\|^2_{L^{\infty}_tL^4_x}+\|v_{\leq N}\|^2_{L^{\infty}_tL^4_x}+\|u_{>N}\|^2_{L^{\infty}_tL^4_x}+\|v_{>N}\|^2_{L^{\infty}_tL^4_x}]\nonumber\\
&\quad+ N^{\frac{5}{3}}[\||\nabla|^{-\frac{2}{3}}u_{>N}\|_{L^2_tL^4_x}+\||\nabla|^{-\frac{2}{3}}v_{>N}\|_{L^2_tL^4_x}]
[\||\nabla|^{\frac{2}{3}}u\|_{L^{\infty}_tL^{\frac{12}{5}}_x}+\||\nabla|^{\frac{2}{3}}v\|_{L^{\infty}_tL^{\frac{12}{5}}_x}]\nonumber\\
&\qquad \qquad \qquad \times [\|u\|_{L^{\infty}_tL^4_x}+\|v\|_{L^{\infty}_tL^4_x}]\nonumber\\
&\lesssim_{(u,v)} N^{\frac{3}{2}}+N^{\frac{5}{3}}\sum_{M>N}M^{-\frac{5}{3}}[\|\nabla u_M\|_{L^2_tL^4_x}+\|\nabla v_M\|_{L^2_tL^4_x}]\nonumber\\
&\lesssim_{(u,v)} N^{\frac{3}{2}}+N^{\frac{5}{3}}N^{-\frac{1}{6}}\lesssim_{(u,v)} N^{\frac{3}{2}}.
\end{align*}
Here all spacetime norms in the estimates above are on $[0,T_{\max})\times \mathbb{R}^4$. (\ref{923x5}) is proved.

Now we are ready to finish the proof of this proposition. Using Bernstein inequality and (\ref{923x5}), we have
\begin{align*}
&\quad \||\nabla|^{-\frac{1}{4}}u\|_{L^{\infty}_tL^2_x}+\||\nabla|^{-\frac{1}{4}}v\|_{L^{\infty}_tL^2_x}\nonumber\\
&\lesssim  \||\nabla|^{-\frac{1}{4}}u_{>1}\|_{L^{\infty}_tL^2_x}+\||\nabla|^{-\frac{1}{4}}u_{\leq 1}\|_{L^{\infty}_tL^2_x}+\||\nabla|^{-\frac{1}{4}}v_{>1}\|_{L^{\infty}_tL^2_x}+\||\nabla|^{-\frac{1}{4}}v_{\leq 1}\|_{L^{\infty}_tL^2_x}\nonumber\\
&\lesssim_{(u,v)}\sum_{N>1}N^{-\frac{5}{4}}+\sum_{N \leq 1}N^{\frac{1}{4}}\lesssim_{(u,v)} 1,
\end{align*}
which means $u, v\in L^{\infty}_t\dot{H}^{-\frac{1}{4}}_x([0,T_{\max})\times \mathbb{R}^4)$.

Fix $t\in [0,T_{\max})$ and let $\eta>0$ be a small constant. By compactness, there exists $c(\eta)>0$ such that
$$
\int_{|\xi|\leq c(\eta)\widetilde{N}(t)}|\xi|^2|\hat{u}(t,\xi)|^2d\xi+\int_{|\xi|\leq c(\eta)\widetilde{N}(t)}|\xi|^2|\hat{v}(t,\xi)|^2d\xi\leq \eta.
$$
Interpolating with $u\in L^{\infty}_t\dot{H}^{-\frac{1}{4}}_x([0,T_{\max})\times \mathbb{R}^4)$ and $v\in L^{\infty}_t\dot{H}^{-\frac{1}{4}}_x([0,T_{\max})\times \mathbb{R}^4)$, we obtain
\begin{align}
\int_{|\xi|\leq c(\eta)\hat{N}(t)}[|\hat{u}(t,\xi)|^2+|\hat{v}(t,\xi)|^2]d\xi\lesssim_{(u,v)} \eta^{\frac{1}{5}}.\label{924s1}
\end{align}
On the other hand,
\begin{align}
&\quad \int_{|\xi|\geq c(\eta)\widetilde{N}(t)}[|\hat{u}(t,\xi)|^2+|\hat{v}(t,\xi)|^2]d\xi\leq [c(\eta)\widetilde{N}(t)]^{-2}\int_{\mathbb{R}^4}|\xi|^2[|\hat{u}(t,\xi)|^2+|\hat{v}(t,\xi)|^2]d\xi\nonumber\\
&\lesssim_{(u,v)} [c(\eta)\widetilde{N}(t)]^{-2}.\label{924s2}
\end{align}
Using Plancherel's theorem, (\ref{924s1}) and (\ref{924s2}), we get
$$
0\leq M(u,v)(t):=\int_{\mathbb{R}^4}[|u(t,x)|^2+|v(t,x)|^2]dx=\int_{\mathbb{R}^4}[|\tilde{u}(t,\xi)|^2+|\hat{v}(t,\xi)|^2]d\xi
\lesssim_{(u,v)} \eta^{\frac{1}{5}}+[c(\eta)\widetilde{N}(t)]^{-2}
$$
for all $t\in [0,T_{\max})$. Recalling (\ref{923x2}) and letting $\eta\rightarrow 0$, we know that $M(u,v)\equiv 0$, which implies that $(u,v)\equiv 0$.
It is a contradiction to $\|u\|_{L^6_{t,x}([0,T_{\max})\times \mathbb{R}^4)}+\|v\|_{L^6_{t,x}([0,T_{\max})\times \mathbb{R}^4)}=+\infty$. Proposition 6.14 is proved.\hfill $\Box$

\subsection{The quasi-solution scenario}
\qquad In this subsection, we preclude the existence of almost periodic solutions as in Proposition 6.6 for which $\int^{T_{\max}}_0\widetilde{N}(t)^{-1}dt=+\infty$.
If there exists such a solution, we will use the weight-coupled interaction Morawetz inequality to deduce a contradiction.

Assume that $(\phi,\psi)$ is a solution of
\begin{align}
i\phi_t+\Delta \phi=\lambda|\psi|^2\phi+\mathcal{F}_1,\quad i\psi_t+\Delta \psi=\mu|\phi|^2\psi+\mathcal{F}_2\label{1126x1}
\end{align}
 when the spatial dimension $d=4$. Similar to (\ref{983}) in the case of the spatial dimension $d=3$, for some weight $a:\mathbb{R}^4\rightarrow \mathbb{R}$, we also define the following weight-coupled Morawetz interaction:\\

\begin{align}
M^{\otimes_2}_a(t)&=2A\int_{\mathbb{R}^4}\int_{\mathbb{R}^4}\widetilde{\nabla} a(x-y)\Im[\bar{\phi}(t,x)\bar{\phi}(t,y)\widetilde{\nabla}(\phi(t,x)\phi(t,y))]dxdy\nonumber\\
&\quad+2B\int_{\mathbb{R}^4}\int_{\mathbb{R}^4}\widetilde{\nabla} a(x-y)\Im[\bar{\psi}(t,x)\bar{\psi}(t,y)\widetilde{\nabla}(\psi(t,x)\psi(t,y))]dxdy\nonumber\\
&\quad+2C\int_{\mathbb{R}^4}\int_{\mathbb{R}^4}\widetilde{\nabla} a(x-y)\Im[\bar{\phi}(t,x)\bar{\psi}(t,y)\widetilde{\nabla}(\phi(t,x)\psi(t,y))]dxdy\nonumber\\
&\quad+2D\int_{\mathbb{R}^4}\int_{\mathbb{R}^4}\widetilde{\nabla} a(x-y)\Im[\bar{\psi}(t,x)\bar{\phi}(t,y)\widetilde{\nabla}(\psi(t,x)\phi(t,y))]dxdy,\label{1126x2}
\end{align}
where $\widetilde{\nabla}=(\nabla_x, \nabla_y)$, $x\in \mathbb{R}^4$ and $y\in \mathbb{R}^4$, $A=4\mu^2$, $B=4\lambda^2$ and $C=D=4\lambda\mu$. If we take $a(x,y)=|x-y|$ in the weight-coupled interaction Morawetz action (\ref{1126x2}), then
\begin{align}
\frac{d}{dt} M(t)^{\otimes}_2&\gtrsim \int_{\mathbb{R}^4}\int_{\mathbb{R}^4}\frac{[|\phi(t,x)|^2+|\psi(t,x)|^2][|\phi(t,y)|^2+|\psi(t,y)|^2]}{|x-y|^3}dxdy\nonumber\\
&\quad +\sum_{m,n=1}^2\int_{\mathbb{R}^4}\int_{\mathbb{R}^4}\Re(\mathcal{F}_m\nabla_x\bar{w}_m-w_m\nabla_x\overline{\mathcal{F}}_m)\cdot\frac{x-y}{|x-y|}dx|w_n(t,y)|^2dy\nonumber\\
&\quad +\sum_{m,n=1}^2\int_{\mathbb{R}^4}\int_{\mathbb{R}^4}\Re(\mathcal{F}_m\nabla_y\bar{w}_m-w_m\nabla_y\overline{\mathcal{F}}_m)\cdot\frac{y-x}{|x-y|}dy|w_n(t,x)|^2dx\nonumber\\
&\quad -\sum_{m,n=1}^2\int_{\mathbb{R}^4}\int_{\mathbb{R}^4}\frac{x-y}{|x-y|}\cdot\Im[\bar{w}_m(t,x)\nabla_xw_m(t,x)]dx
\Im(\overline{\mathcal{F}}_nw_n)(t,y)dy\nonumber\\
&\quad -\sum_{m,n=1}^2\int_{\mathbb{R}^4}\int_{\mathbb{R}^4}\frac{y-x}{|x-y|}\cdot\Im[\bar{w}_m(t,y)\nabla_y w_m(t,y)]dy
\Im(\overline{\mathcal{F}}_nw_n)(t,x)dx.\label{924x1}
\end{align}
Here $w_1$, $w_2$ can be taken as follows: (1) one is $\phi$, another is $\psi$; (2) both are $\phi$; (3) both are $\psi$.
Denote
\begin{align}
&\{\mathcal{F}_m, \phi,\psi\}_{px}:=\Re(\mathcal{F}_m\nabla_x\bar{w}_m-w_m\nabla_x\overline{\mathcal{F}}_m),\quad
\{\mathcal{F}_m, \phi,\psi\}_{py}:=\Re(\mathcal{F}_m\nabla_y\bar{w}_m-w_m\nabla_y\overline{\mathcal{F}}_m),\label{924x4}\\
&\{\mathcal{F}_n, \phi,\psi\}_{ix}:=\Im(\overline{\mathcal{F}}_nw_n)(t,x),\quad \{\mathcal{F}_n, \phi,\psi\}_{iy}:=\Im(\overline{\mathcal{F}}_nw_n)(t,y).\label{924x3}
\end{align}

Integrating (\ref{924x1}) with respect to time, we get

{\bf Proposition 6.15(Weight-coupled interaction Morawetz inequality).} {\it
\begin{align}
&\quad \int_I\int_{\mathbb{R}^4}\int_{\mathbb{R}^4}\frac{[|\phi(t,x)|^2+|\psi(t,x)|^2][|\phi(t,y)|^2+|\psi(t,y)|^2]}{|x-y|^3}dxdydt\nonumber\\
&+\sum_{m,n=1}^2\int_I\int_{\mathbb{R}^4}\int_{\mathbb{R}^4}\{\mathcal{F}_m, \phi,\psi\}_{px}\cdot\frac{x-y}{|x-y|}dx|w_n(t,y)|^2dydt\nonumber\\
&+\sum_{m,n=1}^2\int_I\int_{\mathbb{R}^4}\int_{\mathbb{R}^4}\{\mathcal{F}_m, \phi,\psi\}_{py}\cdot\frac{y-x}{|x-y|}dy|w_n(t,x)|^2dxdt\nonumber\\
&\lesssim [\|\phi\|^3_{L^{\infty}_tL^2_x}+\|\psi\|^3_{L^{\infty}_tL^2_x}][\|\phi\|_{L^{\infty}_t\dot{H}^1_x}+\|\psi\|_{L^{\infty}_t\dot{H}^1_x}]
+\left\{[\|\phi\|_{L^{\infty}_tL^2_x}+\|\psi\|_{L^{\infty}_tL^2_x}]\right.\nonumber\\
&\qquad \qquad \qquad \left.[\|\phi\|_{L^{\infty}_t\dot{H}^1_x}+\|\psi\|_{L^{\infty}_t\dot{H}^1_x}][\|\{\mathcal{F}_1, \phi,\psi\}_{ix}\|_{L^1_{t,x}}+\|\{\mathcal{F}_2, \phi,\psi\}_{ix}\|_{L^1_{t,x}}]\right\}.\label{924x4}
\end{align}
All spacetime norms above are on $I\times \mathbb{R}^4$.
}

Letting $\phi=u_{\geq N}$, $\psi=v_{\geq N}$, $\mathcal{F}_1=P_{\geq N}(|v|^2u)$ and $\mathcal{F}_2=P_{\geq N}(|u|^2v)$ for $N$ small enough such that
the Littlewood-Paley projection captures most of the solution, we will use Proposition 6.15 to prove

{\bf Proposition 6.16(Frequency-localized interaction Morawetz estimate).} {\it
Let $(u,v):[0,T_{\max})\times \mathbb{R}^4\rightarrow \mathbb{C}\times \mathbb{C}$ be an almost periodic solution of (\ref{826x1}) satisfying $\widetilde{N}(t)\equiv N_k\geq 1$ on every characteristic interval $J_k\subset [0,T_{\max})$. Then for any $\eta>0$ there exists $N_0=N_0(\eta)$ such
that
\begin{align}
&\quad \int_I\int_{\mathbb{R}^4}\int_{\mathbb{R}^4}\frac{[|u_{\geq N}(t,x)|^2+|v_{\geq N}(t,x)|^2][|u_{\geq N}(t,y)|^2+|v_{\geq N}(t,y)|^2]}{|x-y|^3}dxdydt\nonumber\\
&\lesssim_{(u,v)} \eta[N^{-3}+\int_I\widetilde{N}(t)^{-1}dt]\label{924x5}
\end{align}
for $N\leq N_0$ and any compact time interval $I\subset [0,T_{\max})$, which is a union of contiguous subintervals $J_k$. Here the implicit constant does not depend on the interval $I$.
}

{\bf Proof:} Let $I$ be a fixed compact interval $I\subset [0,T_{\max})$, which is a union of contiguous subintervals $J_k$, denote $K:=\int_IN(t)^{-1}dt$.
All spacetime norms will be on $I\times \mathbb{R}^4$ in the proof of this proposition below.

Fix $\eta>0$ and chose $N_0=N_0(\eta)$ small enough such that the claims in Corollary 6.13 hold for all $N\leq N_0$. Especially, fixing $N\leq N_0$ and writing $u_{lo}:=u_{\leq N}$, $v_{lo}:=v_{\leq N}$, $u_{hi}:=u_{>N}$ and $v_{hi}:=v_{>N}$, there holds
\begin{align}
&\|\nabla u_{lo}\|_{L^q_tL^r_x}+\|\nabla v_{lo}\|_{L^q_tL^r_x}\lesssim_{(u,v)} \eta(1+N^3K)^{\frac{1}{q}}, \quad {\rm for\ all}\quad \frac{1}{q}+\frac{2}{r}=1,\quad 2\leq q\leq +\infty,\label{924x6}\\
&\| u_{hi}\|_{L^q_tL^r_x}+\|v_{hi}\|_{L^q_tL^r_x}\lesssim_{(u,v)} N^{-1}(1+N^3K)^{\frac{1}{q}}, \quad {\rm for\ all}\quad \frac{1}{q}+\frac{2}{r}=1,\quad
3<q\leq +\infty,\label{924x7}\\
&\|u_{hi}\|_{L^{\infty}_tL^2_x}+\|v_{hi}\|_{L^{\infty}_tL^2_x}\lesssim_{(u,v)}\eta^6N^{-1}.\label{924x8}
\end{align}

Applying Proposition 6.15 with $\phi=u_{hi}$, $\psi=v_{hi}$, $\mathcal{F}_1=P_{hi}f_1(u,v)$ and $\mathcal{F}_2=P_{hi}f_2(u,v)$, and using (\ref{924x8}), we get
\begin{align}
&\quad \int_I\int_{\mathbb{R}^4}\int_{\mathbb{R}^4}\frac{[|u_{hi}(t,x)|^2+|v_{hi}(t,x)|^2][|u_{hi}(t,y)|^2+|v_{hi}(t,y)|^2]}{|x-y|^3}dxdydt\nonumber\\
&\quad+2\sum_{m,n=1}^2\int_{\mathbb{R}^4}\int_{\mathbb{R}^4}\Re(\mathcal{F}_m\nabla_x\bar{w}_m-w_m\nabla_x\overline{\mathcal{F}}_m)
\cdot\frac{x-y}{|x-y|}dx|w_n(t,y)|^2dy\nonumber\\
&\quad+2\sum_{m,n=1}^2\int_{\mathbb{R}^4}\int_{\mathbb{R}^4}\Re(\mathcal{F}_m\nabla_y\bar{w}_m-w_m\nabla_y\overline{\mathcal{F}}_m)
\cdot\frac{y-x}{|x-y|}dy|w_n(t,x)|^2dx\nonumber\\
&\lesssim_{(u,v)} \eta^{18}N^{-3}+\eta^6N^{-1}[\|\{\mathcal{F}_1,u_{hi}, v_{hi}\}_{ix}\|_{L^1_{t,x}}+\|\{\mathcal{F}_2,u_{hi}, v_{hi}\}_{ix}\|_{L^1_{t,x}}].\label{924w1}
\end{align}
Here $w_1$, $w_2$ can be taken as follows: (1) one is $u_{hi}$, another is $v_{hi}$; (2) both are $u_{hi}$; (3) both are $v_{hi}$.

Since
\begin{align*}
&\quad \sum_{m,n=1}^2\Re(\mathcal{F}_m\nabla_x\bar{w}_m-w_m\nabla_x\overline{\mathcal{F}}_m)\nonumber\\
&=\nabla_x [\varnothing(u_{hi}^2v_{hi}^2)+\varnothing(u^2_{hi}v_{hi}v_{lo})+\varnothing(v^2_{hi}u_{hi}u_{lo})+\varnothing(u_{hi}v_{hi}u_{lo}v_{lo})]\nonumber\\
&\quad+\nabla_x [\varnothing(u^2_{hi}v^2_{lo})+\varnothing(u_{hi}u_{lo}v^2_{lo})+\varnothing(u^2_{hi}v^2_{lo})+\varnothing(v_{hi}v_{lo}u^2_{lo})]\nonumber\\
&\quad+ [\varnothing(v_{hi}^2u_{hi})+\varnothing(u_{hi}v_{hi}v_{lo})+\varnothing(u_{hi}v^2_{lo})]\nabla_x u_{lo}\nonumber\\
&\quad+
[\varnothing(u_{hi}^2v_{hi})+\varnothing(v_{hi}u_{hi}u_{lo})+\varnothing(v_{hi}u^2_{lo})]\nabla_x v_{lo}\nonumber\\
&\quad+\nabla_x [\varnothing(u_{hi}P_{lo}f_1(u,v))+\varnothing(v_{hi}P_{lo}f_2(u,v))]+\varnothing(u_{hi}\nabla_x P_{lo}f_1(u,v))+\varnothing(v_{hi}\nabla_x P_{lo}f_2(u,v))
\end{align*}
and
\begin{align*}
\{\mathcal{F}_1,u_{hi}, v_{hi}\}_{ix}&=\varnothing(u_{hi}v^2_{hi}u_{lo})+\varnothing(v_{hi}u^2_{hi}v_{lo})+\varnothing(u^2_{hi}v^2_{lo})
+\varnothing(u_{hi}v_{hi}u_{lo}v_{lo})\\
\{\mathcal{F}_2,u_{hi}, v_{hi}\}_{ix}&=\varnothing(v_{hi}u^2_{hi}v_{lo})+\varnothing(u_{hi}v^2_{hi}u_{lo})+\varnothing(v^2_{hi}u^2_{lo})
+\varnothing(u_{hi}v_{hi}u_{lo}v_{lo}),
\end{align*}
using (\ref{924x6})--(\ref{924x8}), we obtain from (\ref{924w1})
\begin{align}
&\quad\int_I\int_{\mathbb{R}^4}\int_{\mathbb{R}^4}\frac{[|u_{hi}(t,x)|^2+|v_{hi}(t,x)|^2][|u_{hi}(t,y)|^2+|v_{hi}(t,y)|^2]}{|x-y|^3}dxdydt\nonumber\\
&\quad+\int_I\int_{\mathbb{R}^4}\int_{\mathbb{R}^4}\frac{[|u_{hi}(t,x)|^2+|v_{hi}(t,x)|^2]|u_{hi}(t,y)|^2|v_{hi}(t,y)|^2}{|x-y|}dxdydt\nonumber\\
&\lesssim_{(u,v)}\eta^{18}N^{-3}+\eta^6N^{-1}I(u_{hi},v_{hi},u_{lo},v_{lo})+\eta^{12}N^{-2}II(u_{hi},v_{hi},u_{lo},v_{lo})\nonumber\\
&\quad+\int_I\int_{\mathbb{R}^4}\int_{\mathbb{R}^4}\frac{[|u_{hi}(t,x)|^2+|v_{hi}(t,x)|^2]J(u_{hi},v_{hi},u_{lo},v_{lo})}{|x-y|}dxdydt\nonumber\\
&\quad+\int_I\int_{\mathbb{R}^4}\int_{\mathbb{R}^4}\frac{[|u_{hi}(t,x)|^2+|v_{hi}(t,x)|^2][P_{lo}f_1(u,v)||u_{hi}|+|P_{lo}f_2(u,v)||v_{hi}|]}{|x-y|}dxdydt,\label{108x1}
\end{align}
where
\begin{align}
I(u_{hi},v_{hi},u_{lo},v_{lo})&=\|u_{hi}v^2_{hi}u_{lo}\|_{L^1_{t,x}}+\|v_{hi}u^2_{hi}v_{lo}\|_{L^1_{t,x}}+\|u_{hi}v_{hi}u_{lo}v_{lo}\|_{L^1_{t,x}}
+\|u^2_{hi}v^2_{lo}\|_{L^1_{t,x}}\nonumber\\
&\quad+\|v^2_{hi}u^2_{lo}\|_{L^1_{t,x}}+\|u_{hi}P_{hi}f_1(u_{lo},v_{lo})\|_{L^1_{t,x}}+\|v_{hi}P_{hi}f_2(u_{lo},v_{lo})\|_{L^1_{t,x}}\nonumber\\
&\quad+\|u_{hi}P_{lo}f_1(u_{hi},v_{hi})\|_{L^1_{t,x}}+\|v_{hi}P_{lo}f_2(u_{hi},v_{hi})\|_{L^1_{t,x}},\label{108x2}\\
II(u_{hi},v_{hi},u_{lo},v_{lo})&=\|u_{hi}v^2_{hi}\nabla u_{lo}\|_{L^1_{t,x}}+\|u_{hi}v_{hi}v_{lo}\nabla u_{lo}\|_{L^1_{t,x}}+\|u_{hi}v^2_{lo}\nabla u_{lo}\|_{L^1_{t,x}}\nonumber\\
&\quad+\|v_{hi}u^2_{hi}\nabla v_{lo}\|_{L^1_{t,x}}+\|v_{hi}u_{hi}u_{lo}\nabla v_{lo}\|_{L^1_{t,x}}+\|v_{hi}u^2_{lo}\nabla v_{lo}\|_{L^1_{t,x}}\nonumber\\
&\quad+\|u_{hi}\nabla P_{lo}f_1(u,v)\|_{L^1_{t,x}}+\|v_{hi}\nabla P_{lo}f_2(u,v)\|_{L^1_{t,x}},\label{108x3}\\
J(u_{hi},v_{hi},u_{lo},v_{lo})&=\varnothing(u_{hi}u_{lo}v^2_{lo})+\varnothing(v_{hi}v_{lo}u^2_{lo})
+\varnothing(u^2_{hi}v^2_{lo})+\varnothing(v^2_{hi}u^2_{lo})\nonumber\\
&\quad+\varnothing(u^2_{hi}v_{hi}v_{lo})+\varnothing(v^2_{hi}u_{hi}u_{lo})+\varnothing(u_{hi}v_{hi}u_{lo}v_{lo}).\label{108x4}
\end{align}

We estimate the terms through (\ref{108x2})--(\ref{108x4}) below.

Using (\ref{924x6})--(\ref{924x8}) and Sobolev embedding, we obtain
\begin{align*}
\|u_{hi}v^2_{hi}u_{lo}\|_{L^1_{t,x}}&\lesssim \|u_{hi}\|_{L^{\infty}_tL^4_x}\|v_{hi}\|^2_{L^{\frac{7}{2}}_tL^{\frac{14}{5}}_x}\|u_{lo}\|_{L^{\frac{7}{3}}_tL^{28}_x}\lesssim_{(u,v)}\eta N^{-2}(1+N^3K),\\
\|u^2_{hi}v^2_{lo}\|_{L^1_{t,x}}&\lesssim \|u_{hi}\|^2_{L^4_tL^{\frac{8}{3}}_x}\|u_{lo}\|^2_{L^4_tL^8_x}\lesssim_{(u,v)}\eta^2N^{-2}(1+N^3K).
\end{align*}
Similarly,
\begin{align*}
&\|v_{hi}u^2_{hi}v_{lo}\|_{L^1_{t,x}}\lesssim_{(u,v)}\eta N^{-2}(1+N^3K),\quad
\|v^2_{hi}u^2_{lo}\|_{L^1_{t,x}}\lesssim_{(u,v)}\eta^2N^{-2}(1+N^3K), \\
&\|u_{hi}v_{hi}u_{lo}v_{lo}\|_{L^1_{t,x}}\lesssim_{(u,v)}\eta^2N^{-2}(1+N^3K)\quad {\rm because}\ |u_{hi}v_{hi}u_{lo}v_{lo}|\leq |u^2_{hi}v^2_{lo}|+|v^2_{hi}u^2_{lo}|.
\end{align*}
Using Bernstein's inequality, we have
\begin{align*}
&\quad \|u_{hi}P_{hi}f_1(u_{lo},v_{lo})\|_{L^1_{t,x}}\lesssim \|u_{hi}\|_{L^4_tL^{\frac{8}{3}}_x}N^{-1}\|\nabla f_1(u_{lo},v_{lo})\|_{L^{\frac{4}{3}}_tL^{\frac{8}{5}}_x}\\
&\lesssim_{(u,v)}N^{-2}(1+N^3K)^{\frac{1}{4}}[\|\nabla u_{lo}\|_{L^2_tL^4_x}+\|\nabla v_{lo}\|_{L^2_tL^4_x}][\| u_{lo}\|_{L^{\infty}_tL^{\frac{16}{3}}_x}+\| v_{lo}\|_{L^{\infty}_tL^{\frac{16}{3}}_x}]\\
&\lesssim_{(u,v)}\eta^3 N^{-2}(1+N^3K)
\end{align*}
and
\begin{align*}
&\quad \|u_{hi}P_{lo}f_1(u_{hi},v_{hi})\|_{L^1_{t,x}}\lesssim \|u_{hi}\|_{L^{\frac{10}{3}}_tL^{\frac{20}{7}}_x}N^{\frac{7}{5}}\|f_1(u_{hi},v_{hi})\|_{L^{\frac{10}{7}}_tL^1_x}\\
&\lesssim_{(u,v)}N^{\frac{2}{5}}(1+N^3K)^{\frac{3}{10}}[\|u_{hi}\|^{\frac{7}{3}}_{L^{\frac{10}{3}}_tL^{\frac{20}{7}}_x}
+\|v_{hi}\|^{\frac{7}{3}}_{L^{\frac{10}{3}}_tL^{\frac{20}{7}}_x}]
[\|u_{hi}\|^{\frac{2}{3}}_{L^{\infty}_tL^{\frac{40}{11}}_x}+\|v_{hi}\|^{\frac{2}{3}}_{L^{\infty}_tL^{\frac{40}{11}}_x}]\\
&\lesssim_{(u,v)}N^{\frac{2}{5}-\frac{7}{3}}(1+N^3K)[\||\nabla|^{\frac{9}{10}}u_{hi}\|^{\frac{2}{3}}_{L^{\infty}_tL^2_x}
+\||\nabla|^{\frac{9}{10}}v_{hi}\|^{\frac{2}{3}}_{L^{\infty}_tL^2_x}]\\
&\lesssim_{(u,v)}N^{-2}(1+N^3K).
\end{align*}
Similarly,
\begin{align*}
&\|v_{hi}P_{hi}f_2(u_{lo},v_{lo})\|_{L^1_{t,x}}\lesssim_{(u,v)}\eta^3 N^{-2}(1+N^3K),\\
&\|v_{hi}P_{lo}f_2(u_{hi},v_{hi})\|_{L^1_{t,x}}\lesssim_{(u,v)}N^{-2}(1+N^3K).
\end{align*}
(\ref{108x2}) is completed.

Consider the terms in (\ref{108x3}). Using (\ref{924x6})--(\ref{924x8}), Sobolev embedding and Bernstein inequality, we get
\begin{align*}
&\|u_{hi}v^2_{hi}\nabla u_{lo}\|_{L^1_{t,x}}\lesssim \|\nabla u_{lo}\|_{L^{\frac{7}{3}}_tL^{28}_x}\|v_{hi}\|^2_{L^{\frac{7}{2}}_tL^{\frac{4}{5}}_x}\|u_{hi}\|_{L^{\infty}_tL^4_x}\lesssim_{(u,v)}\eta N^{-1}(1+N^3K),\\
&\|u_{hi}v_{hi}v_{lo}\nabla u_{lo}\|_{L^1_{t,x}}\lesssim \|\nabla u_{lo}\|_{L^2_tL^4_x}[\|u_{hi}\|^2_{L^4_tL^{\frac{8}{3}}_x}+\|v_{hi}\|^2_{L^4_tL^{\frac{8}{3}}_x}]\|v_{lo}\|_{L^{\infty}_{t,x}}\lesssim_{(u,v)}\eta^2 N^{-1}(1+N^3K),\\
&\|u_{hi}v^2_{lo}\nabla u_{lo}\|_{L^1_{t,x}}\lesssim \|\nabla u_{lo}\|_{L^2_tL^4_x}\|u_{hi}\|_{L^{\infty}_tL^2_x}\|v_{lo}\|^2_{L^4_tL^8_x}\lesssim_{(u,v)}\eta^9 N^{-1}(1+N^3K).
\end{align*}
Similarly,
\begin{align*}
&\|v_{hi}u^2_{hi}\nabla v_{lo}\|_{L^1_{t,x}}\lesssim_{(u,v)}\eta N^{-1}(1+N^3K),\quad \|v_{hi}u_{hi}u_{lo}\nabla v_{lo}\|_{L^1_{t,x}}\lesssim_{(u,v)}\eta^2 N^{-1}(1+N^3K),\\
&\|v_{hi}u^2_{lo}\nabla v_{lo}\|_{L^1_{t,x}}\lesssim_{(u,v)}\eta^9 N^{-1}(1+N^3K).
\end{align*}
Write $$f_1(u,v)=f_1(u_{lo},v_{lo})+\varnothing(u_{hi}v^2_{hi})+\varnothing(u_{hi}v_{hi}v_{lo})+\varnothing(v_{hi}u_{lo}v_{lo})
+\varnothing(v^2_{hi}u_{lo})+\varnothing(u_{hi}v^2_{lo}).$$
We estimate
\begin{align*}
&\quad \|u_{hi}\nabla P_{lo}f_1(u_{lo},v_{lo})\|_{L^1_{t,x}}\lesssim \|u_{hi}\|_{L^{\infty}_tL^2_x}[\|\nabla u_{lo}\|_{L^2_tL^4_x}+\|\nabla v_{lo}\|_{L^2_tL^4_x}]
[\| u_{lo}\|^2_{L^4_tL^8_x}+\| v_{lo}\|^2_{L^4_tL^8_x}]\\
&\lesssim_{(u,v)} \eta^9N^{-1}(1+N^3K),\\
&\quad \|u_{hi}\nabla P_{lo}\varnothing(u_{hi}v^2_{hi})\|_{L^1_{t,x}}\lesssim \|u_{hi}\|_{L^{\frac{10}{3}}_tL^{\frac{20}{7}}_x}N^{\frac{12}{5}}\|u_{hi}v^2_{hi}\|_{L^{\frac{10}{7}}_tL^1_x}\\
&\lesssim N^{\frac{12}{5}}\|u_{hi}\|_{L^{\frac{10}{3}}_tL^{\frac{20}{7}}_x}[\|u^3_{hi}\|_{L^{\frac{10}{7}}_tL^1_x}+\|v^3_{hi}\|_{L^{\frac{10}{7}}_tL^1_x}]\\
&\lesssim N^{\frac{12}{5}}\|u_{hi}\|_{L^{\frac{10}{3}}_tL^{\frac{20}{7}}_x}
[\|u_{hi}\|^{\frac{7}{3}}_{L^{\frac{10}{3}}_tL^{\frac{20}{7}}_x}\|u_{hi}\|^{\frac{2}{3}}_{L^{\infty}_tL^{\frac{40}{11}}_x}
+\|v_{hi}\|^{\frac{7}{3}}_{L^{\frac{10}{3}}_tL^{\frac{20}{7}}_x}\|v_{hi}\|^{\frac{2}{3}}_{L^{\infty}_tL^{\frac{40}{11}}_x}]\\
&\lesssim_{(u,v)} N^{-1}(1+N^3K),\\
&\quad \|u_{hi}\nabla P_{lo}\varnothing(u_{hi}v_{hi}v_{lo})\|_{L^1_{t,x}}\lesssim N\|u_{hi}\|_{L^{\infty}_tL^4_x}[\|u_{hi}\|^2_{L^{\frac{7}{2}}_tL^{\frac{14}{5}}_x}+\|v_{hi}\|^2_{L^{\frac{7}{2}}_tL^{\frac{14}{5}}_x}]
\|v_{lo}\|_{L^{\frac{7}{3}}_tL^{28}_x}\\
&\lesssim_{(u,v)}\eta N^{-1}(1+N^3K),\\
&\quad \|u_{hi}\nabla P_{lo}\varnothing(v_{hi}u_{lo}v_{lo})\|_{L^1_{t,x}}\lesssim N\|u_{hi}\|_{L^4_tL^{\frac{8}{3}}_x}\|v_{hi}\|_{L^4_tL^{\frac{8}{3}}_x}
[\|u_{lo}\|^2_{L^4_tL^8_x}+\|v_{lo}\|^2_{L^4_tL^8_x}]\\
&\lesssim N[\|u_{hi}\|^2_{L^4_tL^{\frac{8}{3}}_x}+\|v_{hi}\|^2_{L^4_tL^{\frac{8}{3}}_x}]
[\|u_{lo}\|^2_{L^4_tL^8_x}+\|v_{lo}\|^2_{L^4_tL^8_x}]\lesssim_{(u,v)}\eta^2 N^{-1}(1+N^3K),\\
&\quad \|u_{hi}\nabla P_{lo}\varnothing(v^2_{hi}u_{lo})\|_{L^1_{t,x}}\lesssim N\|u_{hi}\|_{L^{\infty}_tL^4_x}\|v_{hi}\|^2_{L^{\frac{7}{2}}_tL^{\frac{14}{5}}_x}
\|u_{lo}\|_{L^{\frac{7}{3}}_tL^{28}_x}\lesssim_{(u,v)}\eta N^{-1}(1+N^3K),\\
&\quad \|u_{hi}\nabla P_{lo}\varnothing(u_{hi}u^2_{lo})\|_{L^1_{t,x}}\lesssim N\|u_{hi}\|^2_{L^4_tL^{\frac{8}{3}}_x}
\|u_{lo}\|^2_{L^4_tL^8_x}\lesssim_{(u,v)}\eta^2 N^{-1}(1+N^3K).
\end{align*}
Similarly, we can obtain the estimate for $\|v_{hi}\nabla P_{lo}f_2(u,v)\|_{L^1_{t,x}}$. (\ref{108x3}) is completed.

Now we consider the contribution of the term containing $J(u_{hi},v_{hi},u_{lo},v_{lo})$. We estimate
\begin{align*}
&\quad\int_I\int_{\mathbb{R}^4}\int_{\mathbb{R}^4}\frac{[|u_{hi}(t,y)|^2+|v_{hi}(t,y)|^2]|u_{hi}(t,x)||u_{lo}(t,x)||v_{lo}(t,x)|^2}{|x-y|}dxdydt\\
&\lesssim
[\|u_{hi}\|^2_{L^{12}_tL^{\frac{24}{11}}_x}+\|v_{hi}\|^2_{L^{12}_tL^{\frac{24}{11}}_x}]\|\frac{1}{|x|}*(|u_{hi}||u_{lo}||v_{lo}|^2)\|_{L^{\frac{6}{5}}_tL^{12}_x}\\
&\lesssim_{(u,v)}N^{-2}(1+N^3K)^{\frac{1}{6}}[\|u_{hi}u^3_{lo}\|_{L^{\frac{6}{5}}_{t,x}}+\|u_{hi}v^3_{lo}\|_{L^{\frac{6}{5}}_{t,x}}]\\
&\lesssim_{(u,v)}N^{-2}(1+N^3K)^{\frac{1}{6}}\|u_{hi}\|_{L^{\infty}_tL^2_x}[\|u_{lo}\|^3_{L^{\frac{18}{5}}_tL^9_x}+\|v_{lo}\|^3_{L^{\frac{18}{5}}_tL^9_x}]
\lesssim_{(u,v)}\eta^9 N^{-3}(1+N^3K).
\end{align*}
Similarly,
\begin{align*}
&\quad\int_I\int_{\mathbb{R}^4}\int_{\mathbb{R}^4}\frac{[|u_{hi}(t,y)|^2+|v_{hi}(t,y)|^2]|v_{hi}(t,x)||v_{lo}(t,x)||u_{lo}(t,x)|^2}{|x-y|}dxdydt\\
&\lesssim_{(u,v)}\eta^9 N^{-3}(1+N^3K).
\end{align*}
Note that $|u^2_{hi}v_{hi}v_{lo}|\lesssim 10^{-6}|u_{hi}|^2|v_{hi}|^2+|u_{hi}|^2|v_{lo}|^2$ and
\begin{align*}
&\quad\int_I\int_{\mathbb{R}^4}\int_{\mathbb{R}^4}\frac{[|u_{hi}(t,y)|^2+|v_{hi}(t,y)|^2]|u_{hi}(t,x)|^2|v_{hi}(t,x)||v_{lo}(t,x)|}{|x-y|}dxdydt\\
&\lesssim 10^{-6}\int_I\int_{\mathbb{R}^4}\int_{\mathbb{R}^4}\frac{[|u_{hi}(t,y)|^2+|v_{hi}(t,y)|^2]|u_{hi}(t,x)|^2|v_{hi}(t,x)|^2}{|x-y|}dxdydt\\
&\quad+\int_I\int_{\mathbb{R}^4}\int_{\mathbb{R}^4}\frac{[|u_{hi}(t,y)|^2+|v_{hi}(t,y)|^2]|u_{hi}(t,x)|^2|v_{lo}(t,x)|^2}{|x-y|}dxdydt.
\end{align*}
Using H\"{o}lder's inequality and Hardy-Littlewood-Sobolev inequality, we have
\begin{align*}
&\quad\int_I\int_{\mathbb{R}^4}\int_{\mathbb{R}^4}\frac{[|u_{hi}(t,y)|^2+|v_{hi}(t,y)|^2]|u_{hi}(t,x)|^2|v_{lo}(t,x)|^2}{|x-y|}dxdydt\\
&\lesssim [\|u_{hi}\|^2_{L^{12}_tL^{\frac{24}{11}}_x}+\|v_{hi}\|^2_{L^{12}_tL^{\frac{24}{11}}_x}]\|\frac{1}{|x|}*(|u_{hi}||u_{hi}v^2_{lo}|)\|_{L^{\frac{6}{5}}_tL^{12}_x}\\
&\lesssim_{(u,v)}N^{-2}(1+N^3K)^{\frac{1}{6}}\||u_{hi}|(|u_{hi}v^2_{lo}|)\|_{L^{\frac{6}{5}}_{t,x}}\\
&\lesssim_{(u,v)}N^{-2}(1+N^3K)^{\frac{1}{6}}\||u_{hi}\|_{L^4_tL^{\frac{8}{3}}_x}\|(|u_{hi}v^2_{lo}|)\|_{L^{\frac{12}{7}}_tL^{\frac{24}{11}}_x}\\
&\lesssim_{(u,v)}N^{-3}(1+N^3K)^{\frac{5}{12}}\|(|u_{hi}v^2_{lo}|)\|_{L^{\frac{12}{7}}_tL^{\frac{24}{11}}_x}\\
&\lesssim_{(u,v)}N^{-3}(1+N^3K)^{\frac{5}{12}}N\|(|u_{hi}v^2_{lo}|)\|_{L^{\frac{12}{7}}_tL^{\frac{24}{17}}_x}\\
&\lesssim_{(u,v)}N^{-3}(1+N^3K)^{\frac{5}{12}}N\|v_{lo}\|^2_{L^4_tL^8_x}\|u_{hi}\|_{L^{12}_tL^{\frac{24}{11}}_x}\\
&\lesssim_{(u,v)}N^{-3}(1+N^3K)^{\frac{5}{12}}\eta^2(1+N^3K)^{\frac{7}{12}}\lesssim_{(u,v)}\eta(N^{-3}+K).
\end{align*}
Similarly,
\begin{align*}
&\quad\int_I\int_{\mathbb{R}^4}\int_{\mathbb{R}^4}\frac{[|u_{hi}(t,y)|^2+|v_{hi}(t,y)|^2]|v_{hi}(t,x)|^2|u_{hi}(t,x)||u_{lo}(t,x)|}{|x-y|}dxdydt\\
&\lesssim 10^{-6}\int_I\int_{\mathbb{R}^4}\int_{\mathbb{R}^4}\frac{[|u_{hi}(t,y)|^2+|v_{hi}(t,y)|^2]|u_{hi}(t,x)|^2|v_{hi}(t,x)|^2}{|x-y|}dxdydt\\
&\quad+\int_I\int_{\mathbb{R}^4}\int_{\mathbb{R}^4}\frac{[|u_{hi}(t,y)|^2+|v_{hi}(t,y)|^2]|v_{hi}(t,x)|^2|u_{lo}(t,x)|^2}{|x-y|}dxdydt.
\end{align*}
and
\begin{align*}
&\quad\int_I\int_{\mathbb{R}^4}\int_{\mathbb{R}^4}\frac{[|u_{hi}(t,y)|^2+|v_{hi}(t,y)|^2]|v_{hi}(t,x)|^2|u_{lo}(t,x)|^2}{|x-y|}dxdydt\lesssim_{(u,v)}\eta(N^{-3}+K),\\
&\quad \int_I\int_{\mathbb{R}^4}\int_{\mathbb{R}^4}\frac{[|u_{hi}(t,y)|^2+|v_{hi}(t,y)|^2]|u_{hi}(t,x)||v_{hi}(t,x)||u_{lo}(t,x)||v_{lo}(t,x)|}{|x-y|}dxdydt\\
&\lesssim_{(u,v)}\eta(N^{-3}+K)\quad{\rm because}\quad |u_{hi}||v_{hi}||u_{lo}||v_{lo}|\lesssim |u_{hi}|^2|v_{lo}|^2+|v_{hi}|^2|u_{lo}|^2.
\end{align*}
The contribution of the term containing $J(u_{hi},v_{hi},u_{lo},v_{lo})$ is completed.

Last, we estimate the terms containing $P_{lo}f_1(u,v)$ and $P_{lo}f_2(u,v)$ respectively. Write
$$
f_1(u,v)=f_1(u_{hi},v_{hi})+\varnothing(u_{lo}v^2_{lo})+\varnothing(u_{hi}v^2_{lo})+\varnothing(v^2_{hi}u_{lo})+\varnothing(u_{hi}v_{hi}v_{lo})
+\varnothing(v_{hi}u_{lo}v_{lo}).
$$
Using H\"{o}lder's inequality and Hardy-Littlewood-Sobolev inequality, we obtain
\begin{align*}
&\quad\int_I\int_{\mathbb{R}^4}\int_{\mathbb{R}^4}\frac{[|u_{hi}(t,y)|^2+|v_{hi}(t,y)|^2]|u_{hi}(t,x)||P_{lo}f_1(u,v)(t,x)|}{|x-y|}dxdydt\\
&\lesssim [\|u_{hi}\|^2_{L^{12}_tL^{\frac{24}{11}}_x}+\|v_{hi}\|^2_{L^{12}_tL^{\frac{24}{11}}_x}]\|\frac{1}{|x|}*(|u_{hi}|||P_{lo}f_1(u,v)(t,x)||)\|_{L^{\frac{6}{5}}_tL^{12}_x}\\
&\lesssim_{(u,v)}N^{-2}(1+N^3K)^{\frac{1}{6}}\||u_{hi}|(|u_{hi}v^2_{lo}|)\|_{L^{\frac{6}{5}}_{t,x}}\\
&\lesssim_{(u,v)}N^{-2}(1+N^3K)^{\frac{1}{6}}\||u_{hi}\|_{L^4_tL^{\frac{8}{3}}_x}\|(||P_{lo}f_1(u,v)(t,x)||)\|_{L^{\frac{12}{7}}_tL^{\frac{24}{11}}_x}\\
&\lesssim_{(u,v)}N^{-3}(1+N^3K)^{\frac{5}{12}}\|(||P_{lo}f_1(u,v)(t,x)||)\|_{L^{\frac{12}{7}}_tL^{\frac{24}{11}}_x}.
\end{align*}
While
\begin{align*}
&\quad \|P_{lo}\varnothing(u_{lo}v^2_{lo})\|_{L^{\frac{12}{7}}_tL^{\frac{24}{11}}_x}\lesssim \|u_{lo}\|_{L^2_tL^{\frac{24}{5}}_x}\|v_{lo}\|^2_{L^4_tL^8_x}\lesssim_{(u,v)} \eta^3(1+N^3K)^{\frac{7}{12}},\\
&\quad \|P_{lo}\varnothing(u_{hi}v^2_{lo})\|_{L^{\frac{12}{7}}_tL^{\frac{24}{11}}_x}\lesssim N\|v^2_{lo}u_{hi}\|_{L^{\frac{12}{7}}_tL^{\frac{24}{17}}_x}\lesssim \|u_{hi}\|_{L^{12}_tL^{\frac{24}{11}}_x}\|v_{lo}\|^2_{L^4_tL^8_x}\lesssim_{(u,v)} \eta^2(1+N^3K)^{\frac{7}{12}},\\
&\quad \|P_{lo}\varnothing(u_{lo}v^2_{hi})\|_{L^{\frac{12}{7}}_tL^{\frac{24}{11}}_x}\lesssim N\|u_{lo}v^2_{hi}\|_{L^{\frac{12}{7}}_tL^{\frac{24}{17}}_x}
\lesssim_{u,v} N\|u_{lo}\|_{L^3_tL^{12}_x}\|v_{hi}\|_{L^4_tL^{\frac{8}{3}}_x}\|v_{hi}\|_{L^{\infty}_tL^4_x}\\&\lesssim_{(u,v)}\eta(1+N^3K)^{\frac{7}{12}}
\end{align*}
and similarly
\begin{align*}
\|P_{lo}\varnothing(u_{hi}v_{hi}v_{lo})\|_{L^{\frac{12}{7}}_tL^{\frac{24}{11}}_x}\lesssim_{(u,v)}\eta(1+N^3K)^{\frac{7}{12}},\quad
\|P_{lo}\varnothing(v_{hi}u_{lo}v_{lo})\|_{L^{\frac{12}{7}}_tL^{\frac{24}{11}}_x}\lesssim_{(u,v)}\eta(1+N^3K)^{\frac{7}{12}}
\end{align*}
because $|u_{hi}v_{hi}v_{lo}|\leq [|u^2_{hi}v_{lo}|+|v^2_{hi}v_{lo}|$ and $|v_{hi}u_{lo}v_{lo}|\leq |v_{hi}u^2_{lo}|+|v_{hi}v^2_{lo}|$.

Consequently,
\begin{align*}
\int_I\int_{\mathbb{R}^4}\int_{\mathbb{R}^4}\frac{[|u_{hi}(t,y)|^2+|v_{hi}(t,y)|^2]|u_{hi}(t,x)||P_{lo}f_1(u,v)(t,x)|}{|x-y|}dxdydt
\lesssim \eta(N^{-3}+K).
\end{align*}
Similarly,
\begin{align*}
\int_I\int_{\mathbb{R}^4}\int_{\mathbb{R}^4}\frac{[|u_{hi}(t,y)|^2+|v_{hi}(t,y)|^2]|u_{hi}(t,x)||P_{lo}f_2(u,v)(t,x)|}{|x-y|}dxdydt
\lesssim \eta(N^{-3}+K).
\end{align*}

Putting all things above together, we obtain (\ref{924x5}).\hfill $\Box$

Now we are ready to preclude the second scenario of Proposition 6.6 and complete the proof of Theorem 7.

{\bf Proposition 6.17(No quasi-solitions).} {\it There are no almost periodic solution $(u,v):[0,T_{\max})\times \mathbb{R}^4\rightarrow \mathbb{C}$ satisfying $\widetilde{N}(t)\equiv N_k\geq 1$ on each characteristic interval $J_k\subset [0,T_{\max})$ and
\begin{align}
\|u\|_{L^6_{t,x}([0,T_{\max})\times \mathbb{R}^4)}+\|v\|_{L^6_{t,x}([0,T_{\max})\times \mathbb{R}^4)}=+\infty,\quad \int_0^{T_{\max}}\widetilde{N}(t)^{-1}dt=+\infty.\label{1011s1}
\end{align}
}

{\bf Proof:} Contradictorily, assume that there exists such a solution $(u,v)$.

Let $\eta>0$ be a small parameter which will be chosen later. By the results of Proposition 6.16, we can find some $N_0=N_0(\eta)$ such that for all $N\leq N_0$ and any compact time interval $I\subset [0,T_{\max})$, which is a union of contiguous subintervals $J_k$,
\begin{align}
&\quad\int_I\int_{\mathbb{R}^4}\int_{\mathbb{R}^4}\frac{[|u_{\geq N}(t,x)|^2+|v_{\geq N}(t,x)|^2][|u_{\geq N}(t,y)|^2+|v_{\geq N}(t,y)|^2]}{|x-y|^3}dxdydt\nonumber\\
&\lesssim_{(u,v)} \eta[N^{-3}+\int_I\widetilde{N}(t)^{-1}dt].\label{1011s2}
\end{align}
Noticing that $\inf_{t\in [0,T_{\max})}\widetilde{N}(t)\geq 1$, we can chose $N_0$ smaller enough and ensure that for all $N\leq N_0$
\begin{align}
\|u_{\leq N}\|_{L^{\infty}_tL^4_x}+\|v_{\leq N}\|_{L^{\infty}_tL^4_x}
+\|\nabla u_{\leq N}\|_{L^{\infty}_tL^2_x}+\|\nabla v_{\leq N}\|_{L^{\infty}_tL^2_x}\leq \eta.\label{1011s3}
\end{align}
Here all the norms are on $([0,T_{\max})\times \mathbb{R}^4)$.

We will show that there exists $C(u,v)>0$ such that
\begin{align}
\widetilde{N}(t)^2\int_{|x-x(t)|\leq \frac{C(u)}{\widetilde{N}(t)}}[|u(t,x)|^2+|v(t,x)|^2]dx\gtrsim_{(u,v)}\frac{1}{C(u,v)}\label{1011s4}
\end{align}
 uniformly for $t\in[0,T_{\max})$. Since $(u(t),v(t)$ is not identically zero, (\ref{1011s4}) is true for each $t$. Furthermore, since $(u(t),v(t))$ is almost periodic, we find that the left hand side of (\ref{1011s4}) is scale invariant and the map $(u(t),v(t))\rightarrow LHS(\ref{1011s4})$ is continuous on $L^4_x$ and $\dot{H}^1_x$, so (\ref{1011s4}) is true uniformly for $t\in [0,T_{\max})$.

Using H\"{o}lder's inequality and (\ref{1011s3}), we get for all $t\in [0,T_{\max})$ and all $N\leq N_0$
\begin{align*}
&\quad\int_{|x-x(t)|\leq \frac{C(u,v)}{\widetilde{N}(t)}}[|u_{\leq N}(t,x)|^2+|v_{\leq N}(t,x)|^2]dx\\
&\lesssim\frac{C^2(u,v)}{\widetilde{N}^2(t)}[\|u_{\leq N}\|^2_{L^{\infty}_tL^4_x([0,T_{\max})\times \mathbb{R}^4)}+\|v_{\leq N}\|^2_{L^{\infty}_tL^4_x([0,T_{\max})\times \mathbb{R}^4)}]\\
&\lesssim_{(u,v)} \eta^2 C^2(u,v)\widetilde{N}(t)^{-2}.
\end{align*}
Putting this and (\ref{1011s4}) together and taking $\eta$ sufficiently small depending on $(u,v)$, we have for all $N\leq N_0$
\begin{align*}
\inf_{t\in [0,T_{\max})}\widetilde{N}^2(t)\int_{|x-x(t)|\leq \frac{C(u,v)}{\widetilde{N}(t)}}[|u_{\geq N}(t,x)|^2+|v_{\geq N}(t,x)|^2]dx\gtrsim_{(u,v)} 1.
\end{align*}
Therefore, we obtain
\begin{align*}
&\quad \int_I\int_{\mathbb{R}^4}\int_{\mathbb{R}^4}\frac{[|u_{\geq N}(t,x)|^2+|v_{\geq N}(t,x)|^2][|u_{\geq N}(t,y)|^2+|v_{\geq N}(t,y)|^2]}{|x-y|^3}dxdydt\\
&\gtrsim \int_I\int\int_{|x-y|\leq \frac{2C(u,v)}{\widetilde{N}(t)}}\left(\frac{\widetilde{N}(t)}{2C(u,v)}\right)^3[|u_{\geq N}(t,x)|^2+|v_{\geq N}(t,x)|^2][|u_{\geq N}(t,y)|^2+|v_{\geq N}(t,y)|^2]dxdydt\\
&\gtrsim \int_I\left(\frac{\widetilde{N}(t)}{2C(u,v)}\right)^3\int_{|x-x(t)|\leq \frac{2C(u,v)}{\widetilde{N}(t)}}[|u_{\geq N}(t,x)|^2+|v_{\geq N}(t,x)|^2]dx\\
&\qquad\qquad \qquad \quad \times \int_{|y-x(t)|\leq \frac{2C(u,v)}{\widetilde{N}(t)}} [|u_{\geq N}(t,y)|^2+|v_{\geq N}(t,y)|^2]dydt\\
&\gtrsim_{(u,v)}\int_I\widetilde{N}(t)^{-1}dt.
\end{align*}
Noticing that (\ref{1011s2}) and choosing $\eta$ small depending on $(u,v)$, we get
\begin{align*}
\int_I\widetilde{N}(t)^{-1}dt\lesssim_{(u,v)} N^{-3}\quad {\rm for \ all}\quad I\in [0,T_{\max})\quad  {\rm and\ all }\quad N\leq N_0.
\end{align*}
Under the assumption of (\ref{1011s1}), if we choose the interval $I$ sufficiently large inside $[0,T_{\max})$, we can derive a contradiction, which completes
the proof of this proposition. \hfill $\Box$

Now we have precluded the existence of the two types of almost periodic solution described in Proposition 6.6 and proved Theorem 7.\hfill$\Box$

\section{$\dot{H}^{s_c}\times \dot{H}^{s_c}$ scattering theory for (\ref{826x1}) in energy-supercritical cases}
\qquad In this section, we consider (\ref{826x1}) in the cases of $\alpha+\beta>2$ when $d=3$ and $\alpha+\beta>0$ when $d\geq 4$.
We will introduce some notations below.

Let
\begin{align*}
s_c=\frac{d}{2}-\frac{2}{\alpha+\beta+2}
\end{align*}
and $\hat{W}^{s,1}(\mathbb{R}^d)$ be the space of functions
$$
\|h\|_{\hat{W}^{s,1}(\mathbb{R}^d)}:=\|\langle\xi\rangle^s\hat{h}\|_{L^1(\mathbb{R}^d)}<+\infty.
$$

Denote $\alpha+\beta+2:=p$. By the results of \cite{Killip20102},  we know that
\begin{align}
&\quad \||\nabla|^{s_c} f\|_{L^{\infty}_tL^2_x(I\times\mathbb{R}^d)}+\||\nabla|^{s_c} f\|_{L^{\frac{2(d+2)}{d}}_{t,x}(I\times\mathbb{R}^d)}+
\|f\|_{L^{\infty}_tL^{\frac{dp}{2}}_x(I\times\mathbb{R}^d)}+\|f\|_{L^{\frac{p(d+2)}{2}}_{t,x}(I\times\mathbb{R}^d)}\nonumber\\
&\lesssim \||\nabla|^{s_c} f\|_{S^0(I)}\quad {\rm for \ all }\quad d\geq 3.\label{12281}
\end{align}

For $a\in \mathbb{R}^+$, define
\begin{equation*}
\chi_{\leq a}=\left\{
\begin{array}{lll}
1,\quad |x|\leq a,\\
{\rm smooth \ function}\quad a\leq |x|\leq \frac{11}{10}a,\\
0,\quad |x|\geq \frac{11}{10}a,
\end{array}\right.
\end{equation*}
$\chi_{\geq a}=1-\chi_{\leq a}$, $\chi_{a\leq \cdot \leq b}=\chi_{\leq b}-\chi_{\leq a}$, $\chi_a=\chi_{\leq 2a}-\chi_{\leq a}$ and
$\chi_{\sim a}=\chi_{\frac{1}{2}\leq \cdot \leq 4a}$.

Let $f\in L^1_{loc}(\mathbb{R}^d)$ be a radial function, we define the incoming component of $f$ as
$$
f_{in}(r)=r^{-\frac{d-1}{2}+2}\int_0^{+\infty}[J(-\rho r)+K(\rho r)]\rho^{d-1}\mathcal{F}f(\rho)d\rho
$$
and the outgoing component of $f$ as
$$
f_{out}(r)=r^{-\frac{d-1}{2}+2}\int_0^{+\infty}[J(-\rho r)-K(\rho r)]\rho^{d-1}\mathcal{F}f(\rho)d\rho.
$$
Here
\begin{align}
&J(r)=\int_0^{\frac{\pi}{2}}e^{2\pi r\sin\theta}\cos^{d-2}\theta d\theta,\quad K(r)=\chi_{\geq 1}(r)[-\frac{1}{2\pi ir}-\frac{d-3}{(2\pi ir)^3}],\\
&\mathcal{F}f(\rho)=\int_0^{+\infty}\int_{-\frac{\pi}{2}}^{\frac{\pi}{2}}e^{-2\pi i\rho r\sin\theta}\cos^{d-2}\theta r^{\frac{3}{2}(d-1)-2}f(r)d\theta dr\quad
{\rm for}\quad  d=3,4,5.
\end{align}

For the radial function $f\in \mathcal{S}(\mathbb{R}^d)$(the Schwartz function space on $\mathbb{R}^d$), define the modified outgoing component of $f$ as
$$
f_+=\frac{1}{2}(1-\chi_{\geq 1})f+\frac{1}{2}(1-P_{\geq 1}\chi_{\geq 1}f+(P_{\geq 1}\chi_{\geq 1} f)_{out}
$$
and the modified incoming component of $f$ as
$$
f_-=\frac{1}{2}(1-\chi_{\geq 1})f+\frac{1}{2}(1-P_{\geq 1}\chi_{\geq 1}f+(P_{\geq 1}\chi_{\geq 1} f)_{in},
$$
 where $P_{\geq 1}$ is the Littlewood-Paley operator which was defined in Section 2. Obviously,
 $$
 f=f_{out}+f_{in},\quad f=f_+ +f_-.
 $$

We need to introduce the working spaces as follows.

Let $\epsilon$ be a fixed small positive constant and $\delta_0$,$\gamma$ be the constants satisfying
\begin{align*}
\frac{1}{\delta_0}=\frac{2d+1}{4d-2}-\frac{2}{d},\quad \gamma=\frac{d-1}{2d-1}-\epsilon.
\end{align*}
Denote  $\infty-=\frac{1}{\epsilon}$ and $2+=2+\epsilon$

 Let $X(I)$, $X_0(I)$, $Y(I)$, $Z_w(I)$ and $Z_z(I)$ for $I\subset \mathbb{R}^+$ respectively be the spaces under the norms
\begin{align*}
\|h\|_{X(I)}&:=\|h\|_{L^{\infty}_t\dot{H}^{s_c}_x(I\times \mathbb{R}^d)}+\|h\|_{L^{2p}_tL^{dp}_x(I\times \mathbb{R}^d)}+
\|h\|_{L^{\frac{(d+2)p}{2}}_{t,x}(I\times \mathbb{R}^d)}+\||\nabla|^{s_c}h\|_{L^{\frac{2(d+2)}{d}}_{t,x}(I\times \mathbb{R}^d)},\\
\|h\|_{X_0(I)}&:=\|h\|_{L^{\infty-}_tL^{2+}_x(I\times \mathbb{R}^d)}+
\|h\|_{L^{\frac{(d+2)p}{2}}_tL^{\frac{2d(d+2)p}{d(d+2)p-8}}_x(I\times \mathbb{R}^d)}+\|h\|_{L^{\frac{2(d+2)}{d}}_{t,x}(I\times \mathbb{R}^d)},\\
\|h\|_{Y(I)}&:=\|h\|_{X(I)}+\|h\|_{X_0(I)},\\
\|h\|_{Z_w(I)}&:=\|h\|_{L^{2p}_tL^{dp}_x(I\times \mathbb{R}^d)}+\sup_{q\in [2, \frac{dp}{2}-\epsilon]}\|h\|_{L^{\infty}_tL^q_x(I\times\mathbb{R}^d)}+\|h\|_{L^{\infty-}_tL^{\frac{dp}{2}}_x(I\times \mathbb{R}^d)}\\
&\qquad +\|h\|_{(L^{\infty}_tH^{s_c}_x+L^2_t\hat{W}^{s_c-\gamma,\sigma_0}_x)(I\times \mathbb{R}^d)},\\
\|h\|_{Z_z(I)}&:=\|h\|_{L^{\frac{2(d+2)}{d}}_{t,x}(I\times \mathbb{R}^d)}+\|h\|_{L^{\frac{(d+2)p}{2}}_{t,x}(I\times \mathbb{R}^d)}+\|h\|_{L^{2p}_tL^{dp}_x(I\times \mathbb{R}^d)}+\sup_{q\in [2+\epsilon, \infty]}\|h\|_{L^{\infty-}_tL^q_x(I\times\mathbb{R}^d)}\\
&\qquad +\sup_{q\in [\frac{2d}{d-2}, \infty]}\||\nabla|^{s_c-\gamma}h\|_{L^2_tL^q_x(I\times\mathbb{R}^d)}
\end{align*}

\subsection{$\dot{H}^{s_c}\times \dot{H}^{s_c}$ scattering theory for (\ref{826x1}) with special radial initial data}
\qquad In this subsection, we will establish the $\dot{H}^{s_c}\times \dot{H}^{s_c}$ scattering theory for the solution of (\ref{826x1}) with special radial initial data when $d=3,4,5$.

We need some lemmas to prove Theorem 8.

Denote $e^{it\Delta} f=f_L$. Let $\psi_1=u-w_{1out,L}-z_{1L}$ and $\psi_2=v-w_{2out,L}-z_{2L}$. It is easy to verify that $\psi_1$ and $\psi_2$ satisfy
\begin{equation*}
\left\{
\begin{array}{lll}
i\partial_t \psi_1+\Delta \psi_1=\lambda|u|^{\alpha}|v|^{\beta+2}u,\quad i\partial_t \psi_2+\Delta \psi_2=\mu|u|^{\alpha+2}|v|^{\beta}v,\quad x\in \mathbb{R}^d,\quad t>0,\\
\psi_1(0,x)=\psi_{10}(x),\quad \psi_2(0,x)=\psi_{20}(x),\quad x\in \mathbb{R}^d.
\end{array}\right.
\end{equation*}
Duhamel formula implies that
\begin{align}
&\psi_1(t)=e^{it\Delta}\psi_{10}-i\lambda\int_0^te^{i(t-s)\Delta}(|u|^{\alpha}|v|^{\beta+2}u)(s)ds,\label{10185}\\ &\psi_2(t)=e^{it\Delta}\psi_{20}-i\mu\int_0^te^{i(t-s)\Delta}(|u|^{\alpha+2}|v|^{\beta}v)(s)ds.\label{1014w1}
\end{align}

{\bf Lemma 7.1.} {\it Assume that $d=3,4,5$, $0\in I$, $\psi_1, \psi_2\in X(I)$ and the assumptions of Theorem 8 hold. Then
\begin{align*}
&\quad \|\int_0^te^{i(t-s)\Delta}(|u|^{\alpha}|v|^{\beta+2}u)(s)ds\|_{X_0(I)}+\|\int_0^te^{i(t-s)\Delta}(|u|^{\alpha+2}|v|^{\beta}v)(s)ds\|_{X_0(I)}\\
&\lesssim [\|\psi_1\|_{X_0(I)}+\|\psi_2\|_{X_0(I)}+\delta_0][\|\psi_1\|^{\alpha+\beta+2}_{X(I)}+\|\psi_2\|^{\alpha+\beta+2}_{X(I)}+\delta_0^{\alpha+\beta+2}].
\end{align*}
}

{\bf Proof:}  Under the assumptions of (\ref{1011w1}) and (\ref{1011w2}), similar to the proofs of Lemma 6.2 and Lemma 6.3 in \cite{Beceanu??}, we can get
 $$\|w_{1out,L}\|_{Z_w(\mathbb{R}^+)}+\|w_{2out,L}\|_{Z_w(\mathbb{R}^+)}\leq \delta_0,\quad \|z_{1L}\|_{Z_z(\mathbb{R})}+\|z_{2L}\|_{Z_z(\mathbb{R})}\leq \delta_0.$$
Especially, for $p=\alpha+\beta+2$,
\begin{align}
&\|w_{1out,L}\|_{L^{\infty-}_tL^{2+}_x}+\|w_{1out,L}\|_{L^{2p}_tL^{dp}_x}
+\|w_{2out,L}\|_{L^{\infty-}_tL^{2+}_x}+\|w_{2out,L}\|_{L^{2p}_tL^{dp}_x}\leq \delta_0,\label{10151}\\
&\|z_{1L}\|_{L^{\infty-}_tL^{2+}_x}+\|z_{1L}\|_{L^{2p}_tL^{dp}_x}+\|z_{2L}\|_{L^{\infty-}_tL^{2+}_x}
+\|z_{2L}\|_{L^{2p}_tL^{dp}_x}\leq \delta_0.\label{10152}
\end{align}
Since $u=\psi_1+w_{1out,L}+z_{1L}$ and $v=\psi_2+w_{2out,L}+z_{2L}$, using (\ref{10151}) and (\ref{10152}), we have
\begin{align*}
&\|u\|_{L^{\infty-}_tL^{2+}_x}\lesssim\|\psi_1\|_{L^{\infty-}_tL^{2+}_x}+\|w_{1out,L}\|_{L^{\infty-}_tL^{2+}_x}+\|z_{1L}\|_{L^{\infty-}_tL^{2+}_x}\lesssim \|\psi_1\|_{X_0(I)}+\delta_0,\\
&\|u\|_{L^{2p}_tL^{dp}_x}\lesssim \|\psi_1\|_{L^{2p}_tL^{dp}_x}+\|w_{1out,L}\|_{L^{2p}_tL^{dp}_x}
+\|z_{1L}\|_{L^{2p}_tL^{dp}_x}\lesssim \|\psi_1\|_{X(I)}+\delta_0,
\end{align*}
and similarly
\begin{align*}
\|v\|_{L^{\infty-}_tL^{2+}_x}\lesssim \|\psi_2\|_{X_0(I)}+\delta_0,\quad \|v\|_{L^{2p}_tL^{dp}_x}\lesssim \|\psi_2\|_{X(I)}+\delta_0.
\end{align*}
All the norms above are on $(\mathbb{R}^+\times \mathbb{R}^d)$ and $p=\alpha+\beta+2$. By Strichartz estimates, we get
\begin{align*}
&\quad\|\int_0^te^{i(t-s)\Delta}(|u|^{\alpha}|v|^{\beta+2}u)(s)ds\|_{X_0(I)}+\|\int_0^te^{i(t-s)\Delta}(|u|^{\alpha+2}|v|^{\beta}v)(s)ds\|_{X_0(I)}\\
&\lesssim \|(|u|^{\alpha}|v|^{\beta+2}u)\|_{L^{2-}_tL^{\frac{2d}{d+2}+}(I\times\mathbb{R}^d)}
+\|(|u|^{\alpha+2}|v|^{\beta}v)\|_{L^{2-}_tL^{\frac{2d}{d+2}+}(I\times\mathbb{R}^d)}\\
&\lesssim [\|u\|_{L^{\infty-}_tL^{2+}_x(I\times\mathbb{R}^d)}+\|v\|_{L^{\infty-}_tL^{2+}_x(I\times\mathbb{R}^d)}][\|u\|^p_{L^{2p}_tL^{dp}_x(I\times\mathbb{R}^d)}
+\|v\|^p_{L^{2p}_tL^{dp}_x(I\times\mathbb{R}^d)}]\\
&\lesssim [\|\psi_1\|_{X_0(I)}+\|\psi_2\|_{X_0(I)}+\delta_0][\|\psi_1\|^p_{X(I)}+\|\psi_2\|^p_{X(I)}+\delta_0^p].
\end{align*}
Lemma 7.1 is proved.\hfill $\Box$

Now we will estimate the nonlinear terms in $X(I)$.

{\bf Lemma 7.2.} {\it Assume that $d=3,4,5$, $0\in I$, $\psi_1, \psi_2\in X(I)$ and the assumptions of Theorem 8 hold. Then
\begin{align*}
&\quad \|\int_0^te^{i(t-s)\Delta}(|u|^{\alpha}|v|^{\beta+2}u)(s)ds\|_{X(I)}+\|\int_0^te^{i(t-s)\Delta}(|u|^{\alpha+2}|v|^{\beta}v)(s)ds\|_{X(I)}\nonumber\\
&\lesssim \|\psi_1\|^{\alpha+\beta+3}_{X(I)}+\|\psi_2\|^{\alpha+\beta+3}_{X(I)}+\delta_0\|\psi_1\|^{\alpha+\beta+2}_{X(I)}+\delta_0\|\psi_2\|^{\alpha+\beta+2}_{X(I)}
+\delta_0^{\alpha+\beta+3}.
\end{align*}
}

{\bf Proof:} Using (\ref{12281}) and Strichartz's inequality in radial case, we get
\begin{align}
&\quad \|\int_0^te^{i(t-s)\Delta}(|u|^{\alpha}|v|^{\beta+2}u)(s)ds\|_{X(I)}+\|\int_0^te^{i(t-s)\Delta}(|u|^{\alpha+2}|v|^{\beta}v)(s)ds\|_{X(I)}\nonumber\\
&\lesssim \||\nabla|^{s_c-\gamma}(|u|^{\alpha}|v|^{\beta+2}u)\|_{L^2_tL_x^{\frac{4d-2}{2d+1}-}(I\times \mathbb{R}^d)}+\||\nabla|^{s_c-\gamma}(|u|^{\alpha+2}|v|^{\beta}v)\|_{L^2_tL_x^{\frac{4d-2}{2d+1}-}(I\times \mathbb{R}^d)}.\label{1015s3}
\end{align}
Here $\gamma=\frac{d-1}{2d-1}-$. Since $u=\psi_1+w_{1out,L}+z_{1L}$, $v=\psi_2+w_{2out,L}+z_{2L}$, $w_{1out,L}=w^I_{1out,L}+w^{II}_{1out,L}$ and
$w_{2out,L}=w^I_{2out,L}+w^{II}_{2out,L}$, using Lemma 7.1, H\"{o}lder's inequality and Young's inequality, we obtain
\begin{align}
&\quad \||\nabla|^{s_c-\gamma}(|u|^{\alpha}|v|^{\beta+2}u)\|_{L^2_tL^{\frac{4d-2}{2d+1}-}}
+\||\nabla|^{s_c-\gamma}(|u|^{\alpha+2}|v|^{\beta}v)\|_{L^2_tL^{\frac{4d-2}{2d+1}-}}\nonumber\\
&\lesssim [\|u\|^{\alpha+\beta+2}_{L^{2(\alpha+\beta+2)}_tL^{d(\alpha+\beta+2)}_x}+\|v\|^{\alpha+\beta+2}_{L^{2(\alpha+\beta+2)}_tL^{d(\alpha+\beta+2)}_x}]
\left\{[\||\nabla|^{s_c-\gamma}\psi_1\|_{L^{\infty}_tL^{q_1}_x}+\||\nabla|^{s_c-\gamma}\psi_2\|_{L^{\infty}_tL^{q_1}_x}]\right.\nonumber\\
&\quad+\left. [\||\nabla|^{s_c-\gamma}w^I_{1out,L}\|_{L^{\infty}_tL^{q_1}_x}+\||\nabla|^{s_c-\gamma}w^I_{2out,L}\|_{L^{\infty}_tL^{q_1}_x}]\right\}\nonumber\\
&\quad+[\|u\|^{\alpha+\beta+2}_{L^{\infty}_tL^{\frac{d(\alpha+\beta+2)}{2}-}_x}+\|v\|^{\alpha+\beta+2}_{L^{\infty}_tL^{\frac{d(\alpha+\beta+2)}{2}-}_x}][\||\nabla|^{s_c}w^{II}_{1out,L}\|_{L^2_tL^{\sigma_0}_x}+\||\nabla|^{s_c}w^{II}_{2out,L}\|_{L^2_tL^{\sigma_0}_x} ]\nonumber\\
&\quad +[\|u\|^{\alpha+\beta+2}_{L^{\infty-}_tL^{\frac{d(\alpha+\beta+2)}{2}}_x}+\|v\|^{\alpha+\beta+2}_{L^{\infty-}_tL^{\frac{d(\alpha+\beta+2)}{2}}_x}][\||\nabla|^{s_c-\gamma}z_{1,L}\|_{L^{2+}_tL^{q_2}_x}+\||\nabla|^{s_c-\gamma}z_{2,L}\|_{L^{2+}_tL^{q_2}_x} ]\nonumber\\
&:=(I)+(II)+(III)+(IV).\label{10181}
\end{align}
Here
$$
q_1=\frac{2d^2-3d+2}{d(2d+1)}+,\quad\sigma_0=\frac{d(4d-2)}{2d^2-7d+4}+,\quad q_2=\frac{d(4d-2)}{2d^2-7d+4}-.
$$

We will estimate (I), (II), (III) and (IV). By interpolation, using H\"{o}lder inequality and Young's inequality,  we have
\begin{align}
(I)&\lesssim [\||\nabla|^{s_c}\psi_1\|_{L^{\infty}_tL^2_x}+\||\nabla|^{s_c}\psi_2\|_{L^{\infty}_tL^2_x} ]\nonumber\\
&\quad \times\sum_{j=1}^2[\|\psi_j\|^{\alpha+\beta+2}_{L^{2(\alpha+\beta+2)}_tL^{d(\alpha+\beta+2)}_x}
+\|w_{jout,L}\|^{\alpha+\beta+2}_{L^{2(\alpha+\beta+2)}_tL^{d(\alpha+\beta+2)}_x}
+\|z_{jout,L}\|^{\alpha+\beta+2}_{L^{2(\alpha+\beta+2)}_tL^{d(\alpha+\beta+2)}_x}] \nonumber\\
&\lesssim [\|\psi_1\|_{X(I)}+\|\psi_2\|_{X(I)}]\nonumber\\
&\quad\times \sum_{j=1}^2[\|\psi_j\|^{\alpha+\beta+2}_{X(I)}
+\|w_{jout,L}\|^{\alpha+\beta+2}_{L^{2(\alpha+\beta+2)}_tL^{d(\alpha+\beta+2)}_x}
+\|z_{jout,L}\|^{\alpha+\beta+2}_{L^{2(\alpha+\beta+2)}_tL^{d(\alpha+\beta+2)}_x}]\nonumber\\
&\lesssim [\|\psi_1\|_{X(I)}+\|\psi_2\|_{X(I)}][\|\psi_1\|^{\alpha+\beta+2}_{X(I)}+\|\psi_2\|^{\alpha+\beta+2}_{X(I)}+\delta_0^{\alpha+\beta+2}]\nonumber\\
&\lesssim \|\psi_1\|^{\alpha+\beta+3}_{X(I)}+\|\psi_2\|^{\alpha+\beta+3}_{X(I)}+\delta_0^{\alpha+\beta+3}.\label{10182}
\end{align}
Similarly,
\begin{align}
(II)&\lesssim [\||\nabla|^{s_c}w^I_{1out,L}\|_{L^{\infty}_tL^2_x}+\||\nabla|^{s_c}w^I_{2out,L}\|_{L^{\infty}_tL^2_x} ]\nonumber\\
&\quad \times\sum_{j=1}^2[\|\psi_j\|^{\alpha+\beta+2}_{L^{2(\alpha+\beta+2)}_tL^{d(\alpha+\beta+2)}_x}
+\|w_{jout,L}\|^{\alpha+\beta+2}_{L^{2(\alpha+\beta+2)}_tL^{d(\alpha+\beta+2)}_x}
+\|z_{jout,L}\|^{\alpha+\beta+2}_{L^{2(\alpha+\beta+2)}_tL^{d(\alpha+\beta+2)}_x}] \nonumber\\
&\lesssim \delta_0\|\psi_1\|^{\alpha+\beta+2}_{X(I)}+\delta_0\|\psi_2\|^{\alpha+\beta+2}_{X(I)}+\delta_0^{\alpha+\beta+3}.\label{10183}
\end{align}
Using H\"{o}lder's inequality and Young's inequality,
\begin{align}
(III)&\lesssim [\||\nabla|^{s_c-\gamma}w^{II}_{1out,L}\|_{L^2_tL^{\sigma_0}_x}+\||\nabla|^{s_c-\gamma}w^{II}_{2out,L}\|_{L^2_tL^{\sigma_0}_x} ]\nonumber\\
&\quad \times\sum_{j=1}^2[\|\psi_j\|^{\alpha+\beta+2}_{L^{\infty}_tL^{\frac{d(\alpha+\beta+2)}{2}-}_x}
+\|w_{jout,L}\|^{\alpha+\beta+2}_{L^{\infty}_tL^{\frac{d(\alpha+\beta+2)}{2}-}_x}
+\|z_{jout,L}\|^{\alpha+\beta+2}_{L^{\infty}_tL^{\frac{d(\alpha+\beta+2)}{2}-}_x}] \nonumber\\
&\lesssim \delta_0\|\psi_1\|^{\alpha+\beta+2}_{X(I)}+\delta_0\|\psi_2\|^{\alpha+\beta+2}_{X(I)}+\delta_0^{\alpha+\beta+3}.\label{10183}
\end{align}

Note that $q_2>\frac{2d}{d-2}$ when $d\geq 3$. Using H\"{o}lder's inequality and Young's inequality, by interpolation,
\begin{align}
(IV)&\lesssim [\||\nabla|^{s_c-\gamma}z_{1,L}\|_{L^{2+}_tL^{q_2}_x}+\||\nabla|^{s_c-\gamma}z_{2,L}\|_{L^{2+}_tL^{q_2}_x} ]\nonumber\\
&\quad \times\sum_{j=1}^2[\|\psi_j\|^{\alpha+\beta+2}_{L^{\infty-}_tL^{\frac{d(\alpha+\beta+2)}{2}}_x}
+\|w_{jout,L}\|^{\alpha+\beta+2}_{L^{\infty-}_tL^{\frac{d(\alpha+\beta+2)}{2}}_x}
+\|z_{jout,L}\|^{\alpha+\beta+2}_{L^{\infty-}_tL^{\frac{d(\alpha+\beta+2)}{2}}_x}] \nonumber\\
&\lesssim \delta_0\|\psi_1\|^{\alpha+\beta+2}_{X(I)}+\delta_0\|\psi_2\|^{\alpha+\beta+2}_{X(I)}+\delta_0^{\alpha+\beta+3}.\label{10184}
\end{align}

Substituting (\ref{10181})--(\ref{10184}) into (\ref{1015s3}), we obtain the conclusion of this lemma.\hfill $\Box$

Now we give the proof of Theorem 8.

{\bf Proof of Theorem 8:} By Duhamel formula (\ref{10185}) and (\ref{1014w1}), for any $I\in \mathbb{R}^+$,
\begin{align*}
&\quad \|\psi_1\|_{Y(I)}+\|\psi_1\|_{Y(I)}\\
&\lesssim \|e^{it\Delta}\psi_{10}\|_{Y(\mathbb{R}^+)}+ \|\int_0^te^{i(t-s)\Delta}(|u|^{\alpha}|v|^{\beta+2}u)(s)ds\|_{X_0(I)}+\|\int_0^te^{i(t-s)\Delta}(|u|^{\alpha}|v|^{\beta+2}u)(s)ds\|_{X(I)}\\ &\quad+\|e^{it\Delta}\psi_{20}\|_{Y(\mathbb{R}^+)}+ \|\int_0^te^{i(t-s)\Delta}(|u|^{\alpha+2}|v|^{\beta}v)(s)ds\|_{X_0(I)}+\|\int_0^te^{i(t-s)\Delta}(|u|^{\alpha+2}|v|^{\beta}v)(s)ds\|_{X(I)}\\
&\lesssim \delta_0+ \|\psi_1\|^{\alpha+\beta+3}_{Y(I)}+\|\psi_2\|^{\alpha+\beta+3}_{Y(I)}+\delta_0\|\psi_1\|^{\alpha+\beta+2}_{Y(I)}+\delta_0\|\psi_2\|^{\alpha+\beta+2}_{Y(I)}
+\delta_0^{\alpha+\beta+2}\|\psi_1\|_{Y(I)}\\
&\qquad \ +\delta_0^{\alpha+\beta+2}\|\psi_2\|_{Y(I)}+\delta_0^{\alpha+\beta+3},
\end{align*}
where the implicit constant is independent on $I$. By the standard continuity argument, we can deduce that $(\psi_1,\psi_2)$ is global existence and
$$
\|\psi_1\|_{Y(\mathbb{R}^+)}+\|\psi_2\|_{Y(\mathbb{R}^+)}\lesssim \delta_0.
$$

Next, we establish the scattering result. Let
\begin{align*}
\psi_{10+}=\psi_{10}-i\lambda\int_0^{+\infty}e^{-is\Delta}(|u|^{\alpha}|v|^{\beta+2}u)(s)ds,\\
\psi_{20+}=\psi_{20}-i\mu\int_0^{+\infty}e^{-is\Delta}(|u|^{\alpha+2}|v|^{\beta}v)(s)ds.
\end{align*}
Then
\begin{align*}
\psi_1-e^{it\Delta}\psi_{10+}=-i\lambda\int_t^{+\infty}e^{i(t-s)\Delta}(|u|^{\alpha}|v|^{\beta+2}u)(s)ds,\\
\psi_2-e^{it\Delta}\psi_{20+}=-i\mu\int_t^{+\infty}e^{i(t-s)\Delta}(|u|^{\alpha+2}|v|^{\beta}v)(s)ds.
\end{align*}
Using Strichartz estimates and the discussion above, we have
\begin{align*}
&\quad \|\psi_1-e^{it\Delta}\psi_{10+}\|_{\dot{H}^{s_c}_x}+\|\psi_2-e^{it\Delta}\psi_{20+}\|_{\dot{H}^{s_c}_x}\\
&\lesssim \|\psi_1\|^{\alpha+\beta+3}_{Y([t,+\infty))}+\|\psi_2\|^{\alpha+\beta+3}_{Y([t,+\infty))}+\|w_{1out,L}\|^{\alpha+\beta+3}_{Z_w([t,+\infty))}
+\|w_{2out,L}\|^{\alpha+\beta+3}_{Z_w([t,+\infty))}\\
&\quad+\|z_{1,L}\|^{\alpha+\beta+3}_{Z_z([t,+\infty))}
+\|z_{2,L}\|^{\alpha+\beta+3}_{Z_z([t,+\infty))}\\
&\quad +\left\{[\|w_{1out,L}\|_{Z_w([t,+\infty))}+\|w_{2out,L}\|_{Z_w([t,+\infty))}+\|z_{1,L}\|_{Z_z([t,+\infty))}+\|z_{2,L}\|_{Z_z([t,+\infty))}]
\right.\\
&\qquad\quad \left. \times [\|\psi_1\|^{\alpha+\beta+2}_{Y([t,+\infty))}+\|\psi_2\|^{\alpha+\beta+2}_{Y([t,+\infty))}] \right\}\displaybreak\\
&\quad +\left\{[\|w_{1out,L}\|_{Z_w([t,+\infty))}+\|w_{2out,L}\|_{Z_w([t,+\infty))}+\|z_{1,L}\|_{Z_z([t,+\infty))}+\|z_{2,L}\|_{Z_z([t,+\infty))}]^{\alpha+\beta+2}
\right.\\
&\qquad \quad \left. \times [\|\psi_1\|_{Y([t,+\infty))}+\|\psi_2\|_{Y([t,+\infty))}] \right\}\\
&\longrightarrow 0\quad {\rm as}\quad t\rightarrow +\infty
\end{align*}
because
\begin{align*}
&\|\psi_1\|_{Y(\mathbb{R}^+)}+\|\psi_2\|_{Y(\mathbb{R}^+)}+\|w_{1out,L}\|_{Z_w(\mathbb{R}^+))}+\|w_{2out,L}\|_{Z_w(\mathbb{R}^+)}
+\|z_{1,L}\|_{Z_z(\mathbb{R}^+)}+\|z_{2,L}\|_{Z_z(\mathbb{R}^+)}\\
&<+\infty.
\end{align*}
Consequently,
$$
\|\psi_1-e^{it\Delta}\psi_{10+}\|_{\dot{H}^{s_c}_x}+\|\psi_2-e^{it\Delta}\psi_{20+}\|_{\dot{H}^{s_c}_x}\rightarrow 0\quad {\rm as }\quad t\rightarrow +\infty.
$$
Letting
$$
u_+=(P_{\geq 1}\chi_{\geq 1}w_1)_{out}+z_1+\psi_{10+},\quad v_+=(P_{\geq 1}\chi_{\geq 1}w_2)_{out}+z_2+\psi_{20+},
$$
we can get
$$
\lim_{t\rightarrow +\infty}[\|u(t)-e^{it\Delta}u_+\|_{\dot{H}^{s_c}_x(\mathbb{R}^d)}+\|v(t)-e^{it\Delta}v_+\|_{\dot{H}^{s_c}_x(\mathbb{R}^d)}]=0.
$$
(\ref{1011w3}) is obtained. Similarly, (\ref{1011w4}) can be obtained. Theorem 8 is proved. \hfill $\Box$

\subsection{$\dot{H}^{s_c}\times \dot{H}^{s_c}$ scattering theory for (\ref{826x1}) with initial data has many bubbles}
\qquad In this subsection, we will establish $\dot{H}^{s_c}\times \dot{H}^{s_c}$ scattering theory for the solution of (\ref{826x1}) with initial data has many bubbles.

{\bf Proof of Theorem 9:} Let $(w_{1k,L}, w_{2k,L})=(e^{it\Delta}w_{1k},e^{it\Delta}w_{2k})$ and
\begin{align*}
w_{1,L}=\sum_{k=0}^{+\infty}w_{1k,L},\quad w_{2,L}=\sum_{k=0}^{+\infty}w_{2k,L}.
\end{align*}
By Strichartz estimate, we have for $\gamma\in [0,s_c)$,
\begin{align}
&\||\nabla|^{\gamma}w_{1,L}\|_{L^{q\gamma}_{t,x}(\mathbb{R}\times \mathbb{R}^d)}\leq \epsilon^{\frac{s_c-\gamma}{d}}\|w_1\|_{\dot{H}^{s_c}_x}
\leq \epsilon^{\frac{s_c-\gamma}{d}-\alpha_0},\label{11011'}\\
&\||\nabla|^{\gamma}w_{2,L}\|_{L^{q\gamma}_{t,x}(\mathbb{R}\times \mathbb{R}^d)}\leq \epsilon^{\frac{s_c-\gamma}{d}}\|w_2\|_{\dot{H}^{s_c}_x}
\leq \epsilon^{\frac{s_c-\gamma}{d}-\alpha_0}.\label{11011}
\end{align}

Let $\varphi=u-w_{1,L}$ and $\psi=v-w_{2,L}$. Then
\begin{equation}
\left\{
\begin{array}{lll}
&i\partial_t \varphi+\Delta \varphi=\lambda|u|^{\alpha}|v|^{\beta+2}u,\quad i\partial_t \psi+\Delta \psi=\mu|u|^{\alpha+2}|v|^{\beta}v,\quad x\in\mathbb{R}^d,\quad t>0,\\
&\varphi(0,x)=\varphi_0,\quad \psi(0,x)=\psi_0, \quad x\in \mathbb{R}^d.
\end{array}\right.
\end{equation}
We introduce the working space and let $X(I)$ for $I\subset \mathbb{R}$ be the space with the following norm
\begin{align*}
\|(\varphi,\psi)\|_{X(I)}&:=\epsilon^{-\frac{s_c}{2d}}\sup_{\gamma \in [0,s_c-\gamma_0]}\||\nabla|^{\gamma}\varphi\|_{L^{q_{\gamma}}_{t,x}}
+\||\nabla|^{s_c}\varphi\|_{L^{\infty}_tL^2_x}+\||\nabla|^{s_c}\varphi\|_{L^{\frac{2(d+2)}{d}}_{t,x}}\\
&\quad+\epsilon^{-\frac{s_c}{2d}}\sup_{\gamma \in [0,s_c-\gamma_0]}\||\nabla|^{\gamma}\psi\|_{L^{q_{\gamma}}_{t,x}}
+\||\nabla|^{s_c}\psi\|_{L^{\infty}_tL^2_x}+\||\nabla|^{s_c}\psi\|_{L^{\frac{2(d+2)}{d}}_{t,x}},
\end{align*}
all norms above are on $I\times \mathbb{R}^d$. Here $\gamma_0$ is a parameter to be determined later while
\begin{align*}
q_{\gamma}=\frac{2(d+2)}{d-2(s_c-\gamma)}.
\end{align*}
Duhamel formulae imply that
\begin{align}
&\varphi(t)=e^{it\Delta} \varphi_0-i\lambda\int_0^te^{i(t-s)\Delta}(|u|^{\alpha}|v|^{\beta+2}u)(s)ds,\label{1101x1}\\
&\psi(t)=e^{it\Delta} \psi_0-i\mu\int_0^te^{i(t-s)\Delta}(|u|^{\alpha+2}|v|^{\beta}v)(s)ds.\label{1101x2}
\end{align}
Using Strichartz estimates, from (\ref{1101x1}) and (\ref{1101x2}), we get
\begin{align}
&\quad\||\nabla|^{\gamma}\varphi\|_{L^{q_{\gamma}}_{t,x}(I\times \mathbb{R}^d)}+\||\nabla|^{\gamma}\psi\|_{L^{q_{\gamma}}_{t,x}(I\times \mathbb{R}^d)}\nonumber\\
&\lesssim \|\varphi_0\|_{\dot{H}^{s_c}(\mathbb{R}^d)}+\|\psi_0\|_{\dot{H}^{s_c}(\mathbb{R}^d)} +\||\nabla|^{\gamma}(|u|^{\alpha}|v|^{\beta+2}u)\|_{L^{q'_1}_tL^{q'_{\gamma}}_x(I\times \mathbb{R}^d)}\nonumber\\
&\quad+\||\nabla|^{\gamma}(|u|^{\alpha+2}|v|^{\beta}v)\|_{L^{q'_1}_tL^{q'_{\gamma}}_x(I\times \mathbb{R}^d)}
\end{align}
for any $I\subset \mathbb{R}$ and $\gamma\in [0,s_c-\gamma_0]$, where $q'_1$ satisfies
$$
\frac{1}{q'_1}=\frac{d+4-2(d+1)(s_c-\gamma)}{2(d+2)}.
$$

Letting
\begin{align*}
q_2=\frac{(d+2)(\alpha+\beta+2)}{2-d(s_c-\gamma)},\quad r_2=\frac{(d+2)(\alpha+\beta+2)}{2[1+(s_c-\gamma)]},
\end{align*}
then
\begin{align}
\frac{1}{q'_{\gamma}}=\frac{1}{q_{\gamma}}+\frac{\alpha+\beta+2}{r_2}.\label{12291}
\end{align}
Using H\"{o}lder's ineuqlity, the fractional Leibniz rules and Young's inequality, we get
\begin{align}
&\quad \||\nabla|^{\gamma}(|u|^{\alpha}|v|^{\beta+2}u)\|_{L^{q'_1}_tL^{q'_{\gamma}}_x(I\times \mathbb{R}^d)}\nonumber\\
&\leq \left\{\int_I\left[\left(\int_{\mathbb{R}^d}||\nabla|^{\gamma}(|u|^{\alpha}u)|^{p_1}dx\right)^{\frac{1}{p_1}}
\left(\int_{\mathbb{R}^d}||v|^{\beta+2})|^{p_2}dx\right)^{\frac{1}{p_2}}\right]^{q'_1}dt\right\}^{\frac{1}{q'_1}}\nonumber\\
&\quad+ \left\{\int_I\left[\left(\int_{\mathbb{R}^d}||u|^{\alpha}u|^{p_3}dx\right)^{\frac{1}{p_3}}\left(\int_{\mathbb{R}^d}||\nabla|^{\gamma}(|v|^{\beta+2})|^{p_4}dx\right)^{\frac{1}{p_4}}
\right]^{q'_1}dt\right\}^{\frac{1}{q'_1}}\nonumber\\
&\leq \left\{\int_I\left[\left(\int_{\mathbb{R}^d}||\nabla|^{\gamma}u|^{q_{\gamma}}dx\right)^{\frac{1}{q_{\gamma}}}
\left(\int_{\mathbb{R}^d}|u|^{r_2}dx\right)^{\frac{1}{p_6}}
\left(\int_{\mathbb{R}^d}||v|^{\beta+2})|^{p_2}dx\right)^{\frac{1}{p_2}}\right]^{q'_1}dt\right\}^{\frac{1}{q'_1}}\nonumber\\
&\quad+ \left\{\int_I\left[\left(\int_{\mathbb{R}^d}||u|^{\alpha}u|^{p_3}dx\right)^{\frac{1}{p_3}}
\left(\int_{\mathbb{R}^d}||\nabla|^{\gamma}v|^{q_{\gamma}}dx\right)^{\frac{1}{q_{\gamma}}}\left(\int_{\mathbb{R}^d}|v|^{r_2}dx\right)^{\frac{1}{p_7}}
\right]^{q'_1}dt\right\}^{\frac{1}{q'_1}}\nonumber\\
&\lesssim [\||\nabla|^{\gamma}u\|_{L^{q_{\gamma}}_{t,x}(I\times\mathbb{R}^d)}+\||\nabla|^{\gamma}v\|_{L^{q_{\gamma}}_{t,x}(I\times\mathbb{R}^d)}]
[\|u\|^{\alpha+\beta+2}_{L^{q_2}_tL^{r_2}_x(I\times\mathbb{R}^d)}+\|v\|^{\alpha+\beta+2}_{L^{q_2}_tL^{r_2}_x(I\times\mathbb{R}^d)}].\label{1101w1}
\end{align}
Here
\begin{align*}
&\frac{1}{q'_{\gamma}}=\frac{1}{p_1}+\frac{1}{p_2},\quad \frac{1}{p_1}=\frac{1}{q_{\gamma}}+\frac{1}{p_6},\quad \frac{1}{p_2}=\frac{\beta+2}{r_2},\quad \frac{1}{p_6}=\frac{\alpha}{r_2},\\
&\frac{1}{q'_{\gamma}}=\frac{1}{p_3}+\frac{1}{p_4},\quad \frac{1}{p_3}=\frac{\alpha+1}{r_2},\quad\frac{1}{p_4}=\frac{1}{q_{\gamma}}+\frac{1}{p_7},\quad \frac{1}{p_7}=\frac{\beta+1}{r_2},
\end{align*}
and recalling (\ref{12291}), if $s_c-\gamma$ is small enough, the $p_j>1$, $j=1,2,...,7$.

Similarly,
\begin{align}
&\quad \||\nabla|^{\gamma}(|u|^{\alpha+2}|v|^{\beta}v)\|_{L^{q'_1}_tL^{q'_{\gamma}}_x(I\times \mathbb{R}^d)}\nonumber\\
&\lesssim [\||\nabla|^{\gamma}u\|_{L^{q_{\gamma}}_{t,x}(I\times\mathbb{R}^d)}+\||\nabla|^{\gamma}v\|_{L^{q_{\gamma}}_{t,x}(I\times\mathbb{R}^d)}]
[\|u\|^{\alpha+\beta+2}_{L^{q_2}_tL^{r_2}_x(I\times\mathbb{R}^d)}+\|v\|^{\alpha+\beta+2}_{L^{q_2}_tL^{r_2}_x(I\times\mathbb{R}^d)}],\label{1101w2}
\end{align}

Choosing $\gamma_0\in [0, s_c)$ such that $s_c-\gamma<s_c-\gamma_0$ and $\theta=\frac{d(d+2)(s_c-\gamma)}{2}<1$, by interpolation and Sobolev embedding, we can get
\begin{align}
&\quad \|u\|_{L^{q_2}_tL^{r_2}_x(I\times \mathbb{R}^d)}\lesssim
\|u\|^{\theta}_{L^{\infty}_tL^{\frac{d(\alpha+\beta+2)}{2}}_x(I\times\mathbb{R}^d)}\|u\|^{1-\theta}_{L^{\frac{(d+2)(\alpha+\beta+2)}{2}}_{t,x}(I\times \mathbb{R}^d)}\nonumber\\
&\lesssim \|u\|^{\theta}_{L^{\infty}_t\dot{H}^{s_c}_x(I\times\mathbb{R}^d)}\|u\|^{1-\theta}_{L^{\frac{(d+2)(\alpha+\beta+2)}{2}}_{t,x}(I\times \mathbb{R}^d)},\label{1101w3}\\
&\|v\|_{L^{q_2}_tL^{r_2}_x(I\times \mathbb{R}^d)}\lesssim \|v\|^{\theta}_{L^{\infty}_t\dot{H}^{s_c}_x(I\times\mathbb{R}^d)}\|v\|^{1-\theta}_{L^{\frac{(d+2)(\alpha+\beta+2)}{2}}_{t,x}(I\times \mathbb{R}^d)}.\label{1101w4}
\end{align}
Under assumption (C9), using (\ref{11011'}) and (\ref{11011}), we have
\begin{align}
&\|u\|_{L^{\infty}_t\dot{H}^{s_c}_x(I\times\mathbb{R}^d)}\lesssim \|w_{1,L}\|_{L^{\infty}_t\dot{H}^{s_c}_x(I\times\mathbb{R}^d)}+\|\varphi\|_{L^{\infty}_t\dot{H}^{s_c}_x(I\times\mathbb{R}^d)}\lesssim \epsilon^{-\alpha_0}+\|\varphi\|_{X(I)},\label{1101w5}\\
&\|u\|_{L^{\frac{(d+2)(\alpha+\beta+2)}{2}}_{t,x}(I\times\mathbb{R}^d)}\lesssim \|w_{1,L}\|_{L^{\frac{(d+2)(\alpha+\beta+2)}{2}}_{t,x}(I\times\mathbb{R}^d)}+\|\varphi\|_{L^{\frac{(d+2)(\alpha+\beta+2)}{2}}_{t,x}(I\times\mathbb{R}^d)}\nonumber\\
&\lesssim \epsilon^{\frac{s_c}{d}-\alpha_0}+\epsilon^{\frac{s_c}{d}}\|\varphi\|_{X(I)},\label{1101w6}\\
&\|v\|_{L^{\infty}_t\dot{H}^{s_c}_x(I\times\mathbb{R}^d)}\lesssim \epsilon^{-\alpha_0}+\|\psi\|_{X(I)},\label{1101w7}\\
&\|v\|_{L^{\frac{(d+2)(\alpha+\beta+2)}{2}}_{t,x}(I\times\mathbb{R}^d)}\lesssim \epsilon^{\frac{s_c}{d}-\alpha_0}+\epsilon^{\frac{s_c}{d}}\|\psi\|_{X(I)}.\label{1101w8}
\end{align}
Putting (\ref{1101w1})--(\ref{1101w8}) together, denoting $\epsilon=(s_c-\gamma)$, we obtain
\begin{align}
&\quad \||\nabla|^{\gamma}(|u|^{\alpha}|v|^{\beta+2}u)\|_{L^{q'_1}_tL^{q'_{\gamma}}_x(I\times \mathbb{R}^d)}+ \||\nabla|^{\gamma}(|u|^{\alpha+2}|v|^{\beta}v)\|_{L^{q'_1}_tL^{q'_{\gamma}}_x(I\times \mathbb{R}^d)}\nonumber\\
&\lesssim \epsilon^{\frac{s_c}{2d}}[\|\varphi\|_{X(I)}+\|\psi\|_{X(I)}]\times \left(\epsilon^{\frac{s_c}{d}(1-\theta)-\alpha_0}+\epsilon^{-\alpha_0\theta+\frac{s_c}{2d}(1-\theta)}[\|\varphi\|^{1-\theta}_{X(I)}+\|\psi\|^{1-\theta}_{X(I)}] \right.\nonumber\\
&\quad\left.+\epsilon^{(\frac{s_c}{d}-\alpha_0)(1-\theta)}[\|\varphi\|^{\theta}_{X(I)}+\|\psi\|^{\theta}_{X(I)}]+\epsilon^{\frac{s_c}{2d}(1-\theta)}
[\|\varphi\|_{X(I)}+\|\psi\|_{X(I)}]\right)^{\alpha+\beta+2}.\label{1101w9}
\end{align}
By Duhamel formulae (\ref{1101x1}), (\ref{1101x2}) and (\ref{1101w9}), using Cauchy-Schwartz inequality, we obtain
\begin{align}
&\quad \epsilon^{-\frac{s_c}{2d}}[\||\nabla|^{\gamma}\varphi\|_{L^{q_{\gamma}}_{t,x}(I\times \mathbb{R}^d)}+\||\nabla|^{\gamma}\psi\|_{L^{q_{\gamma}}_{t,x}(I\times \mathbb{R}^d)}]\nonumber\\
&\lesssim \epsilon^{1-\frac{s_c}{2d}}+[\|\varphi\|_{X(I)}+\|\psi\|_{X(I)}]\times \left(\epsilon^{\frac{s_c}{d}(1-\theta)-\alpha_0}+\epsilon^{-\alpha_0\theta+\frac{s_c}{2d}(1-\theta)}[\|\varphi\|^{1-\theta}_{X(I)}+\|\psi\|^{1-\theta}_{X(I)}] \right.\nonumber\\
&\quad\left.+\epsilon^{(\frac{s_c}{d}-\alpha_0)(1-\theta)}[\|\varphi\|^{\theta}_{X(I)}+\|\psi\|^{\theta}_{X(I)}]+\epsilon^{\frac{s_c}{2d}(1-\theta)}
[\|\varphi\|_{X(I)}+\|\psi\|_{X(I)}]\right)^{\alpha+\beta+2}\nonumber\\
&\lesssim \epsilon^{1-\frac{s_c}{2d}}+\epsilon^{\frac{(\alpha+\beta+2)s_c}{4d}}[\|\varphi\|_{X(I)}+\|\psi\|_{X(I)}]
+[\|\varphi\|_{X(I)}+\|\psi\|_{X(I)}]^{\alpha+\beta+3}\label{11021}
\end{align}
if we choose $\alpha_0$ and $s_c-\gamma_0$ sufficiently small.

Meanwhile, we have
\begin{align*}
&\quad \||\nabla|^{s_c}\varphi\|_{L^{\infty}_tL^2_x(I\times \mathbb{R}^d)}+\|\varphi\|_{L^{\frac{(d+2)(\alpha+\beta+2)}{2}}_{t,x}(I\times \mathbb{R}^d)}+\||\nabla|^{s_c}\psi\|_{L^{\infty}_tL^2_x(I\times \mathbb{R}^d)}+\|\psi\|_{L^{\frac{(d+2)(\alpha+\beta+2)}{2}}_{t,x}(I\times \mathbb{R}^d)}\\
&\lesssim \|\varphi_0\|_{\dot{H}^{s_c}_x(\mathbb{R}^d)}+\|\psi_0\|_{\dot{H}^{s_c}_x(\mathbb{R}^d)}
+\||\nabla|^{s_c}(|u|^{\alpha}|v|^{\beta+2}u)\|_{L^{\frac{2d+4}{d+4}}_{t,x}(I\times\mathbb{R}^d)}
+\||\nabla|^{s_c}(|u|^{\alpha+2}|v|^{\beta}v)\|_{L^{\frac{2d+4}{d+4}}_{t,x}(I\times\mathbb{R}^d)}.
\end{align*}
Recalling (\ref{11011'}) and (\ref{11011}), we get
\begin{align*}
&\||\nabla|^{s_c}u\|_{L^{\frac{2d+4}{d}}_{t,x}(I\times \mathbb{R}^d)}\lesssim \||\nabla|^{s_c}w_{1,L}\|_{L^{\frac{2d+4}{d}}_{t,x}(I\times \mathbb{R}^d)}+\||\nabla|^{s_c}\varphi\|_{L^{\frac{2d+4}{d}}_{t,x}(I\times \mathbb{R}^d)}\lesssim \epsilon^{-\alpha_0}+\|\varphi\|_{X(I)},\\
&\||\nabla|^{s_c}v\|_{L^{\frac{2d+4}{d}}_{t,x}(I\times \mathbb{R}^d)}\lesssim \epsilon^{-\alpha_0}+\|\psi\|_{X(I)},\\ &\|u\|_{L^{\frac{(d+2)(\alpha+\beta+2)}{2}}_{t,x}(I\times \mathbb{R}^d)}+\|v\|_{L^{\frac{(d+2)(\alpha+\beta+2)}{2}}_{t,x}(I\times \mathbb{R}^d)}\lesssim \epsilon^{\frac{s_c}{d}-\alpha_0}+\epsilon^{\frac{s_c}{2d}}[\|\varphi\|_{X(I)}+\|\psi\|_{X(I)}].
\end{align*}
Similar to (\ref{11021}), we can obtain
\begin{align}
&\quad \||\nabla|^{s_c}\varphi\|_{L^{\infty}_tL^2_x(I\times \mathbb{R}^d)}+\|\varphi\|_{L^{\frac{(d+2)(\alpha+\beta+2)}{2}}_{t,x}(I\times \mathbb{R}^d)}+\||\nabla|^{s_c}\psi\|_{L^{\infty}_tL^2_x(I\times \mathbb{R}^d)}+\|\psi\|_{L^{\frac{(d+2)(\alpha+\beta+2)}{2}}_{t,x}(I\times \mathbb{R}^d)}\nonumber\\
&\lesssim \epsilon+\epsilon^{\frac{(\alpha+\beta+2)s_c}{4d}}+\epsilon^{\frac{(\alpha+\beta+2)s_c}{4d}}[\|\varphi\|_{X(I)}+\|\psi\|_{X(I)}]
+[\|\varphi\|_{X(I)}+\|\psi\|_{X(I)}]^{\alpha+\beta+3}.\label{11022}
\end{align}
By (\ref{11021}) and (\ref{11022}), choosing $s_c-\gamma$ small enough, we have for some $\alpha_0>0$,
\begin{align*}
[\|\varphi\|_{X(I)}+\|\psi\|_{X(I)}]\lesssim \epsilon^{\alpha_0}+[\|\varphi\|_{X(I)}+\|\psi\|_{X(I)}]^{\alpha+\beta+3}.
\end{align*}
By the standard bootstrap argument, we can get
\begin{align*}
[\|\varphi\|_{X(I)}+\|\psi\|_{X(I)}]\lesssim \epsilon^{\alpha_0}
\end{align*}
uniformly in interval $I\subset \mathbb{R}$, which implies the global existence of the solution to (\ref{826x1}).

Letting
\begin{align*}
\varphi_{0+}=\varphi_0-i\lambda\int_0^{+\infty}e^{-is\Delta}(|u|^{\alpha}|v|^{\beta+2}u)(s)ds,\\
\psi_{0+}=\psi_0-i\mu\int_0^{+\infty}e^{-is\Delta}(|u|^{\alpha+2}|v|^{\beta}v)(s)ds.
\end{align*}
and
$$
u_+=w_1+\varphi_{0+},\quad v_+=w_2+\psi_{0+},
$$
then using Strichartz estimate, we can get
$$
\lim_{t\rightarrow +\infty}[\|u(t)-e^{it\Delta}u_+\|_{\dot{H}^{s_c}_x(\mathbb{R}^d)}+\|v(t)-e^{it\Delta}v_+\|_{\dot{H}^{s_c}_x(\mathbb{R}^d)}]=0
$$
Similarly, we can prove that
$$
\lim_{t\rightarrow -\infty}[\|u(t)-e^{it\Delta}u_-\|_{\dot{H}^{s_c}_x(\mathbb{R}^d)}+\|v(t)-e^{it\Delta}v_-\|_{\dot{H}^{s_c}_x(\mathbb{R}^d)}]=0.
$$
 Theorem 9 is proved. \hfill $\Box$

\section{Some results on (\ref{826x1}) with focusing nonlinearities}
\qquad In this short section, we give some results on (\ref{826x1}) with focusing nonlinearities and some discussions.

{\bf A. Some results on (\ref{826x1}) with focusing nonlinearities}

If $(u,v)$ is the $H^1\times H^1$-solution of (\ref{826x1}) when $d=3$ with $\lambda<0$, $\mu<0$, $\alpha\geq 0$, $\beta\geq 0$ and $\alpha+\beta\leq 2$, then we can define
\begin{align}
\tilde{E}_w(u,v):=\int_{\mathbb{R}^3}[\frac{\alpha+2}{2|\lambda|}|\nabla u|^2+\frac{\beta+2}{2|\mu|}|\nabla v|^2-|u|^{\alpha+2}|v|^{\beta+2}]dx\label{12061}
\end{align}
and get $\tilde{E}_w(u,v)=\tilde{E}_w(u_0,v_0)$. Using the convexity method in the famous paper \cite{Glassey1977}, we can prove that
the solution $(u,v)$ will blow up in finite time if $\tilde{E}_w(u_0,v_0)\leq 0$ and $\Im\int_{\mathbb{R}^3}[c_1\bar{u}_0x\cdot \nabla u_0+c_2\bar{v}_0x\cdot \nabla v_0]dx<0$. In fact, denoting $(c_1,c_2)=(\frac{\alpha+2}{2|\lambda|},\frac{\beta+2}{2|\mu|})$
\begin{align}
V(t)=\int_{\mathbb{R}^3}|x|^2[c_1|u(t,x)|^2+c_2|v(t,x)|^2]dx\label{12062}
\end{align}
as a function of $t>0$, we have
\begin{align}
&V'(t)=4\Im\int_{\mathbb{R}^3}[c_1\bar{u}x\cdot \nabla u+c_2\bar{v}x\cdot \nabla v]dx,\label{12063}\\
&V''(t)=-2[3(\alpha+\beta)+2]\int_{\mathbb{R}^3}|u|^{\alpha+2}|v|^{\beta+2}dx+8\tilde{E}_w(u_0,v_0).\label{12064}
\end{align}
Now following the standard discussions similar to those in the famous paper \cite{Glassey1977}, it is easy to prove that there exists $T_{\max}>0$ such that
\begin{align}
\lim_{t\rightarrow T_{\max}}\int_{\mathbb{R}^3}(c_1|\nabla u|^2+c_2|\nabla v|^2)dx=+\infty.
\end{align}

Similarly, if $(u,v)$ is the $H^1\times H^1$-solution of (\ref{826x1}) when $d=4$ with $\lambda<0$, $\mu<0$, and $(\alpha,\beta)=(0,0)$, then
the solution $(u,v)$ will blow up in finite time if $\tilde{E}_w(u_0,v_0)\leq 0$ and $\Im\int_{\mathbb{R}^4}[c_1\bar{u}_0x\cdot \nabla u_0+c_2\bar{v}_0x\cdot \nabla v_0]dx<0$.

However, if $\lambda<0$, $\mu<0$, $\alpha+\beta>2$ when $d=3$ and $\alpha+\beta>0$ when $d\geq 4$, there are more complicate phenomena on (\ref{826x1}). Similar to the discussions in Section 7, we can establish $\dot{H}^{s_c}\times \dot{H}^{s_c}$ scattering theory for the global solutions of (\ref{826x1}) under certain assumptions. On the other hand, there exists the solution of (\ref{826x1}) which is endowed with some type of norm and will blow up in finite time.  We would like to discuss it in the next paper.

{\bf B. The generalization of the definition of weighted gradient system of Schr\"{o}dinger equations}

The definition of  weighted(or essential) gradient system of Schr\"{o}dinger equations can be generalized to the following case
\begin{align}
iu_{jt}+\Delta u_j=f_j(|u_1|^2,...,|u_j|^2,...,|u_n|^2)u_j,\quad j=1,2,...,m.\label{1225xw1}
\end{align}
(\ref{1225xw1}) is called as a weighted gradient(or essential gradient) system of Schr\"{o}dinger equations if there exist $a_j\in\mathbb{R}$ and real value function $G(w_1,...,w_j,...,w_m)$ such that $$\frac{\partial G}{\partial w_j}=a_jf_j(w_1,...,w_j,...,w_m),\quad w_j\geq 0,\quad j=1,2,...,m.$$ $(a_1,a_2,...,a_m)$ is also called as the weighted coefficient pair. Especially, if $(a_1, a_2,...,a_m)=(1,1,...,1)$, then the system is a gradient one.


\begin{thebibliography}{20}
\bibitem{Alazard20091}
T. Alazard, R. Carles, Loss of regularity for supercritical nonlinear Schr\"{o}dinger equations, {\it Math.
Ann.}, 343(2009), 397--420.

\bibitem{Alazard20092}
T. Alazard, R. Carles, Supercritical geometric optics for nonlinear Schr\"{o}dinger equations, {\it Arch.
Ration. Mech. Anal.}, 194(2009), 315--347.

\bibitem{Barab1984}
J. E. Barab, Nonexistence of asymptotically free solutions for a nonlinear
Schr\"{o}dinger equation, {\it J. Math. Phys.}, 25(1984), 3270--3273.

\bibitem{Beceanu??}
M. Beceanu, Q. Q. Deng, A. Soffer and Y. F. Wu, Large global solutions for nonlinear Sch\"{o}dinger equations III, energy-supercritical cases,
arXiv:1901.07709v1.

\bibitem{Birnir1996}
B. Birnir, C. Kenig, G. Ponce, N. Svanstedt and L. Vega, On the ill-posedness of the IVP for the generalized Korteweg-de Vries and nonlinear Schr\"{o}dinger equations, {\it J. London Math. Soc.}, 53(1996), 551--559.

\bibitem{Bourgain19991}
J. Bourgain, Global well--posedness of defocusing critical nonlinear Sch\"{o}dinger equations in the radial case, {\it J. Amer. Math. Soc.}, 12(1999):145--171.

\bibitem{Bourgain19992}
J. Bourgain, Global Solutions of Nonlinear Schr\"{o}dinger Equations,
American Mathematical Society Colloquium Publications 46, American
Mathematical Society, Providence, RI, 1999.

\bibitem{Burq2005}
N. Burq, P. G\'{e}rard, N. Tzvetkov, Multilinear eigenfunction estimates and global existence for the
three dimensional nonlinear Schr\"{o}dinger equations, {\it Ann. Sci. \'{E}cole Norm. Sup.}, 38(2005), 255--301.


\bibitem{Carles20071}
 R. Carles, Geometric optics and instability for semi-classical Schr\"{o}dinger equations, {\it Arch. Ration.
Mech. Anal.}, 183(2007), 525--553.

\bibitem{Carles20072}
 R. Carles, On instability for the cubic nonlinear Schr\"{o}dinger equation, {\it C. R. Math. Acad. Sci. Paris},
344(2007), 483--486.

\bibitem{Carles2012}
R. Carles, E. Dumas and C. Sparber, Geometric optics and instability for NLS and Davey-Stewartson models, {\it J. Eur. Math. Soc.(JEMS)}, 14(2012), 1885--1921.

\bibitem{Cassano2015}
 B. Cassano and M. Tarulli, $H^1$-scattering for systems of $N$-defocusing weakly coupled NLS equations in low space dimensions, {\it J. Math. Anal. Appl.}, 430(2015), 528--548.

 \bibitem{Cazenave1990}
 T. Cazenave and F. Weissler, The Cauchy problem for the critical nonlinear Schr\"{o}dinger equation in
$H^s$, Nonlinear Anal., 14(1990), 807--836.

\bibitem{Cazenave1992}
 T. Cazenave and F. Weissler, Rapidly decaying solutions of the nonlinear Schr\"{o}dinger equation, {\it
Comm. Math. Phys.}, 147(1992), 75--100.

\bibitem{Cazenave2003}
 T. Cazenave, Semilinear Schr\"{o}dinger equations, Courant Lecture Notes in
Mathematics 10, New York University, Courant Institute of Mathematical
Sciences, AMS, Providence, RI, 2003.

\bibitem{Cho2013}
Y. Cho and S. Lee, Strichartz estimates in spherical coordinates, {\it Indiana Univ. Math. J.}, 62(2013), 991--1020.

\bibitem{Christ1991}
M. Christ and M. Weinstein, Dispersion of small amplitude solutions of the generalized Korteweg-de Vries equation, {\it J. Funct. Anal.}, 100(1991), 87--109.

\bibitem{Christ2003}
M. Christ, J. Colliander, and T. Tao, Asymptotics, frequency modulation, and low regularity ill-posedness for canonical defocusing equations, {\it Amer. J. Math.}, 125(2003), 1235--1293.

\bibitem{Christ20032}
M. Christ, J. Colliander, and T. Tao, Ill-posedness for nonlinear Sch\"{o}dinger and wave equations, arXiv:math/0311048.

\bibitem{Colliander2004}
J. Colliander, M. Keel, G. Staffilani, H. Takaoka and T. Tao, Global existence
and scattering for rough solutions of a nonlinear Schr\"{o}dinger equation in $R^3$, {\it Commun.
Pure Appl. Math.},  57(2004), 987--1014.

\bibitem{Colliander2008}
J. Colliander, M.  Keel, G. Staffilani, H. Takaoka and T. Tao, Global well--posedness and scattering for the energy--critical nonlinear Schr\"{o}dinger equation in $R^3$, {\it Ann. of Math.}, 167(2008), 767--865.

\bibitem{Dodson2012}
B. Dodson, Global well--posedness and scattering for the defocusing, $L^2$-critical, nonlinear
Schr\"{o}dinger equation when $d\geq 3$, {\it J. Amer. Math. Soc.}, 25(2012), 429--463.

\bibitem{Dodson2015}
B. Dodson, Global well-posedness and scattering for the mass critical nonlinear Schr\"{o}dinger equation with mass below the mass of the ground state, {\it Adv. Math.}, 285(2015), 1589--1618.

\bibitem{Dodson20193}
B. Dodson, Global well--posedness and scattering for the
focusing, cubic nonlinear Schr\"{o}dinger
equation in dimension $d=4$, {\it Ann. Sci. $\tilde{E}$c. Norm. Sup\'{e}r.}, 52(2019), 138--180.

\bibitem{Dodson20190}
B. Dodson, Global well-posedness and scattering for nonlinear Schr\"{o}dinger equations with algebraic nonlinearity when $d=2,3$ and $u_0$ is radial, {\it Camb. J. Math.}, 7(2019), 283--318.

\bibitem{Dodson2019}
B. Dodson, Defocusing nonlinear Schr\"{o}dinger equations. Cambridge Tracts in Mathematics, 217. Cambridge University Press, Cambridge, 2019. xii+242 pp. ISBN: 978-1-108-47208-1.

\bibitem{Dodson2017}
B. Dodson, C. Miao, J. Murphy, and J. Zheng, The defocusing quintic NLS in four space dimensions,
{\it Ann. Inst. H. Poincar\'{e} Anal. Non Lin\'{e}aire}, 34(2017), 759--787.

\bibitem{Duyckaerts2008}
T. Duyckaerts, J. Holmer and S. Roudenko, Scattering for the non--radial 3D cubic nonlinear Schr\"{o}dinger equation, {\it Math. Res. Lett.}, 15(2008),   1233--1250.


\bibitem{Escobedo1991}
M. Escobedo and M. A. Herrero, Boundedness and blow up for
a semilinear reaction--diffusion system, {\it J. Differential Equations}, 89 (1991), 176--202.

\bibitem{Escobedo1995}
M. Escobedo and H. Levine, Critical blowup and global existence numbers for a weakly coupled system
of reaction--diffusion equations, {\it Arch. Rational Mech. Anal.}, 129(1995) 47--100.

\bibitem{Farah2017}
L. G. Farah and A. Pastor, Scattering for a 3D coupled nonlinear Schr\"{o}dinger system, {\it J. Math. Phys.}, 58(2017), 071502,33 pp.

\bibitem{Gao2019}
C. W. Gao, Z. H. Zhao, On scattering for the defocusing high dimensional inter-critical NLS, {\it J. DifferentialEquations}, 267(2019), 6198--6215.

\bibitem{Ginibre1979}
J. Ginibre and G. Velo, On a class of nonlinear Schr\"{o}dinger equations, {\it J. Funct. Anal.}, 32(1979), 1--71.

\bibitem{Ginibre19851}
J. Ginibre and G. Velo, Scattering theory in the energy space for a class of nonlinear Schr\"{o}dinger
equations, {\it J. Math. Pures Appl.} 64(1985), 363--401.

\bibitem{Ginibre19852}
J. Ginibre and G. Velo, Time decay of finite energy solutions of the nonlinear Klein--Gordon
and Schr\"{o}dinger equations, {\it Ann. Inst. H. Poincar\'{e} Phys. Th\'{e}or.}, 43(1985), 399--442.

\bibitem{Ginibre1992}
J. Ginibre and G. Velo, Smoothing properties and retarded estimates for some dispersive evolution equations,
{\it Comm. Math. Phys.}, 144(1992), 163--188.

\bibitem{Glassey1977}
R. T. Glassey, On the blowing up of solutions to the Cauchy
problem for nonlinear Schr\"{o}dinger equations, {\it J. Math. Phys}.,
18(1977), 1794--1797.

\bibitem{Grillakis2000}
M. Grillakis, On nonlinear Schr\"{o}dinger equations, {\it Comm. Partial Differential
Equations}, 25(2000), 1827--1844.

\bibitem{Guo2014}
Z. Guo and Y. Wang, Improved Strichartz estimates for a class of dispersive equations in the radial case
and their applications to nonlinear Schr\"{o}dinger and wave equations, {\it J. Anal. Math.}, 124(2014), 1--38.

\bibitem{Holmer2008}
J. Holmer and S. Roudenko, A sharp condition for scattering of the radial
3D cubic nonlinear Schr\"{o}dinger equation, {\it Comm. Math. Phys.}, 282(2008):435--647.

\bibitem{Kato1987}
T. Kato, On nonlinear Schr\"{o}dinger equations, {\it Ann. Inst. H. Poincar\'e Phys. Th\'{e}or.},
46(1987), 113--129.

\bibitem{Keel1998}
M. Keel and T. Tao, Endpoint Strichartz Estimates, {\it Amer. J. Math.}, 120(1998), 955--980.


\bibitem{Kenig2006}
C. E. Kenig and F. Merle, Global well--posedness, scattering and blow-up for the energy-critical,
focusing, nonlinear Schr\"{o}dinger equation in the radial case, {\it Invent. Math.}, 166(2006), 645--675.

\bibitem{Kenig2008}
C. E. Kenig and F. Merle, Global well--posedness, scattering and blow up for the energy critical,
focusing, nonlinear wave equation, {\it Acta Math.}, 201(2008), 147--212.


\bibitem{Kenig2010}
C. Kenig, F. Merle, Scattering for $\dot{H}^{\frac{1}{2}}$ bounded solutions to the cubic, defocusing NLS in 3 dimensions, {\it
Trans. Amer. Math. Soc.}, 362(2010), 1937--1962.

\bibitem{Kenig20111}
 C. Kenig, F. Merle, Nondispersive radial solutions to energy supercritical non--linear wave equations,
with applications, {\it Amer. J. Math.}, 133(2011), 1029--1065.

\bibitem{Kenig20112}
 C. Kenig, F. Merle, Radial solutions to energy supercritical wave equations in odd dimensions, {\it Discrete
Con.  Dyn. Syst.}, 31(2011), 1365--1381.

\bibitem{Kenig2001}
C. Kenig, G. Ponce and L. Vega, On the ill-posedness of some canonical dispersive equations, {\it Duke Math. J.}, 106(2001), 617--633.

\bibitem{Keraani2001}
S. Keraani, On the defect of compactness for the Strichartz estimates for the Sch\"{o}dinger equations.,
{\it J. Differential Equations}, 175(2001), 353--392.

\bibitem{Keraani2006}
S. Keraani, On the blow up phenomenon of the critical nonlinear
Schr\"{o}dinger equation, {\it J. Funct. Anal.}, 235(2006): 171--192.

\bibitem{Killip2017}
R. Killip, S. Masaki, J. Murphy, and M. Visan, Large data mass--subcritical NLS: critical weighted bounds
imply scattering, {\it NoDEA Nonlinear Differential Equations Appl.}, 24(2017), Paper No. 38, 33pp.

\bibitem{Killip2019}
R. Killip, J. Murphy, and M. Visan, Almost sure scattering for the energy-critical NLS
with radial data below $H^1(\mathbb{R}^4)$, {\it Comm. Partial Differential Equations}, 44(2019), 51--71.


\bibitem{Killip20172}
R. Killip, T. Oh, O. Pocovnicu and M. Visan, Solitons and scattering for the cubic-quintic nonlinear Schr\"{o}dinger equation on $\mathbb{R}^3$, {\it Arch. Ration. Mech. Anal.}, 225(2017), 469--548.

\bibitem{Killip2009}
R. Killip, T. Tao and M. Visan, The cubic nonlinear Schr\"{o}dinger equation
in two dimensions with radial data, {\it J. Eur. Math. Soc.(JEMS)}, 11(2009): 1203--1258.

\bibitem{Killip20101}
R. Killip and M. Visan, The focusing energy--critical nonlinear Schr\"{o}dinger equation in dimensions five
and higher, {\it Amer. J. Math.}, 132(2010), 361--424.

\bibitem{Killip20102}
R. Killip and M. Visan, Energy--supercritical NLS:critical $\dot{H}^s$-bounds imply scattering, {\it Comm. Partial
Differential Equations}, 35(2010), 945--987.

\bibitem{Killip20111}
R. Killip and M. Visan, The radial energy-supercritical nonlinear wave equation in three space
dimensions, {\it Trans. Amer. Math. Soc.}, 363(2011), 3813--3934.

\bibitem{Killip20112}
R. Killip and M. Visan, The radial defocusing energy-supercritical nonlinear wave equation in all space
dimensions, {\it Proc. Amer. Math. Soc.}, 139(2011), 1805--1817.

\bibitem{Killip2012}
R. Killip and M. Visan, Global well-posedness and scattering for the defocusing quintic NLS in three dimensions, {\it Anal. PDE}, 5(2012), 855--885.

\bibitem{Killip2013}
R. Killip and M. Visan, Nonlinear Schr\"{o}dinger equations at critical regularity.  Evolution equations, 325--437, Clay Math. Proc., 17, Amer. Math. Soc., Providence, RI, 2013.


\bibitem{Killip20082}
R. Killip, M. Visan, and X. Zhang, The mass--critical nonlinear Schr\"{o}dinger equation with radial data in dimensions
three and higher, {\it  Analysis and PDE}, 1(2008), 229--266.

\bibitem{Krieger2006}
J. Krieger and W. Schlag, Stable manifolds for all monic supercritical focusing nonlinear Schr\"{o}dinger equations in one dimension,
{\it J. Amer. Math. Soc.}, 19(2006), 815--920.


\bibitem{Merle2015}
F. Merle, P. Rapha\"{e}l and I. Rodnianski, Type II blow up for the energy supercritical NLS, {\it Camb.
J. Math.}, 3(2015), 439--617.


\bibitem{Miao2013}
C. X. Miao, G. X. Xu and L. F. Zhao, The dynamics of the 3D radial NLS with the combined
terms, {\it Commun. Math. Phys.}, 318(2013), 767--808.

\bibitem{Miao2014}
C. X. Miao, J. Murphy and J. Zheng, The defocusing energy--supercritical NLS in four space dimensions,
{\it J. Funct. Anal.}, 267(2014), 1662--1724.

\bibitem{Morawetz1968}
 C. Morawetz, Time decay for the nonlinear Klein-Gordon equations, {\it Proc. Roy. London Soc.
Ser. A}, 306(1968), 291--296.

\bibitem{Murphy2014}
J. Murphy, Intercritical NLS: critical $\dot{H}^s$-bounds imply scattering, {\it SIAM J. Math. Anal.}, 46(2014), 939--997.

\bibitem{Murphy20142}
J. Murphy, The defocusing $\dot{H}^{\frac{1}{2}}$-critical NLS in high dimensions, {\it
Discrete Contin. Dyn. Syst. Series A}, 34(2014):733--748.

\bibitem{Murphy2015}
J. Murphy, The radial defocusing nonlinear Schr\"{o}dinger equation in three space dimensions, {\it Comm.
Partial Differential Equations}, 40(2015), 265--308.

\bibitem{Nakanishi2012}
K. Nakanishi and W. Schlag, Global dynamics above the ground state energy
for the cubic NLS equation in 3D, {\it Calc. Var. Partial Differential Equations}, 44(2012), 1--45.


\bibitem{Roy2011}
T. Roy, Scattering above energy norm of solutions of a loglog energy--supercritical Schr\"{o}dinger equation
with radial data, {\it J. Differential Equations}, 250(2011), 292--319.

\bibitem{Ryckman2007}
E. Ryckman, M. Visan, Global well--posedness and scattering for the defocusing
energy--critical nonlinear Schr\"{o}dinger equation in $\mathbb{R}^{1+4}$, {\it Amer. J. Math.}, 129(2007), 1--60.


\bibitem{Schlag2009}
W. Schlag, Stable manifolds for an orbitally unstable nonlinear Schr\"{o}dinger equation, {\it Ann. of Math.}, 169(2009), 139--227.

\bibitem{Schlag20101}
W. Schlag, A. Soffer and W. Staubach, Decay for the wave and Schr\"{o}dinger evolutions on manifolds with conical ends. I., {\it
Trans. Amer. Math. Soc.}, 362(2010), 19--52.

\bibitem{Schlag20102}
W. Schlag, A. Soffer and W. Staubach, Decay for the wave and Schr\"{o}dinger evolutions on manifolds with conical ends. II., {\it
Trans. Amer. Math. Soc.}, 362(2010), 289--318.

\bibitem{Shao20091}
S. Shao, Maximizers for Strichartz inequalities and Sobolev-Strichartz inequalities for Schr\"{o}dinger equation, {\it Electron. J. Differential
 Equations}, 2009, No. 3, 13pp.

\bibitem{Shao20092}
S. Shao, Sharp linear and bilinear restriction estimates for paraboloids in the cylindrically symmetric
case, {\it Rev. Mat. Iberoam.}, 25(2009), 1127--1168.


\bibitem{Strauss1981}
W. A. Strauss, Nonilinear scattering theory at low energy, {\it J. Funct. Anal.}, 41(1981), 110--133.

\bibitem{Strichartz1977}
R. S. Strichartz, Restriction of Fourier transforms to quadratic surfaces and decay of solutions of wave equations,
{\it Duke Math. J.}, 44(1977), 705--714.

\bibitem{Tao2004}
T. Tao, On the asymptotic behavior of large radial data for a focusing non-linear Schr\"{o}dinger equation,
{\it Dyn. Partial Differ. Equ.}, 1(2004), 1--47.

\bibitem{Tao20051}
T. Tao, Global well--posedness and scattering for the higher--dimensional energy critical
nonlinear Schr\"{o}dinger equation for radial data, {\it  New York J. Math.},
11(2005), 57--80.

\bibitem{Tao2006}
T. Tao, Nonlinear dispersive equations. Local and global analysis. CBMS Regional Conference Series in Mathematics, 106. Published for the Conference Board of the Mathematical Sciences, Washington, DC; by the American Mathematical Society, Providence, RI, 2006. xvi+373 pp. ISBN: 0-8218-4143-2.

\bibitem{Tao2009}
Tao, Terence, Global existence and uniqueness results for weak solutions of the focusing masscritical
nonlinear Schr\"{o}dinger equation, {\it Anal. PDE}, 2(2009), 61--81.

\bibitem{Tao2018}
T. Tao, Finite time blowup for a supercritical defocusing nonlinear Schr\"{o}dinger system, {\it Anal. PDE}, 11(2018), 383--438.

\bibitem{Tao20052}
T. Tao, M. Visan, Stability of energy--critical nonlinear Schr\"{o}dinger equations in
high dimensions, {\it Electron. J. Differetial. Equations}, 118(2005), 1--28.

\bibitem{Tao20071}
T. T. Tao,  M. Visan and X. Y. Zhang, The nonlinear Schr\"{o}dinger equation with combined power--type nonlinearities, {\it Comm. Partial Differential Equations}, 32 (2007), 1281--1343.

\bibitem{Tao20072}
T. Tao, M. Visan and X. Zhang, Global well--posedness and scattering for the mass--critical nonlinear Schr\"{o}dinger
equation for radial data in high dimensions, {\it Duke Math. J.}, 140(2007), 165--202.

\bibitem{Tao2008}
 T. Tao, M. Visan, and X. Zhang, Minimal-mass blowup solutions of the mass-critical NLS, {\it Forum Math}, 20
(2008), 881--919.

\bibitem{Tsutsumi1985}
Y. Tsutsumi, Scattering problem for nonlinear Schr\"{o}dinger equations, {\it Ann. Inst.
H. Poincar\'{e} Phys. Th\'{e}or.}, 43(1985), 321--347.

\bibitem{Tsutsumi1990}
Y. Tsutsumi, Rate of $L^2$ concentration of blow--up solutions for the nonlinear
Schr\"{o}dinger equation with the critical power, {\it Nonlinear Anal.}, 15(1990), 719--
724.

\bibitem{Tsutsumi1984}
Y. Tsutsumi and K. Yajima, The asymptotic behavior of nonlinear Schr\"{o}dinger
equations, {\it Bull. Amer. Math. Soc.}, 11(1984), 186--188.

\bibitem{Visan2006}
M. Visan, The Defocusing Energy--Critical Nonlinear Schr\"{o}dinger Equation in
Dimensions Four and Higher, Ph.D. Thesis, UCLA. (2006)

\bibitem{Visan2007}
M. Visan, The defocusing energy-critical nonlinear Schr\"{o}dinger equation in higher
dimensions, {\it Duke Math. J.}, 138(2007), 281--374.

\bibitem{Visan2011}
M. Visan, Global well--posedness and scattering for the defocusing cubic nonlinear Schr\"{o}dinger equation in four dimensions, {\it
Int. Math. Res. Not. IMRN}, 2012, 1037--1067.

\bibitem{Wang2016}
W. Wang, Energy supercritical nonlinear Schr\"{o}dinger equation: quasiperiodic solutions, {\it Duke Math. J.},
165(2016), 1129--1192.

\bibitem{Xie2013}
J. Xie and D. Fang, Global well--posedness and scattering for the defocusing $\dot{H}^s$--critical NLS, {\it Chin. Ann.
Math. Ser. B}, 34(2013), 801--842.

\bibitem{Xu2016}
Y. S. Xu, Global well-posedness, scattering, and blowup for nonlinear coupled Schr\"{o}dinger equations in $\mathbb{R}^3$, {\it Appl. Anal.}, 95(2016), 483--502.

\bibitem{Zhang2006}
X. Zhang, On Cauchy problem of 3-D energy critical Schr\"{o}dinger equation with
subcritical perturbations, {\it J. Differential Equations}, 230(2006), 422--445.


\end{thebibliography}
\end{document}